\documentclass[nonamelimits,leqno,centertags,a4paper]{amsbook}  
\usepackage{A4} \usepackage{amssymb} \usepackage{latexsym}
\usepackage{upref} \usepackage{amsfonts,amssymb,amsmath}
\usepackage[latin1]{inputenc}  %æøå
\usepackage[english]{babel}

\numberwithin{equation}{chapter}  %%for amsbook %se p.27 i guide
\numberwithin{subsection}{chapter} %%for amsbook
\newtheorem{thm}{Theorem}[chapter] %for amsbook

\swapnumbers

\newtheorem{lemma}[thm]{Lemma}
\newtheorem{cor}[thm]{Corollary}
\newtheorem{prop}[thm]{Proposition}

\newtheorem{defn}[thm]{Definition}

\theoremstyle{definition}
\newtheorem{rmk}[thm]{Remark}
\newtheorem{exmp}[thm]{Example}

\newcommand{\Mdef}[2]{\newcommand{#1}{\relax \ifmmode #2 \else $#2$\fi}}

\newcommand{\set}{\setcounter{equation}{\value{thm}}}
\newcommand{\add}{\addtocounter{thm}{1}}

\newcommand{\pgl}[1]{\PGL(#1,\C)}
\newcommand{\pslr}[1]{\PSL(#1,\R)}
\newcommand{\di}{\mathrm{DI}(4)}
\newcommand{\slr}[1]{\SL(#1,\R)}
\newcommand{\dm}[1]{(\Z/2)^{#1}}

\newcommand{\pcg}{$p\/$-compact group}
\newcommand{\pctg}{$p\/$-compact toral group}
\newcommand{\pct}{$p\/$-compact torus}
\newcommand{\lmn}{elementary abelian $p\/$-group}

\newcommand{\pl}{preferred lift}

\newcommand{\epct}{extended \pct}

\newcommand{\mt}{maximal torus}
\newcommand{\mts}{maximal tori}
\newcommand{\mtn}{maximal torus normalizer}

\newcommand{\cd}{\operatorname{cd}}

\newcommand{\gen}[1]{\left\langle {#1} \right\rangle}
\newcommand{\op}{\operatorname{op}}
\newcommand{\ch}[1]{{\Check #1}}
\newcommand{\order}[1]{\lvert #1 \rvert}

\newcommand{\he}{homotopy equivalence}
\newcommand{\m}{morphism}
\newcommand{\ses}{short exact sequence}
\newcommand{\func}[3]{\mbox{$#1\colon #2\rightarrow #3$}}

\newcommand{\twocg}{$2\/$-compact group}
\newcommand{\twoctg}{$2\/$-compact toral group}

\newcommand{\twoct}{$2$-compact torus}
\newcommand{\twocts}{$2$-compact tori}
\newcommand{\etwoct}{extended \twoct}
\newcommand{\etwocts}{extended \twocts}
\newcommand{\mtnp}{\mtn\ pair}
\newcommand{\lmntwo}{elementary abelian $2$-group}
\newcommand{\AM}{\operatorname{AM}}

\newcommand{\N}{N}  %% put % in front of \N from defs2e.tex 
\newcommand{\T}{T}
\newcommand{\W}{W}
\newcommand{\Ze}{Z}

\newcommand{\A}{{\mathbf A}}

\newcommand{\I}{{\mathbf I}}

\newcommand{\C}{{\mathbf C}}
\newcommand{\Ha}{{\mathbf H}}
\newcommand{\Z}{{\mathbf Z}}
\newcommand{\Q}{{\mathbf Q}}
\newcommand{\R}{{\mathbf R}}
\newcommand{\F}{{\mathbf F}}
\newcommand{\Ab}{{\mathbf{Ab}}}

\newcommand{\Topspaces}{{\mathbf{Top}}}

\newcommand{\map}{\operatorname{map}}

\newcommand{\Rep}{\operatorname{Rep}}

\newcommand{\Hom}{\operatorname{Hom}}

\newcommand{\Aut}{\operatorname{Aut}}
\newcommand{\End}{\operatorname{End}}

\newcommand{\Ob}{\operatorname{Ob}}

\newcommand{\Out}{\operatorname{Out}}
\newcommand{\GL}{\mathrm{GL}}
\newcommand{\PSL}{\mathrm{PSL}}

\newcommand{\St}{\operatorname{St}}
\newcommand{\im}{\operatorname{im}}
\newcommand{\diag}{\operatorname{diag}}

\newcommand{\tr}{\operatorname{tr}}

\newcommand{\hocolim}{\operatorname{hocolim}}

\newcommand{\U}{\mathrm{U}}
\newcommand{\SU}{\mathrm{SU}}
\newcommand{\SO}{\mathrm{SO}}
\newcommand{\Or}{\mathrm{O}}
\newcommand{\Pin}{\mathrm{pin}}
\newcommand{\Spin}{\mathrm{Spin}}
\newcommand{\Symp}{\mathrm{Sp}}
\newcommand{\Ffour}  {{\mathrm{F}_4}}
\newcommand{\Gtwo}{{\mathrm{G}_2}}
\newcommand{\E}{{\mathrm{E}}}

\newcommand{\B}{{\mathrm B}}

\newcommand{\SL}{\mathrm{SL}}
\newcommand{\PGL}{\operatorname{PGL}}

\newcounter{opg}
\newenvironment{simplelist}%
{\begin{list}
    {(\arabic{opg})}{\usecounter{opg}%
    \setlength{\leftmargin}{6pt}%
    \setlength{\itemsep}{0.1\itemsep}%
    \setlength{\parsep}{0\parsep}
    \setlength{\topsep}{0\topsep}}}%
{\end{list}}

\input xypic
\xyoption{all}

\CompileMatrices
\newdir{ >}{{}*!/-5pt/\dir{>}} %\ar@{ >->} gives a monomorphism arrow

\begin{document}

%%%%% Heading  %%%%%%
\title{$N$-determined $2\/$-compact groups}
\author[J.M. Møller]{Jesper M.\ Møller}
\address{Matematisk Institut\\
  Universitetsparken 5\\
  DK--2100 København} 
\email{moller@math.ku.dk}
%\urladdr{\url{http://www.math.ku.dk/~moller}} 
\date{\today}
\begin{abstract}
  We first formulate a general scheme for the classification of \twocg
  s in terms of \mtnp s. Applying this scheme, we show that all simple
  \twocg s are $N$-determined. We also compute auto\m\ groups in many
  cases. As an application we confirm the splitting conjecture
  formulated by Dwyer and Wilkerson.
\end{abstract}

\maketitle
\tableofcontents

\chapter{Introduction}
\label{sec:intro}

A $p$-compact group, where $p$ is a prime number, is a $p$-complete
space $BX$ whose loop space $X=\Omega BX$ has finite mod $p$ singular
cohomology. If $G$ is a Lie group and $\pi_0(G)$ is a finite $p$-group
then the $p$-completed classifying space of $G$ is a $p$-compact
group. The Sullivan spheres $(BS^{2n-1})_p^{\wedge}$, $n \vert (p-1)$,
or, more generally, the Clark--Ewing spaces \cite{ce} are also
examples of \pcg s.  These homotopy Lie groups were defined and
explored by W.G.\ Dwyer and C.W.\ Wilkerson in a series of papers
\cite{dw:fixpt, dw:center, dw:split, dw:new}. (Consult the survey
articles \cite{dwyer:survey, lannes:rep, no:survey, jmm:survey} for a
quick overview.) They show that any \pcg\ $BX$ has a \mt\ $BT \to BX$
and a Weyl group acting on the \mt . The Borel construction for this
action is the total space $BN$ of a fibration \cite[9.8]{dw:fixpt}
\begin{equation*}
  BT \to BN \to BW
\end{equation*}
whose fibre is the \mt\ and whose base space is the classifying space
of the Weyl group. By construction the mono\m\ $BT \to BX$ extends to
a mono\m\ $BN \to BX$. (Strictly speaking, $BN$ is in general not a
\pcg\ as its fundamental group may not be a finite $p$-group; instead,
$BN$ is an example of an \epct\ \cite[3.12]{dw:center}.)  In case $G$
is a compact connected Lie group, the \mtn s in the Lie and in the
\pcg\ sense are essentially identical.  The (discrete approximation to
the) \mtn\ of the Sullivan sphere is the semi-direct product
$\Z/p^{\infty} \rtimes \Z/n$ for the action of the cyclic group $\Z/n
< \Z/(p-1)=\Z_p^{\times} = \Aut(\Z/p^{\infty})$ on $\Z/p^{\infty}$.
It is a conjecture, suggested by the analogous situation for connected
compact Lie groups \cite{cww} (and some nonconnected ones
\cite{hammerli:thesis, hammerli:remarks}), that $BN$ determines $BX$.
This classification conjecture has actually been verified for odd
primes \cite{jmm:deter,jmm:ndet,agmv:efamily}.  For $p=2$, however,
only scattered results are known. 
%First, we have to be
%clear about what it means for two \twocg s to have the same \mtn .

The first obstacle for a classification of \twocg s in terms of their
\mtn s is, in contrast to the odd $p$ case, that the \mtn\ does not
retain information about component groups . For instance, the
nonconnected \twocg\ $B\Or(2n)$ and the connected \twocg\ $B\SO(2n+1)$
have identical \mtn s. Thus the \mtn\ in itself is simply not a strong
enough \twocg\ invariant.  That is why we replace \mtn s by \mtnp s.
%One way to store information about
%component groups is to replace \mtn s by \mtn\ {\em pairs}.

Consider a \twocg\ $BX$ with identity component $BX_0 \to BX$ and
component group $\pi_0(X)=X/X_0$. Then there is a map of fibrations
\begin{equation*}
  \xymatrix{
  BX_0 \ar[d] & BN_0 \ar[l]_-{Bj_0} \ar[d] \\
  BX \ar[d] & BN \ar[l]_-{Bj} \ar[d] \\
  B\pi_0(X) & B\pi \ar[l]_-{\cong} }  
\end{equation*}
where the top two horizontal arrows are \mtn s and the bottom
horizontal arrow is an iso\m\ between $\pi=N/N_0$ and the component
group. We say that $(BN,BN_0) \to (BX_0,BX)$ is a {\em \mtnp\/} for
the \twocg\ $BX$ (\ref{defn:mtnp}).  Accordingly, two \twocg s, $BX$
and $BX'$, have the same \mtnp\ if there exists a commutative diagram
\begin{equation*}
  \xymatrix{
    BX_0 \ar[d] & 
    BN_0 \ar[l]_{Bj_0} \ar[r]^{Bj_0'} \ar[d] & BX_0' \ar[d]\\
    BX \ar[d] & BN  \ar[l]_{Bj} \ar[r]^{Bj'} \ar[d] & BX' \ar[d] \\
    B\pi_0(X) & B\pi \ar[l]_-{\cong} \ar[r]^-{\cong} & B\pi_0(X')}
\end{equation*}
where the horizontal maps are \mtn s of \twocg s.

%Let $BN$ be an \etwoct\  with a normal maximal rank subgroup $BN_0
%\to BN$; this simply means that there is a fibration sequence of the
%form $BN_0 \to BN \to B\pi$ where $\pi=N/N_0$ is a finite group. Consider a
%\twocg\ $BX$ with identity component $BX_0 \to BX$ and component group
%$\pi_0(X)=X/X_0$.  We shall say that the pair $(BN,BN_0)$ is a {\em
%  \mtnp\/} for the \twocg\ $BX$ (\ref{defn:mtnp}) if there
%exists a map of fibrations
%\begin{equation*}
%  \xymatrix{
%  BX_0 \ar[d] & BN_0 \ar[l]_-{Bj_0} \ar[d] \\
%  BX \ar[d] & BN \ar[l]_-{Bj} \ar[d] \\
%  B\pi_0(X) & B\pi \ar[l]_-{\cong} }  
%\end{equation*}
%where the horizontal maps are \mtn s of \twocg s (so that the map
%between the base spaces is a \he ).  In particular, the \mtnp\ 
%determines the component group in that $\pi_0(X) \cong N/N_0$.

%According to this
%definition, two \twocg s, $BX$ and $BX'$ have the same \mtnp\ if there
%exists a commutative diagram
%\begin{equation*}
%  \xymatrix{
%    BX_0 \ar[d] & 
%    BN_0 \ar[l]_{Bj_0} \ar[r]^{Bj_0'} \ar[d] & BX_0' \ar[d]\\
%    BX \ar[d] & BN  \ar[l]_{Bj} \ar[r]^{Bj'} \ar[d] & BX' \ar[d] \\
%    B\pi_0(X) & B\pi \ar[l]_-{\cong} \ar[r]^-{\cong} & B\pi_0(X')}
%\end{equation*}
%where the horizontal maps are \mtn s of \twocg s.

We can now give a precise formulation of the classification conjecture
for \twocg s.
We shall say that a \twocg\ $BX$ with \mtnp\ $(BN,BN_0) \to (BX_0,BX)$
is {\em totally $N$-determined\/} (\ref{defn:det}) if
\begin{enumerate}
\item auto\m s of $BX$ are determined by their restrictions to $BN$, and
\item for any other connected \twocg\ $BX'$ with the same \mtnp\ as
  $BX$ there exist an iso\m\ \func{Bf}{BX}{BX'} and an
  auto\m\ \func{B\alpha}{BN}{BN}, inducing the identity map on
  homotopy groups, making the diagram
\begin{equation*}
  \xymatrix{
      BN \ar[d] \ar[r]^{B\alpha}_{\cong} & BN \ar[d] \\
      BX \ar[r]^{\cong}_{Bf} & BX' }
\end{equation*}
commutative  
\end{enumerate}
We shall say that $BX$ is
{\em uniquely $N$-determined\/} if in addition the auto\m s of $BX$
are determined by their effect on the (two nontrivial) homotopy groups
of $BN$. (See Lemma~\ref{defcons} for a justification of the
terminology.)
The role of $B\alpha$ is to compensate for the auto\m s of $BN$ that
do not extend to auto\m s of $BX$; such auto\m s do exist when $p=2$
whereas they do not occur at odd primes. 

Here is the main result of this paper. The first item confirms a
conjecture from \cite[5.1, 5.2]{dwyer:survey}.

\begin{thm}\label{thm:conclusion}
  \begin{enumerate}
  \item  Any connected \twocg\ is isomorpic to $BG \times B\di^t$ for some
  compact connected Lie group $G$ and some integer $t \geq 0$. 
\item  All connected \twocg s are uniquely $N$-determined.
\item All LHS
  (\ref{LHSinit}) \twocg s are totally $N$-determined. 
  \end{enumerate}
\end{thm}

The proof of this theorem, following the inductive principle of
\cite[9.1]{dw:center}, has two stages. The first stage consists in
the reduction of the problem to case of connected, simple and
centerless \twocg s. The next stage is an inductive case-by-case
checking of the simple \twocg s.

The classification conjecture for \twocg s can be reduced to the case
of connected, centerless, and simple \twocg s because $N$-determinism
is to a large extent hereditary. It turns out
(Chapter~\ref{cha:ndet}.\ref{sec:red}) that, under certain extra
conditions, a \twocg\ with an $N$-determined identity component, is
itself $N$-determined. Also, a connected \twocg\ whose adjoint form,
the quotient by the center, is $N$-determined, is itself
$N$-determined. Since any connected \twocg\ with trivial center is a
product of simple \twocg s \cite{dw:split,no:split} and products of
$N$-determined \twocg s are $N$-determined, this reduces the problem
to the case of a connected, simple \twocg\ with trivial center. This
is the first stage in the classification procedure.

Furthermore, any connected and centerless \twocg\ can
be decomposed into a homotopy colimit of a system of \twocg\ of
smaller dimension \cite{dw:codec} and, under certain hypotheses
(\ref{ndetauto}, \ref{indstepalt}), $N$-determinism is hereditary also
under homotopy colimits. Thus there is a theoretical possibility of
proving the classification conjecture by induction over the
dimension. This is the second stage in the classification procedure.

In this paper we go through the classification procedure for the
\twocg s associated to the classical matrix Lie groups of the infinite
$A$-, $B$-, $C$- and $D$-families, for the exceptional compact Lie
groups $\Gtwo$, $\Ffour$, $E_6$, $PE_7$ and $E_8$, and for the exotic
\twocg\ $\di$ from \cite{dw:new}. The main difficulty is to show that
the obstruction groups of \ref{ndetauto} and \ref{indstepalt} vanish.
For this we use the computer algebra program {\em magma}.

The main results are the theorems below (whose proofs are in
Chapter~\ref{chp:di4andf4} and in Chapter~\ref{cha:proofs}).

Any connected Lie group $G$ has an associated \twocg\ obtained as the
$2$-completion of the classifying space of $G$. We shall denote the
\twocg\ associated to $G$ by $G$ also. As \twocg s, $\SL(2n+1,\R)$ and
$\SO(2n+1)$, for instance, are synonyms because their classifying spaces
are homotopy equivalent.

\begin{thm}\label{thm:afam}
  The simple \twocg\ $\pgl{n+1}$, $n \geq 1$, is uniquely
  $\N$-determined, and its auto\m\ group is 
  \begin{equation*}
    \Aut(\pgl{n+1}) = 
    \begin{cases}
      \Z^{\times} \backslash  \Z_2^{\times} & n=1 \\
       \Z_2^{\times} & n>1 
    \end{cases}
  \end{equation*}
\end{thm}

\begin{thm}\label{mainthm}
  The simple \twocg s $\pslr{2n}$, $n \geq 4$, $\SL(2n+1,\R)$, $n \geq
  2$, and $\PGL(n,\Ha)$, $n \geq 3$,
  are uniquely $N$-determined for all $n \geq 1$. Their auto\m\ groups
  are 
  \begin{equation*}
    \Aut(\PSL(2n,\R)) = 
    \begin{cases}
       \Z^{\times} \backslash  \Z_2^{\times} \times \Sigma_3  & n=4 \\
      \Z^{\times} \backslash  \Z_2^{\times} \times \gen{c_1} 
       & \textmd{$n> 4$ even} \\
      \Z_2^{\times} & \textmd{$n>4$ odd}
    \end{cases} 
  \end{equation*}
  where $\gen{c_1}$ is a group of order two, $\Aut(\SL(2n+1,\R))
  = \Z^{\times} \backslash \Z_2^{\times}$ for $n \geq 2$, and
  $\Aut(\PGL(n,\Ha)) = \Z^{\times} \backslash \Z_2^{\times}$ for $n
  \geq 3$.
  %Are they cohomologically unique? What
  %are the auto\m\ groups?
\end{thm}

\begin{thm} \cite[1.3]{antonio:bg2} \label{thm:g2}
  The simple \twocg\ $\Gtwo$ is uniquely ${\N}$-determined and its auto\m\ 
  group is $\Aut(\Gtwo) =\Z^{\times} \backslash \Z_2^{\times} \times
  C_2$.
\end{thm}

\begin{thm}\label{thm:di4}
  The simple \twocg\ $\di$ is uniquely $N$-determined and its auto\m\ 
  group is $\Aut(\di)=\Z^{\times} \backslash \Z_2^{\times}$.
\end{thm}

\begin{thm}\label{thm:f4}
  The simple \twocg\ $\Ffour$ is uniquely $N$-determined and its
  auto\m\ group is $\Aut(\Ffour) = \Z^{\times} \backslash
  \Z_2^{\times}$. 
\end{thm}

\begin{thm}\label{thm:efam}
  The simple \twocg s $E_6$, $PE_7$ and $E_8$ are uniquely $N$-determined.
\end{thm}

The methods are not limited to simple nor even to connected \twocg s.
Here are two examples of the type of consequences that can obtained
for more general \twocg s.

\begin{cor}\label{cor:afam} \cite[1.9]{mn:maxtorus}
  The \twocg\ $\GL(n,\C)$ is uniquely $\N$-determined and its auto\m\ 
  group is
  \begin{equation*}
 \Aut(\GL(n,\C)) = 
 \begin{cases}
   \Z^{\times} \backslash
  \Aut_{\Z_2\Sigma_2}(\Z_2^2) & n=2 \\
 \Aut_{\Z_2\Sigma_n}(\Z_2^n) & n>2
 \end{cases}
  \end{equation*}
\end{cor}

\begin{cor}\label{cor:glnR}
  The \twocg\ 
  $\GL(n,\R)$ is totally $N$-determined for all $n \geq 2$ and its
  auto\m\ group is 
  \begin{equation*}
    \Aut(\GL(n,\R)) = 
    \begin{cases}
      \Z^{\times} \backslash \Z_2^{\times} & \text{$n \geq 3$ odd} \\
      \Z_2^{\times} & n=2 \\
     \Z_2^{\times} \times \gen{\delta} & \text{$n \equiv 2 \bmod 4$, $n
        > 2$} \\
      \Z^{\times} \backslash \Z_2^{\times} \times \gen{c_1} \times
      \gen{\delta}  &  \text{$n \equiv 0 \bmod 4$}
    \end{cases}
  \end{equation*}
  where $\gen{\delta}$ and $\gen{c_1}$ are subgroups of order two.
\end{cor}

Related uniqueness results can be found in the papers
\cite{no:uniqueness, no:o(n), no:DI(4), antonio:bg2, antonio:f4p2,
  antonio:spn} by 
Dietrich Notbohm, Antonio Viruel and Ale{\v{s}} Vavpeti{\v{c}},
respectively.

\chapter{$\N$-determined $2\/$-compact groups}
\label{cha:ndet}

This chapter contains the fundamental definitions and the first
general results. Whereas \pcg s are determined by their \mtn s
\cite{jmm:ndet, agmv:efamily} when $p>2$, a finer invariant is needed
for \twocg s as there are examples (\ref{onson}) of distinct \twocg s
with identical \mtn s.

\section{Maximal torus normalizer pairs}
\label{sec:Z2refl}\add
Let ${\N}_0 \to {\N}$ be a maximal rank normal mono\m\ between two
extended $2$-compact tori, meaning simply that there exists a \ses\ 
\cite[3.2]{dw:fixpt} of loop spaces ${\N}_0 \to {\N} \to \pi$ for some
finite group $\pi$. For a \twocg , $X$, let $(X,X_0)$ be the
pair consistsing of $X$ and its identity component $X_0$. Then there
is a \ses\ $X_0 \to X \to \pi_0(X)$ of loop spaces where
$\pi_0(X)=X/X_0$ is a finite $2$-group, the component group of $X$.

\begin{defn}\label{defn:mtnp}
If there exists a \m\ of loop space \ses s \cite[2.1]{dw:center}
\begin{equation*}
  \xymatrix{
       {\N}_0 \ar[r] \ar[d]_{j_0} &
       {\N}   \ar[r] \ar[d]^{j} &
       {\pi} \ar[d]^{\cong} \\
       X_0 \ar[r] & X \ar[r] & {\pi_0(X)} }
\end{equation*}
where \func{j_0}{\N_0}{X_0} and \func{j}{\N}{X} are \mtn s
\cite[9.8]{dw:fixpt}, and $\pi \to \pi_0(X)$  an iso\m\ of finite
$2$-groups, then we say that $(N,N_0)$ is a \mtnp\ for $(X,X_0)$.
\end{defn}

A \mtnp\ for $X$ determines the \mt\ ${\T}(X)$, isomorphic to the
identity component of $N$, the Weyl groups, ${\W}(X)=\pi_0(N)$ and
${\W}(X_0)=\pi_0(N_0)$, of $X$ and $X_0$, the component group
$\pi_0(X) = N/N_0 ={\W}(X)/{\W}(X_0)$ \cite[3.8]{mn:center}, and
\cite[7.5]{dw:center} the center ${\Ze}(X_0) \to X_0$ of $X_0$
\cite{mn:center,dw:center}.

\begin{exmp}\label{onson}
  (1) Since $\GL(2,\R) \wr \Sigma_n = N(\SL(2n+1,\R)) \subset
  \GL(2n,\R) \subsetneq \SL(2n+1,\R)$, $\GL(2n,\R)$ and $\SL(2n+1,\R)$
  (Chp~\ref{sec:slodd})
  have the same \mtn .  Their \mtnp s are distinct, however, as their
  component groups are distinct. 
% for $\SL(2n+1,\R)$ is connected and $\Or(2n)$ disconnected.\\
%$(N,N_0)(\Or(2n))=\big(N(\Or(2n)),N(\SO(2n))\big)$ and
%$N(\SO(2n+1))=\big(N(\Or(2n)),N(\Or(2n))\big)$ are distinct \mtnp
%s. 

  \noindent
  (2) More generally \cite{hammerli:thesis}, let $G$ be any compact
  connected Lie group and $\N(G)$ its \mtn . If $\N(G)$ is not maximal,
  there exists a compact Lie group $H$ such that $\N(G) \subseteq H
  \subsetneq G$. The two compact Lie groups, $G$ and $H$, have
  isomorphic \mtn s but distinct \mtn\ pairs as $H$ is
  nonconnected \cite{bo9}. 
  %\cite[Exercise 6, p.\ 113]{bo9}.

  \noindent
  3.\ The Weyl groups for $\SL(2n+1,\R)$ (Chp~\ref{sec:slodd}) and
  $\GL(n,\Ha)$, $n \geq 3$, (Chp~\ref{sec:spn}) are isomorphic as
  reflection groups but $\N(\SL(2n+1,\R))$ is a split and
  $\N(\GL(n,\Ha))$ a nonsplit extension \cite{cww, matthey:second} of
  the Weyl group by the \mt . Thus connected \twocg s can not be
  classified by their Weyl group alone.
\end{exmp}

\setcounter{subsection}{\value{thm}}
\subsection{The Adams--Mahmud homo\m}
\label{sec:AM}\add

For a \twocg\ (or \etwoct\ \cite[3.12]{dw:center}) $X$, we let
$\End(X)=[BX,*;BX,*]$ denote the monoid of pointed homotopy classes of
self-maps of $BX$. The {\em auto\m\ group\/} $\Aut(X) \subseteq
[BX,*;BX,*]$ of $X$ is the group of invertible elements in $\End(X)$
and the {\em outer auto\m\ group\/} $\Out(X)= \pi_0(X) \backslash
\Aut(X) \subseteq [BX;BX]$ is the group of conjugacy classes (free
homotopy classes \cite[2.1]{dw:center}) of auto\m s of $X$.

Let $X$ be a \twocg\ with \mtnp\ $(\N,\N_0)$.  Turn the \mtn\ 
\func{Bj}{B\N}{BX} into a fibration. Any auto\m\ \func{f}{X}{X} of the
\twocg\ $X$ restricts to an auto\m\ \func{\AM(f)}{\N}{\N} of the \mtn
, unique up to the action of the Weyl group ${\W}(X_0)=\pi_1(X/\N)$
\cite[3.8, 5.6.(1)]{mn:center} of the identity component $X_0$ of $X$,
such that the diagram
\begin{equation*}
  \xymatrix@C=50pt{B\N \ar[r]^{B(\AM(f))} \ar[d]_{Bj} & B\N \ar[d]^{Bj} \\
            BX \ar[r]_{Bf} & BX}  
\end{equation*}
commutes up to based homotopy \cite[\S3]{jmm:deter}. The Adams--Mahmud
homo\m\ is the
resulting homo\m\ 
\set\begin{equation}\label{AM}
  \func{\AM}{\Aut(X)}{{\W}(X_0) \backslash \Aut(\N)}
\end{equation}\add
of auto\m\ groups. 

The auto\m\ group of $\N$ sits\set\ \cite[5.2]{jmm:ext} in a
\ses\ ($\ch{\T(X)}$ is the discrete approximation \cite[6.5]{dw:fixpt}
to $\T(X)$)
\begin{equation}\label{eq:H1TW} 
 0 \to
  H^1({\W}(X);\ch{{\T}}(X)) \to \Aut(\N) \xrightarrow{\pi_*} 
  \Aut({\W}(X),\ch{{\T}}(X),e(X)) \to 1 
\end{equation} 
where\add\ the normal subgroup to the left consists of all auto\m s of
$\N$ that induce the identity on homotopy groups and the group to the
right consists of all pairs
$(\alpha,\theta)\in\Aut({\W}(X))\times\Aut(\ch{{\T}}(X))$ such that
$\theta$ is $\alpha$-linear and the induced auto\m\ 
$H^2(\alpha^{-1},\theta)$ \cite[6.7.6]{weibel} preserves the extension
class $e(X)\in H^2({\W}(X);\ch{{\T}}(X))$.  The image of ${\W}(X_0)$
in $\Aut(\N)$ does not intersect the subgroup $
H^1({\W}(X);\ch{{\T}}(X))$ (as ${\W}(X_0)$ is represented faithfully
in $\Aut(\ch{{\T}}(X))$ \cite[9.7]{dw:fixpt}) so there is an
induced\set\ \ses\ 
\begin{equation}\label{autNW0} 0 \to
  H^1({\W}(X);\ch{{\T}}(X)) \to {\W}(X_0) \backslash\Aut(\N) \xrightarrow{\pi_*}
  {\W}(X_0) \backslash \Aut({\W}(X),\ch{{\T}}(X),e(X)) \to 1
\end{equation}\add
whose middle term is the target of the Adams--Mahmud homo\m . In
particular, if $X$ is {\em connected}, this\set\ \ses\ 
\begin{equation}\label{outN} 0 \to
  H^1({\W}(X);\ch{{\T}}(X)) \to \Out(\N) \xrightarrow{\pi_*}
  {\W}(X) \backslash \Aut({\W}(X),\ch{{\T}}(X),e(X)) \to 1
\end{equation}\add
has the group $\Out(\N)={\W}(X) \backslash \Aut(\N)$ of outer auto\m s
of $\N$ as its middle term. The group $ \Aut({\W}(X),\ch{{\T}}(X),0)$
may also be described as the normalizer $\N_{\GL({ L}(X))}({\W}(X))$
of ${\W}(X)$ in $\GL({ L}(X))$ where ${ L}(X) = \pi_2(B{\T}(X))$. This
group evidently fits into an exact sequence \cite[\S2]{jmm:red}\set
\begin{equation}\label{autwNgl}
  1 \to {\Ze}({\W}(X)) \backslash \Aut_{\Z_2{\W}(X)}({  L}(X)) \to 
  {\W}(X) \backslash \N_{\GL({  L}(X))}({\W}(X)) \to 
  \Out_{\mathrm{tr}}({\W}(X))
\end{equation}\add
where $\Aut_{\Z_2{\W}(X)}({ L}(X)) = \Z_2^{\times}$ if $X$ is simple
by Schur's lemma, and $ \Out_{\mathrm{tr}}({\W}(X))$ is the group of
outer auto\m s of $\W(X)$ that preserve the trace taken in ${L}(X)$.

\setcounter{subsection}{\value{thm}}
\subsection{Totally $\N$-determined \twocg s}
\label{sec:defns}\add
We are now ready to formulate the concept of $\N$-determinism that will
be used in this paper. The extra complications compared to the odd $p$
case \cite[7.1]{jmm:deter} stem from the fact that $H^1({\W};\ch{{\T}})$,
the first cohomology group of the Weyl group with coefficients in the
discrete \mt , is trivial for any connected \pcg\ when $p$ is odd
\cite{kksa:split} but when $p=2$ it may very well be nontrivial
\cite{hms:first}.

\begin{defn}\label{defn:det}
  Let $X$ be a \twocg\ with \mtnp\ $(\N,\N_0) \xrightarrow{(j,j_0)}
  (X,X_0)$ (\ref{defn:mtnp}).  
  \begin{enumerate}
  \item \label{det1} $X$ {\em has $\N$-determined
      ($\pi_*(\N)$-determined) auto\m s\/} if
    \begin{equation*}
    \func{\mathrm{AM}}{\Aut(X)}{{\W}_0\backslash\Aut(\N)} \quad
      \left(\func{\pi_*
      \circ
      \mathrm{AM}}{\Aut(X)}{{\W}_0\backslash\Aut({\W},\ch{{\T}},e)}
                                            \right)
    \end{equation*}
     is injective.
  \item \label{det2} $X$ is {\em $\N$-determined} if, for any other
    \twocg\ $X'$ with \mtnp\ $(\N,\N_0) \xrightarrow{(j',j'_0)}
    (X',X'_0)$, there exist an iso\m\ \func{f}{X}{X'} and an auto\m\ 
    $\alpha \in H^1({\W};\ch{{\T}}) \subset {\W}_0 \backslash \Aut(\N)$ such that
    the diagram 
        \set\begin{equation}\label{dia:BaBf} \xymatrix{
        {B\N}\ar[d]_{Bj}\ar[r]^{B\alpha}_{\cong} & {B\N}\ar[d]^{Bj'} \\
        {BX}\ar[r]^{\cong}_{Bf} & {BX'} }
    \end{equation}\add
    commutes up to based homotopy.
  \end{enumerate}
\end{defn}

Furthermore, we say that 
\begin{itemize}
\item  $X$ is {\em totally $\N$-determined} if $X$ has $\N$-determined
   auto\m s and is $\N$-determined,
\item $X$ is {\em uniquely\/} $\N$-determined if $X$ is totally
  $N$-determined and $X$ has $\pi_*(N)$-determined 
auto\m s.
\end{itemize}
If $X$ is a totally $\N$-determined \twocg\ then
\begin{equation*}
\text{$X$ is uniquely $\N$-determined} \iff 
                H^1({\W};\ch{{\T}}) \cap \Aut(X)=0 
\end{equation*}
as we see from the \ses\ (\ref{autNW0}).

\begin{lemma}\label{defcons}
  Let $X$ be a \twocg\ as in Definition~\ref{defn:det}. 
  \begin{enumerate}
  \item \label{defcons1} $X$ has $\N$-determined auto\m s if and only
    if for any $\alpha\in {\W}_0 \backslash \Aut(\N)$ with
    $\pi_*(B\alpha)=\mathrm{Id}$ and for any \twocg\ $X'$ as in
    \ref{defn:det}.(\ref{det2}) there is at most one iso\m\ 
    \func{f}{X}{X'} such that diagram (\ref{dia:BaBf}) commutes up to
    based homotopy.
  \item \label{defcons2} $X$ has $\pi_*(\N)$-determined auto\m s if
    and only if
    for any given $X'$ as in \ref{defn:det}.(\ref{det2}), diagram
    (\ref{dia:BaBf}) has at most one solution $(f,\alpha)$ with
    $\pi_*(B\alpha)=\mathrm{Id}$ .
  \end{enumerate}
\end{lemma}
%%%%%%%%%%%%%old version%%%%%%%%%%%%%%%%%%%%%%%%%
%\begin{lemma}\label{defcons}
%  Let $X$ be a \twocg\ with \mtnp\ $(\N,\N_0) \xrightarrow{(j,j_0)}
%  (X,X_0)$.  
%  \begin{enumerate}
%  \item \label{defcons1} $X$ has $\N$-determined auto\m s if and only
%    if for fixed $\alpha\in {\W}_0 \backslash \Aut(\N)$ with
%    $\pi_*(B\alpha)=\mathrm{Id}$ there is at most one iso\m\ 
%    \func{f}{X}{X'} such that diagram (\ref{dia:BaBf}) commutes up to
%    based homotopy (for any given $X'$ as in
%    \ref{defn:det}.(\ref{det2})).
%  \item \label{defcons2} $X$ has $\pi_*(\N)$-determined auto\m s if
%    and only if diagram (\ref{dia:BaBf}) has at most one solution
%    $(f,\alpha)$ with $\pi_*(B\alpha)=\mathrm{Id}$ (for any given $X'$
%    as in \ref{defn:det}.(\ref{det2})).
%  \end{enumerate}
%\end{lemma}
%%%%%%%%%%%%%%%%%%%%%%%%%%%%%%%%%%%%%%%%%%%%
\begin{proof}
  (1) Suppose that $X$ has $\N$-determined auto\m s. Let
  $(f_1,\alpha)$ and $(f_2,\alpha)$ be two solutions to diagram
  (\ref{dia:BaBf}) with the same $\alpha \in
  H^1({\W};\ch{{\T}})\subset {\W}(X_0)\backslash \Aut(\N)$. Then
  $\AM(f_2^{-1}f_1)$ is the identity of $\W(X_0) \backslash \Aut(\N)$
  and since \func{\AM}{\Aut(X)}{{\W}(X_0)\backslash \Aut(\N)} is
  injective, $f_1=f_2$. For the converse, take $B\alpha$ to be the
  identity of $BN$ and take $X'$ to be $X$. Then the assumption is
  precisely that $\mathrm{AM}$ is injective.
 % Conversely, if there exist two solutions to
%  diagram (\ref{dia:BaBf}), $(f_1,\alpha)$ and $(f_2,\alpha)$, with
%  the same $\alpha \in H^1({\W};\ch{{\T}})\subset {\W}(X_0)\backslash
%  \Aut(\N)$ but with $f_1 \neq f_2$, then
%  $\AM(f_2^{-1}f_1)=\alpha^{-1}\alpha$ is the identity element of
%  $\W(X_0) \backslash \Aut(\N)$, but $f_2^{-1}f_1\in\Aut(X)$ is not
%  the identity so $\AM$ is not injective and $X$ does not have
%  $\N$-determined auto\m s.

  \noindent (2)
  Suppose that $X$ has $\pi_*(\N)$-determined auto\m s and let
  $(f_1,\alpha_1)$ and $(f_2,\alpha_2)$ be two solutions to diagram
  (\ref{dia:BaBf}). Then $\AM(f_2^{-1}f_1)=\alpha_2^{-1}\alpha_1 \in
  \AM(\Aut(X)) \cap H^1({\W}(X);\ch{{\T}}(X))$ and this intersection
  is trivial by hypothesis.  Thus $\AM(f_2^{-1}f_1)=1$ and $f_2=f_1$
  as $\AM$ is injective.  If $X$ does not have $\pi_*(\N)$-determined
  auto\m s, then $\AM(f)$ lies in $H^1({\W}(X);\ch{{\T}}(X)) \subset
  {\W}(X_0)\backslash\Aut({\N})$ for some nontrivial $f\in\Aut(X)$ so
  that $(f,\AM(f))$ and $(1,0)$ are two solutions to diagram
  (\ref{dia:BaBf}) with
  $X'=X$ and $j'=j$.
\end{proof}

\begin{exmp}
  For (the \twocg\ associated to) a connected Lie group $G$, the
  cohomology group $H^1({\W}(G);\ch{{\T}}(G))$ is always an \lmntwo\ 
  \cite[1.1]{matthey:normalizers} (\ref{kerneltheta}). For instance,
  this first cohomology group has order two for $G=\PGL(4,\C)$
  \cite[Appendix B]{matthey:second}. Let $\alpha$ be an iso\m\ of
  $\N(\PGL(4,\C))$ representing the nontrivial element of
  $H^1(W;\ch{T})$.  The unique solution
  (\ref{defcons}.(\ref{defcons2})) to diagram (\ref{dia:BaBf}) is
\begin{equation*}
  \xymatrix@C=35pt{
    {\N(\PGL(4,\C))} \ar[d]_j \ar[r]^{\alpha} &  {\N(\PGL(4,\C))}
    \ar[d]^{j'}\\
    {\PGL(4,\C)} \ar@{=}[r] &  {\PGL(4,\C)} }
\end{equation*}
when we use the \m s $j$, induced by an inclusion of Lie groups, and
$j'=j\alpha$ for \mtn s. This example demonstrates that, in contrast
with the $p$ odd case \cite[7.1]{jmm:deter} \cite{jmm:ndet,
  agmv:efamily}, diagram (\ref{dia:BaBf}) can not always be solved with
$\alpha$ the identity.
%$\PGL(4,\C)$ is a uniquely but not strongly
%$\N$-determined \twocg .
\end{exmp}

\begin{lemma}\label{underT}
  Let $X$ be a connected \twocg\ with \mtn\ \func{j}{\N}{X} and \mt\ ${\T}
  \hookrightarrow \N \xrightarrow{j} X$. 
  \begin{enumerate}
  \item  \label{underT1} $X$
  is $\N$-determined if and only if for any other
  connected \twocg\ $X'$ with \mtn\ \func{j'}{\N}{X'} there exists a
   \m\ \func{f}{X}{X'} such that
  \set\begin{equation}\label{dia:underT}
    \xymatrix{
              & {BT}\ar[dl]_{Bj \vert {BT}} \ar[dr]^{Bj' \vert {BT}} \\
              BX \ar[rr]_{Bf} && {BX'} }
  \end{equation}\add
  commutes up to conjugacy.
\item \label{underT2} $X$
  is uniquely $\N$-determined if and only if for any other
  connected \twocg\ $X'$ with \mtn\ \func{Bj'}{BN}{BX'} there exists a
  unique \m\ \func{Bf}{BX}{BX'} such that (\ref{dia:underT}) 
  %\begin{equation*}
%    \xymatrix{
%              & {\T}\ar[dl]_{j \vert {\T}} \ar[dr]^{j' \vert {\T}} \\
%              X \ar[rr]_f && {X'} }
%  \end{equation*}
  commutes up to homotopy.
  \end{enumerate}
\end{lemma}
\begin{proof}  %%Fri Aug 29 14:34:58 CEST 2003
  
  (1) Suppose that the connected \twocg\ $X$ is $\N$-determined and
  let $X'$ be another connected \twocg\ with the same \mtn . Then $X$
  and $X'$ have the same \mtnp , $(\N,\N)$, and therefore diagram
  (\ref{dia:BaBf}) admits a solution $(f,\alpha)$ such that
  $\pi_*(B\alpha)$ is the identity. In particular, $\pi_2(B\alpha)$ is
  the identity of $\pi_2(B{\T})$ which means that $B\alpha$ restricts
  to the identity on the identity component $B\T$ of $B\N$.
  
  Conversely, under the existence assumption of point (1), we shall
  show that $X$ is $\N$-determined. Let $X'$ be another \twocg\ with
  the same \mtnp\ as $X$. Since the \mtnp\ informs about component
  groups (see the remarks just below \ref{defn:mtnp}), $X'$ is connected.
  By assumption, there exists a \m , in fact \cite[5.6]{dw:split}
  \cite[3.11]{jmm:normax} an iso\m , \func{Bf}{BX}{BX'} under $B{\T}$.
  Let \func{B\alpha}{B\N}{B\N}, $B\alpha \in \Out(\N)={\W} \backslash
  \Aut(\N)$, be the restriction of $Bf$ to $B\N$ \cite[\S3]{jmm:ndet}
  so that
  \begin{equation*}
    \xymatrix{
      B\N \ar[d]_{Bj} \ar[r]^-{B\alpha}_-{\cong} & B\N \ar[d]^{Bj'} \\
      BX \ar[r]^-{\cong}_-{Bf} & BX'}
  \end{equation*}
  commutes up to based homotopy as in the definition of the
  Adams--Mahmud homo\m\ (\S\ref{sec:AM}). The further restriction of $B\alpha$
  to the \mt\ $B{\T}$ agrees with the restriction of $Bf$ to $B{\T}$,
  the identity of $B{\T}$, up to the action of a Weyl group element $w
  \in \W$ because ${\W}\backslash [B{\T},B{\T}]=[B{\T},BX']$
  \cite[3.4]{jmm:rip} \cite[3.4]{dw:split}. Since $\pi_1(B\N)={\W}$ is
  faithfully represented in $\pi_2(B{\T})$ for the connected \twocg\ 
  $X'$ \cite[9.7]{dw:fixpt}, it follows that $\pi_1(B\alpha)$ is
  conjugation by $w$. Thus $B\alpha$ belongs (\ref{outN}) to the
  subgroup $H^1({\W};\ch{{\T}})$ of $\Out(\N)$ so that $(f,\alpha)$
  is a legitimate solution to diagram (\ref{dia:BaBf}).
  
  (2) Suppose that $X$ is uniquely $\N$-determined and let $X'$ be
  another connected \twocg\ with the same \mtn\ as $X$.  From point
  (1) we already know that there exists at least one iso\m\ 
  \func{f}{X}{X'} under ${\T}$.  Suppose \func{f_1,f_2}{X}{X'} are two
  such iso\m s under ${\T}$. Then $f_2^{-1}f_1$ is an auto\m\ of $X$
  under ${\T}$, i.e.\ $\pi_*(\B\!\AM(f_2^{-1}f_1)) \in \W \backslash
  \Aut(\W,\T)$ is the identity. As $\pi_* \circ \AM$ is injective,
  $f_2^{-1}f_1$ is the identity of $X$, so $f_1=f_2$.

  Conversely, under the existence and uniqueness assumption of point
  (2), we shall show that $X$ is uniquely $\N$-determined. By point
  (1), $X$ is $\N$-determined, so we only need to show that $\pi_*
  \circ \AM$ is injective. Let \func{f}{X}{X} be an auto\m\ of $X$
  such that $\pi_*(\B\!\AM(f)) \in \W \backslash \Aut(\W,\T)$ is the
  identity. Since $\B\!\AM(f)$ is determined only up to conjugacy, we
  may assume that $\pi_*(\B\!\AM(f))$ {\em is\/} the identity of
  $\pi_*(B\N)$. In particular, $\pi_2(\B\!\AM(f))$ is the identity of
  $\pi_2(\B\T)$ meaning that $f$ is an auto\m\ under ${\T}$. The
  identity of $X$ is also an auto\m\ under ${\T}$, so $f$ is the
  identity auto\m\ of $X$ by the uniqueness hypothesis. This shows
  that $\pi_* \circ \AM$ is injective.
\end{proof}

\begin{lemma}\label{lemma:autX}
  Let $X$ be a connected \twocg\ with \mtn\ $N \to X$.
  \begin{enumerate}
  \item $\Out(\N) = H^1({\W}(X);\ch{{\T}}(X)) \cdot \AM(\Aut(X))$ if
    $X$ is $N$-determined.
  \item $\Out(\N) \cong H^1({\W}(X);\ch{{\T}}(X)) \rtimes \Aut(X)$ and
    $\Aut(X) \cong {\W}(X) \backslash \Aut({\W}(X),\ch{{\T}}(X),e(X))$
    if $X$ is uniquely $N$-determined. The group
    $\Aut({\W}(X),\ch{{\T}}(X),e(X))$ is a subgroup of
    $N_{\GL(L(X))}(W(X))$ (\ref{autwNgl}) and isomorphic to this group
    if $e(X)=0$.
  \end{enumerate}
\end{lemma}
\begin{proof}
  (1) For any $\beta \in \Out(N)$ there exist an auto\m\ $\alpha \in
  H^1(W(X);\ch{T}(X)) \subset \Out(N)$ and an auto\m\ $f \in \Aut(X)$
  such that the diagram
  \begin{equation*}
    \xymatrix{
       BN \ar[r]^{B\alpha} \ar[d]_{Bj} &
       BN \ar[d]^{Bj \circ B\beta} \\
       BX \ar[r]_{Bf} & BX}
  \end{equation*}
  commutes up to homotopy (\ref{defn:det}.(\ref{det2})). Thus
  $\AM(f)=\beta\alpha$ in $\Out(N)$ (\S\ref{sec:AM}).

  \noindent
  (2) If the connected \twocg\ $X$ is uniquely $N$-determined, then
  there is commutative diagram
  \begin{equation*}
    \xymatrix{
    0 \ar[r] &
  H^1({\W}(X);\ch{{\T}}(X)) \ar[r] &
  {\Out(\N)} \ar[r]^-{\pi_*} &
  {\W}(X) \backslash \Aut({\W}(X),\ch{{\T}}(X),e(X)) \ar[r] &
  1 \\
  && {\Aut(X)} \ar@{^(->}[u]^{\AM} \ar@{^(->}[ur] }
  \end{equation*}
  where the top row is the \ses\ (\ref{outN}). The composite homo\m\ 
  $\pi_* \circ \AM$ is injectice by assumption
  (\ref{defn:det}.(\ref{det1})). It is surjective since $\Out(N)$ is
  generated by $H^1(W(X);\ch{T}(X))$ by item (1) of this lemma. Thus
  $\pi_* \circ \AM$ is an iso\m\ and $\AM$ is a splitting of the \ses\ 
  (\ref{outN}).
\end{proof}
  
As evidence of the conjecture that all connected \twocg s are uniquely
$N$-determined we note that all compact connected Lie groups have
$\pi_*(\N)$-determined auto\m s \cite[2.5]{jmo:selfho} and satisfy the
above two formulas for auto\m\ groups \cite[3.10]{hammerli:thesis}.

%\begin{proof}
%  The \m\ \func{f}{X}{X'} in the above commutative diagram is in fact
%  an iso\m\ \cite[5.6]{dw:split} \cite[3.11]{jmm:normax}. The
%  assumption of the lemma that $f$ be a \m\ under $T$ means (use
%  $W\backslash [BT,BT]=[BT,BX]$ \cite[3.4]{jmm:rip}
%  \cite[3.4]{dw:split}) that $f$ admits a restriction $N(f)$ to $N$
%  which is the identity on $T$, i.e.\ such that
%  \begin{equation*}
%    \xymatrix{
%    BT \ar@{=}[d]\ar[r] & BN \ar[d]_{BN(f)} \ar[r]^{Bj} &
%    BX\ar[d]^{Bf} \\
%    BT \ar[r] & BN \ar[r]_{Bj'} & BX' }
%  \end{equation*}
%  is homotopy commutative. But then also \func{\pi_0N(f)}{W}{W} is the
%  identity map for  $W$ is faithfully represented as a group of
%  operators on $T$ \cite[9.7]{dw:fixpt}. Thus $\pi_*(BN(f))$ is the
%  identity auto\m\ of $\pi_*(BN)$.  
  
%  Assume that the iso\m\ $f$ exists and is uniquely determined. In
%  particular, the identity of $X$ is the only auto\m\ under $T$. That
%  $f \in \Aut(X)$ is a map under $T$ means precisely that $\AM(f) \in
%  H^1(W;\ch{T})$. Thus $X$ is uniquely $N$-determined. Suppose,
%  conversely, that $X$ has this property and let \func{f_0,f_1}{X}{X'}
%  be two iso\m s under $T$. Then $f_1^{-1}f_0 \in \Aut(X)$ is an
%  iso\m\ under $T$ so equals the identity.
%\end{proof}

With a view to the situation for possibly nonconnected \twocg s, let
$\Aut(N,N_0)$ denote the subgroup of $\Aut(N)$ consissting of all
auto\m\ $\phi\in\Aut(N)$ such that $\pi_0(\phi)$ takes  $\pi_0(N_0)$
to  itself inducing an iso\m\ 
\begin{equation*}
  \xymatrix{
    {\pi_0(N_0)} \ar[r] \ar[d]^{\cong} & 
    {\pi_0(N)} \ar[d]^{\pi_0(\phi)}_{\cong} \ar[r] & {\pi}\ar[d]^{\cong} \\
     {\pi_0(N_0)} \ar[r] &  {\pi_0(N)}  \ar[r] & {\pi} }
\end{equation*}
of \ses s. Since $H^1(W;\ch{T})$ is contained in $\Aut(N,N_0)$, there
are \ses s similar to (\ref{eq:H1TW}) and (\ref{autNW0}) except that
$\Aut(W,\ch{T},e)$ has been replaced by its subgroup
$\Aut(\ch{T},W,W_0,e)$ of all $(\alpha,\theta) \in \Aut(W,\ch{T},e)$
for which $\alpha(W_0)=W_0$. (If $N=N_0$, then $\Aut(N)=\Aut(N,N_0)$.)
Observe that the Adams--Mahmud homo\m\ for a nonconnected \twocg\ 
actually takes values in the subgroup $W(X_0)\backslash \Aut(N,N_0)$
of $W(X_0)\backslash \Aut(N)$.

\begin{lemma}
  \label{lemma:autXnoncon} Let $X$ be a \twocg\ with \mtnp\ $(N,N_0)
  \to (X,X_0)$.
  \begin{enumerate}
  \item $W(X_0)\backslash \Aut(N,N_0) = H^1(W(X);\ch{T}(X)) \cdot
    \AM(\Aut(X))$ if $X$ is $N$-determined.
  \item $W(X_0)\backslash \Aut(N,N_0) \cong H^1(W(X);\ch{T}(X))
      \rtimes_{H^1(\pi_0(X);\ch{Z}(X_0))} \Aut(X)$ if $X$ is totally
    $N$-determined.
  \end{enumerate}
\end{lemma}
\begin{proof}
  The first item is proved like the first item in
  \ref{lemma:autX}. The claim of the second item is that 
  \begin{equation*}
    \xymatrix{
          H^1(\pi_0(X);\ch{Z}(X_0)) \ar[d] \ar@{^(->}[r] &
          {\Aut(X)} \ar@{^(->}[d]^{\AM} \\
          H^1(W(X);\ch{T}(X)) \ar@{^(->}[r] &
          {W(X_0) \backslash \Aut(N,N_0)} }
  \end{equation*}
  is a push-out diagram. This is proved in \ref{idcomp} (allowing
  ourselves to refer ahead!).
\end{proof}

\begin{rmk}\label{unbasedOut}
  When the \twocg\ $X$ has $\N$-determined auto\m s, also the unbased
  Adams--Mahmud homo\m\ 
  \begin{equation*}
  \Out(X)= \pi_0(X) \backslash \Aut(X) \to
  \Out(\N) = \pi_0(\N) \backslash \Aut(\N) = \pi_0(X) \backslash
  W(X_0) \backslash \Aut(N)    
  \end{equation*}
  is injective.
\end{rmk}

%%%%%%%%%%%%%%%%%%%%%%%%%%%%%%
\setcounter{subsection}{\value{thm}}

\subsection{Regular \twocg s}
\label{subsec:regular}
For a {\em connected\/} \twocg\ $X$ with \mt\ ${\T} \to X$ and Weyl group
${\W}$, let
\set\begin{equation}
  \label{eq:theta}
   \func{\theta=\theta(X)}{\Hom({\W},\ch{{\T}}^{{\W}}) =
  H^1({\W};\ch{{\T}}^{{\W}})}{H^1({\W};\ch{{\T}})} 
\end{equation}\add
be the homo\m\ induced by the inclusion $\ch{{\T}}^{{\W}} \hookrightarrow
\ch{{\T}}$. Following \cite[5.3]{hms:first} we say that $X$ is {\em
  regular\/} if (\ref{eq:theta}) is surjective. See 
\cite{matthey:normalizers} for a thorough investigation of $\theta$.

\begin{lemma} \label{kerneltheta}
  \cite{matthey:normalizers} Let $X$ be the connected
  \twocg\ associated to a connected Lie group. Assume that $X$
  contains no direct factors isomorphic to an odd orthogonal group
  $\SO(2n+1)$, $n \geq 1$.  Consider the homo\m\ $\theta=\theta(X)$
  (\ref{eq:theta}) associated to $X$.
  %\func{\theta}{\Hom({\W},\ch{{\T}}^{\W})}{H^1({\W};\ch{{\T}})} for $X$.
  \begin{enumerate}
  \item $\Hom({\W},\ch{{\T}}^{\W})$ and $H^1({\W};\ch{{\T}})$ are $\F_2$-vector
    spaces, and the kernel of $\theta$, consisting of those homo\m s
    ${\W} \to \ch{{\T}}^{{\W}}$ that are principal crossed homo\m s ${\W} \to
    \ch{{\T}}$, is an $\F_2$-vector space of dimension equal to the
    number of direct factors of $PX$ isomorphic to $\SO(2n+1)$, $n
    \geq 1$. \label{kerneltheta.1}
  \item Suppose that the projective group $PX$ contains no direct
    factors isomorphic to an odd orthogonal group $\SO(2n+1)$, $n \geq
    1$, $\mathrm{PSU}(4)$, $\mathrm{PSp}(3)$, $\mathrm{PSp}(4)$, or
    $\mathrm{PS0}(8)$.  Then $X$ is regular.
    \label{kerneltheta.2}
  \end{enumerate}
\end{lemma}
\begin{proof} \eqref{kerneltheta.1}
  $\Hom({\W},\ch{{\T}})$ and its subgroup $\Hom({\W},\ch{{\T}}^{\W})$
  are \lmntwo s since the abelianization ${\W}_{\mathrm{ab}}$ of
  ${\W}$ is an \lmntwo\ of finite rank. The cohomology group
  $H^1({\W};\ch{{\T}})$ is isomorphic to $H^2({\W};  L \otimes
  \Z_2)$ where $  L$ is the fundamental group of the Lie group
  maximal torus of the Lie group underlying the \twocg\ $X$.
  Homological algebra shows that $H^2({\W};  L \otimes \Z_2)
  \cong H^2({\W};  L) \otimes \Z_2$ where $H^2({\W};  L)$
  is an \lmntwo\ by \cite[1.1]{matthey:normalizers}.  The injection
  $\ch{{\T}}^{\W} \to \ch{{\T}}$ of ${\W}$-modules gives a coefficient
  group long exact sequence
  \begin{equation*}
    0 \to \left(\ch{{\T}}/\ch{{\T}}^{\W}\right)^{\W} \to 
                   \Hom({\W},\ch{{\T}}^{\W}) 
    \xrightarrow{\theta} H^1({\W};\ch{{\T}}) \to 
    H^1({\W};\ch{{\T}}/\ch{{\T}}^{\W}) \to 
              H^2({\W};\ch{{\T}}^{\W}) \to \cdots
  \end{equation*}
  in cohomology. Thus the kernel of $\theta$ is isomorphic to $
  \left(\ch{{\T}}/\ch{{\T}}^{\W}\right)^{\W}$ in general. If $X$ is
  without direct factors isomorphic to $\SO(2n+1)$, then
  $\ch{{\T}}^{\W}$ is the center of $X$, $\ch{{\T}}/\ch{{\T}}^{\W}$ is
  the \mt\ of the adjoint \twocg\ $PX$, and
  $\left(\ch{{\T}}/\ch{{\T}}^{\W}\right)^{\W}$ is isomorphic to
  $(\Z/2)^s$ where $s$ is the number of direct factors isomorphic to
  an odd special orthogonal group in the adjoint \twocg\ $PX$
  \cite[1.6]{matthey:normalizers} \cite[4.6, 4.7]{mn:center}. (See
  \ref{rmk:T/TW} for the general case.)

  \noindent \eqref{kerneltheta.2} The discrete \mt\ of $PX=X/{\Ze}(X)$ is
  $\ch{{\T}}(PX)=\ch{{\T}}/\ch{{\T}}^{\W}$ for
  $\ch{{\Ze}}(X)=\ch{{\T}}^{\W}$ as $X$ contains no direct factors
  isomorphic to an odd orthogonal group.  The projective group $PX=
  \prod G_i$ splits as a product of simple and centerfree compact Lie
  groups $G_i$ all of which satisfy $\ch{{\T}}^{{\W}}(G_i)=0$ since
  they are not odd orthogonal groups.  Therefore
  $H^1({\W};\ch{{\T}}/\ch{{\T}}^{\W})= H^1(\prod {\W}(G_i); \prod
  \ch{{\T}}(G_i)) = \prod H^1({\W}(G_i);\ch{{\T}}(G_i))$ and these
  cohomology groups are trivial except in the excluded cases
  \cite{hms:first}. By the above exact sequence, $\theta$ is
  surjective.
\end{proof}

For a compact connected Lie group $X$, let $s(X)$ denote the number of
direct factors of $X$ isomorphic to $\SO(2n+1)$ with $n \geq 1$. (Keep
the low degree identifications (\ref{eq:lowdegree}) in
mind.)
\begin{lemma}\label{lemma:kerntheta}
  Let $X$ be a compact connected Lie group and $PX$ its adjoint form.
  The kernel of $\theta(X) \colon H^1(W;\Check{Z})(X) \to
  H^1(W;\Check{T})(X)$ is an $\F_2$-vector space of dimension
  $s(PX)-s(X)$.
\end{lemma}
\begin{proof}
  In the exact sequence
  \begin{equation*}
    0 \to \Check{Z} \to  \Check{T}^W \to \left( \Check{T}/ \Check{Z}
    \right)^W \to H^1(W;\Check{Z}) \to H^1(W;\Check{T})    
  \end{equation*}
  induced from the inclusion $\Check{Z} \to \Check{T}$ of $W$-modules,
  the fixed point groups $\Check{T}^W = \Check{Z}(X) \times 2^{s(X)}$
  and $\left(\Check{T}/\Check{Z}\right)^W = \Check{Z}(X/Z) \times
  2^{s(X/Z)} =2^{s(X/Z)}$ \cite[1.6]{matthey:normalizers}.
\end{proof}

%% intended for machine computation of theta
%\begin{rmk} 
%  The commutative diagram
%  \begin{equation*}
%    \xymatrix{
%      & H^1(W;L) \ar[d] \\
%      H^1(W;t^W) \ar[d] \ar[r] & H^1(W;t) \ar[d] \\
%      H^1(W;\ch{T}^W) \ar[r]^{\theta} & H^1(W;\ch{T}) } 
%  \end{equation*}
%  may be useful for computations of $\theta$ as we often have
%  $t^W=\ch{T}^W$. The column is exact.
%\end{rmk}

\begin{lemma}
  \label{regquot} Let $X$ be a connected \twocg\ with \mt\ ${\T} \to X$
  and Weyl group ${\W}$, and let ${\Ze} \to {\T} \to X$ be a central mono\m . If
  $X$ is regular and $H^2({\W};\ch{{\Ze}}) \to H^2({\W};\ch{{\T}})$ is injective,
  then the quotient \twocg\ $X/{\Ze}$ is regular.
\end{lemma}
\begin{proof}
  Since the hypothesis implies that $H^1({\W};\ch{{\T}}) \to
  H^1({\W};\ch{{\T}}/\ch{{\Ze}})$ is surjective, the claim follows from the
  commutative square
  \begin{equation*}
    \xymatrix@C=40pt{
      {\Hom({\W},\ch{{\T}}^{\W})} \ar@{->>}[d]_{\theta(X)} \ar[r] &
      {\Hom({\W},(\ch{{\T}}/\ch{{\Ze}})^{\W})} \ar[d]^{\theta(X/{\Ze})} \\
      H^1({\W};\ch{{\T}}) \ar@{->>}[r] & H^1({\W};\ch{{\T}}/\ch{{\Ze}}) }
  \end{equation*}
  induced by the projection $\ch{{\T}} \to \ch{{\T}}/\ch{{\Ze}}$ of
  ${\W}$-modules \cite[4.6]{mn:center}.
\end{proof}

\begin{exmp}\label{exmp:glmCreg}
  (1) $\GL(m,\C)$ is regular for all $m \geq 1$. For $m=1$, this is
  obvious. For $m>2$, the restriction homo\m\ ($\ch{S}={\Z}/2^{\infty}$)
  \begin{equation*}
    \Hom(\Sigma_m,\ch{S}) = H^1(\Sigma_m;\ch{S})
    \xrightarrow{\mathrm{res}=\theta(\GL(m,\C))}
    H^1(\Sigma_m;\ch{S}^m) \stackrel{\text{Shapiro}}{=}
    H^1(\Sigma_{m-1};\ch{S}) = \Hom(\Sigma_{m-1},\ch{S})
  \end{equation*}
   is bijective and for $m=2$ it is surjective. It now follows
   \cite[5.7]{hms:first} that all products $\prod \GL(m_j,\C)$ are
   regular. 

   \noindent
   (2) $\PGL(m,\C)$, $2 \leq m$, is regular for $m \neq 4$ since
   (\ref{regquot})
   \begin{equation*}
     \Hom(H_2(\Sigma_m),\ch{S}) = H^2(\Sigma_m;\ch{S})
     \xrightarrow{\mathrm{res}}
     H^2(\Sigma_m;\ch{S}^m) \stackrel{\text{Shapiro}}{=}
     H^2(\Sigma_{m-1};\ch{S}) =  \Hom(H_2(\Sigma_{m-1}),\ch{S})
   \end{equation*}
   is then an iso\m . The \twocg\ $\PGL(4,\C)$ is not
   regular as $H^1({\W};\ch{{\T}})=\Z/2$ is nontrivial while the
   discrete center $\ch{{\T}}^{{\W}}$ is trivial.
\end{exmp}

\begin{rmk}\label{rmk:T/TW}
  If $X=\SO(2n+1)$, $n \geq 1$, then $\ch{{\T}}^{\W}=\Z/2$,
  ${\W}_{\mathrm{ab}}$ is $\Z/2$ for $n=1$ and $(\Z/2)^2$ for $n \geq 2$,
  \func{\theta}{\Hom({\W},\ch{{\T}}^{\W})}{H^1({\W};\ch{{\T}})} is surjective
  \cite[5.5]{hms:first}, and $H^1({\W};\ch{{\T}})$ is trivial for $n=1$,
  $\Z/2$ for $n=2$, and $(\Z/2)^2$ for $n>2$ \cite[Main Theorem,
  5.5]{hms:first}. Thus the kernel of $\theta$ is
  \begin{equation*}
    \left( \ch{{\T}}/\ch{{\T}}^{\W} \right)^{\W} = 
    \begin{cases}
      \Z/2 & n=1,2 \\ 0 & n>2
    \end{cases}
  \end{equation*}
  In general, write the connected Lie group $X=X_1 \times X_2$ where
  $X_1$ is the product of all direct factors of $X$ isomorphic to
  $\SO(2n+1)$ for some $n \geq 1$ and $X_2$ is without such direct
  factors. Then
  \begin{equation*}
     \left( \ch{{\T}}/\ch{{\T}}^{\W} \right)^{\W} =
     \left( \ch{{\T}}_1/\ch{{\T}}_1^{{\W}_1} \right)^{{\W}_1} \times  
     \left( \ch{{\T}}_2/\ch{{\T}}_2^{{\W}_2} \right)^{{\W}_2} =
     \left( \Z/2 \right)^{s_{\scriptscriptstyle{\leq 2}}(X)} \times
     \left( \Z/2 \right)^{s(PX_2)}
  \end{equation*}
  where $s_{\scriptstyle{\leq 2}}(X)$ is the number of direct factors of $X$
  isomorphic to $\SO(3)$ or $\SO(5)$ and $s(PX_2)$ is the number of
  direct factors of $PX_2$ isomorphic to $\SO(2n+1)$ for some $n \geq
  1$.  
\end{rmk}

\setcounter{subsection}{\value{thm}}
\subsection{LHS \twocg s}\add\label{subsec:LHS}
Let $\N_0 \to \N$ be maximal rank normal mono\m\ between two \etwocts ,
i.e.\ a commutative diagram with rows and columns that are \ses s of
loop spaces \cite[3.2]{dw:fixpt}
\begin{equation*}
  \xymatrix{
 {{\T}} \ar@{=}[r] \ar[d] & {{\T}} \ar[r] \ar[d] & {\{1\}} \ar[d] \\
 {\N_0} \ar[d] \ar[r] & {\N} \ar[d] \ar[r] & {{\W}/{\W}_0} \ar@{=}[d] \\
 {\W}_0 \ar[r] & {\W} \ar[r] & {\W}/{\W}_0 }
\end{equation*}
where ${\T}$ is a \twoct\ and ${\W}_0=\pi_0(\N_0)$ a normal subgroup of the
finite group ${\W}=\pi_0(\N)$. The $5$-term exact sequence
\begin{equation*}
  0 \to H^1({\W}/{\W}_0;\ch{{\T}}^{{\W}_0}) 
  \xrightarrow{\mathrm{inf}} H^1({\W};\ch{{\T}}) 
  \xrightarrow{\mathrm{res}}{H^1({\W}_0;\ch{{\T}})^{{\W}/{\W}_0}} 
  \xrightarrow{d_2} H^2({\W}/{\W}_0;\ch{{\T}}^{{\W}_0})
  \xrightarrow{\mathrm{inf}} H^2({\W};\ch{{\T}})
\end{equation*}
is part of the Lyndon--Hochschild--Serre spectral sequence
\cite{hochserre} converging to $H^*({\W};\ch{{\T}})$.

\begin{defn}\label{LHSinit}
  A \twocg\ with \mtnp\ $(\N,\N_0)$ is LHS if the restriction homo\m\
  \func{\mathrm{res}}{H^1({\W};\ch{{\T}})}{H^1({\W}_0;\ch{{\T}})^{{\W}/{\W}_0}} is
    surjective. 
\end{defn}

Thus $X$ is LHS if and only if the initial segment of the
Lyndon--Hochschild--Serre spectral sequence
\begin{equation*}
   0 \to H^1({\W}/{\W}_0;\ch{{\T}}^{{\W}_0})
  \xrightarrow{\mathrm{inf}} H^1({\W};\ch{{\T}})
  \xrightarrow{\mathrm{res}}{H^1({\W}_0;\ch{{\T}})^{{\W}/{\W}_0}} \to 0
\end{equation*}
is exact. If $\ch{{\T}}^{{\W}_0}=0$ or ${\W}={\W}_0 \times
{\W}/{\W}_0$ is a direct product, then $X$ is LHS. Note that the Weyl
group of a compact Lie group $G$ is always the semi-direct product
$W(G)=W(G_0) \rtimes \pi_0(G)$ for the action of the component group
$\pi_0(G)$ on the Weyl group $W(G_0)$ of the identity component
\cite[\S2.5]{hammerli:thesis}. (In fact, It is not so easy to find a
nonconnected compact Lie group $G$ for which the extension $G_0 \to G
\to G/G_0=\pi$ is nonsplit \cite{hammerli:remarks}.)

\begin{lemma}
    \label{lhscrit1}
    Let $W=W(X)$ be the Weyl group of the \twocg\ $X$, $W_0=W(X_0)$
    the Weyl group of the identity component, and $\pi=W/W_0$ the
    component group \cite[3.8]{mn:center} of $X$. If ${\W} = {\W}_0
    \rtimes \pi$ is a semi-direct product and
    \begin{equation*}
    \func{\theta(X_0)^{\pi}}{\Hom({\W}_0,\ch{{\T}}^{{\W}_0})^{\pi}}
    {H^1({\W}_0;\ch{{\T}})^{\pi}}      
    \end{equation*}
    is surjective, then $X$ is LHS.
   % \begin{equation*}
%      H^1(W_0;\ch{T})^{\pi} \subset \im \left(
%        \Hom(W_0,\ch{T}^{W_0}) \xrightarrow{\theta}
%         H^1(W_0;\ch{T}) \right)
%    \end{equation*}
%\func{\theta^{\pi}}{\Hom(W_0,\ch{T}^{W_0})^{\pi}}{H^1(W_0;\ch{T})^{\pi}} is
%%then $X$ is LHS.
    %$H^1(W_0;\ch{T})^{\pi}$ is contained in the image of $\theta$,
%    then $X$ is LHS.
\end{lemma}
\begin{proof}
  Assume that the group $G=H \rtimes Q$ is the semi-direct product for
  a group action $Q \to \Aut(H)$, and let $A$ be a $G$-module. We
  show that the image of the restriction homo\m\ 
  \func{\mathrm{res}}{H^1(G;A)}{H^1(H;A)^Q} contains the image of
  \func{\theta^Q}{\Hom(H,A^H)^Q}{H^1(H;A)^Q}.  Let $\phi \in
  \Hom(H,A^H)^Q$ be a $Q$-equivariant homo\m\ of $H$ into the fixed
  point module $A^H$. Then $\theta(\phi) \in H^1(H;A)^Q$ is the
  cohomology class represented by the crossed homo\m\ 
  \func{\phi}{H}{A^H \subset A}. If we define \func{\overline{\phi}}{H
    \rtimes Q}{A} by $\overline{\phi}(nq)=\phi(n)$, $n \in H$, $q \in
  Q$, then
  \begin{multline*}
      \overline{\phi}(n_1q_1n_2q_2) =
      \overline{\phi}(n_1(q_1n_2q_1^{-1})q_1q_2)=
      \overline{\phi}(n_1(q_1 \cdot n_2)q_1q_2) \stackrel{\mathrm{def}}{=}
      \phi(n_1(q_1 \cdot n_2)) =
      \phi(n_1) + \phi(q_1 \cdot n_2) \\ =
      \phi(n_1) + q_1\phi(n_2)
  \end{multline*}
  and also 
  \begin{equation*}
      \overline{\phi}(n_1q_1) + n_1q_1\overline{\phi}(n_2q_2)
      \stackrel{\mathrm{def}}{=} 
      \phi(n_1) +  n_1q_1\phi(n_2) =
      \phi(n_1) +  q_1\phi(n_2) 
  \end{equation*}
  as $q_1\phi(n_2) \in A^H$. This shows that the crossed homo\m\ 
  $\phi$ defined on $H$ extends to a crossed homo\m\ $\overline{\phi}$
  defined on $G=H\rtimes Q$. (I do not know if the LHS spectral
  sequence differential \func{d_2}{H^1(H;A)^Q}{H^2(Q;A^H)} is always
  trivial for a semi-direct product $H \rtimes Q$ of finite groups.)
\end{proof}

The next example demonstrates that condition \ref{lhscrit1} is not
necessary.

\begin{exmp}\label{exmp:sl2C}

  \noindent
  (1) $X=\PGL(6,\R) = \PSL(6,\R) \rtimes C_2 $ does not satisfy the
  condition of \ref{lhscrit1} for $H¹({\W}_0;\ch{{\T}})=\Z/2$ \cite[Main
  Theorem]{hms:first} while
  $\ch{{\T}}^{{\W}_0}=\ch{{\Ze}}(X_0)=0$. Nevertheless, $X$ is LHS because
  also $H¹({\W};\ch{{\T}})=\Z/2$ (computer computation).

  \noindent
  (2) $X=\PGL(8,\R) = \PSL(8,\R) \rtimes C_2 $ does not satisfy the
  condition of \ref{lhscrit1} for $H¹({\W}_0;\ch{{\T}})=\Z/2 \oplus \Z/2$
  \cite[Main Theorem]{hms:first} while $\ch{{\T}}^{{\W}_0}=\ch{{\Ze}}(X_0)=0$.
  Nevertheless, $X$ is LHS because $H¹({\W};\ch{{\T}})=\Z/2$ and the outer
  auto\m\ group $C_2$ acts nontrivially on $H¹({\W}_0;\ch{{\T}})$ (computer
  computation).

  \noindent
  (3) When $X_0=\SL(2,\C)$, the Weyl group ${\W}_0=\Sigma_2$ has order
  two, the center $\ch{{\Ze}}=\ch{{\T}}^{{\W}_0}$ also has order two, and
  $H¹({\W}_0;\ch{{\T}})=0$ is trivial, so the homo\m\ $\theta(X_0)$
  is trivial as well, of
  course.  Indeed, the nontrivial homo\m\ ${\W}_0 \to \ch{{\Ze}} \subset
  \ch{{\T}}$ is the principal crossed homo\m\ corresponding to the
  element $\diag(i,-i)$ of the \mt . More generally, the direct
  product $X_0^r=\SL(2,\C)^r$ is regular \cite[5.7]{hms:first}, 
  has Weyl group ${\W}_0^r$, center $\ch{{\Ze}}^r$, and 
  \ref{kerneltheta}.(\ref{kerneltheta.1}) identifies the kernel of
  $\theta(X_0)$ enabling us to  conclude that
  \set\begin{equation}\label{H1W0}
    H¹({\W};\ch{{\T}})(X_0^r) = 
    \frac{\Hom({\W}_0^r,\ch{{\Ze}}^r)}{\Hom({\W}_0,\ch{{\Ze}})^r}
  \end{equation}\add
  is an $\F_2$-vector space of dimension $r²-r$ as in
  \cite[5.8]{hms:first}. Let $X=X_0 \rtimes C_2$ be the semi-direct
  product for the nontrivial outer auto\m\ of $X_0$. The component
  group $C_2^r$ of $X^r$ acts trivially on (\ref{H1W0}) and as
  $H¹({\W};\ch{{\T}})(X^r)$ has dimension $2r²-r$ (by induction) and
  $H¹(C_2^r;\ch{{\Ze}}^r)$ dimension $r²$, the direct product $X^r$ is LHS
  for all $r \geq 1$.

   \noindent
  (4) When $X=\slr{4}$, the Weyl group $W=\gen{\sigma,c_1c_2} = \Z/2
  \times \Z/2$ is elementary abelian generated by $\sigma=
    \begin{pmatrix}
      0 & E \\ E & 0
    \end{pmatrix}$ and $c_1c_2=\diag(-1,1,-1,1)$. The center 
    $\ch{Z} = \ch{T}^W = \gen{\diag(-1,-1,-1,-1)}=\Z/2$ has order
    $2$, and the first cohomology group $H¹(W;\ch{T})=0$ is trivial,
    so the homo\m\ \func{\theta}{\Hom(W;\ch{T}^W)}{H¹(W;\ch{T})} is
    also trivial, of course.  Indeed, the principal homo\m\ 
    $\varphi(w) = (w \cdot t) \cdot t^{-1} \colon W \to \ch{T}$, is
    the first coordinate function $W \to \ch{Z}(X)$ when
    $t=\diag(-E,E)$ and the second coordinate function when
    $t=\diag(I,I)$. The outer auto\m , conjugation with
    $D=\diag(-1,1,1,1) \in \GL(4,\R)$, sends $\sigma$ to
    $\sigma^D=\sigma (c_1c_2)$ and $c_1c_2$ to itself.
    
    More generally, when $X^r$ is a product of $r$ copies of
    $\slr{4}$, the Weyl group $W^r$ is a product of $r$ copies
    $W=W(\slr{4})=\Z/2 \times \Z/2$, the center $\ch{Z}(X)=\ch{Z}^r$
    is a product of $r$ copies of $\ch{Z}=\ch{Z}(\slr{4})=\Z/2$ and as
    \func{\theta}{\Hom(W^r,\ch{Z}^r)}{H¹(W;\ch{T})(X^r)}
    is surjective \cite[5.5, 5.7]{hms:first}, the first cohomology
    group
    \begin{equation*}
      H¹(W;\ch{T})(X^r)= \frac{\Hom(W^r,\ch{Z}^r)}{\Hom(W,\ch{Z})^r}
    \end{equation*}
    has dimension $2r²-2r$ over $\F_2$
    (\ref{kerneltheta}). The component group
    $\pi_0(\GL(4,\R)^r)=C_2^r$ acts on this $\F_2$-vector space such
    that the space of fixed vectors has dimension $r^2-r$. By
    induction we see that $H¹(W;\ch{T})(\GL(4,\R)^r)$ is an
    $\F_2$-vector space of dimension $2r²-r$ and clearly
    $H¹(C_2^r;\ch{Z}^r)$ has dimension $r^2$. Thus $\GL(4,\R)^r$ is
    LHS for all $r \geq 1$.

    \noindent (5) 
    The homo\m s $\theta$ is surjective for $\SL(2n,\R)$ for all $n
    \geq 1$ \cite[Main Theorem, 5.4]{hms:first} and
    $H^1(W;\ch{T})(\SL(2n,\R))=0$ for $n=1,2$ and
    $H^1(W;\ch{T})(\SL(2n,\R))=\Z/2$ for $n \geq 3$. Hence
    $\GL(2n,\R)$ is LHS for all $n \geq 1$ by
    \ref{lhscrit1}.
 \end{exmp}

%The $\F_2$-vector space $H^1(W;\ch{T})$ can be computed from the \ses\
%\begin{equation*}
%  0 \to \ch{T}^W \otimes \Z/2 \to H^1(W;t) \to H^1(W;\ch{T}) \to 0
%\end{equation*}
%induced from $0 \to t \to \ch{T} \xrightarrow{\cdot 2} \ch{T} \to 0$;
%note that $\ch{T}^W \otimes \Z/2=0$ if $\ch{T}^W$ is a $2$-discrete
%torus.

I do not know any examples of \twocg s that are not LHS.

The coefficient group \ses\ $0 \to L \to L \otimes \Q \to \ch{T} \to
0$ gives the exact sequence
\begin{equation*}
  0 \to H^0(W;L) \to H^0(W;L\otimes\Q) \to H^0(W;\ch{T}) \to
        H^1(W;L) \to 0
\end{equation*}
form which we see that 
\set\begin{equation}\label{eq:hiwt}
  H^i(W;\ch{T}) = 
  \begin{cases}
    H^0(W;L\otimes\Q) /H^0(W;L) \oplus H^1(W;L) & i=0 \\
    H^{i+1}(W;L) & i>0
  \end{cases}
\end{equation}\add

\setcounter{subsection}{\value{thm}}
\subsection{The center of the \mtn}
\label{sec:centermtn}\add
We need criteria to ensure that the center of the \twocg\ $X$ agrees
with the center of its \mtn . (This is automatic when $p>2$
\cite[3.4]{jmm:deter} but not when $p=2$ \cite[\S7]{dw:center}.)
\begin{prop}\label{prop:centerofN}
  Let $X$ be a \twocg\ with identity component $X_0$. If
  ${\Ze}(X_0)={\Ze}(\N(X_0))$ and $X_0$ has $\N$-determined auto\m s,
  then ${\Ze}(X)={\Ze}(\N(X))$.
\end{prop}
\begin{proof}
  This is proved in \cite[4.12]{jmm:ndet} for \pcg s where $p$ is odd.
  If we replace the assumption that $p$ is odd by the assumption that
  ${\Ze}(X_0)={\Ze}(\N(X_0))$ (which always holds when $p>2$
  \cite[7.1]{dw:center}), then the same proof works also for \twocg s.
\end{proof}

Assume now that $X$ is a {\em connected\/} \twocg .  Then
$\ch{{\Ze}}(\N(X)) = \ch{{\T}(X)}^{{\W}(X)}$ and there is an injection $\ch{{\Ze}}(X)
\hookrightarrow \ch{{\Ze}}(\N(X))$ which is not necessarily an iso\m\ 
\cite[\S7]{dw:center}.

Inspection shows that ${\Ze}(G)={\Ze}\N(G)$ for any {\em simply
  connected\/} compact Lie group $G$; see \cite[1.4]{dw:coxeter} for a
conceptual proof of this fact. In fact, ${\Ze}(G)={\Ze}\N(G)$ for any
connected compact Lie group $G$ containing no direct factors
isomorphic to $\SO(2n+1)$ \cite[1.6]{matthey:normalizers}.

%\begin{rmk}\label{rmk:ZLie}
%  Inspection shows that ${\Ze}(G)={\Ze}\N(G)$ for any {\em simply connected\/}
%  compact Lie group $G$; see \cite[1.4]{dw:coxeter} for a conceptual
%  proof of this fact. In fact, ${\Ze}(G)={\Ze}\N(G)$ for any compact connected
%  Lie group $G$ containing no direct factors isomorphic to $\SO(2n+1)$
%  \cite[1.6]{matthey:normalizers}.
%\end{rmk}

Let ${\Ze} \to \N(X)$ be a central mono\m\ such that also the
composition ${\Ze} \to \N(X) \to X$ is central. Under these
assumptions, the quotient loop spaces $N(X)/Z$ and $X/Z$ can be
defined \cite[2.8]{dw:center}.  The action map \cite[8.6]{dw:fixpt}
$B{\Ze} \times B\N(X) \to B\N(X)$ induces an action $[B\N(X),B{\Ze}]
\times \Out(\N(X)) \to \Out(\N(X))$ of the group $[B\N(X),B{\Ze}]
\cong H^1(\ch{\N}(X);\ch{{\Ze}})$ on the set $ \Out(\N(X))$. Let
$[B\N(X),B{\Ze}]_{(1)}$ denote the isotropy subgroup at $(1) \in
\Out(\N(X))$.
\begin{lemma}
  If ${\Ze}(X)={\Ze}(\N(X))$ and $[B\N(X),B{\Ze}]_{(1)}=0$, 
  then ${\Ze}(X/{\Ze})={\Ze}\N(X/{\Ze})$.
\end{lemma}
\begin{proof}
  Using \cite[4.6.4]{mn:center}, the assumption of the lemma, and
  %\ref{quotcenter} below, 
  \cite[5.11]{jmm:ndet},
  we get ${\Ze}(X/{\Ze}) = {\Ze}(X)/{\Ze} = {\Ze}(\N(X))/{\Ze} =
  {\Ze}(\N(X)/{\Ze})={\Ze}\N(X/{\Ze})$.
\end{proof}

\section{Reduction to the connected, centerless (simple) case}\label{sec:red}

In this section we reduce the general classification problem first to
the connected case and next to the connected and centerless case. We
first show (\ref{lemma:autononcon}, \ref{redtoconnected}) that if $X$
is any nonconnected \twocg\ with identity component $X_0$ then
\begin{equation*}
  \left.
    \begin{array}[l]{l}
      \text{$X_0$ is uniquely ${\N}$-determined} \\
      \text{$X$ is LHS} \\
      \text{$H^i({\W}/{\W}_0;\ch{{\Ze}}(X_0)) \to 
             H^i({\W}/{\W}_0;\ch{{\T}}^{{\W}_0})$ is injective for $i=1,2$}
    \end{array}
    \right\} \Longrightarrow
    \text{$X$ is totally ${\N}$-determined}
\end{equation*}
The $H^i$-injectivity conditions holds when $X_0$ is a connected Lie
group \cite[1.6]{matthey:normalizers} or equals $\mathrm{DI}(4)$
\cite{dw:new,no:DI(4)}. To see this observe that the condition
obviously holds when $\ch{{\Ze}}(X_0)= \ch{{\T}}(X_0)^{{\W}_0}$ or
$\ch{{\Ze}}(X_0)$ is trivial.  We shall later see that any connected
\twocg\ splits as a product of a compact connected Lie group and a
finite number of $\mathrm{DI}(4)$, and from this it follows that this
condition is always satisfied.  Indeed, let $X_0=G' \times G'' \times
\mathrm{DI}(4)^s$ where $G'$ is a connected compact Lie group with no
direct factors isomorphic to $\SO(2n+1)$, $G''$ is a direct product of
$\SO(2n+1)$s, and $s \geq 0$. The $\pi_0(X)$-equivariant group homo\m\ 
$\ch{Z}(G')=\ch{Z}(X_0) \to \ch{T}(G')^{W(G')} \times
\ch{T}(G'')^{W(G'')}$ has a left inverse since it takes $\ch{Z}(G')$
isomorphically to the $\pi_0(X)$-subgroup $\{1\} \times
\ch{T}(G'')^{W(G'')}$ of the left hand side. The induced map on
cohomology therefore also has a left inverse. However, it is not at
present clear if all nonconnected \twocg s are LHS.

Next we consider a connected \twocg\ $X$ with adjoint form $PX=X/Z(X)$
\cite[2.8]{dw:center}
and show (\ref{autondetcenter}, \ref{ndetcenter}) that
\begin{equation*}
  \text{$PX$ is uniquely ${\N}$-determined} \Longrightarrow
  \text{$X$ is uniquely ${\N}$-determined}
\end{equation*}
This reduces in principle the problem to the connected and centerless
case. One can go a little further since connected, centerless \twocg s
split into products of simple factors \cite{dw:split,no:split}. We show 
(\ref{prodauto}, \ref{ndetprod}) that
\begin{equation*}
  \text{$X_1$ and $X_2$ are uniquely ${\N}$-determined} \Longrightarrow
  \text{$X_1 \times X_2$ is uniquely ${\N}$-determined}
\end{equation*}
when $X_1$ and $X_2$ are connected.  Therefore it suffices to show
that all {\em connected, centerless and simple\/} \twocg s are
uniquely ${\N}$-determined. It is already known that all connected
compact Lie groups as well as $\mathrm{DI(4)}$ have
$\pi_*({\N})$-determined auto\m s \cite{jmo:selfho,no:DI(4)}.

Let $X$ be a \twocg\ with \mtnp\ $({\N},{\N}_0)(X)=({\N},{\N}_0)$.

\begin{lemma}\label{lemma:autononcon}\cite[4.2]{jmm:deter}
Suppose that $X_0$ has ${\N}$-determined auto\m s. Then
\begin{equation*}
  \text{$X$ has ${\N}$-determined auto\m s} \iff
       \text{$H^1({\W}/{\W}_0;\ch{{\Ze}}(X_0)) \to
   H^1({\W}/{\W}_0;\ch{{\T}}^{{\W}_0})$ is injective}    
\end{equation*}
\end{lemma}
\begin{proof}
  The restriction of $\AM$ to the subgroup
  $H^1({\W}/{\W}_0;\ch{{\Ze}}(X_0)) \subset \Aut(X)$ is the homo\m\ 
  \set\begin{equation}\label{AMres}
    H^1({\W}/{\W}_0;\ch{{\Ze}}(X_0)) \to  
    H^1({\W}/{\W}_0;\ch{{\T}}^{{\W}_0}) \xrightarrow{\mathrm{inf}}
    H^1({\W};\ch{{\T}})
  \end{equation}\add
  where $\mathrm{inf}$ is the inflation mono\m . If the first homo\m\ 
  has a nontrivial kernel, $X$ does not have ${\N}$-determined auto\m
  s.  Conversely, assume that the first homo\m\ is injective, and let
  $f \in \Aut(X)$ be an auto\m\ such that $\AM(f) \in {\W}_0
  \backslash \Aut({\N})$ is the identity. Then $\AM(f_0) \in {\W}_0
  \backslash \Aut({\N}_0)$ and $\pi_0(f)$ equal the respective
  identity maps. Since $X_0$ has ${\N}$-determined auto\m s by
  assumption, $f_0$ is the identity. Thus $f$ belongs to the subgroup
  $ H^1({\W}/{\W}_0;\ch{{\Ze}}(X_0))$ of $\Aut(X)$ \cite[5.2]{jmm:ext}
  where $\AM$ is injective, so $f$ is the identity auto\m\ of $X$. (The
  description of the kernel in the short exact sequence of
  \cite[5.2]{jmm:ext} holds for all \pcg s, not just those with 
  a completely reducible identity component.)
\end{proof} 

\begin{lemma}\label{idcomp}
  Suppose that $X$ has ${\N}$-determined auto\m s and that $X_0$ has
  $\pi_*({\N})$-determined auto\m s. Then $\Aut(X) \cap H^1({\W};\ch{{\T}})=
  H^1({\W}/{\W}_0;\ch{{\Ze}}(X_0))$ so that
  \begin{equation*}
    \text{$X$ has $\pi_*({\N})$-determined auto\m s} \iff
    H^1({\W}/{\W}_0;\ch{{\Ze}}(X_0)) = 0
  \end{equation*}
\end{lemma}
\begin{proof}
  Let $f \in \Aut(X)$ be an auto\m\ such that $\pi_*\AM(f)$ is the
  identity. Then also $\pi_*\AM(f_0)$ and $\pi_0(f)$ equal the
  respective identity maps. Since $X_0$ is assumed to have
  $\pi_*({\N})$-determined auto\m s, $f_0$ is the identity. Thus $f$
  belongs to the subgroup $ H^1({\W}/{\W}_0;\ch{{\Ze}}(X_0))$ of
  $\Aut(X)$ \cite[5.2]{jmm:ext}. This shows that $\Aut(X) \cap
  H^1({\W};\ch{{\T}}) \subset H^1({\W}/{\W}_0;\ch{{\Ze}}(X_0))$. The
  opposite inclusion is immediate from (\ref{AMres}).
\end{proof}

%\begin{lemma}\label{lemma:autononcon}\cite[4.2]{jmm:deter}
%  Suppose that 
%  \begin{enumerate}
%  \item $X_0$ has $N$-determined auto\m s
%  \item $H^1\big(W/W_0;\ch{Z}(X_0)\big) \to
%    H^1\big(W/W_0;Z(\ch{N_0})\big)$ is injective
%  \end{enumerate}
%  Then $X$ has $N$-determined auto\m s.
%\end{lemma}

\begin{lemma}\label{autondetcenter}\cite[4.8]{jmm:deter}
  Suppose that $X$ is connected.
  If the adjoint form $PX=X/{\Ze}(X)$ has $\pi_*({\N})$-determined auto\m
  s, so does $X$.
\end{lemma}
\begin{proof}
  If $f \in \Aut(X)$ is an auto\m\ under ${\T}(X)$, the induced auto\m\ 
  $Pf \in \Aut(PX)$ is an auto\m\ under ${\T}(PX)$, hence equals the
  identity, and the induced auto\m\ ${\Ze}(f) \in \Aut({\Ze}X)$ is also the
  identity since the center ${\Ze}X \to X$ factors through the \mt\ ${\T}(X)
  \to X$ \cite[7.5]{dw:center} \cite[4.3]{mn:center}. But then $f$
  itself is the identity for $\Aut(X)$ embeds into $\Aut(PX) \times
  \Aut({\Ze}X)$ \cite[4.3]{jmm:rip}.
\end{proof}

\begin{lemma}\label{prodauto}\cite[9.4]{jmm:ndet}
  If the two connected \twocg s $X_1$ and $X_2$ have ${\N}$-determined
  (resp.\ $\pi_*({\N})$-determined) auto\m s, so does the product $X_1
  \times X_2$.
\end{lemma}
\begin{proof}
  Since the statement concerning ${\N}$-determined auto\m s is proved in
  \cite[9.4]{jmm:ndet} we deal here only with the case of
  $\pi_*({\N})$-determined auto\m s.  Let $f$ be an auto\m\ under ${\T}_1
  \times {\T}_2$ of the product \twocg\ $X_1 \times X_2$. Then
  \begin{align*}
    f_1 &\colon X_1 \to X_1 \times X_2 \xrightarrow{f} X_1 \times X_2
    \to X_1 \\
    f_2 &\colon X_2 \to X_1 \times X_2 \xrightarrow{f} X_1 \times X_2
    \to X_2
  \end{align*}
  are endo\m s under the \mts\ and therefore conjugate to the
  respective identity maps. But $f$ is \cite[9.3]{jmm:ndet} in fact
  conjugate to the product \m\ $(f_1,f_2)$ which is the identity.
\end{proof}

\begin{lemma}\label{redtoconnected} (Cf \cite[7.8]{jmm:deter})
  Suppose that
  \begin{enumerate}
  \item $X_0$ is uniquely ${\N}$-determined.
  \item $X$ is LHS. \label{redtoconnected.2}
  \item $H^2\big({\W}/{\W}_0,\ch{{\Ze}}(X_0)\big) \to
    H^2\big({\W}/{\W}_0,\ch{{\T}}^{{\W}_0}\big)$ is injective.
    \label{redtoconnected.3}
  \end{enumerate}
  Then $X$ is ${\N}$-determined. 
  %If also $H^1(W_0;\ch{T})^{W/W_0} = 0$
%  and the conditions of \ref{lemma:autononcon} are satisfied then
%  $X$ is strongly $N$-determined.
\end{lemma}

\begin{proof}
  Let $X'$ be another \twocg\ with \mtnp\ $({\N},{\N}_0)$. The
  assumption on the identity component $X_0$ means (\ref{underT}) that
  there exists an iso\m\ \func{f_0}{X_0}{X_0'} under ${\T}$. For any
  $\xi \in {\W}/{\W}_0 = {\N}/{\N}_0 = X/X_0 = X'/X_0'$, the iso\m\ 
  $\xi f_0 \xi^{-1}$ is also an iso\m\ under ${\T}$ and thus $\xi f_0
  = f_0 \xi$ as $X_0$ is uniquely ${\N}$-determined. By the second
  assumption, the auto\m\ \func{\alpha_0 = \AM(f_0)}{{\N}_0}{{\N}_0}
  with $\pi_*(B\alpha_0)=\mathrm{Id}$ extends to an iso\m\ 
  \func{\alpha}{{\N}}{{\N}} with $\pi_*(B\alpha)=\mathrm{Id}$.

  The situation is now as shown in the commutative diagram
  \begin{equation*}
    \xymatrix{
      BX_0 \ar[d]\ar[d] \ar@/^1.5pc/[rrr]^{Bf_0} &
      BN_0 \ar[d]\ar[l]^{Bj_0} \ar[r]_{B\alpha_0} &
      BN_0 \ar[d]\ar[r]_{Bj_0'} &
      BX_0' \ar[d] \\
      BX \ar[d] &
      BN \ar[d]\ar[l]_{Bj} \ar[r]^{B\alpha} &
      BN \ar[d]\ar[r]^{Bj'} &
      BX'\ar[d] \\
      B\pi_0(X) &
      B(W/W_0) \ar[l]_{\cong} &
      B(W/W_0) \ar@{=}[l]\ar[r]^{\cong} &
      B\pi_0(X') }
  \end{equation*}
  
  Our aim is to find an iso\m\ \func{f}{X}{X'} to fill in the based
  homotopy commutative diagram
\begin{equation*}
    \xymatrix{
         {BX_0}\ar[d]\ar[r]^{Bf_0}_{\cong} & {BX'_0}\ar[d] \\
         {BX} \ar@{ .>}[r]\ar[d] & {BX'}\ar[d] \\
         {B\pi_0(X)}\ar[r]_{\cong} & {B\pi_0(X')} }
\end{equation*}
where the iso\m\ between the base \twocg s is given by the iso\m s
$\pi_0(X) \leftarrow {\N}/{\N}_0 \rightarrow \pi_0(X')$.  Since $f_0$
is $\pi_0(X)$-equivariant up to homotopy, $\map(BX_0,BX_0')_{Bf_0}$
is a $\pi_0(X)$-space in the sense that there exists a fibration
\begin{equation*}
  \map(BX_0,BX_0';Bf_0) \to  \map(BX_0,BX_0';Bf_0)_{h\pi_0(X)} \to
  B\pi_0(X) 
\end{equation*}
over $B\pi_0(X)$ with $\map(BX_0,BX_0')_{Bf_0}$, here written as
$\map(BX_0,BX_0';Bf_0)$, as fibre. The space of sections of this
fibration, $\map(BX_0,BX_0')_{Bf_0}^{h\pi_0(X)}$, is a space of fibre
maps of $BX$ to $BX'$.  This fibration sits to the left in the
commutative diagram
\begin{equation*}
  \xymatrix@C=10pt{
 {\map(BX_0,BX_0';Bf_0)} \ar[r] \ar[d] &
 {\map(BN_0,BX_0';B(j_0'\alpha))} \ar[d] &
 {\map(BN_0,BN_0;B\alpha_0)} \ar[d] \ar[l]_-{\simeq} \\
 {\map(BX_0,BX_0';Bf_0)_{h\pi_0(X)}} \ar[r] \ar[d] &
 {\map(BN_0,BX_0';B(j_0'\alpha))_{h(W/W_0)}}  \ar[d] &
 {\map(BN_0,BN_0;B\alpha_0)_{h(W/W_0)}} \ar[d] \ar[l] \\
 B\pi_0(X) &
 B(W/W_0) \ar[l] \ar@{=}[r] &
 B(W/W_0) }
\end{equation*}
where the columns are fibrations and the horizontal maps are defined
as composition with $Bj$ and $Bj'$, respectively. The fibre map from
the right column to the central one is actually a fibre homotopy
equivalence because the centralizer of the \mt\ in $X_0'$ and in
$\N_0$ are isomorphic in that they are both isomorphic to the \mt .

The middle fibration admits a section corresponding to the fibrewise
map $Bj' \circ B\alpha$. But then the left fibration also admits a
section: The obstruction to a section of the left fibration is a
cohomology class in $H^2(\pi_0(X); \ch{\Ze}(X_0))$. Since the middle
fibration does admit a section, this obstruction class is in the
kernel of the coefficient group homo\m\ $H^2(\pi_0(X); \ch{\Ze}(X_0))
\to H^2\big({\W}/{\W}_0,\ch{{\T}}^{{\W}_0}\big)$. But the assumption
is that this is an injection and therefore the obstruction must
vanish. (We are here tacitly replacing the three fibrations above by
their fibrewise discrete approximations \cite[4.3]{jmm:ext}.)

%Composition with $BX \xleftarrow{Bj} B{\N}
%\xrightarrow{Bj'} BX'$ gives maps
%  \begin{multline*}
%    \map(BX_0,BX_0';Bf_0)^{h{\W}/{\W}_0} \xrightarrow{Bj^*}
%    \map(B{\N}_0,BX_0';{B(j_0'\alpha)})^{h{\W}/{\W}_0} \\
%    \xleftarrow[\simeq]{Bj'_*} \map(B{\N}_0,B{\N}_0;{B\alpha_0})^{h{\W}/{\W}_0} 
%  \end{multline*}
%  of homotopy fixed point spaces. The space to the right is non-empty
%  for it contains the iso\m\ \func{B\alpha}{B{\N}}{B{\N}}.  Using
%  obstruction theory and assumption (\ref{redtoconnected.3}), we see
%  that also the homotopy fixed point space to the left is non-empty;
A section of the left fibration corresponds to a \m\ 
\func{Bf}{BX}{BX'} under the iso\m\ \func{Bf_0}{BX_0}{BX'_0} and over
$B\pi_0(X) \xrightarrow{\cong}B\pi_0(X')$ such that $Bf \circ Bj$ and
$Bj \circ B\alpha$ are homotopic over $B({\N}/{\N}_0) \rightarrow
B\pi_0(X')$.  But since the fibre $BX'_0$ of $BX' \to B\pi_0(X')$ is
simply connected this means that $Bf \circ Bj$ and $Bj \circ B\alpha$
are based homotopic maps $B{\N} \to BX'$.
\end{proof}

%\begin{cor}\label{cor:redtoconnected}
%  Suppose that 
%  \begin{enumerate}
%  \item $X_0$ is uniquely $N$-determined and completely reducible,
%  \item \func{\mathrm{res}}{H^1(W;\ch{T})}{H^1(W_0;\ch{T})^{W/W_0}} is
%  an iso\m , \label{cor:redtoconnected2}
%  \item $\ch{Z}(X_0)=Z(\ch{N}_0)$.
%  \end{enumerate}
%Then $X$ is uniquely $N$-determined.
%\end{cor}
%\begin{proof}
%Combine \ref{idcomp} and \ref{redtoconnected} and use
%the $5$-term exact sequence
%\begin{equation*}
%  0 \to H^1(W/W_0;\ch{T}^{W_0}) 
%  \xrightarrow{\mathrm{inf}} H^1(W;\ch{T}) 
%  \xrightarrow{\mathrm{res}}{H^1(W_0;\ch{T})^{W/W_0}} 
%  \xrightarrow{d_2} H^2(W/W_0;\ch{T}^{W_0})
%  \xrightarrow{\mathrm{inf}} H^2(W;\ch{T})
%\end{equation*}
%from the Lyndon--Hochschild--Serre spectral sequence
%\cite{hochserre}. 
%\end{proof}

\begin{exmp}\label{ex:toral}
  1.\ Any \twoct\ ${\T}$ is uniquely ${\N}$-determined for if
  \func{j}{{\T}}{X} is the \mtn\ for the connected \twocg\ $X$, then $j$
  is an iso\m .  Indeed, $H^*(B{\T};\Q_2) \cong H^*(BX;\Q_2)$
  \cite[9.7.(3)]{dw:fixpt} and the connected space $X/{\T}$ has
  cohomological dimension $\cd_{\F_2}(X/{\T})=0$ \cite[4.5,
  5.6]{dw:center} so is a point.  \\
  2.\ Any \twoctg\ $G$ is totally
  ${\N}$-determined: $G$ clearly has ${\N}$-determined auto\m s as
  $G$ is its own \mtn .
  If the \twocg\ $X$ has the same \mtnp\ $(G,{\T})$ as $G$, then $X$ is a
  \twoctg\ and \func{j'}{G}{X} is an iso\m .
  %so
%  \begin{equation*}
%    \xymatrix{
%      G \ar[d]_j^{\cong} \ar[r]^{\alpha}_{\cong} &
%      G \ar[d]^{j'}_{\cong} \\
%      G \ar[r]^{\cong}_{j'\alpha j^{-1}} & X }
%  \end{equation*}
%  solves diagram (\ref{dia:BaBf}) for any $\alpha\in
%  H^1(\pi_0(G);\ch{T})$.  Any \twoctg\ $G$ is clearly strongly totally
%  $N$-determined \eqref{redtoconnected}.  Thus $G$ is uniquely
%  $N$-determined if and only if this cohomology group vanishes.
  $G$ is uniquely ${\N}$-determined if and only if $H^1(\pi_0(G);\ch{{\T}})
  = 0$. In
  particular, $\GL(2,\R)$ is uniquely ${\N}$-determined.\\
 \end{exmp}

%\begin{lemma}
%  Let $X$ be a connected \twocg\ and $Z \to X$ its center. Suppose
%  that $X/Z$ is strongly $N$-determined or $N$-determined in such a
%  way that we can always find $\alpha/Z$ in the image of the
%  coefficient homo\m\ $H^1(W(X);\ch{T}(X)) \to
%  H^1(W(X);\ch{T}(X)/\ch{Z})$. Then $X$ is (strongly) $N$-determined.
%\end{lemma}
\begin{lemma}\label{ndetcenter} (Cf \cite[7.10]{jmm:deter})
    Suppose that $X$ is connected.
  If the adjoint form $PX=X/{\Ze}(X)$ is ${\N}$-determined, so is
  $X$. 
\end{lemma}
\begin{proof}
  Let \func{j}{{\N}}{X} be the \mtn\ for $X$ and \func{j'}{{\N}}{X'} the
  \mtn\ for some other connected \twocg\ $X'$. It  suffices
  \eqref{underT} to find a
  \m\ \func{f}{X}{X'} under the \mts\ $X \xleftarrow{i} {\T}
  \xrightarrow{i'} X'$. The $2$-discrete center $\ch{{\Ze}}$ of $X$ and
  $X'$ is contained in the the $2$-discrete \mt\ $\ch{{\T}}$
  \cite[7.5]{dw:center}.  Factoring out \cite[8.3]{dw:fixpt} these
  central mono\m s we obtain the commutative diagram
  \begin{equation*}
    \xymatrix@C=40pt{
       B\ch{X} \ar[d] & B\ch{{\T}} \ar[l]_{Bi}\ar[r]^{Bi'}\ar[d] 
       & B\ch{X}' \ar[d] \\
       B(X/{\Ze})  \ar@/_1.5pc/[rr]|-{B(f/{\Ze})} & B({\T}/{\Ze}) \ar[l]_{B(i/{\Ze})}
       \ar[r]^{B(i'/{\Ze})} & B(X'/{\Ze}) }
  \end{equation*}
  where the vertical maps are fibrations with fibre $B\ch{{\Ze}}$, the
  total spaces, such as $B\ch{X}$, are the fibre-wise discrete
  approximations, and \func{f/{\Ze}}{X/{\Ze}}{X'/{\Ze}} is the iso\m\ under ${\T}/{\Ze}$
  that exists because $X/{\Ze}$ is ${\N}$-determined. Construct the fibration
  \begin{equation*}
    \map(B\ch{{\Ze}},B\ch{{\Ze}};B1) \to B\ch{{\Ze}}_{h(X/{\Ze})} \to B(X/{\Ze})
  \end{equation*}
  whose sections are maps $BX \to BX'$ over $B(f/{\Ze})$ and under
  $B\ch{{\Ze}}$. There are two other such fibrations related to this one
  as shown in the commutative diagram 
  \begin{equation*}
    \xymatrix@C=45pt{
        {\map(B\ch{{\Ze}},B\ch{{\Ze}};B1)} \ar@{=}[r] \ar[d] &
        {\map(B\ch{{\Ze}},B\ch{{\Ze}};B1)} \ar@{=}[r] \ar[d] &
        {\map(B\ch{{\Ze}},B\ch{{\Ze}};B1)}  \ar[d] \\
        B\ch{{\Ze}}_{h(X/{\Ze})} \ar[d] & 
        B\ch{{\Ze}}_{h({\T}/{\Ze})} \ar[d] \ar[l] \ar[r]^{Bi^*}  &
        B\ch{{\Ze}}_{h({\T}/{\Ze})} \ar[d] \\
        B(X/{\Ze}) & B({\T}/{\Ze})\ar[l]^{B(i/{\Ze})} \ar@{=}[r] & B({\T}/{\Ze}) }
  \end{equation*}
  where the middle fibration is the   pull-back along $B(i/{\Ze})$ of the
  left fibration and the fibre over $b \in B({\T}/{\Ze})$ of the right
  fibration consists of one component of
  the space of maps of the fibre $B\ch{{\T}}_b$ over $b$ into the fibre
  $B\ch{X}'_{B(i'/{\Ze})(b)}$ over $B(i'/{\Ze})(b)$. The fibre equivalence
  $Bi^*$ is induced by \func{Bi}{B\ch{{\T}}}{B\ch{X}}. The middle
  fibration has a section $u'$ such that $Bi^* \circ u'$ is the
  section \func{Bi'}{B\ch{{\T}}}{B\ch{X}'} of the right fibration. We now
  have fibre maps 
  \begin{equation*}
    \xymatrix{
          X/{\T} \ar[d] \ar[rr]^{u \vert X/{\T}} && B\ch{{\Ze}}\ar[d] \\
    B\ch{{\T}} \ar[dr]_{B(i/{\Ze})} \ar[rr] && 
               B\ch{{\Ze}}_{h(X/{\Ze})}\ar[dl] \\
          & B(X/{\Ze}) }
  \end{equation*}
  where $u$ is the composition of $u'$ and
  $B\ch{{\Ze}}_{h({\T}/{\Ze})} \to B\ch{{\Ze}}_{h(X/{\Ze})}$. The
  canonical map, given by constants, $B\ch{{\Ze}} \to
  \map(X/{\T},B\ch{{\Ze}})$ is a \he\ since $X/{\T}$ is simply
  connected \cite[5.6]{mn:center} and hence a version
  \cite[6.6]{jmm:deter} of the Zabrodsky lemma implies that $u=v \circ
  B(i/{\Ze})$ for some section
  \func{v}{B(X/{\Ze})}{B\ch{{\Ze}}_{h(X/{\Ze})}} of the left
  fibration. The section $v$ is, after fibre-wise completion, a fibre
  map $BX \to BX'$ under $B{\T}$.
\end{proof}

Let \func{j_1}{{\N}_1}{X_1} and \func{j_2}{{\N}_2}{X_2} be \mtn s for the
connected \twocg s $X_1$ and $X_2$ and suppose that $X'$ is some
connected \twocg\ that admits a \mtn\ of the form \func{j'}{{\N}_1 \times
  {\N}_2}{X'}.  The Splitting Theorem \cite[1.4]{dw:split}, or more
explicitly in the form of \cite[5.5]{no:split}, says that there exist
\twocg s $X_1'$ and $X_2'$ and an iso\m\ $X_1' \times X_2' \to X'$
such that
\begin{equation*}
  \xymatrix{
    & B{\N}_1 \times B{\N}_2 \ar[dl]_{Bj_1' \times Bj_2'\;\;} \ar[dr]^{Bj'} \\
    BX_1' \times BX_2' \ar[rr]^{\simeq} && BX' }
\end{equation*}
commutes where \func{j_1'}{{\N}_1}{X_1'} and \func{j_2'}{{\N}_2}{X_2'} are
\mtn s. The following lemma is an immediate consequence
\begin{lemma}\label{ndetprod}
  The product of two ${\N}$-determined connected \twocg s is
  ${\N}$-determined.  
\end{lemma}
\begin{proof}
  Since $X_1$, $X_2$ are ${\N}$-determined there exist iso\m s
  \func{f_1}{X_1}{X_1'}, \func{f_2}{X_2}{X_2'} and auto\m s $\alpha_1
  \in H^1({\W}_1;\ch{{\T}}_1) \subset \Out({\N}_1)$, $\alpha_2 \in
  H^1({\W}_2;\ch{{\T}}_2) \subset \Out({\N}_2)$ such that
  \begin{equation*}
    \xymatrix@C=45pt{
      B{\N}_1 \times B{\N}_2 \ar[d]_-{Bj_1 \times Bj_2} \ar[r]^{B\alpha_1
        \times B\alpha_2} &
      B{\N}_1 \times B{\N}_2 \ar[d]_{Bj_1' \times Bj_2'} \ar[dr]^{Bj'} \\
      BX_1 \times BX_2 \ar[r]_-{Bf_1 \times Bf_2} &
      BX_1' \times BX_2' \ar[r]^-{\simeq} &
      BX' }
  \end{equation*}
  commutes up to based homotopy.
\end{proof}

\section{$\N$-determined connected, centerless $2\/$-compact groups}
\label{sec:conncenterless}

In this section we formulate inductive criteria that, at least in
favorable cases, can be used to show total ${\N}$-determinacy for 
connected, centerless (simple) \twocg s $X$. The key tool is the
homology decomposition \cite[8.1]{dw:center}
\set\begin{equation}\label{hodecomp}
  \hocolim_{\A(X)^{\op}} BC_X \to BX
\end{equation}\add
of $BX$ in terms of centralizers of elementary abelian subgroups.
Since $X$ has no center, the cohomological dimension of each
centralizer $C_X(V,\nu)$ is smaller than the cohomological dimension
of $X$. As part of an inductive argument we will therefore assume that
all centralizers are totally ${\N}$-determined and formulate criteria
(\ref{ndetauto}, \ref{indstepalt}) that imply that also $X$ is totally
${\N}$-determined. 

\begin{defn}\cite[\S8]{dw:center}
  \label{defn:AX}
   The objects of the Quillen category $\A(X)$  are
conjugacy classes of mono\m s \func{\nu}{V}{X} of nontrivial \lmntwo s
into $X$; the \m s \func{\alpha}{(V_1,\nu_1)}{(V_2,\nu_2)} are
injective group homo\m s \func{\alpha}{V_1}{V_2} such that $\nu_1$ and
$\nu_2\alpha$ are conjugate mono\m s $V_1 \to X$. We shall write
$\A(X)(V_1,V_2)$ for the set of \m\ $V_1 \to V_2$ and $\A(X)(V)$ for
the group of all endo\m s (which are all iso\m s) of $V$. 
\end{defn}

The  functor 
\set\begin{equation}
  \label{eq:defnBCX}
  \func{BC_X}{\A(X)^{\textmd{op}}}{\Topspaces} \qquad \textmd{(topological
spaces)}
\end{equation}\add
 takes an object $(V,\nu)$ of the Quillen category $\A(X)$ to
its centralizer $BC_X(V,\nu)=\map(BV,BX)_{B\nu}$.  The covariant
functor 
\set\begin{equation}
  \label{eq:defnpiiBZC}
\func{\pi_i(B{\Ze}C_X)}{\A(X)}{\Ab} \qquad \textmd{(abelian groups) }
\end{equation}\add
takes $(V,\nu)$ into the abelian homotopy group
$\pi_i(\map(BC_X(V,\nu),BX),{e(\nu)})$ based at the evaluation map
\func{e(\nu)}{BC_X(V,\nu)}{BX}. The space $\map(BC_X(V,\nu),BX)$ is
homotopy equivalent to $B\Ze C_X(V,\nu)$ \cite{dwk:centric}).

\begin{lemma}\label{ndetauto}\cite[4.9]{jmm:deter}
Suppose that $X$ is connected and centerless. If
\begin{enumerate}
\item \label{ndetauto1} $C_X(L,\lambda)$ has ${\N}$-determined
  (resp.\ $\pi_*({\N})$-determined) auto\m s for each rank $1$ object
  $(L,\lambda)$ of $\A(X)$ and
\item \label{ndetauto2}
  ${\operatorname{lim}}^1({\A(X)};\pi_1(B{\Ze}C_X))=0=
  {\operatorname{lim}}^2({\A(X)};\pi_2(B{\Ze}C_X))$
\end{enumerate}
Then $X$ has ${\N}$-determined (resp.\ $\pi_*({\N})$-determined) auto\m s.
\end{lemma}
\begin{proof}
  Suppose first that each line centralizer has $\pi_*({\N})$-determined
  auto\m s.
  Let \func{f}{X}{X} be an auto\m\ under the \mt\ ${\T} \to X$. Since any
  mono\m\ \func{\lambda}{L}{X}, $L=\Z/2$, factors through the \mt ,
  the commutative diagram
  \begin{equation*}
    \xymatrix@R=10pt{
    && {\N} \ar[dd]^{\AM(f)} \ar[r] & X \ar[dd]^f \\
    L \ar[r]^{\lambda^{\T}} & {\T} \ar[ur] \ar[dr] \\
    && {\N} \ar[r] & X}
  \end{equation*}
  shows that $f\lambda = \lambda$ and gives a commutative diagram
  \begin{equation*}
    \xymatrix@R=10pt{
    & C_{\N}(L) \ar[dd]^{C_{\AM(f)}(L)} \ar[r] & C_X(L) \ar[dd]^{C_f(L)}
    \\
    {\T} \ar[ur] \ar[dr] \\
    & C_{\N}(L) \ar[r] & C_X(L) }
  \end{equation*}
  of auto\m s under ${\T}$. Thus $\AM(C_f(L)) = C_{\AM(f)}(L) \colon
  C_{\N}(L) \to C_{\N}(L)$.  Now, $\pi_*(C_{\N}(L))$ is a subgroup of
  $\pi_*({\N})$ (for $\pi_1(C_{\N}(L)) = \pi_1({\N})$ and $\pi_0(C_{\N}(L)) =
  {\W}(X)(L)$ is \cite[7.6]{dw:center} \cite[3.2.(1)]{jmm:rip} the
  stabilizer subgroup at $L < \ch{{\T}}$ for the action of ${\W}(X)$ on
  $\ch{{\T}}$) so $\pi_*( C_{\AM(f)}(L)) =1$ and $C_f(L) \simeq
  1_{C_X(L)}$ since $C_X(L)$ has $\pi_*({\N})$-determined auto\m s. For
  any other object $(V,\nu)$ of $\A(X)$ of rank $>1$, choose a line
  $L$ in $V$. Since the mono\m\ \func{\nu}{V}{X} canonically factors
  through $C_X(L)$ \cite[8.2]{dw:fixpt} \cite[3.18]{jmm:ndet}, the
  commutative diagram
  \begin{equation*}
    \xymatrix@R=10pt{
    && X \ar[dd]^f \\
    V \ar[r] \ar@/^/[urr]^{\nu}  \ar@/_/[drr]_{\nu} & C_X(L) \ar[ur] \ar[dr] \\
    && X}
  \end{equation*}
  shows that $f\nu=\nu$ and the induced diagram
  \begin{equation*}
    \xymatrix@R=10pt{
     & C_X(V) \ar[dd]^{C_f(V)} \\
     C_{C_X(L)}(V) \ar[ur]^{\cong} \ar[dr]_{\cong} \\
     & C_X(V) }
  \end{equation*}
  that \func{C_f(V)}{C_X(V)}{C_X(V)} is conjugate to the
  identity.  The third %corrected from ``second''
  assumption of the lemma assures that there are no
  obstructions to conjugating $f$ to the the identity now that we 
  know that the restriction of  $f$ to each of the
  centralizers is conjugate to the identity, see
  \cite[4.9]{jmm:deter}. 

  Suppose next that each line centralizer has ${\N}$-determined auto\m
  s. Let \func{f}{X}{X} be an auto\m\ such that the diagram
  \begin{equation*}
    \xymatrix@R=10pt{
      & X \ar[dd]^f \\
      {\N} \ar[ur] \ar[dr] \\
      & X}
  \end{equation*}
  commutes up to conjugacy. For each line $L$ in ${\T}$, the induced
  diagram 
  \begin{equation*}
    \xymatrix@R=10pt{
      & C_X(L) \ar[dd]^{C_f(L)} \\
      C_{\N}(L) \ar[ur] \ar[dr] \\
      & C_X(L)}
  \end{equation*}
  also commutes up to conjugacy. By assumption, this means
  (\ref{unbasedOut}) that the induced auto\m s $C_f(L)$ of line
  centralizers are conjugate to the identity. As above, this implies
  that the induced map \func{C_f(V)}{C_X(V)}{C_X(V)} is conjugate to
  the identity for any object $(V,\nu)$ of the Quillen category for
  $X$ and that $f$ is conjugate to the identity.
\end{proof}

Consider next an \etwoct\ ${\N}$ and two connected,
centerless \twocg s $X$ and $X'$ both having ${\N}$ as their \mtn\ 
\set\begin{equation}\label{eq:situation} \xymatrix@C=30pt{ X & {\N}
    \ar[r]^-{j'} \ar[l]_-{j} & X'}
\end{equation}\add
Our task is (\ref{underT}.\ref{underT1}) to construct an iso\m\ $X \to
X'$ under the \mt .

\begin{defn}
  \label{defn:toralAX}
  An object $(V,\nu)$ of $\A(X)$ is {\em toral\/} if the mono\m\ 
  \func{\nu}{V}{X} factors through the \mt\ ${\T} \to X$. Let
  $\A(X)^{\leq t}$ denote the full subcategory of toral objects, and
  $\A(X)^{\leq t}_{\leq 2}$ the full subcategory of toral objects of
  rank $\leq 2$.
\end{defn}

  For each toral object $(V,\nu)$ of $\A(X)^{\leq t}$, let
\func{\nu^{\N}}{V}{{\N}} be the unique \pl\ \cite[4.10]{jmm:normax} of $\nu$
(which factors through the identity component of ${\N}$) and let
$(V,\nu')$ be the toral object of $\A(X')$ defined by \func{\nu' = j
  \circ \nu^{\N}}{V}{X'} as in the commutative diagram
\begin{equation*}
  \xymatrix{
    & V \ar[dl]_{\nu} \ar[d]^{\nu^{\N}} \ar[dr]^{\nu'} \\
    X & {\N} \ar[l]^j \ar[r]_{j'} & X'}
\end{equation*}
The functor $\A(X)^{\leq t} \to \A(X')^{\leq t}$ that takes the object
$(V,\nu)$ to the object $(V,\nu')$ and is the identity on \m s is an
equivalence of toral Quillen categories \cite[2.8]{jmm:ndet}.

\begin{thm}\label{indstepalt} (Cf \cite[3.8]{jmm:ndet})
  In the situation of (\ref{eq:situation}), assume the following:
  \begin{enumerate}
  \item \label{indstepalt3} The centralizer $C_X(V,\nu)$ of any
    $(V,\nu) \in \Ob(\A(X)^{\leq t}_{\leq 2})$ has ${\N}$-determined
    auto\m s.
  \item \label{indstepalt1} There exists a self-homotopy equivalence
    $\alpha \in H^1({\W};\ch{{\T}}) \subseteq \Out({\N})$ such that for every
    object $(L,\lambda) \in \Ob(\A(X)^{\leq t}_{\leq 1})$ the diagram 
    \begin{equation*}
      \xymatrix@C=50pt{
          C_{\N}(L,\lambda^{\N}) \ar[d]_{j \vert C_{\N}(\lambda^{\N})}
          \ar[r]^{\alpha \vert C_{\N}(\lambda^{\N})} &
          C_{\N}(L,\lambda^{\N}) \ar[d]^{j' \vert C_{\N}(\lambda^{\N})} \\
          C_X(L,\lambda) \ar[r]_{f_{\lambda}} & C_{X'}(L,\lambda') }
    \end{equation*}
    commutes for some iso\m\ $f_{\lambda}$.
    \item \label{indstepalt4} For any nontoral rank two object
    $(V,\nu)$ of $\A(X)$ the composite mono\m\ 
    \begin{equation*} \nu_L' \colon
      \xymatrix@C=40pt{
        V \ar[r]^-{\overline{\nu}(L)} &
        C_X(L,\nu\vert L) \ar[r]^-{f_{\nu\vert L}}_-{\cong} &
        C_{X'}(L,(\nu\vert L)') \ar[r]^-{\mathrm{res}} &
        X'}
    \end{equation*}
    and the induced iso\m\ 
    \func{f_{\nu,L}}{C_X(V,\nu)}{C_{X'}(V,\nu_L')} defined by the
    commutative diagram
    \begin{equation*}
      \xymatrix@C=40pt{
        C_{ C_X(L,\nu\vert L)}(V,\overline{\nu}(L)) \ar[d]_-{\cong}
        \ar[r]^-{C_{f_{\nu\vert L}}} &
        C_{ C_{X'}(L,(\nu\vert L)')}(V,f_{\nu\vert
        L}\circ\overline{\nu}(L)) \ar[d]^{\cong} \\
        C_X(V,\nu) \ar[r]_-{f_{\nu,L}} &
        C_{X'}(V,{\nu}'_L) } 
    \end{equation*}
    do not depend on the choice of line $L<V$. (See \ref{rmk:canfact}
    for the definition of the canonical factorization
    $\overline{\nu}(L)$.) 
  \item \label{indstepalt5}
    ${\operatorname{lim}}^2({\A(X)};\pi_1(B{\Ze}C_X))=0=
    {\operatorname{lim}}^3({\A(X)};\pi_2(B{\Ze}C_X))$.
  \end{enumerate}
  Then there exists an iso\m\ \func{f}{X}{X'} under ${\T}$
  (\ref{underT}). 
\end{thm}

\begin{proof}
  The idea is that the iso\m s
  \func{f_{\lambda}}{C_X(\lambda)}{C_{X'}(\lambda')} on the line
  centralizers restrict to iso\m s
  \func{f_{\nu}}{C_X(\nu)}{C_{X'}(\nu')} for all centralizers in the
  $\F_2$-homology decomposition (\ref{hodecomp}) of $BX$. These
  locally defined iso\m s combine to a globally defined iso\m\ $BX \to
  BX'$.

  First observe that the iso\m s $f_{\lambda}$ on the line
  centralizers are uniquely determined by the cohomology class
  $\alpha\in H^1({\W};\ch{{\T}})$ (\ref{defcons}.(\ref{defcons1})).
  
  Let now $(V,\nu)$ be a rank two object of $\A(X)$ and $L$ a line in
  the plane $V$. If $(V,\nu)$ is {\em toral}, define
  \func{f_{\nu}}{C_{X}(V,\nu)}{C_{X'}(V,\nu')} to be the iso\m\ 
  induced by \func{f_{\nu\vert L}}{C_X(L,\nu\vert
    L)}{C_{X'}(L,(\nu\vert L)'}. Since $f_{\nu}$ is an iso\m\ under
  $\alpha \vert C_{\N}(V,\nu^{\N})$ it does not depend on the choice of $L$
  in $V$ (\ref{defcons}.(\ref{defcons1})). 
  If $(V,\nu)$ is {\em nontoral}, define
  $\nu'$ to be $\nu_L'$ and define
  \func{f_{\nu}}{C_{X}(V,\nu)}{C_{X'}(V,\nu')} to be $f_{\nu,L}$. By
  assumption \ref{indstepalt}.(\ref{indstepalt4}), the mono\m\ $\nu'$
  and the iso\m\ $f_{\nu,L}$ are independent of the choice of $L$.

  This construction respects \m s in $\A(X)$. Consider first, for
  instance, a \m\ \func{\beta}{(L_1,\lambda_1)}{(L_2,\lambda_2)}
  between two lines in $X$. Then $\lambda_1=\lambda_2\beta$ and
  $\lambda_1^{\N}=\lambda_2^{\N}\beta$. The commutative diagram of iso\m s
  \begin{equation*}
    \xymatrix{
      C_X(\lambda_1) \ar[ddd]_{f_{\lambda_1}} &&& 
      C_X(\lambda_2) \ar[ddd]^{f_{\lambda_2}} \ar[lll]_{C_X(\beta)} \\
      & C_{\N}(\lambda_1^{\N}) \ar[d]_{\alpha\vert
        C_{\N}(\lambda_1^{\N})} \ar[ul] &
      C_{\N}(\lambda_2^{\N}) 
      \ar[d]^{\alpha\vert C_{\N}(\lambda_2^{\N})} \ar[l]_{C_N(\beta)} 
        \ar[ur] \\
      & C_{\N}(\lambda_1^{\N}) \ar[dl] &
      C_{\N}(\lambda_2^{\N}) \ar[l]^{C_N(\beta)} \ar[dr] \\
      C_{X'}(\lambda_1') &&&
      C_{X'}(\lambda_2') \ar[lll]^{C_{X'}(\beta)} }
  \end{equation*}
  shows that $C_{X'}(\beta)^{-1} \circ f_{\lambda_1} \circ C_{X}(\beta) =
  f_{\lambda_2}$ for they are both iso\m\ under $C_N(\beta)^{-1}
  \circ \alpha\vert C_{\N}(\lambda_1^{\N}) \circ C_N(\beta) =
  \alpha\vert C_{\N}(\lambda_2^{\N})$. Second, by the very definition
  of $f_{\nu}$, the diagram
  \begin{equation*}
    \xymatrix{
      C_X(V,\nu) \ar[d] \ar[r]^-{f_{\nu}} &
      C_{X'}(V,\nu') \ar[d] \\
      C_X(L,\nu\vert L) \ar[r]_-{f_{\nu\vert L}} &
      C_{X'}(L,(\nu\vert L)') }
  \end{equation*}
  commutes whenever $L<V$ and $(V,\nu)$ is (toral or nontoral) rank
  $2$ object of $\A(X)$.

  We have now defined natural iso\m s
  \func{f_{\nu}}{C_X(V,\nu)}{C_{X'}(V,\nu')} for all objects
  $(V,\nu)\in\Ob(\A(X))$ of rank $\leq 2$. For any other object
  $(E,\varepsilon)$ of $\A(X)$, choose a line $L<E$ and proceed as
  for toral rank $2$ objects. That is, define
  \func{\varepsilon'}{E}{X'} to be the mono\m\
  \begin{equation*}
    \xymatrix@1@C=40pt{
      E \ar[r]^-{\overline{\varepsilon}(L)} & 
      C_X(E,\varepsilon\vert L) \ar[r]^-{f_{\varepsilon\vert L}} &
      C_{X'}(E,(\varepsilon\vert L)') \ar[r]^-{\mathrm{res}} &
      X' }
  \end{equation*}
  and define
  \func{f_{\varepsilon}}{C_X(E,\varepsilon)}{C_{X'}(E,\varepsilon')}
  to be the iso\m\
  \begin{equation*}
    \xymatrix@C=45pt{
      C_{C_X(E,\varepsilon\vert L)}(\overline{\varepsilon}(L))
      \ar[r]^-{(f_{\varepsilon\vert L})_*} \ar[d]_{\cong} &
      C_{C_{X'}(E,(\varepsilon\vert L)')}(f_{\varepsilon\vert L} \circ
      \overline{\varepsilon}(L)) \ar[d]^{\cong} \\
      C_X(E,\varepsilon) \ar[r]_{f_{\varepsilon}} &
      C_{X'}(E,\varepsilon') }
  \end{equation*}
  induced by $f_{\varepsilon\vert L}$. If $L_1$ and $L_2$ are two
  distinct lines in E, let $P=\gen{L_1,L_2}$ be the plane generated by
  them. Then the commutative diagram 
  \begin{equation*}
  \xymatrix@C=45pt{
    & C_X(L_1,\varepsilon\vert L_1) 
    \ar[r]^-{f_{\varepsilon\vert L_1}}_-{\cong} &
      C_{X'}(L_1,(\varepsilon\vert L_1)') \ar[dr]^-{\mathrm{res}} \\
    P \ar[ur]^-{\overline{\varepsilon}(L_1)} 
    \ar[r]^(.5){\overline{\varepsilon}(P)}
      \ar[dr]_-{\overline{\varepsilon}(L_2)} &
      C_X(P,\varepsilon\vert P) 
    \ar[r]_-{f_{\varepsilon\vert P}}^-{\cong} \ar[u] \ar[d] &
      C_{X'}(P,(\varepsilon\vert P)') \ar[u] \ar[d]
      \ar[r]^(.5){\mathrm{res}} & X' \\
    & C_X(L_2,\varepsilon\vert L_2) 
   \ar[r]_-{f_{\varepsilon\vert L_2}}^-{\cong} &
      C_{X'}(L_2,(\varepsilon\vert L_2)')\ar[ur]_-{\mathrm{res}} }
  \end{equation*}
  shows that neither $(E,\varepsilon')\in\Ob(\A(X'))$ nor the iso\m\ 
  $f_{\varepsilon}$ depend on the choice of line in $E$.  Thus we have
  constructed a collection of centric \cite{dwk:centric} maps
  \set\begin{equation}\label{eq:collection} \B C_X(V,\nu) \to \B X',
    \quad (V,\nu) \in \Ob(\A(X)),
\end{equation}\add
that are homotopy invariant under $\A(X)$-\m s.  The vanishing
(\ref{indstepalt}.(\ref{indstepalt5})) of the
obstruction groups means \cite{w:obs} that these homotopy
$\A(X)$-invariant maps can be realized by a map
  \begin{equation*}
    \B f \colon \B X \xleftarrow{\simeq} \hocolim \B C_X \to \B X'
  \end{equation*}
  such that $f \circ \mathrm{res} = \mathrm{res} \circ f_{\nu}$ for
  all $(V,\nu)\in\Ob(\A(X))$. In particular, $f$ is a map under ${\T}$
  and an iso\m\ (\ref{underT}).
\end{proof}

\setcounter{subsection}{\value{thm}}
\subsection{Verification of condition
  \ref{indstepalt}.(\ref{indstepalt1})} 
\label{sec:cond2}\add
Let $\A(X)^{\leq t}$ be the toral part of the Quillen category and let
\func{H^1({\W}_0;\ch{{\T}})^{{\W}/{\W}_0}}{\A(X)^{\leq t}}{\Ab} 
the functor with
%%A functor because the LHSss is natural (because the Serre ss is)
value $H^1({\W}(C_X(V,\nu)_0);\ch{{\T}})^{\pi_0C_X(V,\nu)}$ on the object
$(V,\nu)$. If the \twocg\ $C$ satisfies the conditions of
Lemma~\ref{redtoconnected} and $\ch{{\Ze}}(C_0)=\ch{{\T}}(C_0)^{{\W}(C_0)}$ we
say that $C$ satisfies the the conditions of
Lemma~\ref{redtoconnected} {\em in the  strong sense}.
\begin{lemma}\label{lemma:altcond2}
  Suppose that
  \begin{itemize}
  \item  The centralizers $C_X(V,\nu)$ of all
    $(V,\nu) \in \Ob(\A(X)^{\leq t}_{\leq 2})$ satisfy the conditions
    of Lemma~\ref{redtoconnected} in the strong sense, 
  \item $H^1({\W};\ch{{\T}}) \to \lim^0(\A(X)^{\leq t}_{\leq 2};
    H^1({\W}_0;\ch{{\T}})^{{\W}/{\W}_0})$ is surjective
  \end{itemize}
  Then conditions \ref{indstepalt}.(\ref{indstepalt3}) and
  \ref{indstepalt}.(\ref{indstepalt1}) are satisfied.
\end{lemma}
\begin{proof}
  Let $(V,\nu)$ be an object of $\A(X)^{\leq t}$ of rank $\leq 2$.
  Since (\ref{redtoconnected}) the centralizer $C_X(V,\nu)$ is
  ${\N}$-determined there is a solution $(f(V,\nu),\alpha(V,\nu))$ to
  the iso\m\ problem
  \begin{equation*}
    \xymatrix{
      C_{\N}(V,\nu^{\N}) \ar[d] \ar[r]^{\alpha(V,\nu)} &  
      C_{\N}(V,\nu^{\N}) \ar[d] \\
      C_X(V,\nu) \ar[r]_{f(V,\nu)} &  C_{X'}(V,\nu') }
  \end{equation*}
  and the set of all solutions is (\ref{lemma:autononcon},
  \ref{idcomp}) a $H^1({\W}/{\W}_0;\ch{{\T}}^{{\W}_0})(C_X(V,\nu))$-coset. Let
  \begin{equation*}
    \overline{\alpha}(V,\nu) \in H^1({\W}_0;\ch{{\T}}^{{\W}/{\W}_0})(C_X(V,\nu))
  \end{equation*}
  be the restriction of any solution $\alpha(V,\nu) \in
  H^1({\W};\ch{{\T}})(C_X(V,\nu))$ to the above iso\m\ problem. Then
  \set\begin{equation}\label{limalpha}
    \{ \overline{\alpha}(V,\nu) \}_{(V,\nu) \in \Ob(\A(X)^{\leq
      t}_{\leq 2}} \in 
     \makebox{$\lim^0(\A(X)^{\leq t}_{\leq 2};
    H^1({\W}_0;\ch{{\T}})^{{\W}/{\W}_0})$}
  \end{equation}\add
  because the restriction of a solution is a solution. By assumption,
  there is an element $\alpha \in H^1({\W};\ch{{\T}})$ that maps to
  (\ref{limalpha}) and $\alpha$ satisfies
  \ref{indstepalt}.(\ref{indstepalt1}).
\end{proof}

In case $H^1({\W};\ch{{\T}})=0$, the second point reduces to
$\lim^0(\A(X)^{\leq t}_{\leq 2};
H^1({\W}_0;\ch{{\T}})^{{\W}/{\W}_0})=0$. Alternatively, if 
$\lim^1(\A(X)^{\leq t}_{\leq 2}; H^1({\W}/{\W}_0;\ch{{\T}}^{{\W}_0}))=0$, then the
\ses s (\ref{LHSinit}) for $C_X(V,\nu)$, $(V,\nu) \in
\Ob(\A(X)^{\leq t}_{\leq 2})$, will produce a \ses\
\begin{multline*}
   \mbox{$0 \to \lim^0\big(\A(X)^{\leq t}_{\leq 2},
     H^1({\W}/{\W}_0;\ch{{\T}}^{{\W}_0})\big) 
  \to \lim^0\big(\A(X)^{\leq t}_{\leq 2}, H^1({\W};\ch{{\T}})\big)$} \\
   \mbox{$
  \to \lim^0\big(\A(X)^{\leq t}_{\leq 2}, H^1({\W}_0;\ch{{\T}})^{{\W}/{\W}_0}\big)
  \to 0$}, 
\end{multline*}
in the limit. Since $H^1({\W};\ch{{\T}})$ is isomorphic to the middle
term by \cite[8.1]{dw:new}, it  maps onto the third term.

\setcounter{subsection}{\value{thm}}
\subsection{Verification of condition
  \ref{indstepalt}.(\ref{indstepalt4})} 
\label{sec:cond4}\add
In this subsection we assume that conditions
\ref{indstepalt}.(\ref{indstepalt3}) and
\ref{indstepalt}.(\ref{indstepalt1}) are satisfied.  The following
observations can sometimes be useful in the verification of condition
\ref{indstepalt}.(\ref{indstepalt4}).

Let $(V,\nu)$ be a nontoral rank two object of $\A(X)$ and $L<V$ a
rank one subgroup. The commtutative diagram
\set\begin{equation}\label{dia:altdefnuL}
    \xymatrix@C=40pt@R=20pt{
      &
      N \ar[r]^{\alpha} &
      N \ar@/^/[dr]^{j'} \\
      V \ar@/^/[ur]^{\nu_L^N} \ar[r]^(.45){\overline{\nu_L^N}(L)}
      \ar@/_/[dr]_{\overline{\nu}(L)} &
      C_N(L,\nu_L^N \vert L) \ar[d] \ar[r]^{C_{\alpha}} \ar[u] &
      C_N(L,\nu_L^N \vert L) \ar[d] \ar[u] &
      X' \\
      &
      C_X(L,\nu\vert L) \ar[r]_-{f_{\nu\vert L}} &
      C_{X'}(L,(\nu\vert L)') \ar@/_/[ur]_-{\mathrm{res}} }
  \end{equation}\add
shows that $\nu_L'$, which is defined to be $\mathrm{res} \circ
f_{\nu\vert L} \circ \overline{\nu}(L)$, is equal to the composite
$\nu_L' = j' \circ \alpha \circ \nu_L^N$. Moreover, we see by taking
the centralizer of $\overline{\nu}(L)$ that
\set\begin{equation}\label{dia:undernuL}
\xymatrix{ & {V} \ar[dl]_-{\overline{\nu}(V)}
    \ar[dr]^-{\overline{\nu}_L'(V)} \\
    C_X(V,\nu) \ar[rr]_-{f_{\nu,L}}^-{\cong} &&
    C_{X'}(V,\nu') }
\end{equation}\add
commutes.
%%%%%%%%%%%%%%%%
%We are in particular assuming the existence of iso\m s
%\func{f_{\lambda}}{C_X(L,\lambda)}{C_{X'}(L,\lambda')} as in
%\ref{indstepalt}.(\ref{indstepalt1}) for all rank one objects
%$(L,\lambda)$ of $\A(X)$. Note that $f_{\lambda}$ is an iso\m\ under
%$L$ in the sense that
%\begin{equation*}
%  \xymatrix{
%     & L
%     \ar[dl]_{\overline{\lambda}(L)}\ar[dr]^{\overline{\lambda}'(L)}
%     \\
%     C_X(L,\lambda) \ar[rr]^{\cong}_{f_{\lambda}} &&
%     C_{X'}(L,\lambda') }
%\end{equation*}
%commutes. The induced diagram
%\set\begin{equation}\label{dia:undernuL}
%\xymatrix{ & {V} \ar[dl]_-{\overline{\nu}(V)}
%    \ar[dr]^-{\overline{\nu}_L'(V)} \\
%    C_X(V,\nu) \ar[rr]_-{f_{\nu,L}}^-{\cong} &&
%    C_{X'}(V,\nu') }
%\end{equation}\add
%will then also commute for any nontoral rank two object $(V,\nu)$ of
%$\A(X)$. Moreover, for any line $L < V$, the commutative diagram
%\set\begin{equation}\label{dia:altdefnuL}
%    \xymatrix@C=40pt@R=20pt{
%      &
%      N \ar[r]^{\alpha} &
%      N \ar@/^/[dr]^{j'} \\
%      V \ar@/^/[ur]^{\nu_L^N} \ar[r]^(.45){\overline{\nu_L^N}(L)}
%      \ar@/_/[dr]_{\overline{\nu}(L)} &
%      C_N(L,\nu_L^N \vert L) \ar[d] \ar[r]^{C_{\alpha}} \ar[u] &
%      C_N(L,\nu_L^N \vert L) \ar[d] \ar[u] &
%      X' \\
%      &
%      C_X(L,\nu\vert L) \ar[r]_-{f_{\nu\vert L}} &
%      C_{X'}(L,(\nu\vert L)') \ar@/_/[ur]_-{\mathrm{res}} }
%  \end{equation}\add
%  shows that $\nu_L'= \mathrm{res} \circ f_{\nu\vert L} \circ
%  \overline{\nu}(L) = j' \circ \alpha \circ \nu_L^N$. 
%%%%%%%%%%%%%%%%%%%%%%%%%%%%%%%

  We are looking for criteria that ensure that \func{\nu'_L}{V}{X'} is
  independent of the choice of $L < V$.

\begin{lemma}\label{lemma:C3cond}
  Let $(V,\nu)$ be a nontoral rank two object of $\A(X)$ and $L<V$ a
  line in $V$. Write $C_3$ for the Sylow $3$-subgroup of $\GL(V)$.
  Suppose that
  \begin{enumerate}
  \item $C_3 \subseteq \A(X)(V,\nu) \cap  \A(X')(V,\nu_L')$
    \label{lemma:C3cond.1} 
  \item \func{f_{\nu,L}}{C_X(V,\nu)}{C_{X'}(V,\nu_L')} is
    $C_3$-equivariant \label{lemma:C3cond2}
  \end{enumerate}
  Then  condition \ref{indstepalt}.(\ref{indstepalt4}) is satisfied.
\end{lemma}
\begin{proof}
  Let $\beta$ be an auto\m\ of $V$. For general reasons, $\nu_L^N\beta
  = (\nu\beta)^N_{\beta^{-1}L}$ and the
  diagram
  \begin{equation*}
    \xymatrix@C=40pt{
      C_X(V,\nu) \ar[d]_{C_X(\beta)}^{\cong} \ar[r]^-{f_{\nu,L}} &
      C_{X'}(V,\nu'_L) \ar[d]^{C_{X'}(\beta)}_{\cong} \\
      C_X(V,\nu\beta)  \ar[r]_-{f_{\nu\beta,\beta^{-1}L}} &
      C_{X'}(V,\nu'_L\beta) }
  \end{equation*}
  commutes. Now, if $\beta \in \A(X)(V,\nu) \cap \A(X')(V,\nu'_L)$, then
  $\nu\beta=\nu$, $\nu'_L\beta=\nu'_L$, and $f_{\nu\beta,\beta^{-1}L} =
  f_{\nu,\beta^{-1}L}$ so that $f_{\nu,\beta^{-1}L} = C_{X'}(\beta) \circ
  f_{\nu,L} \circ C_X(\beta)^{-1}$ according to the above diagram. If
  also $f_{\nu,L}$ commutes with the action of $\beta$, we conclude
  that $f_{\nu,L} = f_{\nu,\beta^{-1}L}$.
  %%%%%%%%%%%%%%%%%%%%%%%%%%%%%%%
 % Let $L_1$ and $L_2=L$ be lines in $V$. Choose an auto\m\ $\alpha \in
%  C_3$ of $(V,\nu)$ that takes $L_1$ to $L_2$. Then
%  $\nu'_{L_2}\alpha=\nu'_{L_1}$ and $f_{\nu,L_1}=C_X(\alpha)\circ
%  f_{\nu,L_2} \circ C_X(\alpha)^{-1}$ (\ref{rmk:canfact}).
\end{proof}

The following lemma assures that condition
\ref{lemma:C3cond}.(\ref{lemma:C3cond.1}) holds.

\begin{lemma}\label{help1}
  Let $L$ and $V$ denote \lmntwo s of rank one and two, respectively.
  Suppose that
  \begin{enumerate}
  \item There is (up to conjugacy) a unique mono\m\
    \func{\lambda}{L}{X} with nonconnected centralizer
  \item  There is (up to conjugacy) a unique nontoral mono\m\
    \func{\nu}{V}{X}
  \end{enumerate}
  Then the same holds for $X'$, and
  $\A(X)(V,\nu)=\GL(V)=\A(X')(V,\nu')$ for the unique nontoral rank
  two objects $(V,\nu )$ of $\A(X)$ and $(V,\nu')$ of $\A(X')$.
\end{lemma}
\begin{proof}
  Let \func{\nu'}{V}{X'} be a nontoral mono\m\ and \func{i}{L}{V} an
  inclusion. Then $(L,\nu' i)=(L,\lambda')$ for $C_{X'}(L,\nu'i)$ is
  nonconnected so that $\nu'i$ and $\lambda$ must correspond under the
  bijection $\A(X)^{\leq t} \to \A(X')^{\leq t}$ between toral
  categories. 
  Moreover, the diagram \set\begin{equation}\label{dia:D1} \xymatrix{
      X & C_X(L,\lambda) \ar[l]_-{\mathrm{res}}
      \ar[rr]^-{f_{\lambda}}_-{\cong} && C_{X'}(L,\lambda')
      \ar[r]^-{\mathrm{res}} &
      X' \\
      && {V} \ar@/^/[ull]^-{\nu} \ar[ul]_{\overline{\nu}(L)}
      \ar[ur]^{\overline{\nu}'(L)} \ar@/_/[urr]_-{\nu'} }
  \end{equation}\add
  is commutative. To see this, observe that $(V,\mathrm{res} \circ
  f_{\lambda}^{-1} \circ \overline{\nu}'(L))$ is a nontoral rank two
  object of $\A(X)$ (its centralizer is isomorphic to
  $C_{C_{X'}(L,\lambda')}(V,\overline{\nu}'(L))=C_{X'}(V,\nu')$) so
  that $(V,\nu)=(V,\mathrm{res} \circ
  f_{\lambda}^{-1} \circ \overline{\nu}'(L))$ by uniqueness of
  $(V,\nu)$. Also, we see from the commutative diagram
  \begin{equation*}
    \xymatrix@C=35pt{
      && L \ar@/_1pc/[dll]_{\lambda} 
               \ar[dl]_{\overline{\lambda}(L)} \ar[d]^i \\
      X & C_X(L,\lambda) \ar[l]^(.6){\mathrm{res}} &
      V \ar@<-1ex>[l]_-{\overline{\nu}(L)}
      \ar@<1ex>[l]^-{f_{\lambda}^{\-1} \circ \overline{\nu}'(L)} }
  \end{equation*}
  that $\overline{\nu}(L)= f_{\lambda}^{-1} \circ \overline{\nu}'(L)$
  by uniqueness of canonical factorizations under $L$
  \cite[3.9]{jmm:trep}. We conclude that $\nu' = \mathrm{res} \circ
  \overline{\nu}'(L) = \mathrm{res} \circ f_{\lambda} \circ
  \overline{\nu}(L)$. This means (\ref{dia:altdefnuL}) that
  $\nu'=\nu'_L$ for any choice of line $L < V$. Since thus $\nu'$ is
  unique up to conjugacy, $\nu'\beta=\nu'$ for any auto\m\ $\beta$ of
  $V$. 
\end{proof}

Note in connection with the verification of condition
\ref{lemma:C3cond}.(\ref{lemma:C3cond2}), that if
\ref{lemma:C3cond}.(\ref{lemma:C3cond.1}) is satisfied so that
$\nu'_L=\nu'$ is independent of $L$, then
(\ref{dia:undernuL}) shows that $f_{\nu,L}$ is a map under $V$ in the sense
that \set\begin{equation}\label{dia:triangle} \xymatrix{ & {V}
    \ar[dl]_-{\overline{\nu}(V)}
    \ar[dr]^-{\overline{\nu}'(V)} \\
    C_X(V,\nu) \ar[rr]_-{f_{\nu,L}}^-{\cong} && C_{X'}(V,\nu') }
\end{equation}\add
commutes. Since the canonical mono\m s, ${\overline{\nu}(V)}$ and
${\overline{\nu}'(V)}$, are $\GL(V)$-equivariant, the restriction of
$f_{\nu,L}$ to $V$ is $C_3$-equivariant.

%%%%%%%%%%%%%from dfam 260202 %%%%%%%%%%%%%%

For any nontoral object (not necessarily of rank two) $(V,\nu)$ of
$\A(X)$ and any rank one subgroup $L \subset V$, let
\func{\nu_L^N}{V}{N} be a preferred lift of $\nu$ such that $\nu_L^N
\vert L$ is the \pl\ of $\nu\vert L$, ie $\nu_L^N\vert L = (\nu \vert
L)^N$. (It is always possible to extend a \pl\ given on the subgroup
$L$ to a \pl\ defined on all of $V$ but a \pl\ defined on $V$ may not
restrict to a \pl\ on $L$ \cite[4.9]{jmm:normax}.)  Also, define
\func{\nu'_L}{V}{X'} and \func{f_{\nu,
    L}}{C_X(V,\nu)}{C_{X'}(V,\nu'_L)} as in
\ref{indstepalt}.(\ref{indstepalt4}).

\begin{lemma}\label{lemma:uniquenu}
  %Suppose that the first two conditions,
  %\ref{indstepalt}.(\ref{indstepalt3}) and
  %\ref{indstepalt}.(\ref{indstepalt1}), hold.  
  Let \func{\nu}{V}{X} be
  any nontoral  object of $\A(X)$.
  \begin{enumerate}
  \item If the centralizer of $\nu$ has a nontrivial identity
    component, then \func{\nu_L'}{V}{X'} is independent up to conjugacy
    of the choice of $L \subset V$, and $\nu'_L=j' \circ \alpha \circ
    \nu_L^N$. \label{lemma:uniquenu1}
  \item \label{lemma:uniquenu2} If also there exist a \twoct\ 
    $T_{\nu}$ and iso\m s $T_{\nu} \to C_N(V,\nu_L^N)_0$ such that the
    composites $T_{\nu} \to C_N(V,\nu_L^N)_0 \to T$ are independent up
    to conjugacy of $L<V$,
  % diagrams
%    \begin{equation*}
%      \xymatrix{
%        T_{\nu} \ar[r] \ar[d] & T \ar[d]\\
%        C_N(\nu_L^N) \ar[r] & N}
%    \end{equation*}
%  commute, 
    then \func{f_{\nu,L}}{C_X(V,\nu)}{C_{X'}(V,\nu')} are iso\m s under
    the \mt\ $T_{\nu}$ for all $L<V$.
  \end{enumerate}
\end{lemma}
\begin{proof}
  (1) Just as in (\ref{dia:altdefnuL}) we see that $\nu_L'=
  \mathrm{res} \circ f_{\nu\vert L} \circ \overline{\nu}(L) = j' \circ
  \alpha \circ \nu_L^N$.
  The hypothesis implies that there exists \cite[5.4, 7.3]{dw:fixpt} a \m\ 
  \func{\phi}{L_1 \times V}{X} extending \func{\nu}{V}{X} whose
  adjoint ${L_1} \to {C_X(\nu)}$ factors through the identity
  component of $C_X(\nu)$. Let $L_1 \to C_N(V,\nu_L^N)$ be the \pl\ of
  $L_1 \to C_X(V,\nu)$ as in the commutative diagram
  \begin{equation*}
    \xymatrix{
      &
      C_N(V,\nu_L^N) \ar[r]^-{\mathrm{res}} \ar[d] &
      N \ar[d]^j \\
      L_1 \ar@/^/[ur] \ar[r] &
      C_X(V,\nu) \ar[r]_-{\mathrm{res}} &
      X }
  \end{equation*}
  This \pl\ will factor through the identity component of
  $C_N(\nu_L^N)$ (and hence its composition with $C_N(\nu_L^N) \to N$
  will factor through the identity component of $N$) since $L_1 \to
  C_X(V,\nu)$ factors through the identity component of $C_X(\nu)$
  \cite[4.10]{jmm:normax}. Let \func{\phi_L^N}{L_1 \times V}{N} be the
  adjoint of the \pl\ $L_1 \to C_N(\nu_L^N)$.  Then $\phi_L^N \vert
  L_1 \colon L_1 \to N$ factors through the identity component of $N$
  (the \mt ) so it is \cite[4.10]{jmm:normax} the \pl\ of  $\phi \vert L_1
  \colon L_1 \to X$.  In particular, $\phi_L^N
  \vert L_1 = (\phi \vert L_1)^N$ does not depend on the choice of
  $L$.
  
  The adjoints, \func{\phi_2^N}{V}{C_N(\phi_L^N\vert L_1)} and
  \func{\phi_2}{V}{C_X(\phi \vert L_1)}, of $\phi_L^N$ and $\phi$,
  respectively, with respect to the second factor, give a commutative
  diagram
  \begin{equation*}
    \xymatrix@C=40pt@R=20pt{
      &
      C_X(L_1,\phi \vert L_1) \ar[r]^-{f_{\phi \vert L_1}} &
      C_{X'}(L_1,(\phi \vert L_1)') \ar@/^/[dr]^-{\mathrm{res}} \\
      V \ar@/^/[ur]^-{\phi_2} \ar[r]^(.45){\phi_2^N} \ar@/_/[dr]_{\nu_L^N} &
      C_N(L_1,\phi_L^N \vert L_1) \ar[d] \ar[r]^-{C_{\alpha}} \ar[u] &
      C_N(L_1,\phi_L^N \vert L_1) \ar[d] \ar[u] &
      X' \\
      &
      N \ar[r]_-{\alpha} & 
      N \ar@/_/[ur]_-{j'} }
  \end{equation*}
  We conclude that
  $\nu_L' = j' \circ \alpha \circ \nu_L^N = \mathrm{res} \circ
  f_{\phi \vert L_1} \circ \phi_2 \colon V \to X'$ is independent of the
  choice of $L < V$.
  
  \noindent
  (2) The upper square in the diagram
  \begin{equation*}
    \xymatrix{
      T_{\nu} \ar@{=}[r] \ar[d] & T_{\nu} \ar[d] \\
      C_N(V,\nu_L^N) \ar[d] \ar[r]^-{\alpha} & 
      C_N(V,\alpha\nu_L^N) \ar[d] \\
      C_X(V,\nu) \ar[r]_-{f_{\nu,L}} &  C_{X'}(V,\nu') }
  \end{equation*}
  commutes because $\alpha$ restricts to the identity on the identity
  component $T$ of $N$ and hence also on $T_{\nu}$. That the lower
  square is commutative is consequence of the commutative diagram
  \begin{equation*}
    \xymatrix{
      & C_N(L,\nu_L^N \vert L) \ar[d] \ar[r]^{C_{\alpha}} &
      C_N(L,\alpha\nu_L^N) \ar[d]  \\
      V \ar@/^/[ur]^-{\overline{\nu}_L^N(L)} \ar[r]_-{\overline{\nu}(L)} &
      C_X(L,\nu\vert L) \ar[r]_-{f_{\nu\vert L}} &
      C_{X'}(L,(\nu\vert L)') }
  \end{equation*}
  where $\overline{\nu}_L^N(L)$ and $\overline{\nu}(L)$ are 
    the canonical factorizations (\ref{rmk:canfact}). 
   %(Well, this is
   % part of previuos diagram so maybe no need to write it again?)
\end{proof}

Let \func{\mu}{U}{X} be a nontrivial \lmntwo\ and \func{\mu}{U}{X} a
mono\m\ whose centralizer $C_X(U,\mu)$ has nontrivial identity
component. Suppose that $U$ contains a nontrivial subgroup $V < U$
such that the restriction of $\mu$ to $V$ is nontoral.  Choose a rank
one subgroup $L \subset V \subset U$. We may choose the \pl s
$\mu_L^N$ and $(\mu\vert V)_L^N$ such that $\mu_L^N \vert V =
(\mu\vert V)_L^N$.  Since $C_X(U,\mu)$ has nontrivial identity
component, the conjugacy classes of the mono\m s $\mu'=\mu'_L$ and
$(\mu \vert V)'_L=\mu' \vert L$ are independent of the choice of $L$
by \ref{lemma:uniquenu}.(\ref{lemma:uniquenu1}).  Then there is a
commutative diagram \set\begin{equation}\label{dia:Utriangle}
  \xymatrix{
    & U \ar[dl]_{\overline{\mu}(V)} \ar[dr]^{\overline{\mu}'(V)} \\
    C_X(V, \mu\vert V) \ar[rr]_{f_{\mu \vert V,L }} && C_{X'}(V,
    \mu'\vert V) }
\end{equation}\add
similar to (\ref{dia:triangle}). 
%%%%%%%%%%%%%%%%%%%%%%%%%%%%%%%%%%%%%%
%\begin{equation*}
%  \xymatrix{
%    & U \ar[dl]_{\overline{\mu}(U)} \ar[dr]^{\overline{\mu'}(U)} \\
%    C_X(U,\mu) \ar[rr]_{f_{\mu,L }} \ar[d] && C_{X'}(U,\mu') \ar[d] \\
%    C_X(V, \mu\vert V)  \ar[rr]_{f_{\mu \vert V,L }} 
%    &&  C_{X'}(V, \mu'\vert V) }
%\end{equation*}
%where the slanted arrows are canonical lifts. The square commutes
%because $\overline{\mu \vert V}(L)=\overline{\mu}(L) \vert V$ by
%(\ref{dia:canfact1}). 
%The induced commutative
%diagram on the level of component groups
%\set\begin{equation}\label{dia:componentgr}
%  \xymatrix{
%  & U  \ar[dl]_{\pi_0(\overline{\mu}(V))} 
%         \ar[dr]^{\pi_0(\overline{\mu'}(V))} \\
%   \pi_0(C_X(V,\mu\vert V))  \ar[rr]_{\pi_0(f_{\mu\vert V,L })}^{\cong} &&
%   \pi_0(C_{X'}(V,\mu' \vert V)) } 
%\end{equation}\add
%shows that $\pi_0(f_{\mu \vert V,L})$ is an iso\m\ under $U$.

%%%%%%%%%%%%%%%%from dfam 260202%%%%%%%%%%%%%%%%%%%%%%%%%

\setcounter{subsection}{\value{thm}}
\subsection{Canonical factorizations} 
\label{rmk:canfact}\add
  Let \func{\nu}{V}{X} be a mono\m\ from an \lmn\ to the \pcg\ $X$.
  The canonical factorization of $\nu$ through its centralizer is the
  central mono\m\ \func{\overline{\nu}(V)}{V}{C_X(V,\nu)} whose
  adjoint is $V \times V \xrightarrow{+} V \xrightarrow{\nu} X$
  \cite[8.2]{dw:fixpt}. If \func{\alpha}{(V_1,\nu_1)}{(V_2,\nu_2)} is
  a \m\ in $\A(X)$ then the canonical factorizations are related by a
  commutative diagram \set\begin{equation}\label{dia:canfact1}
    \xymatrix@C=40pt{ V_1 \ar[d]_-{\alpha}
      \ar[r]^-{\overline{\nu}_1(V_1)} & C_X(V_1,\nu_1)
      \ar[r]^-{\mathrm{res}} &
      X \\
      V_2 \ar[r]_-{\overline{\nu}_2(V_2)} & C_X(V_2,\nu_2)
      \ar[u]^-{C_X(\alpha)} \ar[r]_-{\mathrm{res}} & X \ar@{=}[u] }
  \end{equation}\add
  and we shall write $\overline{\nu}_2(V_1) \colon V_2 \to
  C_X(V_1,\nu_1)$ for $C_X(\alpha)\circ \overline{\nu}_2(V_2)$ and
  call it the canonical factorization of $\nu_2$ through the
  centralizer of $\nu_1$. The induced diagram
  \set\begin{equation}\label{dia:canfact2}
    \xymatrix@C=45pt{
      C_{C_X(V_2,\nu_2)}(V_2,\overline{\nu}_2(V_2))
      \ar[r]_-{\cong}^-{C_{C_X(\alpha)}} \ar[d]_{\cong} &
      C_{C_X(V_1,\nu_1)}(V_2,\overline{\nu}_2(V_1))
      \ar[r]^-{C_{C_X(V_1,\nu_1)}(\alpha)} &
      C_{C_X(V_1,\nu_1)}(V_1,\overline{\nu}_1(V_1)) \ar[d]^{\cong} \\
      C_X(V_2,\nu_2) \ar[rr]_-{C_X(\alpha)} &&  C_X(V_1,\nu_1)  } 
  \end{equation}\add
  is a factorization of $C_X(\alpha)$.

\section{An exact functor}
\label{sec:limAWt}

Let $W$ be a finite group, $p$ a prime, and \func{\rho}{W}{\GL(t)} a
representation of $W$ in an $\F_p$-vector space $t$ of finite
dimension.  For any nontrivial subgroup $V \subset t$, let
\begin{equation*}
  {\W}(V)=\{ w\in {\W} \mid \forall v \in V \colon wv=v\}
\end{equation*}
be the subgroup of elements of $W$ that act as the identity on $V$.
For any two nontrivial subgroups $V_1,V_2 \subset t$, let
\begin{equation*}
 \overline{{\W}}(V_1,V_2) =\{ w \in {\W} \vert wV_1 \subset
  V_2\} 
\end{equation*}
be the transporter set. (Even though suppresed in the notation, these
set depend on the representation $\rho$.)

Suppose that we are given also a $\Z_pW$-module $L$. 

\begin{defn}\cite[2.2]{jmm:ndet}\label{defn:AWt}
  $\A(\rho,t)$ is the category whose objects are nontrivial subspaces
  of $V$ and whose \m s are group homo\m s induced by the
  ${\W}$-action. The functor \func{L_i}{\A(\rho,t)}{\Ab} is the
  functor that takes the object $V \subset t$ to $H^i(W(V);L)$ and the
  \m\ \func{w}{V_1}{V_2} to $\xymatrix@1{ H^i(W(V_1);L) \ar[r]^-{w^*}
    & H^i(W(V_1)^w;L) \ar[r]^-{\mathrm{res}}& H^i(W(V_2);L) }$ where
  $\mathrm{res}$ is restriction and $w^*$ induced from conjugation
  with $w \in W$. 
\end{defn}

The category $\A(\rho,t)$ depends only on the image of $W$ in $\GL(t)$
but the functor $L_i$ depends on the actual representation. The \m\
set in   $\A(\rho,t)$ is the set of orbits
\begin{equation*}
  \A(\rho,t)(V_1,V_2) =  \overline{{\W}}(V_1,V_2)/W(V_1)
\end{equation*}
for the action of the group $W(V_1)$ on the set $
\overline{{\W}}(V_1,V_2)$.  We shall often write $\A(W,t)$ for
$\A(\rho,t)$ when the representation $\rho$ is clear from the context
and $\A(\W,t)(V)$ will be used as an abbreviation for the endo\m\ 
group $\A(\W,t)(V,V) = \overline{W}(V,V)/W(V)$.

\begin{lemma}\label{dw:limits}\cite[8.1]{dw:new}
  $L_i$ is an exact functor with limit $H^i(W;L)$:
  \begin{equation*}
    \lim^j(\A(W,t),L_i) = 
    \begin{cases}
      H^i(W;L) & j=0 \\
      0        & j>0
    \end{cases}
  \end{equation*}
\end{lemma}
\begin{proof}
  The proof of \cite[8.1]{dw:new} also applies to this slightly
  different setting where the action of $W$ on the $\F_p$-vector space
  $t$ may not be faithful and $L$ is a $\Z_pW$-module (and not an
  $\F_pW$-module).

  Another possibility is to use the ideas of \cite{jm}.  It suffices
  to show that the category $\A(W,t)$ satisfies (the dual of) the
  conditions of \cite[5.16]{jm} and that $L_*$ is a proto-Mackey
  functor. Define $L^* \colon \A(W,t) \to \Ab$ to be the contravariant
  functor that agrees with $L_*$ on objects but takes the $\A(W,t)$-\m\
  \func{w}{E_0}{E_1} to the group homo\m\
  \begin{equation*}
    \xymatrix@1@C=40pt{
    H^*(W(E_0);L) &  
    H^*(W(E_0)^w;L) \ar[l]_-{(w^{-1})^*} &
    H^*(W(E_1);L) \ar[l]_-{\mathrm{tr}} }
  \end{equation*}
  where $\mathrm{tr}$ is transfer.
  To prove the existence of coproducts and push-outs in the
  multiplicative extension $\A(W,t)_{\prod}$ we follow \cite[6.3]{jm}.
  Let $E_0,E_1,E_2$ be elementary abelian subgroups of $t$ where
  $E_0\subset E_1$ and there is a \m\ $E_0 \rightarrow E_2$
  represented by an element $w\in\overline{W}(E_0,E_2)\subset W$.
  ($E_0$ is possibly empty to allow for the construction of
  coproducts.) Each coset $gW(E_1) \in W(E_0)/W(E_1)$ has an
  associated special diagram
  \begin{equation*}
    {\xymatrix{ E_0 \ar[d]_w \ar@{^{(}->}[r] & E_1\ar[d]^{wg} \\
                E_2 \ar@{^{(}->}[r] & E_2+wgE_1 }}
  \end{equation*}
  where we note that $W(E_2+wgE_1)=W(E_2) \cap W(E_1)^{wg}$. This
  construction determines a bijection between the double coset
  $w^{-1}W(E_2)w\backslash W(E_0) / W(E_1)$ and the set of iso\m\
  classes of special diagrams, cf.\ \cite[7.3]{jm}, and therefore 
  \begin{equation*}
    {\xymatrix{ E_0 \ar[d]_w \ar@{^{(}->}[r] & E_1\ar[d]^{\prod wg} \\
                 E_2 \ar@{^{(}->}[r] & 
          {\prod}(E_2+wgE_1) }}
  \end{equation*}
  where the product is taken over all $g \in w^{-1}W(E_2)w\backslash
  W(E_0) / W(E_1)$, is a push-out diagram in $\A(W,t)_{\prod}$
  \cite[6.3]{jm}. By \cite[5.13]{jm}, we need  to show that the diagram
  \begin{equation*}
    {\xymatrix@C=90pt{H^*(W(E_0);L) \ar[d]_{L_*(w)} &
               H^*(W(E_1);L) \ar[d]^{\prod L_*(wg)}
               \ar[l]_{L^*(E_0\subseteq E_1)} \\
               H^*(W(E_2);L) & 
               {\prod H^*(W(E_2) \cap W(E_1)^{wg};L)}
               \ar[l]^-{\sum L^*(E_2\subseteq E_2+gwE_1)} }}
  \end{equation*}
  commutes. But this is precisely the content of the Cartan--Eilenberg
  double coset formula relating the restriction and transfer homo\m s
  in group cohomology \cite{cartaneilenberg}
  \cite[4.2.6]{evens:groupcohomology}.
  
  The restriction homo\m\ $H^*(W;L) \rightarrow \lim^0(\A(W,t);L_*)$
  is injective since $t$ contains an elementary abelian subgroup
  $E\subset t$ such that the index of $W(E)$ in $W$ is prime to $p$.
  To show surjectivity, we use the argument from the proof of
  \cite[7.2]{jm}.
\end{proof}

\chapter{The $A$-family}
\label{cha:afam}

The $A$-family consists of the matrix groups 
\begin{equation*}
  \pgl{n+1}=\frac{\GL(n+1,\C)}{\GL(1,\C)}, \qquad n \geq 1,
\end{equation*}
where $\GL(n+1,\C)$ is the Lie group of complex $(n+1) \times (n+1)$
matrices with center $\GL(1,\C)$ consisting of scalar matrices.  The
\mtn\ for $\pgl{n+1}$ is
\begin{equation*}
  \N(\pgl{n+1}) = \frac{\GL(1,\C)^{n+1}}{\GL(1,\C)} \rtimes
  \Sigma_{n+1} 
\end{equation*}
where $\Sigma_{n+1}=\W(\pgl{n+1}) \subset \PGL(n+1,\C)$ is the Weyl
group of permutation matrices. It is known \cite{hms:first,
  matthey:normalizers} that
\begin{equation}\label{eq:afamH0H1}
  H^0(\W;\ch{\T})=
  \begin{cases}
    \Z/2 & n=1 \\
     0   & n>1
  \end{cases}, \qquad
  H^1(\W;\ch{\T})=
  \begin{cases}
    \Z/2 & n=3 \\
     0   & n \neq 3
  \end{cases}
\end{equation}\add
for $\pgl{n+1}$. For all $n$, $\PGL(n+1,\C)=\PSL(n+1,\C)$. When $n+1$
is odd,  $\PGL(n+1,\C)=\PSL(n+1,\C)=\SL(n+1,\C)$ as \twocg s.

\section{The structure of $\PGL(n+1,\C)$}
\label{sec:pglnc}

In this and the following section we use the results of
Chapter~\ref{cha:ndet} to show that the \twocg s $\pgl{n+1}$, $n \geq
1$,  are
uniquely ${\N}$-determined. This section provides the information
about the Quillen category needed for the calculation (\ref{lim=0}) of
the higher limit obstruction groups from \ref{ndetauto}
and \ref{indstepalt}.

\setcounter{subsection}{\value{thm}}
\subsection{The toral subcategory of $\A(\PGL(n+1,\C))$}
\label{sec:Apglnct}\add

We consider the full subcategory of $\A(\PGL(n+1,\C))$ generated by
the toral nontrivial \lmntwo s in $\PGL(n+1,\C)$,
$\A(\PGL(n+1,\C))^{\leq t}$ (\ref{defn:toralAX}) .

%Define $\A(\Sigma_{n+1},t(\PGL(n+1,\C)))$ to be the category whose
%objects are nontrivial subspaces of $t(\PGL(n+1,\C))$, the
%$\F_2\Sigma_{n+1}$-module of elements of order $\leq 2$ in
%$\ch{T}((\PGL(n+1,\C))$, and whose \m s $V_1 \to V_2$ are the group
%homo\m s induced by the $\Sigma_{n+1}$-action; i.e.\ the set of orbits
%$\overline{\Sigma}_{n+1}(V_1,V_2)/\Sigma_{n+1}(V_1)$ for the action of
%the group $\Sigma_{n+1}(V_1)=\{\sigma\in\Sigma_{n+1}\vert\sigma v=v
%\text{\ for all\ } v \in V_1\}$ on the set
%$\overline{\Sigma}_{n+1}(V_1,V_2) = \{\sigma\in\Sigma_{n+1}\vert\sigma
% (V_1)\subseteq V_2\}$.

\begin{lemma}\label{tpglnC}
  The mono\m\ \func{\nu}{V}{\PGL(n+1,\C)} is toral if and only if it
  lifts to a \m\ $ V \to \GL(n+1,\C)$.  If $n+1$ is odd, all objects
  of $\A(\PGL(n+1,\C))$ are toral.
\end{lemma}
\begin{proof}
  Any mono\m\ $V \to \GL(n+1,\C) \to \PGL(n+1,\C)$ is toral since it
  is toral already in $\GL(n+1,\C)$ by complex representation theory.
  Conversely, any toral mono\m\ $V \to \GL(1,\C)^{n+1}/\GL(1,\C)
  \subset \PGL(n+1,\C)$ lifts to $\GL(1,\C)$ since $\GL(1,\C)$ is
  divisible. When $n+1$ is odd, $\PGL(n+1,\C)=\SL(n+1,\C) \subset
  \GL(n+1,\C)$ as \twocg s so all mono\m s $V \to \PGL(n+1,\C)$ are
  toral.
%%%%%%%%%%%%%%%%%%%%%%%%%%%%%%%%%%%%%%%%%% 
 % All objects of $\A(\GL(n+1,\C))$ are toral by complex representation
%  theory. Any mono\m\ $V \to (\C^{\times})^{n+1}/\C^{\times}$ lifts to
%  $(\C^{\times})^{n+1}$ since $\C^{\times}$ is a divisible abelian
%  group. If $n$ is even, $\PGL(n+1,\C)=\SL(n+1,\C)$ as \twocg s and all
%  mono\m s $V \to \SL(n+1,\C)$ are toral by complex representation
%  theory.
\end{proof}

Let
\begin{equation*}
  e_i = \diag(+1,\ldots,+1,-1,+1,\ldots,+1) \in \GL(n+1,\C), \qquad 1
  \leq i \leq n+1,
\end{equation*}
be the diagonal matrix with $-1$ in position $i$ and $+1$ at all other
positions. The maximal toral \lmntwo s
\begin{align*}
  &\Delta_{n+1} = \gen{e_1, \ldots ,e_{n+1}} =
  \gen{\diag(\pm 1, \ldots, \pm 1)} \cong (\Z/2)^{n+1} \subset
  \GL(n+1,\C), \\
  &P\Delta_{n+1} = \gen{e_1, \ldots ,e_{n+1}}/\gen{e_1\cdots e_{n+1}}
  \cong  (\Z/2)^{n} \subset
  \PGL(n+1,\C),
\end{align*}
have Quillen auto\m\ groups $\Sigma_{n+1} \cong
\A(\GL(n+1,\C))(\Delta_{n+1}) \cong \A(\PGL(n+1,\C))(P\Delta_{n+1})$.

\begin{lemma}  \label{lemma:afamtoralcat}
  The inclusion functors
  \begin{equation*} 
  \A(\Sigma_{n+1},\Delta_{n+1}) \to \A(\GL(n+1,\C)), \qquad
  \A(\Sigma_{n+1},P\Delta_{n+1}) \to
  \A(\PGL(n+1,\C))^{\leq t}    
  \end{equation*}
   are equivalence of categories.
\end{lemma}
\begin{proof}
  This is a general fact; the first part of \cite[2.8]{jmm:ndet} 
  also holds for the case $p=2$. However, it may be more illustrative
  to prove the lemma directly in this special case. 
  
  By complex represention theory, any nontrivial \lmntwo\ in
  $\GL(n+1,\C)$ is conjugate to a subgroup of $\Delta_{n+1}$ and
  $\A(\GL(n+1,\C))(\Delta_{n+1})=\Sigma_{n+1}$. Thus there is a
  faithful inclusion functor $\A(\Sigma_{n+1},\Delta_{n+1}) \to
  \A(\GL(n+1,\C))$ which is surjective on the sets of iso\m\ classes
  of objects. It remains to show that this functor is full. Since any
  \m\ in the category $\A(\GL(n+1,\C))$ is an iso\m\ followed by an
  inclusion, it is enough to show that any conjugation induced iso\m\ 
  $V_1 \to V_2$ between nontrivial subgroups $V_1,V_2 \subset
  \Delta_{n+1}$ is actually induced from conjugation by an element of
  $N(\GL(n+1,\C))$. But this is well-known fact from Lie group theory
  easily derived from eg \cite[IV.2.5]{brockerdieck}.

  Any toral nontrivial \lmntwo\ in $\PGL(n+1,\C)$ is the image of a
  \lmntwo\ in $\GL(n+1,\C)$ and hence conjugate to subgroup of
  $P\Delta_{n+1}$. Since any $\A(\PGL(n+1,\C))$-\m\ between subgroups
  of $P\Delta_{n+1}$ are induced from conjugation with an element of
  $N(\PGL(n+1,C))$, it follows that $\A(\Sigma_{n+1},P\Delta_{n+1})
  \to \A(\PGL(n+1,\C))^{\leq t}$ is an equivalence of categories.
\end{proof}

For any partition $n+1=i_0+i_1+ \cdots + i_r$ of $n+1$ into a sum of
$r$ positive integers, let $(\pm 1)^{i_0}(\pm 1)^{i_1} \cdots (\pm
1)^{i_r}$ denote the diagonal matrix
\begin{equation*}
  \diag(\overbrace{\pm 1, \ldots \pm 1}^{i_0},
        \overbrace{\pm 1, \ldots \pm 1}^{i_1}, \ldots ,
        \overbrace{\pm 1, \ldots \pm 1}^{i_r})
\end{equation*}
in $\GL(n+1,\C)$.

For any partition $(i_0,i_1)$ of $n+1=i_0+i_1$ into a sum of two
positive integers $i_0 \geq i_1 \geq 1$, let $L[i_0,i_1] \subset
\PGL(n+1,\C)$ be the image in $\PGL(n+1,\C)$ of the \lmntwo\ 
\begin{equation*}
 L[i_0,i_1]^*=\big\langle
 (+1)^{i_0}(-1)^{i_1}, (-1)^{n+1}
 \big\rangle 
\end{equation*}
in $\GL(n+1,\C)$. The centralizer of $L[i_0,i_1]$ is
\set\begin{equation}
  \label{eq:afamCL}
  C_{\PGL(n+1,\C)}L[i_0,i_1]=
  \begin{cases}
    \frac{\GL(i_0,\C)^2}{\GL(1,\C)}\rtimes C_2 & i_0=i_1\\
    \frac{\GL(i_0,\C) \times \GL(i_1,\C)}{\GL(1,\C)} & i_0>i_1
  \end{cases}
\end{equation}\add
where the action of
\begin{equation*}
  C_2 = \left\langle
      \begin{pmatrix}
        0&E\\E&0
      \end{pmatrix}\right\rangle
\end{equation*}
interchanges the two $\GL(i_0,\C)$-factors.  The center of the
centralizer of $L[i_0,i_1]$ is 
\set\begin{equation}
  \label{eq:afamZCL}
   \Ze C_{\PGL(n+1,\C)}L[i_0,i_1]=
   \begin{cases}
     L[i_0,i_1] & i_0=i_1 \\
     \frac{\GL(1,\C) \times \GL(1,\C)}{\GL(1,\C)} & i_0>i_1
   \end{cases}
\end{equation}\add

For any partition $(i_0,i_1,i_2)$ of $n+1=i_0+i_1+i_2$ into a sum of
three positive integers $i_0 \geq i_1 \geq i_2 \geq 1$ let
$P[i_0,i_1,i_2] \subset \PGL(n+1,\C)$ be the image in
$\PGL(n+1,\C)$ of the \lmntwo\ 
\begin{equation*}
  P[i_0,i_1,i_2]^*= \gen{ 
  (+1)^{i_0}(-1)^{i_1}(+1)^{i_2}, 
  (+1)^{i_0}(+1)^{i_1}(-1)^{i_2}, (-1)^{n+1} }
\end{equation*}
in $\GL(n+1,\C)$. The centralizer of $P[i_0,i_1,i_2]$ is
\set\begin{equation}
  \label{eq:afamCP3}
  C_{\PGL(n+1,\C)} P[i_0,i_1,i_2] = \frac{\GL(i_0,\C) \times
    \GL(i_1,\C) \times \GL(i_2,\C)}{\GL(1,\C)}
\end{equation}\add
so that the center of the centralizer is 
\set\begin{equation}
  \label{eq:afamZCP3}
   \Ze C_{\PGL(n+1,\C)} P[i_0,i_1,i_2] = \frac{\GL(1,\C) \times
    \GL(1,\C) \times \GL(1,\C)}{\GL(1,\C)}
\end{equation}\add
connected.

For any partition $(i_0,i_1,i_2,i_3)$ of $n+1$ into a sum
$n+1=i_0+i_1+i_2+i_3$ of $n+1$ into a sum of four positive integers
$i_0 \geq i_1 \geq i_2 \geq i_3 \geq 1$ let $P[i_0,i_1,i_2,i_3]
\subset \PGL(n+1,\C)$ be the image in $\PGL(n+1,\C)$ of the \lmntwo\ 
\begin{equation*}
  P[i_0,i_1,i_2,i_3]^*= \gen{ 
  (+1)^{i_0}(-1)^{i_1}(+1)^{i_2}(-1)^{i_3}, 
  (+1)^{i_0}(+1)^{i_1}(-1)^{i_2}(-1)^{i_3}, (-1)^{n+1} }
 % \big\langle
%  \diag(\overbrace{1,\ldots,1}^{i_0},\overbrace{-1,\ldots,-1}^{i_1},
%        \overbrace{1,\ldots,1}^{i_2},\overbrace{-1,\ldots,-1}^{i_3}), \\
%  \diag(\overbrace{1,\ldots,1}^{i_0},\overbrace{1,\ldots,1}^{i_1},
%        \overbrace{-1,\ldots,-1}^{i_2},\overbrace{-1,\ldots,-1}^{i_3}),
%   -E\big\rangle
\end{equation*} 
The centralizer of $P[i_0,i_1,i_2,i_3]$ is
\begin{equation}
  \label{eq:afamCP}
  C_{\PGL(n+1,\C)}P[i_0,i_1,i_2,i_3]=
  \begin{cases}
    \frac{\GL(i_0,\C)^4}{\GL(1,\C)} \rtimes (C_2 \times C_2)  
     & i_0=i_1=i_2=i_3 \\
  \frac{\GL(i_0,\C)^2 \times \GL(i_2,\C)^2}{\GL(1,\C)} \rtimes C_2 
     & i_0=i_1>i_2=i_3 \\
  \frac{\GL(i_0,\C) \times \GL(i_1,\C) \times \GL(i_2,\C) \times
    \GL(i_3,\C)}{\GL(1,\C)} & \textmd{otherwise}
\end{cases}
\end{equation}
where 
\begin{equation*}
  C_2 \times C_2 = \gen{
    \begin{pmatrix}
      0&E&0&0\\E&0&0&0\\0&0&0&E\\0&0&E&0
    \end{pmatrix},
    \begin{pmatrix}
      0&0&E&0\\0&0&0&E\\E&0&0&0\\0&E&0&0
    \end{pmatrix}}, \qquad
   C_2= \gen{
 \begin{pmatrix}
      0&E&0&0\\E&0&0&0\\0&0&0&E\\0&0&E&0
    \end{pmatrix}}
\end{equation*}
The center of the centralizer of $P[i_0,i_1,i_2,i_3]$ is
\begin{multline}
   \label{eq:afamZCP}
   \Ze C_{\PGL(n+1,\C)}P[i_0,i_1,i_2,i_3]= \\
   \begin{cases}
   P[i_0,i_1,i_2,i_3] & i_0=i_1=i_2=i_3 \\
  \frac{\GL(1,\C) \times \GL(1,\C)}{\GL(1,\C)} \times
  \gen{(+1)^{i_0}(-1)^{i_1}(+1)^{i_2}(-1)^{i_3}} 
  & i_0=i_1>i_2=i_3 \\
   \frac{\GL(1,\C) \times \GL(1,\C) \times \GL(1,\C) \times
    \GL(1,\C)}{\GL(1,\C)} & \textmd{otherwise}
  \end{cases}
\end{multline}

We collect the information about the toral subcategory that we shall
need later on in the folllowing proposition.  Let $P(m,k)$ denote the
number of partitions of $m$ into sums of $k$ natural integers.

\begin{prop}\label{prop:afamtoral}
  The category $\A(\PGL(n+1,\C))$ contains precisely 
  \begin{itemize}
  \item $P(n+1,2)$ iso\m\ classes of toral rank one objects
    represented by the lines $L[i_0,i_1]$.
  \item $P(n+1,3)+P(n+1,4)$ iso\m\ classes of toral rank two objects
    represented by the planes $P[i_0,i_1,i_2]$ and
    $P[i_0,i_1,i_2,i_3]$.
  \end{itemize}
  The centralizers of these objects are listed in (\ref{eq:afamCL}),
  (\ref{eq:afamCP3}), and (\ref{eq:afamCP}).
\end{prop}

The auto\m\ groups are easily computed using complex representation
theory because 
\begin{equation*}
  \A(\GL(n+1,\C))(P[i_0,i_1,i_2,i_3]^*) \to
  \A(\PGL(n+1,\C))(P[i_0,i_1,i_2,i_3]) 
\end{equation*}
is surjective  (as in \ref{Qauto}). One finds that
\begin{equation*}
  \A(\PGL(n+1,\C))P[i_0,i_1,i_2,i_3]=
  \begin{cases}
    \GL(2,\F_2) & \ell(i_0,i_1,i_2,i_3) \geq 3 \\
    C_2 & \ell(i_0,i_1,i_2,i_3) = 2 \\
    \{1\} & \ell(i_0,i_1,i_2,i_3) =1 \\
  \end{cases}
\end{equation*}
where $\ell(i_0,i_1,i_2,i_3)=\max_{1 \leq j \leq 4}\#\{k \mid
i_k=i_j\}$ is the maximal number of repetitions in the sequence
$(i_0,i_1,i_2,i_3)$. This formula also
holds for the objects $P[i_0,i_1,i_2]$ when interpreted as
$P[i_0,i_1,i_2,0]$.

\section{Centralizers of objects of $\A(\pgl{n+1})^{\leq t}_{\leq 2}$
  are LHS} 
\label{sec:afamlhs}

In this section we check that all toral objects of rank $\leq 2$ have
LHS (\ref{cha:ndet}.\S\ref{subsec:LHS}) centralizers.

%A \twocg\ $X$ is LHS \cite[2.22]{jmm:2cgs} if the restriction map
%$H¹(W;\ch{T}) \to H^1(W_0;\ch{T})^{W/W_0}$ is surjective or,
%equivalently, if the initial segment
%\begin{equation*}
%  0 \to H^1(\pi;\ch{T}^{W_0}) \to H^1(W;\ch{T}) \to
%  H^1(W_0;\ch{T})^{W/W_0} \to 0
%\end{equation*}
%of the Lyndon--Hochschild--Serre spectral sequence for the extension
%$W_0 \to W \to W/W_0$ is exact. Here, $\ch{T}$ is the discrete \mt\ of
%$X$, $W$ is the Weyl group, and $W_0$ the Weyl group of the identity
%component.

%The objective of this section is to verify that the centralizers of
%all toral \lmntwo s of rank $\leq 2$ are LHS. (No nonLHS \twocg s are
%presently known.) This will pave the way for a later application of
%\cite[3.6]{jmm:2cgs}.

\begin{lemma}\label{ex:LHS} 
  The centralizers of the objects of $\A(\pgl{n+1})^{\leq t}_{\leq
    2}$, 
  \begin{enumerate}
  \item $\frac{\GL(i,\C)^2}{\GL(1,\C)} \rtimes C_2$
    (\ref{eq:afamCL}),
  \item  $\frac{\GL(i,\C)^4}{\GL(1,\C)} \rtimes (C_2 \times C_2)$
    (\ref{eq:afamCP}), 
  \item $\frac{\GL(i_0,\C)^2 \times \GL(i_2,\C)^2}{\GL(1,\C)} \rtimes
    C_2$ (\ref{eq:afamCP})
  \end{enumerate}
  are LHS.
\end{lemma}
\begin{proof}
  (1) Let $X=\frac{\GL(i,\C)^2}{\GL(1,\C)} \rtimes C_2$, $i \geq 1$,
  where the $C_2$-action switches the two $\GL(i,\C)$-factors.  For
  $i=1$, $X$ is a \twoctg , hence LHS. 
  For
  $i=2$ explicit computer computation yields 
  \begin{center}
    \begin{tabular}[c]{|c||c|c|c|c|} \hline 
      $\frac{\GL(i,\C)^2}{\GL(1,\C)} \rtimes C_2$ &
      $H^1(\pi;\ch{{\T}}^{{\W}_0})$ & $H^1({\W};\ch{{\T}})$ &  
      $H^1({\W}_0;\ch{{\T}})^{\pi}$ &
      $H^1({\W}_0;\ch{{\T}})$ \\ \hline \hline
      $i=2$ & $0$ & $\Z/2$ & $\Z/2$ & $\Z/2$ \\ \hline
    \end{tabular}
  \end{center}
  so $X$ is manifestly LHS in this case (even though $X_0$ is not
  regular).  For $i>2$, $\theta(X_0)$ is bijective and thus $X$ is LHS
  by \ref{lhscrit1}. ($\theta(X_0)$ is injective by
  \ref{kerneltheta}.(\ref{kerneltheta.1}) and surjective by
  \ref{kerneltheta}.(\ref{kerneltheta.2}) for $i \neq 4$ and for $i=4$
  by inspection or by
  \ref{regquot} and \ref{exmp:glmCreg} for all $i>2$.)
  
 % for $i>2$ and $i \neq 4$ by \ref{kerneltheta}; for $i=4$,
%  $\theta(X_0)$ is an injective homo\m\ from
%  $\Hom(W_0,\ch{T}^{W_0})=(\Z/2)^2$ to $H^1(W_0;\ch{T})=(\Z/2)^2$
%  (computer computation).)
  
  \noindent
  (2) Let $X= \frac{\GL(i,\C)^4}{\GL(1,\C)} \rtimes (C_2\times
  C_2)$, $i \geq 1$, where $C_2 \times C_2 = \gen{(12)(34),(13)(24)}$
  permutes the four $\GL(i,\C)$-factors. For $i=1$, $X$ is a \twoctg ,
  hence LHS.  For
  $i=2$ explicit computer computation yields 
  \begin{center}
    \begin{tabular}[c]{|c||c|c|c|c|} \hline 
      $\frac{\GL(i,\C)^4}{\GL(1,\C)} \rtimes (C_2 \times C_2)$ &
      $H^1(\pi;\ch{{\T}}^{{\W}_0})$ & $H^1({\W};\ch{{\T}})$ &  
                 $H^1({\W}_0;\ch{{\T}})^{\pi}$ &
      $H^1({\W}_0;\ch{{\T}})$ \\ \hline \hline
      $i=2$ & $\Z/2$ & $(\Z/2)^3$ & $(\Z/2)^2$ & $(\Z/2)^8$ \\ \hline
    \end{tabular}
  \end{center}
  so $X$ is manifestly LHS in this case. (Alternatively, observe that
  $X_0$ is regular (\ref{exmp:glmCreg}, \ref{regquot}), the kernel of
  $\theta(X_0)$ is $(\Z/2)^4$, and $\theta(X_0)^{\pi}$ is surjective
  because $H^1(C_2 \times C_2;(\Z/2)^4)=0$ for the regular
  representation.)  For $i>2$, we see as in \ref{ex:LHS}.(1) above
  that $\theta(X_0)$ is bijective and hence $X$ is LHS by
  \ref{lhscrit1}.

  \noindent
  (3) Let $X=\frac{(\GL(i_0,\C)\times\GL(i_2,\C))^2}{\GL(1,\C)}\rtimes
  C_2$, $1 \leq i_0 < i_2$, where $C_2$ switches the two identical
  factors.  Using \ref{regquot} and \ref{exmp:glmCreg} we see (details
  omitted) that $X_0$ is regular. By
  \ref{kerneltheta}.(\ref{kerneltheta.1}), $\theta(X_0)$ is in fact
  bijective except when $i_0$ or $i_2$ is $2$.  In those cases, the
  kernel of $\theta(X_0)$ is $(\Z/2)^2$ and $\theta(X_0)^{C_2}$ is
  surjective as $H^1(C_2;(\Z/2)^2)=0$ for the regular representation.
  Therefore $X$ is LHS by \ref{lhscrit1}.
\end{proof}

\section{Limits over the Quillen category of $\pgl{n+1}$}
\label{sec:afamlim}

In this section we show that the problem of computing the higher
limits of the functors $\pi_i(B\Ze C_{\pgl{n+1}})$, $i=1,2$,
(\ref{eq:defnpiiBZC}) is concentrated on
the nontoral objects of the Quillen category. 

\begin{lemma}\label{lemma:pijB}\cite[2.8]{jmm:ndet}
  Let $V \subset P\Delta_{n+1}$ be a nontrivial subgroup
  representing an object of
  $\A(\Sigma_{n+1},P\Delta_{n+1})=\A(\PGL(n+1,\C))^{\leq t}$
  (\ref{lemma:afamtoralcat}). Then
  \begin{equation*}
    \ch{Z}C_{\PGL(n+1,\C)}(V) = \ch{T}^{\Sigma_{n+1}(V)}
  \end{equation*}
  where $\ch{T} = \ch{T}(\PGL(n,\C))$ is the discrete approximation
  \cite[\S3]{dw:center} to the \mt\ of $\PGL(n+1,\C)$ and
  $\Sigma_{n+1}(V)$ is the point-wise stabilizer subgroup
  (\ref{defn:AWt}).
 % \begin{equation*}
%    \pi_i(\B
%    {\Ze}C_{\PGL(n+1,\C)})(V)=H^{2-i}(\Sigma_{n+1}(V);{L}),
%    \quad i=1,2,
%  \end{equation*}
%  where ${L}=\pi_2(B\T(\PGL(n+1,\C)))$ and $\Sigma_{n+1}(V)$
%  is the point-wise stabilizer subgroup (\ref{defn:AWt}). The
%  cohomology groups $H^{2-i}(W;L)$, $i=1,2$, are trivial for $n>1$.  
\end{lemma}
\begin{proof}
Let \func{\nu^*}{V}{T(\GL(n+1,\C))} be a lift to $\GL(n+1,\C)$ of the
inclusion homo\m\ of $V$ into $T(\PGL(n+1,\C))$. Then 
\begin{equation*}
  C_{\GL(n+1,\C)}(\nu^*V)= \prod_{\rho \in V^{\vee}} \GL(i_{\rho},\C),
  \qquad \Sigma_{n+1}(\nu^*V) = \prod_{\rho \in
    V^{\vee}}\Sigma_{i_{\rho}} 
\end{equation*}
where \func{i}{V^{\vee}}{\Z} records the multiplicity of each linear
character $\rho \in V^{\vee}$ in the representation $\nu^*$. Using
\cite[5.11]{jmm:ndet} and \ref{lemma:WVWVast} we get that
\begin{equation*}
  C_{\PGL(n+1,\C)}(V) =  \frac{ C_{\GL(n+1,\C)}(\nu^*V) }{\GL(1,\C)}
  \rtimes V^{\vee}_{\nu^*}, \qquad
  \Sigma_{n+1}(V) = \Sigma_{n+1}(\nu^*V) \rtimes V^{\vee}_{\nu^*}
\end{equation*}
where $V^{\vee}_{\nu^*} = \{ \zeta \in  V^{\vee}= \Hom(V,\GL(1,\C))  \mid
\forall \rho \in V^{\vee} \colon i_{\zeta\rho} = i_{\rho}\}$. The
semi-direct products are obtained because the elements of
$V^{\vee}_{\nu^*}$ can be effectuated by permutations from
$\Sigma_{n+1}$ that fix $V \subset \PGL(n+1,\C)$ point-wise. The
discrete approximation \cite[\S3]{dw:center} to the center of the
centralizer is therefore
\begin{align*}
  \ch{Z} C_{\PGL(n+1,\C)}(V) 
  &= \ch{Z} \left( \frac{\prod \GL(i_{\rho},\C)}{\GL(1,\C)} \rtimes
    V^{\vee}_{\nu^*} \right) 
  \stackrel{\textmd{(\ref{cor:cent})}}{=} 
  \ch{Z} \left( \frac{\prod \GL(i_{\rho},\C)}{\GL(1,\C)}
  \right)^{V^{\vee}_{\nu^*}}  \\
  &=\left( \frac{\prod \ch{Z}\GL(i_{\rho},\C)}{\GL(1,\C)} 
  \right)^{V^{\vee}_{\nu^*}}
  =\left( \frac{\ch{T}(\GL(n+1,\C))^{\Sigma_{n+1}(\nu^*V)}}{\GL(1,\C)} 
  \right)^{V^{\vee}_{\nu^*}} \\
  &=\left( \ch{T}(\PGL(n+1,\C))^{\Sigma_{n+1}(\nu^*V)}
  \right)^{V^{\vee}_{\nu^*}} 
  = \ch{T}(\PGL(n+1,\C))^{\Sigma_{n+1}(V)}
\end{align*}
where the penultimate equlity sign is justified by the fact that 
$H^1(\Sigma_{n+1}(\nu^*V);\GL(1,\C)) \to
 H^1(\Sigma_{n+1}(\nu^*V);\ch{T}(\GL(n+1,\C)))$ is injective.
\end{proof}

Define $\pi_i(\B {\Ze}C_{\PGL(n+1,\C)})_{\not \leq t}$ to be the
subfunctor of $\pi_i(\B {\Ze}C_{\PGL(n+1,\C)})$
(\ref{eq:defnpiiBZC}) that vanishes on all toral
objects and is unchanged on all nontoral objects of the Quillen
category. This means that \set\begin{equation}\label{defn:nontoralsub}
  \pi_i(\B{\Ze}C_{\PGL(n+1,\C)})_{\not \leq t}(V,\nu) =
   \begin{cases}
     0 & \text{$(V,\nu)$ is toral} \\
      \pi_i(\B{\Ze}C_{\PGL(n+1,\C)})(V,\nu) & \text{$(V,\nu)$ is
        nontoral} 
   \end{cases}
  \end{equation}\add
for all objects $(V,\nu)$ of $\A(\pgl{n+1})$. The reason for
introducing this subfunctor is that in the computation of the higher
limits, we can ignore the toral objects.

\begin{cor}\label{torallimforpglnC}
  When $n>1$, 
  \begin{equation*}
    \lim^*(\A(\PGL(n+1,\C));\pi_i(BZC)_{\nleq t}) \cong
    \lim^*(\A(\PGL(n+1,\C));\pi_i(BZC)), \quad i=1,2,
  \end{equation*}
  where $\pi_i(BZC)=\pi_i(\B {\Ze}C_{\PGL(n+1,\C)})$
  (\ref{eq:defnpiiBZC}). 
\end{cor}
\begin{proof}
The result of Lemma~\ref{lemma:pijB} is (\ref{eq:hiwt}) equivalent to 
\begin{equation*}
  \pi_i(BZC)(V) = H^{2-i}(\Sigma_{n+1}(V);L) , \qquad V \subset
  P\Delta_{n+1},
\end{equation*}
where $L$ is the $\Z_2\Sigma_{n+1}$-module $\pi_2BT(\PGL(n+1,\C))$ and
therefore (\ref{dw:limits})
 \begin{equation*}
\lim^j(\A(\PGL(n+1,\C))^{\leq t},
\pi_i(\B{\Ze}C)) =
\begin{cases}
  H^{2-i}(\Sigma_{n+1}, L) & j=0 \\
  0 & j>0
\end{cases}
\end{equation*}
where the cohomology groups $H^{2-i}(\Sigma_{n+1}; L)$,
$i=1,2$, are trivial for $n>1$ (\ref{eq:afamH0H1}).

Since the quotient functor $\pi_i(\B{\Ze}C)/\pi_i(\B{\Ze}C)_{\not\leq
  t}$ vanishes on all nontoral objects
\begin{equation*}
  \lim^j(\A(\PGL(n+1,\C)),
  \pi_i(\B{\Ze}C)/\pi_i(\B{\Ze}C)_{\not\leq t})
  \stackrel{\textmd{\cite[13.12]{jmm:ndet}}}{\cong}
  \lim^j(\A(\PGL(n+1,\C))^{\leq t},\pi_i(\B{\Ze}C))
\end{equation*}
We conclude that $\lim^*(\A(\PGL(n+1,\C)),
\pi_i(\B{\Ze}C)/\pi_i(\B{\Ze}C)_{\not\leq t})=0$.  The long exact
coefficient functor sequence for higher limits now shows that
$\lim^*(\A(\PGL(n+1,\C)),\pi_i(\B{\Ze}C)_{\not\leq t})$ and
$\lim^*(\A(\PGL(n+1,\C)), \pi_i(\B{\Ze}C))$ are isomorphic.
\end{proof}

%Thus the problem of computing the higher limits of
%the functor $\pi_i(\B{\Ze}C_{\PGL(n+1,\C)})$ is concentrated on the
%nontoral objects of the Quillen category. This does not mean, for
%reasons that will become clear later (\ref{nonTpglnC}), that the toral
%objects are irrelevant, so we shall now continue with their
%classification.

\section{The category $\A(\PGL(n+1,\C))^{[\ , \ ] \neq 0}$}
\label{sec:nontoralapglnc}
%For \twocg\ $X$, let $\A(X)_{\nleq t}$ denote
%the full subcategory of $\A(X)$ on all nontoral objects and their
%sub-objects.  
%We determine the nontoral objects of the category
%$\A(\PGL(n+1,\C))$.

For any nontrivial \lmntwo \ $V$ in $\PGL(n+1,\C)$, let \func{[\; , \;
  ]}{V \times V}{\F_2} be the symplectic bilinear form
\cite[II.9.1]{huppert:I} given by $[u\C^{\times},v\C^{\times}]=r$ if
$[u,v]=(-E)^r$ where $u,v \in \GL(n+1,\C)$ are such that
$u\C^{\times},\, v\C^{\times} \in V$. (The elements $[u,v]$ and $u^2$
lie in the center $\C^{\times}$ of $\GL(n+1,\C)$ so that
$E=[u^2,v]=[u,v]^u[u,v]=[u,v]^2$ and thus $[u,v] \in \C^{\times}$ has
order $2$. Therefore $[u,v] = [u,v]^{-1} = [v,u]$.)

\begin{lemma}\label{altform}
  $\text{$V$ in $\PGL(n+1,\C)$ is toral} \iff [V,V] =0$ 
\end{lemma}
\begin{proof}
  Let $e_i\C^{\times}$, $1\leq i \leq d$, be a basis for $V$. Since
  $\C^{\times}$ is divisible, we can assume that each $e_i \in
  \GL(n+1,\C)$ has order $2$. If $[V,V]=0$, these $e_i$s commute and
  span a lift to $\GL(n+1,\C)$ of $V \subseteq \PGL(n+1,\C)$.
\end{proof}

An extra special $2$-group is of {\em positive type\/} if
it is isomorphic to a central product of dihedral groups $D_8$ of order
$8$ \cite[p 145--146]{robinson:groups}.
\begin{lemma}\label{griesslemma}
  \cite[3.1]{griess:elem} \cite[5.4]{jmm:ndet} Let
  \func{\nu}{V}{\PGL(n,\C)} be a nontoral mono\m\ of a nontrivial
  elementary abelian $2$-group $V$ into $\PGL(n+1,\C)$. Then there
  exists a \m\ of \ses s of groups
  \begin{equation*}
    \xymatrix{
     1 \ar[r] & {\Ze}(P) \ar[r]\ar@{^{(}->}[d] & PE \ar[r]\ar@{^{(}->}[d] & V
     \ar[r]\ar@{^{(}->}[d]^{\nu} & 1 \\
     1 \ar[r] & {\C}^{\times} \ar[r] & {\GL(n+1,\C)} \ar[r] &
     {\PGL(n+1,\C)} \ar[r] & 1 } 
  \end{equation*}
  where $PE$ is the direct product of an extra special $2$-group
  $P\subseteq\GL(n+1,\C)$ of positive type and an elementary abelian
  $2$-group $E\subseteq\GL(n+1,\C)$ with $ P \cap E =\{1\}=[P,E]$.
\end{lemma}

%Write $\C^{n+1}=\C^{2^d}\otimes\C^m$ for some $d>0$ and some $m\geq
%0$. Let the extra-special $2$-group $2_+^{1+2d}$ act faithfully on the
%first factor of the tensor product and let the (possibly trivial)
%\lmntwo\ $E$ act faithfully on the second factor such that no
%nontrivial element of $E$ acts as scalar multiplication. This makes
%$\C^{n+1}$ a  $\C[2_+^{1+2d} \times E]$-module. The image
%of the group $2_+^{1+2d} \times E \subseteq \GL(n+1,\C)$ in
%$\PGL(n+1,\C)$ is a nontoral \lmntwo\ (\ref{altform}) and any
%nontoral \lmntwo\ in $\PGL(n+1,\C)$ has this form (\ref{griesslemma}).

Let $G=\gen{P,E,i}= P \circ C_4 \times E$ be the group generated by
$E$ and the central product $P \circ C_4$ of $P$ and the cyclic group
$C_4=\gen{i} \subseteq \C^{\times}$ with $C_2=\gen{-E}$ amalgamated.
The image of $G$ in $\PGL(n+1,\C)$ is $V$. 

Let $\A(\GL(n,\C))(G)$ be the subgroup, isomorphic to
$N_{\GL(n+1,\C)}(G)/G\cdot C_{\GL(n+1,\C)}(G)$, of $\Out(G)$
consisting of all outer auto\m s of $G$ induced from conjugation in
$\GL(n+1,\C)$ \cite[5.8]{jmm:ndet}.  In other words,
$\A(\GL(n,\C))(G)=\Out_{\textmd{tr}}(G)$ is the group of trace
preserving outer auto\m s of $G$.
%%%%%%%%%%%%%%%%%%%%%%%%%%%%%%%%%%%%%%%
%% do we need this?
%and $q(v\C^{\times})=v^2$, $v\in G$,
%is a quadratic form on $V$ such that $q(u\C^{\times}+v\C^{\times})=
%q(u\C^{\times})+q(v\C^{\times})+[u\C^{\times},v\C^{\times}]$ for all
%$u\C^{\times},v\C^{\times}\in V$.
%and there is a symmetric bilinear form
%\func{f}{V\times V} {\F_2} given by
%$f(u,v)=[\overline{u},\overline{v}] \in Z(P)=\gen{-1}$ where
%$\overline{u},\overline{v}\in G$ are lifts of $u,v \in V$.

\begin{lemma}\label{Qauto}
  $\A(\GL(n+1,\C))(G) \to \A(\PGL(n+1,\C))(V)$ is surjective.
\end{lemma}
\begin{proof}
  Suppose that $B\in\GL(n+1,\C)$ is such that $V^{B\C^{\times}}=V$.
  Then $G^B \subseteq G \cdot \C^{\times}$: For any $g \in G$ there
  exist $h \in G$ and $z \in \C^{\times}$ such that $g^B=hz$. But
  since $G$ has exponent $4$, $z^4=1$ so $z\in C_4$ and $g^B \in G$.
 % This shows that $\A(\GL(n+1,\C))(G,G) \to \A(\PGL(n+1,\C))(V,V)$ is
%  surjective. 
  %Now $\A(\GL(n+1,\C))(G,G)$ is a subgroup of $\Out(G)$.
%  When $E$ is trivial, this group was determined in \cite[p.\ 
%  403]{griess:autos} to be $\Z/2 \times \Aut(V,f)$. We now jump to the
%  conclusion of the lemma.
\end{proof}

%%%%%%%%%%%% begin from nontoral.tex %%%%%%%%%%%%%%%%%%%%%%

A monomorphic conjugacy class \func{\nu}{V}{\PGL(n+1,\C)} is said to
be a $(2d+r,r)$ object of $\A(\PGL(n+1,\C))$ if the underlying
symplectic vector space of $(V,\nu)$ is isomorphic to $V=H^d \times
V^{\perp}$ where $H$ denotes the symplectic plane over $\F_2$ and
$\dim_{\F_2}V^{\perp}=r$ \cite[II.9.6]{huppert:I} (so that
$\dim_{\F_2}V=2d+r$). An $(r,r)$ object is the same thing as an
$r$-dimensional toral object.  We write $\Symp(V)$ or $\Symp(2d+r,r)$
(abbreviated to $\Symp(2d)$ if $r=0$) for the group of linear auto\m s
of $V$ that preserve the symplectic form.

\begin{cor}\label{nonTpglnC}
Suppose that $n+1=2^dm$ for some natural numbers $d\geq 1$ and $m \geq 1$.
\begin{enumerate}
\item There is up to iso\m\ a unique $(2d,0)$ object $H^d$ of
$\A(\PGL(n+1,\C))$, and
\begin{equation*}
 \A(\PGL(n+1,\C))(H^d) =\Symp(2d), \quad
 C_{\PGL(n+1,\C)}(H^d) = H^d \times \PGL(m,\C)
\end{equation*}
for this object.
\item \label{nonTpglnC2} Isomorphism classes of
$(2d+r,r)$, $r>0$, objects $V$ of $\A(\PGL(2^dm,\C))$
correspond bijectively to
isomorphism classes of $(r,r)$ objects $V^{\perp}$ of
$\A(\PGL(m,\C))$, and
\begin{align*}
& \A(\PGL(2^dm,\C))(V) = \begin{pmatrix}
\Symp(2d) & 0 \\
* & \A(\PGL(m,\C))(V^{\perp})
\end{pmatrix} \\
& C_{\PGL(2^dm,\C)}(V) = V/V^{\perp} \times C_{\PGL(m,\C)}(V^{\perp})
\end{align*}
for these objects.
\end{enumerate}
\end{cor}
\begin{proof}
1.\ The group $2_+^{1+2d}\circ 4$ has \cite[7.5]{huppert:char} $2^{1+2d}$
characters of degree $1$ and $2$ irreducible characters of degree
$2^d$ (interchanged by the action of $\Out(2_+^{1+2d}\circ 4) \cong
\Symp(2d) \times \Aut(C_4)$ \cite[pp.\ 403--404]{griess:autos}) given by
\begin{equation*}
  \chi_{\lambda}(g)=
  \begin{cases}
    2^d\lambda(g)  & g\in C_4 \\
    0 & g\not\in C_4
  \end{cases}
\end{equation*}
where \func{\lambda}{C_4}{\C^{\times}} is an injective group homo\m\
($\lambda(i)=\pm i$). The linear characters vanish on the derived
group $2=[2_+^{1+2d}\circ 4,2_+^{1+2d}\circ 4]$ but the irreducible
characters of degree $2^d$ do not. Thus the only faithful representations of
$2_+^{1+2d}\circ 4$ with central centers are multiples
$m\chi_{\lambda}$ of $\chi_{\lambda}$ for a fixed $\lambda$. Phrased
slightly differently, $\GL(m2^d,\C)$ contains up to conjugacy a unique
subgroup with central center isomorphic to $2_+^{1+2d}\circ 4$. For
this group and its image $H^d$ in $\PGL(2^dm,\C)$ we have
\begin{gather*}
  \A(\GL(m2^d,\C))(2_+^{1+2d}\circ 4,2_+^{1+2d}\circ 4) \cong
 \Symp(2d) \cong \A(\PGL(m2^d,\C))(H^d,H^d) \\
  C_{\GL(m2^d,\C)}(2_+^{1+2d} \circ 4) \cong \GL(m,\C), \quad
  C_{\PGL(m2^d,\C)}(H^d) \cong H^d \times \PGL(m,\C)
\end{gather*}
where the last iso\m\ is a consequence of \cite[5.9]{jmm:ndet}. \\
2.\ The $(2d+r,r)$ object $(V,\nu)$ of $\A(\PGL(2^dm,\C))$ and the
$(r,0)$ object $(V^{\perp},\nu^{\perp})$ of $\A(\PGL(m,\C))$
correspond to each other iff
there is an $m$-dimensional representation
\func{\mu}{V^{\perp}}{\GL(m,\C)} such that $\C^{2^d} \otimes \mu$ is a
lift of $\nu \vert V^{\perp}$ and $\mu$ a lift of $\nu^{\perp}$.
According to \ref{griesslemma} any lift of $\nu \vert V^{\perp}$ has
this form for some $\mu$ uniquely determined up to the action of
$(V^{\perp})^{\vee}$.

We use \ref{Qauto}  to calculate the Quillen auto\m\ group of a
$(2d+r,r)$ object $H^d \times V^{\perp}$ of $\A(\PGL(2^dm,\C))$. Let
$H^d \times V^{\perp}$ be covered by the group $P \circ C_4 \times
V^{\perp}$ as in \ref{griesslemma}. Let $\alpha$ be an auto\m\ of $P
\circ C_4$, let $\beta$ be any homo\m\ of the form $P \circ C_4 \to
H^d \to V^{\perp}$, and let $\gamma$ be any Quillen auto\m\ of
$(V^{\perp},\nu^{\perp})$. Choose a homo\m\ $\zeta_1 \colon P \circ
C_4 \to H^d \times C_4/C_2 \to C_4$ such that
$\lambda(\zeta_1(x)\alpha(x)) = \lambda(x)$ for all $x \in C_4$ and a
homo\m\ \func{\zeta_2}{V^{\perp}}{C_4} such that
$\lambda(\zeta_2(v)) \mu(\gamma(v)) = \mu(v)$ for all $v \in
V^{\perp}$. Then the auto\m\ of $P \circ C_4$ that takes $(x,v)$ to
$(\zeta_1(x)\zeta_2(v)\alpha(x),\beta(x)+\gamma(v))$ preserves the
trace of $\chi_{\lambda} \# \mu$ and therefore the auto\m\ induced
on the quotient is a Quillen auto\m\ of $H^d \times V^{\perp}$.
Conversely, any auto\m\ of $P \circ C_4 \times V^{\perp}$ takes the
center $C_4 \times V^{\perp}$ isomorphically to itself and hence
it is of the form $(x,v) \to (\zeta(x,v)\alpha(x),
\beta(x)+\gamma(v))$ for some auto\m\ $\alpha$ of $P \circ C_4$, some
homo\m\ \func{\beta}{P \circ C_4}{V^{\perp}} vanishing on $C_4$, and
some homo\m\ \func{\zeta}{P \circ C_4 \times V^{\perp}}{C_4}. Such an
auto\m\ preserves the trace of $\chi_{\lambda} \# \mu$ iff $\lambda(\zeta(x,v)\alpha(x)) =
\mu(\gamma(v))$ for all $(x,v) \in {\Ze}(P \circ C_4 \times V^{\perp}) = C_4 \times
V^{\perp}$. But this means that the induced auto\m\ of $H^d \times
V^{\perp}$ is of the stated form.
\end{proof}

We conclude that the nontoral objects of
$\A(\PGL(2m,\C))$ of rank $\leq 4$ are
\begin{itemize}
\item One $(2,0)$ object $H$, $\A(\pgl{2m})(V)=\Symp(2)$, 
\item $P(m,2)$ $(3,1)$ objects $V$, $\A(\pgl{2m})(V)=\Symp(3,1)$,
\item $P(m,3)+P(m,4)$ $(4,2)$ objects $E$, $\A(\pgl{2m})(E)=\begin{pmatrix}
\Symp(2) & 0 \\ * & \A(\pgl{m})(E^{\perp})
\end{pmatrix}$ where $\A(\pgl{m})(E^{\perp}) =1$, $C_2$, or $\GL(E^{\perp})$, 
\item One $(4,0)$ object $H^2$ if $m$ is even,
  $\A(\PGL(m,\C))(H^2)=\Symp(4)$. 
\end{itemize}
This information will be needed in the next section as input for
Oliver's cochain complex \cite{bob:steinberg} for computing higher
limits of the functors $\pi_i(\B {\Ze}C_{\PGL(n+1,\C)})_{\nleq t}$.

\section{Higher limits of the functor $\pi_iB\Ze C_{\pgl{n+1}}$ on
$\A(\PGL(n+1,\C))^{[\ , \ ] \neq 0}$} 
\label{sec:Afam}

We compute the higher limits from \ref{ndetauto}.(\ref{ndetauto2}) and
\ref{indstepalt}.(\ref{indstepalt5}) by means of
\ref{torallimforpglnC} and Oliver's cochain complex
\cite{bob:steinberg}.
%Since \ref{torallimforpglnC} is not valid for
%$\pgl{2}$ we first consider this case separately.

\begin{lemma}\label{lim=0} The low degree higher limits  of the functors
  $\pi_i B{\Ze}C_{\pgl{n+1}}$, 
  $i=1,2$, are:
  \begin{enumerate}
  \item $\lim^j(\A(\pgl{n+1}),\pi_1B{\Ze}C_{\pgl{n+1}}) = 0$ for $j=1,2$,
  \item $\lim^j(\A(\pgl{n+1}),\pi_2B{\Ze}C_{\pgl{n+1}}) = 0$ for $j=2,3$,
  \end{enumerate}
for all $n \geq 1$.
\end{lemma}

For any \lmntwo\ $E$ in $\pgl{n+1}$ we shall write
\begin{equation*}
  [E] = \Hom_{\A(\pgl{n+1})(E)}(\St(E),\pi_1(B\Ze C_{\pgl{n+1}}(E)))
\end{equation*}
for the $\F_2$-vector space of $\F_2\A(\pgl{n+1})(E)$-module homo\m s
from the Steinberg representation $\St(E)$ over $\F_2$ of $\GL(E)$ to
$\pi_1(B\Ze C_{\pgl{n+1}}(E)$. 

Oliver's cochain complex for computing the first limits of the functor
$\pi_1(B\Ze C_{\pgl{n+1}})_{\not\leq t}$ has the form
\set\begin{equation}
  \label{bobccc}
  0 \to [H] \xrightarrow{d^1} 
  \prod_{1 \leq i \leq [m/2]}[H\#L[m-i,i]]  \xrightarrow{d^2}
  [H\#P[1,1,m-2] \times \prod_{2<i<[m/2]}[H\#P[1,i-1,m-i]]
\end{equation}\add
where we only list some of the nontoral rank four objects. Here,
\begin{align*}
  &[H] = \Hom_{\Symp(2)}(\St(H),H) \cong \F_2 \\
  &[H\#L[m-i,i]] = \Hom_{\Symp(3,1)}(\St(V),V) \cong \F_2, 
     \quad V=H\#L[m-i,i], \\
  &[H\#P[1,1,m-2]] = \Hom_{\begin{pmatrix} \Symp(2) & 0 \\ * & C_2
\end{pmatrix}}
(\St(E_2),E_2/E_2^{\perp}) \cong \F_2, \quad E_2=H\#P[1,1,m-2], \\
  &[H\#P[1,i-1,m-i]] =  \Hom_{\begin{pmatrix} \Symp(2) & 0 \\ * & 1
\end{pmatrix}}
(\St(E_i),E_i) \cong \F_2 \times \F_2, 
 \quad E_i=H\#P[1,i-1,m-i]
\end{align*}
where $ 2<i\leq [m/2]$ in the last line. The dimensions of these
spaces were found using the computer algebra program {\em magma}. It
suffices to show that the first differential $d^1$ is injective and
that the second differential $d^2$ has rank $[m/2]-1$. 
%% /home/moller/manus/2cgs/magma/steinberg/delta2m5.prg to
%% delta2m11.prg and delta1.prg

Let $H=\F_2e_1+\F_2e_2$ be a $2$-dimensional vector space over $\F_2$
with basis $\{e_1,e_2\}$ and symplectic
inner product matrix
\begin{equation*}
  \begin{pmatrix}
    0&1\\1&0
  \end{pmatrix}
\end{equation*}
Let $\F_2[1]$ be the $3$-dimensional $\F_2$-vector space on all length
zero flags $[L]$ of nontrivial proper subspaces $L \subset H$. The
Steinberg module $\St(H)$ for $H$ is the $2$-dimensional
$\F_2\GL(H)$-module that is the kernel for the augmentation
\func{d}{\F_2[1]}{\F_2} given by $d[L]=1$ for all $L$. Let
\func{f}{\St(H)}{H} be the restriction to $\St(H)$ of the
$\F_2\GL(H)$-module homo\m\ \func{\overline{f}}{\F_2[1]}{H} given by
$\overline{f}[L]=L$.

Let $V=\F_2e_1 + \F_2e_2 + \F_2e_3$ be a $3$-dimensional vector space
over $\F_2$ with basis $\{e_1,e_2,e_3\}$ and (degenerate) symplectic
inner product matrix
\begin{equation*}
\begin{pmatrix}
0 & 1 & 0 \\
1 & 0 & 0 \\
0 & 0 & 0
\end{pmatrix}
\end{equation*}
Let $\F_2[1]$ be the $21$-dimensional $\F_2$-vector space on all
length one flags $[P>L]$ and $\F_2[0]$ the $14$-dimensional
$\F_2$-vector space on all length zero flags, $[P]$ or $[L]$, of
non-trivial and proper subspaces of $V$.  The Steinberg module
$\St(V)$ over $\F_2$ for $V$ is the $2^3=8$-dimensional kernel of the
linear map \func{d}{\F_2[1]}{\F_2[0]} given by $d[P>L] = [P]+[L]$.
Define \func{df}{\St(V)}{V} to be the
restriction to $\St(V)$ of the linear map
\func{\overline{df}}{\F_2[1]}{V} given by
\set\begin{equation}\label{eq:defdfafam}
\overline{df}[P>L] =
\begin{cases}
L & P \cap P^{\perp}=\{0\} \\
0 & \text{otherwise}
\end{cases}
\end{equation}\add
on the basis vectors.

Let $E=\F_2e_1 + \F_2e_2 + \F_2e_3 + \F_2e_4$ be a $4$-dimensional
vector space over $\F_2$ with basis $\{e_1,e_2,e_3,e_4\}$ and
(degenerate) symplectic inner product matrix
\begin{equation*}
\begin{pmatrix}
0 & 1 & 0 & 0\\
1 & 0 & 0 & 0\\
0 & 0 & 0 & 0\\
0 & 0 & 0 & 0\\
\end{pmatrix}
\end{equation*}
Let $\F_2[2]$ be the $315$-dimensional $\F_2$-vector space on all
length two flags $[V>P>L]$ and $\F_2[1]$ the also $315$-dimensional
$\F_2$-vector space on all length one flags, $[P>L]$ or $[V>L]$ or
$[V>P]$, of non-trivial, proper subspaces of $E$.  The Steinberg
module $\St(E)$ over $\F_2$ for $E$ is the $2^6=64$-dimensional kernel
of the linear map \func{d}{\F_2[2]}{\F_2[1]} given by $d[V>P>L] =
[P>L]+[V>L]+[V>P]$. Define
\func{F_1=\overline{F}_1\vert\St(E)}{\St(E)}{E} as the restriction to
$\St(E)$ of the linear map \func{\overline{F}_1}{\F_2[2]}{E} with
values 
\set\begin{equation}\label{defF1}
  \overline{F}_1[V>P>L]=\begin{cases}
    L &  P \cap P^{\perp}=0, V \cap V^{\perp}=\F_2e_3 \\
    0 & \text{otherwise}
\end{cases}
\end{equation}\add
on the basis elements. 
%%/home/moller/manus/2cgs/magma/steinberg/delta2m8.prg
%%OrbitsOfSpaces(Aut,3);
Define 
\func{F_2=\overline{F}_2\vert\St(E)}{\St(E)}{E} similarly
but replace
the condition $V \cap V^{\perp}=\F_2e_3$ by $V \cap
V^{\perp}=\F_2e_4$. The linear maps $F_1$ and $F_2$ are $
\begin{pmatrix}
  \Symp(2) & 0 \\ * & 1
\end{pmatrix}$-equivariant because this group preserves the symplectic
inner product on $E$ and preserves $V^{\perp}= \F_2\gen{e_3,e_4}$
pointwise.

\begin{lemma}\label{lemma:fF1F2basis}
Let $f$ and $F_1,F_2$ be the linear maps defined above.
  \begin{enumerate}
  \item The vector $f$ is a basis for $[H]$.
   % \begin{equation*}
%      \Hom_{\Symp(2)}(\St(H),H) \cong \F_2
%    \end{equation*}
  \item The vector $df$ is a basis for
    $[H\#L[m-i,i]]$, $1 \leq i \leq [m/2]$.
%\begin{equation*}
%\Hom_{\begin{pmatrix}
%\Symp(2) & 0 \\ * & 1
%\end{pmatrix}}
%(\St(V),V) \cong
%\Hom_{\begin{pmatrix}
%\Symp(2) & 0 \\ * & 1
%\end{pmatrix}}
%(\St(V),V/V^{\perp}) \cong \F_2
%\end{equation*}
\item The vector $F_2$ is a basis for $[H\#P[1,1,m-2]]$.
\item The set $\{F_1,F_2\}$ is a basis for
  $[H\#P[1,i-1,m-i]]$, $2<i\leq [m/2]$.
%\begin{equation*}
%\Hom_{\begin{pmatrix}
%\Symp(2) & 0 \\ * & 1
%\end{pmatrix}}
%(\St(E),E) \cong
%\Hom_{\begin{pmatrix}
%\Symp(2) & 0 \\ * & 1
%\end{pmatrix}}
%(\St(E),E/E^{\perp}) \cong \F_2^2
%\end{equation*}
The sum $F_1+F_2$ is the linear map defined as in (\ref{defF1}) but with 
condition $V \cap V^{\perp}=\F_2e_3$ replaced by $V \cap
V^{\perp}=\F_2(e_3+e_4)$.
  \end{enumerate}
\end{lemma}

\begin{proof}
This can be directly verified by machine computation.
\end{proof}

\begin{proof}[Proof of Lemma~\ref{lim=0}]
  Since we already know that these higher limits vanish when $n+1$ is
  odd (\ref{tpglnC}, \ref{torallimforpglnC}) we can assume that 
  $n+1=2m$ is even.  
  
  \noindent (1) 
  See \ref{pgl2C} for the case $m=1$ and assume now that $m \geq
  2$. The image in $[H\#L[m-i,i]]$ of $f \in [H]$ is
  \begin{equation*}
  df_{L[m-i,i]}[P>L]=
  \begin{cases}
    L & P=H \\
    0 & \textmd{otherwise}
  \end{cases}
  \end{equation*}
  which equals $df$ (\ref{eq:defdfafam}). For $1<i \leq  [m/2]$, let
  \begin{equation*}
     ddf_{L[m-i,i]}[V>P>L]=
     \begin{cases}
       L & V=H\#L[m-i,i],\; P=H \\
       0 & \textmd{otherwise}
     \end{cases}
  \end{equation*}
  The object $H\#P[1,1,m-2]$ receives \m s from $H\#L[m-1,1]$ and
  (when $m>2$) $H\#L[m-2,2]$. Using a computer program one easily
  checks that $ddf_{L[m-1,1]} = F_2 = ddf_{L[m-2,2]}$ in
  $[H\#P[1,1,m-2]]$.  The object $H\#P[1,i-1,m-i]$ receives \m s from
  $H\#L[m-1,1]$, $H\#L[m-i+1,i-1]$, and $H\#L[m-i,i]$. Using a
  computer program one easily checks that $ddf_{L[m-i+1,i-1]}=F_1$,
  $ddf_{L[m-i,i]}=F_1$, and $ddf_{L[m-1,1]}=F_1+F_2$ in
  $[H\#P[1,i-1,m-i]$.  For $m=2$ or $m=3$, the cochain complexes
  (\ref{bobccc}) take the form
  \begin{align*}
      0 \to &[H] \xrightarrow{d^1} [H\#L[1,1]] \xrightarrow{d^2} 0 \\
      0 \to &[H] \xrightarrow{d^1} [H\#L[2,1]] \xrightarrow{d^2}
      [H\#P[1,1,1]] 
  \end{align*}
  where $d^1$ is an iso\m . For $m \geq 4$, and with our choice of
  basis (\ref{lemma:fF1F2basis}), the matrix for the differential
  $d^1$ is the injective $(1 \times [m/2])$-matrix
\begin{equation*}
  \begin{pmatrix}
    1 & 1 & \cdots & 1
  \end{pmatrix}
\end{equation*}
and the matrix for $d^2$ (or rather, the components of $d^2$ shown in
\ref{bobccc}) is the $([m/2] \times (2[m/2]-3))$-matrix 
\begin{center}
  \begin{tabular}[c]{r|cccc}
& $[H\#P[1,1,8]]$ & $[H\#P[1,2,7]]$ & $[H\#P[1,3,6]]$ &
$[H\#P[1,4,5]]$  \\ \hline
$[H\#L[9,1]]$ & $
\begin{pmatrix}
  1
\end{pmatrix}
$ & $
\begin{pmatrix}
  1 & 1
\end{pmatrix}$ &
$\begin{pmatrix}
  1 & 1
\end{pmatrix}$ &
$\begin{pmatrix}
  1 & 1
\end{pmatrix}$ \\
$[H\#L[8,2]]$ & $
\begin{pmatrix}
  1
\end{pmatrix}
$ & $
\begin{pmatrix}
  1 & 0
\end{pmatrix}$ & {} & {} \\
$[H\#L[7,3]]$ & {} &
$
\begin{pmatrix}
  0 & 1
\end{pmatrix}$ &
$
\begin{pmatrix}
  1 & 0
\end{pmatrix}$ & {} \\
$[H\#L[6,4]]$ & {} & {} &
$
\begin{pmatrix}
  0 & 1
\end{pmatrix}$ &
$
\begin{pmatrix}
  1 & 0
\end{pmatrix}$ \\
$[H\#L[5,5]]$ & {} & {} & {} &
$
\begin{pmatrix}
  0 & 1
\end{pmatrix}$
\end{tabular}
\end{center}
(shown here for $m=10$) of rank $[m/2]-1$.

\noindent
(2) Oliver's cochain complex for computing these higher limits over
$\A(\pgl{2m})$ involve the $\Z_2$-modules
(\ref{nonTpglnC}.(\ref{nonTpglnC2})) 
\begin{equation*}
  \Hom_{
    \begin{pmatrix}
      \Symp(2) & 0 \\ * & \A(\pgl{m})(E^{\perp})
    \end{pmatrix}}
   (\St(E), \pi_2(B{\Ze}C_{\pgl{2m}}(E^{\perp}))), \quad \dim_{\F_2}E=3,4,
\end{equation*}
that are submodules of finite products of $\Z_2$-modules of the form
\begin{equation*}
  \Hom_{
    \begin{pmatrix}
      \Symp(2) & 0 \\ * & 1
    \end{pmatrix}}
   (\St(E), \Z_2),  \quad \dim_{\F_2}E=3,4,
\end{equation*}
where the action on $\Z_2$ is trivial. According to the computer
program {\em magma}, these latter modules are trivial.
\end{proof}

%\begin{lemma}\label{lemma:1and2forA}
%  Suppose that $n+1=2m \geq 2$ is even.
%  \begin{enumerate}
%  \item There is a unique mono\m\ conjugacy class
%    \func{\lambda}{\Z/2}{\PGL(n+1,\C)} with disconnected centralizer.
%    The centralizer of this mono\m\ is ${\GL(m,\C)}^2/{\C^{\times}}
%    \rtimes \Z/2$ \label{1and2forA.1}
%  \item There is a unique mono\m\ conjugacy class
%    \func{\nu}{(\Z/2)²}{\PGL(n+1,\C)} such that $\nu$ is nontoral.
%    The centralizer of this mono\m\ is $(\Z/2)^2 \times \PGL(m,\C)$
%    and the Quillen auto\m\ group is $\Aut((\Z/2)²)$.
%  \end{enumerate}
%\end{lemma}
%\begin{proof}
%  This is contained in \ref{prop:afamtoral} and \ref{nonTpglnC}.
% % Use that any mono\m\ of $\Z/2$ into $\PGL(n+1,\C)$ lifts to
%%  $\mu \colon \Z/2 \to \GL(n+1,\C)$. The only possibility is that $\mu
%%  = m \cdot \mathrm{reg}$ is a direct sum of regular
%%  representations. The result for nontoral rank $2$ objects in
%%  $\A(\PGL(n+1,\C))$ is a special case of \ref{nonTpglnC}.
%\end{proof}

\chapter{The $D$-family}
\label{sec:dfam}

Let $\GL(2n,\R)$, $n\geq 1$, be the matrix group of $2n \times 2n$
real matrices and $\SL(2n,\R)$ the closed subgroup of matrices with
determinant $1$.
The $D$-family is the infinite family of matrix groups 
\begin{equation*}
  \PSL(2n,\R) = \frac{\SL(2n,\R)}{\gen{-E}}, \quad n \geq 4,
\end{equation*}
with trivial center. These groups also exist for $n=1,2,3$; however,
$\pslr{2}=\{1\}$ is the trivial group, and $\PSL(4,\R)=\PGL(2,\C)^2$,
$\PSL(6,\R)=\PGL(4,\C)$ are already known to be uniquely
$N$-determined (\ref{thm:afam}). 

The \mt , the \mtn\ of $\GL(2n,\R)$, $\SL(2n,\R)$,
and $\PSL(2n,\R)$ are\set
\begin{equation}\label{eq:dfamTNW}
  \begin{split}
  &T(\GL(2n,\R))=\SL(2,\R)^n,  \qquad\;\;\;
  N(\GL(2n,\R))=\GL(2,\R) \wr \Sigma_n
  %%\;\; W(\GL(2n,\R))=W(\GL(2,\R)) \wr \Sigma_n 
  \\
  &T(\SL(2n,\R))=\SL(2,\R)^n, \qquad\quad N(\SL(2n,\R))=\SL(2n,\R) \cap
  N(\GL(2n,\R)) \\
  &T(\PSL(2n,\R))=\frac{\SL(2,\R)^n}{\gen{-E}}, \qquad
   N(\PSL(2n,\R))=\frac{N(\SL(2n,\R))}{\gen{-E}}
\end{split}
\end{equation}\add
In all three cases, the \mtn\ is the semi-direct product for the
action of the Weyl group\set
\begin{equation}\label{eq:dfamW}
  \begin{split}
  &W(\GL(2n,\R))= 
  \Sigma_2\wr\Sigma_n, \quad \Sigma_2=W(\GL(2,\R))=\big\langle 
  \begin{pmatrix}
    0&1\\1&0
  \end{pmatrix}\big\rangle, \\
  &W(\SL(2n,\R))= 
  A_{2n} \cap (\Sigma_2\wr\Sigma_n) =
  W(\PSL(2n,\R))
 \end{split}
\end{equation}\add
on the \mt .
%where the Weyl group $W(\GL(2,\R))=\pi_0\GL(2,\R)$ has order two.
%The corresponding invariants for $\SL(2n,\R)$ are
%\begin{equation*}
%  T(\SL(2n,\R))=\SL(2,\R)^n, \quad N(\SL(2n,\R))=\SL(2n,\R) \cap
%  N(\GL(2n,\R))
%\end{equation*}
%The Weyl group $W(\SL(2n,\R))$ is the subgroup
%$S(C_2^n)\rtimes\Sigma_n$  of
%$W(\GL(2n,\R))=C_2 \wr \Sigma_n$ generated by 
%$S(C_2^n)=\{(x_1, \ldots ,x_n) \in C_2^n \mid x_1 \cdots x_n=1\}$,
%$C_2=\{\pm 1\}$, and $\Sigma_n$. 
It is known that for $n\geq 3$ \cite{bo7-8, hms:first, matthey:second,
  matthey:normalizers}\set
\begin{equation}\label{eq:dfamH01WT}
   H^0(W;\ch{T})(\PSL(2n,\R)) = 0, \qquad 
   H^1(W;\ch{T})(\PSL(2n,\R)) = 
   \begin{cases}
     \Z/2 & n=3 \\
     \Z/2 \times \Z/2 & n=4 \\
     0 & n>4
   \end{cases}
\end{equation}\add
for these projective groups. (The group of outer Lie auto\m s of the
Lie group $\PSL(8,\R)$, isomorphic to $\Sigma_3$, is faithfully
represented in $H^1(W;\ch{T})(\PSL(8,\R))$.) %% po8.prg.
%and that
%\begin{equation*}
%  \Out(\PSL(2n,\R)) = 
%  \begin{cases}
%    \Sigma_3 & n=4 \\
%    C_2 & n > 4
%  \end{cases}
%\end{equation*}
%The outer involution is induced by conjugation with a matrix of
%determinant $-1$. (The fixed point group varies with which matrix we
%take \cite[2.18]{griess:elem}.  What about the homotopy fixed point
%group?) It seems that $\Out(\pslr{8})$ is faithfully represented as a
%group of auto\m s of $H¹(W;\ch{T})$ -- see magma program po8.prg.

The Lie groups 
\begin{equation*}
  \GL(2n,\R) = \SL(2n,\R) \rtimes \gen{D},  \quad
  \PGL(2n,\R) = \PSL(2n,\R) \rtimes \gen{D\gen{-E}}
\end{equation*}
are the semi-direct products of their identity components with the
order two subgroup generated by the matrix $D=\diag(-1,1,\ldots ,1)$
(or any other order two matrix with negative determinant) and
conjugation with $D$ induces an outer auto\m\ of the Lie groups
$\SL(2n,\R)$ and $\PSL(2n,\R)$.

\section{The structure of $\pslr{2n}$}
\label{sec:structure}
In this section we investigate the Quillen category $\A(\pslr{2n})$
(\ref{defn:AX}) for $\pslr{2n}$ (and related \twocg s $\SL(2n,\R)$,
$\GL(2n,\R)$, $\PGL(2n,\R)$).
%%%%%%%%%%%%% Thu Feb 17 13:36:21 CET 2005%%%%%%%%%%%%%%%%%%%%%%
%Let $\A(\pslr{2n})$ be the category whose objects are$\A(\pslr{2n})$ all nontrivial
%\lmntwo s $V \subset \pslr{2n}$  and whose \m s are commutative
%diagrams 
%\begin{equation*}
%  \xymatrix@C=30pt{
%    V_1 \ar@{^(->}[r] \ar[d] & {\pslr{2n}} \ar[d]^M \\
%    V_2 \ar@{^(->}[r]  & {\pslr{2n}} }
%\end{equation*}
%where the arrow marked $M$ is conjugation with some matrix $M \in
%\PSL(2n,\R)$. The set of isomorphism classes of objects of the small
%category $\A(\pslr{2n})$ is the set of conjugacy classes of nontrivial
%elementary abelian $2$-subgroups of $\pslr{2n}$. The categories
%$\A(\SL(2n,\R))$ and $\A(\GL(2n,\R)$ are defined similarly.
%%%%%%%%%%%%%%%%%%%%%%%%%%%%%%%%%%%%%%%%%%%%%%%%%%%%%%%%%%%%%%%%%%%%

Consider the \lmntwo s\set
\begin{equation}\label{eq:dfamt}
  \begin{split}
    & t(\SL(2n,\R)) = t(\GL(2n,\R)) =
   \gen{e_1,\ldots ,e_n} \subset \SL(2n,\R) \subset \GL(2n,\R) \\
     &{\Delta_{2n}} =\gen{e_1,\ldots ,e_n,c_1,\ldots ,   c_n} 
  = \gen{\diag(\pm 1, \ldots , \pm 1)}
  \cong (\Z/2)^{2n} \subset \GL(2n,\R) \\
  &P\Delta_{2n}=\Delta_{2n}/\gen{e_1\cdots e_n} \cong (\Z/2)^{2n-1}
  \subset \PGL(2n,\R)\qquad (e_1\cdots e_n=-E) \\
  &S{\Delta_{2n}} =\gen{e_1,\ldots ,e_n,c_1c_2,\ldots , c_1c_n} = \SL(2n,\R)
  \cap {\Delta_{2n}}
  \cong
  (\Z/2)^{2n-1} \subset \SL(2n,\R) \\
  &PS{\Delta_{2n}} = S{\Delta_{2n}}/\gen{e_1\cdots e_n} \cong
  (\Z/2)^{2n-2} \subset \pslr{2n} \qquad (e_1\cdots e_n=-E) \\
   & t(\PSL(2n,\R)) = t(\PGL(2n,\R)) = \gen{I,e_1,\ldots
     ,e_n}/\gen{e_1 \ldots e_n} \subset \PSL(2n,\R) \subset
   \PGL(2n,\R) \\
   &  Pt(\SL(2n,\R)) = Pt(\GL(2n,\R)) =
   \gen{e_1,\ldots ,e_n}/\gen{e_1 \cdots e_n} \subset \SL(2n,\R) 
   \subset \GL(2n,\R) 
  \end{split}\add
\end{equation}
%%%%%%
%\marginpar{$(S\Delta_{2n})^D$ different??}
%%%%%%%%%%%%%%%%%%%%%071002%%%%%%%%%%%%%%%%%%%%%%%%%%%%%
where\set
\begin{equation}\label{eq:dfamIc}
  \begin{split}
  e_j &= \diag\left(\begin{pmatrix} 1&0\\0&1 \end{pmatrix},\ldots
    ,\begin{pmatrix} -1& \phantom{-}0\\ \phantom{-}0&-1 \end{pmatrix},\ldots ,
 \begin{pmatrix} 1&0\\0&1 \end{pmatrix} \right) \in \SL(2n,\R),\quad 1
\leq j \leq n  \\
I &=\diag\left(
  \begin{pmatrix}
    0 & -1 \\ 1 & 0
  \end{pmatrix}, \ldots ,
  \begin{pmatrix}
    0 & -1 \\ 1 & 0
  \end{pmatrix}\right) \in \slr{2n}, \\
  c_j &= \diag\left(\begin{pmatrix} 1&0\\0&1 \end{pmatrix},\ldots
    ,\begin{pmatrix} -1&\phantom{-}0\\\phantom{+}0&+1 \end{pmatrix},\ldots ,
   \begin{pmatrix} 1&0\\0&1 \end{pmatrix} \right) \in \GL(2n,\R),
 \quad 1 \leq j \leq n
\end{split}  
\end{equation}\add
%\marginpar{\bf Sign change in $c_j$!!}  
The matrices $e_j$ and $c_j$
have order two and commute with each other while $Ie_j=e_jI$,
$Ic_j=e_jc_jI$, and $I^2=e_1 \cdots e_n=-E$. 

The representation of the  Weyl groups
\set\begin{align}
  \label{WgroupsGL}
    W(\GL(2n,\R)) &= \gen{c_1,\ldots , c_n} \rtimes \Sigma_n =
    \Sigma_2\wr\Sigma_n, \\
    \label{WgroupsSL}
   W(\SL(2n,\R)) &= \gen{c_1c_2,\ldots ,c_1c_n} \rtimes \Sigma_n =
    A_{2n} \cap (\Sigma_2 \wr \Sigma_n)
\end{align}\add\add
on the maximal toral \lmntwo\ $t(\SL(2n,\R))=t(\GL(2n,\R))$ is trivial
on the subgroup $\gen{c_1,\ldots, c_n}=\Sigma_2^n$ while $\Sigma_n
\subset \GL(n,\C) \subset \SL(2n,\R)$ permutes the $n$ basis vectors
$e_1,\ldots ,e_n$ of $t(\SL(2n,\R))=t(\GL(2n,\R))$.
Let $V$ be a nontrivial \lmntwo\ in
$\PGL(2n,\R)$ and ${V^*}$ its inverse image in $\GL(2n,\R)$.  Let
\func{q}{V}{\F_2=\{0,1\}} be the function and \func{[\; , \;]}{V
  \times V}{\F_2=\{0,1\}} the bilinear map given by
${v^*}^2=(-E)^{q(v)}$ and $[{v}_1^*,{v}_2^*]= (-E)^{[v_1,v_2]}$ where
${v^*},{v}_1^*,{v}_2^*\in\SL(2n,\R)$ are preimages of
$v,v_1,v_2\in\PSL(2n,\R)$, respectively. The equations
\begin{equation*}
  [v_1,v_2]=[v_2,v_1], \quad [v,v]=0 ,\quad
  q(v_1+v_2)= q(v_1) + q(v_2) + [v_1,v_2]
\end{equation*}
show that $q$ is the quadratic function associated to the symplectic
bilinear form $[\; , \;]$ \cite[p.\ 356]{huppert:I}.  The
bilinear form is the deviation from linearity of the quadratic
function.  Define $V^{\perp} \supset R(V)$ to be the subgroups
\begin{equation*}
  V^{\perp} = 
  \{ v \in V \mid [v,V]=0 \} \supset
  \{ v \in V^{\perp} \mid q(v)=0 \} = R(V)
\end{equation*}
of $V$. Since $q$ is a group homo\m\ on $V^{\perp}$, the subgroup
$R(V)$ is either all of $V^{\perp}$ or a subgroup of index $2$.

In the following we write $G \circ H$ for the product of the groups
$G$ and $H$ with a common central subgroup amalgamated. The subgroup
$\mho_1(V^*)$ is generated by all squares of elements of $V^*$
\cite[III.10.4]{huppert:I}.

\begin{lemma}\label{lemma:Vast}
  Let $V$ be a nontrivial \lmntwo\ in $\PGL(2n,\R)$. %% $\pslr{2n}$. 
  The preimage $V^*$
  in $\GL(2n,\R)$ %%$\SL(2n,\R)$ 
  is
  \begin{equation*}
    V^* =
    \begin{cases}
      C_2 \times V     & q(V)=0 \\
      C_4 \circ V   & [V,V]=0,\; q(V) \neq 0 \\ %% C_4 \times R(V) 
      P \times R(V)    & [V,V] \neq 0,\; q(V^{\perp})=0\\
      (C_4 \circ P) \times R(V) &[V,V] \neq 0,\; q(V^{\perp}) \neq 0\\ 
    \end{cases}
  \end{equation*}
  where $C_2=\gen{-E} \subset C_4  \subset \SL(2n,\R)$,
  $P=2^{1+2d}_{\pm}$ is extraspecial, $C_4 \circ P$ is generalized
  extraspecial with center of order $4$, and $\mho_1(V^*)
  \subset \gen{-E}$.
 % and the first factor group in
%  each case has central center in $\SL(2n,\R)$. {\bf Only applies to
%    the factor groups - not the central product} $\mho_1(V^*) \subset
%  \gen{-E}$ is that what we want?
\end{lemma}
\begin{proof}
  %The central group extension $1 \to C_2 \to V^* \to V \to 1$ is
 % classified \cite[IV.3]{maclane} by the cohomology class represented
  %by the bilinear form $[\; , \; ] \in \F_2[V^{\vee}] = H^*(BV;\F_2)$.
%%%%%%%%%%%%%%%%%%%%%%%%%%%%%%%%
  As long as the bilinear form is trivial, $[V,V]=0$,
  $V^*$ is  abelian and the structure theorem for finitely generated
  abelian groups applies.  Assume that the bilinear form does not
  completely vanish, $[V,V] \neq 0$. Then $V^*$ is nonabelian with
  commutator subgroup $[V^*,V^*]=C_2$. Write $V = U \times R(V)$ for
  some nontrivial subgroup $U$ complementary to $R(V)$. Then
  $V^{\perp} = V^{\perp} \cap (U \times R(V)) = (V^{\perp} \cap U)
  \times R(V)$ and $q(V^{\perp}) = q(V^{\perp} \cap U)$. If $U^*$
  denotes the preimage of $U$, we have $V^* = U^* (C_2 \times R(V)) =
  U^* \times R(V)$ as the preimage of $R(V)$, $C_2 \times R(V)$, is
  central in $V^*$.  The commutator subgroup $[U^*,U^*] =
  [U^*R(V),U^*R(V)] = [V^*,V^*] = C_2$ and the center $Z(U^*)$ is the
  preimage of $V^{\perp} \cap U$.  If $q(V^{\perp}) =0$,
  $R(V)=V^{\perp}$ and $V^{\perp} \cap U = R(V) \cap U$ is trivial so
  $Z(U^*)=C_2$ and $U^*=P$ is extraspecial.  If $q(V^{\perp}) \neq 0$,
  $R(V)$ has index $2$ in $V^{\perp}$, $V^{\perp} \cap U$ has order
  $2$, and $q(V^{\perp} \cap U) \neq 0$ so that $Z(U^*)$ contains an
  element of order $4$. Therefore $Z(U^*)=C_4$ and $U^*$ is
  generalized extraspecial. There are two isomorphism classes of such
  groups but only $U^*=C_4 \circ D_8 \circ \cdots D_8=C_4 \circ P$ has
  elementary abelian abelianization \cite[Ex.\ 8, p.\ 
  146]{robinson:groups}.
\end{proof}

For instance, the preimage of the maximal toral \lmntwo\ 
$t(\PSL(2n,\R)$ of $\PSL(2n,\R)$ is the abelian group\set
\begin{equation}
  \label{eq:tPSL2nR} 
   t(\PSL(2n,\R))^*=\gen{I,e_1,\ldots ,e_n}, 
\end{equation}\add
generated by $I$ and $t(\SL(2n,\R))$.

\begin{cor}\label{cor:toralnontoral}
  Let $V$ be a nontrivial \lmntwo\ in $\pslr{2n}$. If 
  \begin{list}{}{}
  \item[{$q(V)=0$, $[V,V]=0$}:] $\text{$V$ is toral in $\pslr{2n}$ }
    \iff \text{$V^*=C_2 \times V$ is toral in $\SL(2n,\R)$}$
  \item[{$q(V) \neq 0$, $[V,V]=0$}:] $V$ is toral
  \item[{$q(V) \neq 0$, $[V,V] \neq 0$}:] $V$ is nontoral
  \end{list}
\end{cor}
\begin{proof}
  We have 
  \begin{equation*}
    \text{$V$ is toral} \iff
    V \subset t(\pslr{2n}) \iff
    V^* \subset  t(\pslr{2n})^* 
  \end{equation*}
  where the symbol \lq$\subset$\rq\ reads \lq is subconjugate to\rq .
  In the first case of the corollary, the preimage $V^*$ contains no
  elements of order $4$ so that
  \begin{equation*}
     V^* \subset  t(\pslr{2n})^* \iff
     V^* \subset t(\SL(2n,\R))
  \end{equation*}
  as $ t(\SL(2n,\R))$ consists of the elements of order $\leq 2$ in
  $t(\pslr{2n})^*$. In the second case, $V^* = C_4 \times R(V)$ so
  that $R(V) \subset C_{\SL(2n,\R)}(I) = \GL(n,\C)$. But any complex
  representation  of the \lmntwo\ $R(V)$ is toral, so $R(V) \subset
  t(\GL(n,\C)) = t(\SL(2n,\R))$ and $V^* \subset
  \gen{C_4,t(\SL(2n,\R))} = t(\pslr{2n})^*$. In the third case, the
  nonabelian group $V^*$ can not be a subgroup of the abelian group
  $t(\pslr{2n})^*$.  
\end{proof}

%%%%%%%%%%%%%%%%%%%%%%%%%% 30.09.02 %%%%%%%%%%%%%%%%%%%%%%%%%

\begin{lemma} \label{lemma:Vconj} %%25.08.03
  Let $V_1$ and $V_2$ be \lmntwo s in $\pslr{2n}$. Then
  \begin{equation*}
    \text{$V_1$ and $V_2$ are conjugate in $\pslr{2n}$} \iff
    \text{$V_1^*$ and $V_2^*$ are conjugate in $\slr{2n}$}
  \end{equation*}
  where $V_1^*,V_2^* \subset \slr{2n}$ are the preimages.
\end{lemma}
\begin{proof}
  This is clear.
\end{proof}

Write $\A(\PGL(2n,\R))^{q=0}$ and $\A(\PGL(2n,\R))^{\leq t,q=0}$ for
the full subcategories of $\A(\PGL(2n,\R))$ generated by all \lmntwo s
$V \subset \PGL(2n,\R)$ with trivial quadratic function $q$,
respectively, all toral \lmntwo s $V \subset \PGL(2n,\R)$ with trivial
quadratic function $q$. Define  $\A(\PSL(2n,\R))^{q=0}$ and
$\A(\PSL(2n,\R))^{\leq t,q=0}$ similarly as full subcategories of
$\A(\PSL(2n,\R))$.

\begin{lemma}\label{cor:equivcat}
  Write $\GL$ for $\GL(2n,\R)$, $\SL$ for $\SL(2n,\R)$, and $\PSL$ for
  $\PSL(2n,\R)$. The
  inclusion functors
  \begin{alignat*}{2}
  &\A(\Sigma_{2n},{\Delta_{2n}}) \to \A(\GL) &
  &\A(\Sigma_{2n},S{\Delta_{2n}}) \to \A(\SL) \\
  % &\A(\Sigma_2\wr\Sigma_n, t(\GL)) \to  \A(\GL)^{\leq t}&
  & &
  \quad &
  \A(W(\SL),t(\SL)) \to \A(\SL)^{\leq t} \\
  &\A(\Sigma_{2n},P{\Delta_{2n}}) \to \A(\PGL)^{q=0} & 
  &\A(\Sigma_{2n},PS{\Delta_{2n}}) \to \A(\PSL)^{q=0} \\
  % &\A(\Sigma_2\wr\Sigma_n, t(\PGL)) \to  \A(\PGL)^{\leq t}&
  &&
  \quad & \A(W(\PSL),t(\PSL)) 
  \to \A(\PSL)^{\leq t} \\
  % &\A(\Sigma_2\wr\Sigma_n, Pt(\GL)) \to \A(\PGL)^{\leq t,q=0}& 
  &&
  \quad & \A(W(\PSL),Pt(\SL)) \to \A(\PSL)^{\leq
    t, q=0} 
  \end{alignat*}
  are equivalences of categories. In particular, $\A(\SL)$ and
  $\A(\PSL)$ are full subcategories of $\A(\GL)$ and $\A(\PGL)$,
  respectively.  (See \ref{defn:AWt} for the meaning of
  $\A(\Sigma_{2n},\Delta_{2n})$.)
\end{lemma}
\begin{proof}
  By real representation theory any nontrivial \lmntwo\ of
  $\GL(2n,\R)$ is conjugate to a subgroup $V$ of $\Delta_{2n}$ and
  \begin{equation*}
    C_{\GL(2n,\R)}(V)= \prod_{\rho\in V^{\vee}} \GL(i_{\rho},\R)
  \end{equation*}
  where \func{i}{V^{\vee}}{\Z} records the multiplicity of $\rho\in
  V^{\vee}$ in the representation $V \subset \Delta_{2n} \subset
  \GL(2n,\R)$. Observe that $\Delta_{2n}$ is the maximal \lmntwo\ in
  $C_{\GL(2n,\R)}(V)$.  (For any $i \geq 1$, $\GL(i,\R)$ contains the
  subgroup $\Delta_i$, consisting of diagonal matrices with $\pm 1$ in
  the diagonal, as a maximal \lmntwo .) Therefore, by the standard
  argument from \cite[IV.2.5]{brockerdieck}, used also in
  \ref{lemma:afamtoralcat}, any group homo\m\ between two
  nontrivial subgroups of $\Delta_{2n}$ induced by conjugation
  with a matrix from $\GL(2n,\R)$, is in fact induced by conjugation
  with a matrix from $N_{\GL(2n,\R)}(\Delta_{2n}) = \Delta_{2n}
  \rtimes \Sigma_{2n}$ \cite[Lemmma 3]{bob:stubborn}. Thus
  the inclusion functor $\A(\Sigma_{2n},\Delta_{2n}) \to
  \A(\GL(2n,\R))$ is a category equivalence.

   Any nontrivial \lmntwo\ $V \subset \PGL(2n,\R)$ with $q(V)=0$ is
  conjugate to a subgroup of $P\Delta_{2n}$ since $V^*$, the preimage
  in $\GL(2n,\R)$, is conjugate to subgroup of $\Delta_{2n}$.  Let
  $V_1,V_2$ be two nontrivial subgroups of $P\Delta_{2n}$. From the
  commutative diagram of \m\ sets
  \begin{equation*}
    \xymatrix{
    {\A(\Sigma_{2n},\Delta_{2n})(V_1^*,V_2^*)} \ar@{->>}[d] \ar@{=}[r] &
    {\A(\GL(2n,\R))(V_1^*,V_2^*)}  \ar@{->>}[d] \\
    {\A(\Sigma_{2n},P\Delta_{2n})(V_1,V_2)} \ar@{^(->}[r] &
    {\A(\PGL(2n,\R))^{q=0}(V_1,V_2)}}
  \end{equation*}
  we see that the the bottom horizontal arrow is a bijection. This
  implies that $\A(\Sigma_{2n},P\Delta_{2n}) \to
  \A(\PGL(2n,\R))^{q=0}$ is an equivalence of categories.

 % In particular for the subgroup $t(\GL) \subset \Delta_{2n}$
%  (\ref{eq:dfamt}) we have
%  that  (see \ref{defn:AWt} for notation),
%  \begin{equation*}
%    \overline{\Sigma}_{2n}(t(\GL),t(\GL))
%    \stackrel{\textmd{\cite[Lemma 3]{bob:stubborn}}}{=}
%  \Sigma_2\wr\Sigma_n, \qquad
%    \Sigma_{2n}(t(\GL))=\Sigma_2^n
%  \end{equation*}
%  so that $\A(\Sigma_2\wr\Sigma_n,t(\GL))\to\A(\GL(2n,\R))^{\leq t}$
%  is an equivalence of categories by an argument similar to the one
%  just given.
  
  Any nontrivial \lmntwo\ in $\SL(2n,\R)$ is conjugate in $\GL(2n,\R)$
  to a subgroup of $\SL(2n,\R) \cap \Delta_{2n}=S\Delta_{2n}$
  (\ref{eq:dfamt}).  The Quillen category of $\SL(2n,\R)$ is a full
  subcategory of the Quillen catgory of $\GL(2n,\R)$ since
  $C_{\GL(2n,\R)}(V) \not\subset \SL(2n,\R)$ for all objects $V$ of
  $\A(\SL(2n,\R))$.  Thus the inclusion functor
  $\A(\Sigma_{2n},S\Delta_{2n}) \to \A(\SL(2n,\R))$ is an equivalence
  of categories.

  Any toral \lmntwo\ in $\SL(2n,\R)$ is conjugate to a subgroup of
  $t(\SL(2n,\R))$ by its very definition (\ref{defn:toralAX}).  Any
  \m\ between two nontrivial subgroups of $t(\SL(2n,\R))$ induced by
  conjugation with a matrix from $\SL(2n,\R)$, is in fact induced by
  conjugation with a matrix from $N(\SL(2n,\R))$ and hence from
  $W(\SL(2n,\R))$ \cite[IV.2.5]{brockerdieck}. Thus
  $\A(W(\SL),t(\SL))\to\A(\SL(2n,\R))^{\leq t}$ is a category
  equivalence. The same argument can be used to identify the toral
  subcategory for $\PSL(2n,\R)$ (and it is actually a general fact that
  the inclusion functor $\A(W(X),t(X)) \to \A(X)^{\leq t}$ is an
  equivalence of categories where $t(X) \to X$ is the maximal toral
  elementary abelian $p$-group in the connected \pcg\ $X$
  \cite[2.8]{jmm:ndet}).
 
  Any nontrivial toral \lmntwo\ $V \subset \PSL(2n,\R)$ with $q(V)=0$
  is conjugate to a subgroup of $Pt(\SL)$ (\ref{eq:dfamt}) since
  $V^*$, the preimage (\ref{lemma:Vast}) in $\GL(2n,\R)$, is conjugate
  to subgroup of $t(\SL) \subset t(\PSL)^*$ (\ref{eq:tPSL2nR}). As
  $\A(\PSL)^{\leq t, q=0}$ is a full subcategory of $\A(\PSL)^{\leq t}
  = \A(W(\PSL),t(\PSL))$, this means that
  $\A(W(\PSL), Pt(\SL)) \to \A(\PSL)^{\leq t,q=0}$ is a
  category equivalence.

\end{proof}

We now specialize to full subcategory $\A(\PSL(2n,\R))^{\leq t}_{\leq
  2}$ (\ref{defn:toralAX}).

\begin{prop}\label{prop:Apslr2nleqt}
The chart
\begin{center}
\begin{tabular}{|c||c|c|c|c|}\hline
  \raisebox{-12pt}{$\A(\PSL(2n,\R))^{\leq t}_{\leq 2}$} 
  & \multicolumn{2}{c|}{\raisebox{-4pt}{Lines}} &
  \multicolumn{2}{c|}{\raisebox{-4pt}{Planes}} \\ \cline{2-5} 
  & $q=0$ & $q \neq 0$ &  $q=0$ & $q \neq 0$ \\ \hline\hline
  $n$ even & $n/2$ & $2$ & $P(n,3)+P(n,4)$ & $n/2+[n/4]$ \\ \hline
  $n$ odd  & $[n/2]$ & $1$ & $P(n,3)+P(n,4)$ & $[n/2]$ \\ \hline 
\end{tabular}
\end{center}
gives the number of iso\m\ classes of toral objects of rank $1$ and
$2$ in $\A(\PSL(2n,\R))$.  

When $n$ is even, the $\frac{n}{2}$ toral lines with $q=0$ are
$L(2i,2n-2i)$, $1 \leq i \leq \frac{n}{2}$, and the two toral lines
with $q \neq 0$ are $I$ and $I^D$. The toral planes with $q=0$ are the
planes $P(2i_0,2i_1,2i_2,0)$ where $(i_0,i_1,i_2)$ is a partition of
$n$ into three natural numbers, $P(2i_0,2i_1,2i_2,2i_3)$ where
$(i_0,i_1,i_2,i_3)$ is a partition of $n$ into four natural numbers,
and the toral planes with $q \neq 0$ are $I\#L(i,n-i)$, $1 \leq i \leq
\frac{n}{2}$, and $I\#L(i,n-i)^D$ for even $i$. 

When $n$ is odd, the $\left[\frac{n}{2}\right]$ toral lines with $q=0$
are $L(2i,2n-2i)$, $1 \leq i \leq \left[\frac{n}{2}\right]$, and the
toral line with $q \neq 0$ is $I$. The toral planes
with $q=0$ are the planes $P(2i_0,2i_1,2i_2,0)$ where $(i_0,i_1,i_2)$
is a partition of $n$ into three natural numbers,
$P(2i_0,2i_1,2i_2,2i_3)$ where $(i_0,i_1,i_2,i_3)$ is a partition of
$n$ into four natural numbers, and the toral planes with $q \neq 0$
are $I\#L(i,n-i)$, $1 \leq i \leq \left[\frac{n}{2}\right]$.

In (\ref{CLqeq0}) and (\ref{CLneq0}) we list the centralizers of the
rank one objects and in (\ref{CVqeq0}) and (\ref{CVqneq0}) the
centralizers of the rank two objects.
\end{prop}

Proposition~\ref{prop:Apslr2nleqt} is the conclusion of the following
considerations.

%\section{The category $\A(\pslr{2n})^{\leq t}_{\leq 2}$.}
%\label{sec:toralsubcat}

%In this section we concentrate on the full subcategory of
%$\A(\pslr{2n})$ generated by toral objects of rank one or two.  The
%relevant properties of $\A(\pslr{2n})^{\leq t}_{\leq 2}$ are listed
%below.

%Let $e_1$, $e_2$ be the standard generators of $\Z/2 \oplus \Z/2$. Put
%$e_3=e_1+e_2$ and let $e_0$ be the trivial element. Let
%\func{\rho_0,\rho_1,\rho_2,\rho_3}{\Z/2 \oplus \Z/2}{\R^{\times}} be
%the homo\m s given by
%\begin{center}
%  \begin{tabular}[r]{|l||c|c|c|c|} \hline
%    {} & $e_0$ & $e_1$ & $e_2$ & $e_3$ \\ \hline \hline
%    $\rho_0$ & $+1$ & $+1$ & $+1$ & $+1$ \\ \hline 
%    $\rho_1$ & $+1$ & $-1$ & $+1$ & $-1$ \\ \hline 
%    $\rho_2$ & $+1$ & $+1$ & $-1$ & $-1$ \\ \hline 
%    $\rho_3$ & $+1$ & $-1$ & $-1$ & $+1$ \\ \hline 
%  \end{tabular}
%\end{center}

For any partition $i=(i_0,i_1)$ of $n=i_0+i_1$ into a sum of two
positive integers $i_0\geq i_1 \geq 1$ let $L[i]=L[2i_0,2i_1] \subset
t(\SL(2n,\R)) \subset \SL(2n,\R)$ be the toral subgroup generated by
\begin{equation*}
  \diag(\overbrace{+E,\ldots,+E}^{i_0},\overbrace{-E,\ldots,-E}^{i_1})
  %%= (2i_0\rho_0 + 2i_1\rho_1)(e_1)
\end{equation*}
Then the centralizer (of the image in $\PSL(2n,\R)$) of this subgroup  is 
\set\begin{equation}\label{CLqeq0}
  C_{\PSL(2n,\R)}L[2i_0,2i_1]=
  \begin{cases}
    \frac{\SL(2i_0,\R)\times\SL(2i_i,\R)}{\gen{-E}}\rtimes\gen{\diag(D_1,D_2)} 
    & i_0 \neq i_1 \\
   \frac{\SL(2i_0,\R)^2}{\gen{-E}}\rtimes
            \gen{\diag(D_1,D_2),
              \begin{pmatrix} O&E\\E&0 \end{pmatrix}}  & i_0 = i_1
  \end{cases}
\end{equation}\add
where $D_j=\diag(-1,1,\ldots,1)\in\GL(2i_j,\R)$ are matrices of
determinant $-1$. The diagonal matrix $\diag(D_1,D_2)$ acts on the
identity component of the centralizer by the outer action on both
factors.  In the second case, which only occurs when $n=2i_0$ is even,
the matrix $\begin{pmatrix} O&E\\E&0 \end{pmatrix}$ acts by permuting
the factors.

The element $I \in t(\PSL(2n,\R))^* \subset \SL(2n,\R)$ of order four
generates an order two toral subgroup of $\pslr{2n}$ with centralizer
\cite[5.11]{jmm:ndet} 
\set\begin{equation}\label{CLneq0}
  C_{\PSL(2n,\R)}(I) =
  \begin{cases}
    {\GL(n,\C)}/\!{\gen{-E}} & \text{$n$ odd} \\
    {\GL(n,\C)}/\!{\gen{-E}} \rtimes \gen{c_1 \cdots c_n} &  \text{$n$ even}
  \end{cases}
\end{equation}\add
where, in the even case, the component group acts on the
identity component through the unstable Adams operation
$\psi^{-1}$. The nontrivial outer auto\m\ of $\pslr{2n}$ takes $I$ to
$I^D$ where $I \neq I^D$ if and only if $n$ is even
(\ref{exmp:IDQ}.(\ref{exmp:IDQ4})). 

For any partition $i=(i_0,i_1,i_2,0)$ of $n=i_0+i_1+i_2$ into a sum of
three positive integers $i_0\geq i_1 \geq i_2 >0$ or any partition
$i=(i_0,i_1,i_2,i_3)$ of $n=i_0+i_1+i_2+i_3$ into a sum of four
positive integers $i_0\geq i_1 \geq i_2 \geq i_3 >0$ let $P[i] =
P[2i_0,2i_1,2i_2,2i_3] \subset t(\SL(2n,\R)) \subset \SL(2n,\R)$ be
the subgroup generated by the two elements
\begin{align*}
 &\diag(\overbrace{+E, \ldots, +E}^{i_0},
       \overbrace{-E, \ldots ,-E}^{i_1},
       \overbrace{+E, \ldots ,+E}^{i_2},
       \overbrace{-E, \ldots ,-E}^{i_3}) 
      %% =(i_0\rho_0+i_1\rho_1+i_2\rho_2+i_3\rho_3)(e_1)   
                    \\
  &\diag(\overbrace{+E, \ldots, +E}^{i_0},
        \overbrace{+E, \ldots ,+E}^{i_1},
        \overbrace{-E, \ldots ,-E}^{i_2},
        \overbrace{-E, \ldots ,-E}^{i_3}) 
       %%=(i_0\rho_0+i_1\rho_1+i_2\rho_2+i_3\rho_3)(e_2)
\end{align*}
%%%%%%%%%%%%%%%%%%%%
%%%%%%%%% action on the \rho s Tue Jan 25 16:07:31 CET 2005
%%%%%%%%%%%%%%%%%%%%%%%%%%%%%%% 
%The action of
%$V^{\vee}=\{\rho_0,\rho_1,\rho_2,\rho_3\}$ on
%$i=(2i_0,2i_1,2i_2,2i_3)$ is given by
%\begin{align*}
% % \rho_0 \cdot i =(2i_0,2i_1,2i_2,2i_3), \quad
%  &\rho_1 \cdot i =(2i_1,2i_0,2i_3,2i_2) &
%  &\rho_2 \cdot i =(2i_2,2i_3,2i_0,2i_1)  &
%  &\rho_3 \cdot i =(2i_3,2i_2,2i_1,2i_0) \\
%&\xymatrix@1@C=25pt{
%{\bullet} \ar@{<->}@/_/[r] & {\bullet} & 
%                {\bullet}\ar@{<->}@/^/[r] & {\bullet}  } &  
%&\xymatrix@1@C=25pt{
%{\bullet} \ar@{<->}@/_/[rr] & {\bullet} \ar@{<->}@/^/[rr]
%         & {\bullet} & {\bullet}  } & 
%&\xymatrix@1@C=25pt{
%{\bullet} \ar@{<->}@/_/[rrr] & {\bullet} \ar@{<->}@/^/[r]
%         & {\bullet} & {\bullet}  } 
%\end{align*}
%%%action on the \rho s
%The action of $\GL(V)\cong\Sigma_3$ on $i=(i_0,i_1,i_2,i_3)$ fixes
%$i_0$ and permutes the other three indices $i_1,i_2,i_3$. Since the
%two actions together generate all permutations of the four symbols
%$i_0, i_1, i_2, i_3$, there are $P(n,3)+P(n,4)$ iso\m\ classes of such
%objects where $P(n,k)$ is the number of $k$-integer partitions of the
%integer $n$, namely $P(n,3)$ objects of the form $(i_0,i_1,i_2,0)$
%with $i_0+i_1+i_2=n$ and $i_0 \geq i_1 \geq i_2>0$, and $P(n,4)$
%objects of the form $(i_0,i_1,i_2,i_3)$ with $i_0+i_1+i_2+i_3=n$ and
%$i_0 \geq i_1 \geq i_2 \geq i_3>0$. 
%%% (the equalizer condition seems to be redundant in rank $\leq 2$.)
%%%%%%%%%%%%%%%%%%%
%%%%%%%%%% END action of the \rho s Tue Jan 25 16:07:31 CET 2005
%%%%%%%%%%%%%%%%%%%%%%%%%%%%%%%%%
The centralizers in $\pslr{2n}$ are 
\set\begin{equation} \label{CVqeq0} C_{\pslr{2n}}P(i)=
  \begin{cases}
    \frac{\SL(2i_0,\R)² \times \SL(2i_2,\R)²}{\gen{-E,-E,-E,-E}} \rtimes 
    \left( \ker\big( C_2^{S(i)} \to C_2 \big) \rtimes
    \Z/2 \right) & i=(2i_0,2i_0,2i_2,2i_2) \\
    \frac{\SL(2i_0,\R)^4}{\gen{-E,-E,-E,-E}} \rtimes
    \left( \ker\big( C_2^{S(i)} \to C_2 \big) \rtimes 
    (\Z/2 \times \Z/2) \right) & i=(2i_0,2i_0,2i_0,2i_0) \\
 \frac{\prod_{S(i)} \SL(2i_j,\R)}{\gen{-E}} \rtimes 
   \ker\big(C_2^{S(i)} \to C_2 \big)  & 
    \text{otherwise}
  \end{cases}
\end{equation}\add
where $\ker\big( C_2^{S(i)} \to C_2 \big) = \gen{ \diag(D_1,D_2,E,E),
  \diag(D_1,E,D_3,E), \diag(D_1,E,E,D_4)}$ (when $\# S(i)=4$) is
generated by diagonal matrices, $D_j=\diag(-1,1,\ldots
,1)\in\GL(2i_j,\R)$, and
\begin{equation*}
  \Z/2 = \gen{\begin{pmatrix}
      0 & E & 0 & 0 \\
      E & 0 & 0 & 0 \\
      0 & 0 & 0 & E \\
      0 & 0 & E & 0
    \end{pmatrix}}, \quad
  \Z/2 \times \Z/2 = \gen{
    \begin{pmatrix}
      0 & E & 0 & 0 \\
      E & 0 & 0 & 0 \\
      0 & 0 & 0 & E \\
      0 & 0 & E & 0
    \end{pmatrix},
    \begin{pmatrix}
      0 & 0 & E & 0 \\
      0 & 0 & 0 & E \\
      E & 0 & 0 & 0 \\
      0 & E & 0 & 0
    \end{pmatrix}}
\end{equation*}
are generated by block permutation matrices. 
(The component group of
the first line is $C_2 \times D_8$; the component group of second line
is extra special of order $32$ isomorphic to $D_8 \circ D_8$.) 

For any partition $i=(i_0,i_1)$ of $n=i_0+i_1$ into a sum of two
positive integers $i_0\geq i_1 >0$ let $I\#L[i_0,i_1] \subset
\PSL(2n,\R)$ be the \lmntwo\ that is the quotient of
\begin{equation*}
  (I\#L[i_0,i_1])^*=\big\langle I,
  \diag(\overbrace{+E,\ldots,
    +E}^{i_0},\overbrace{-E,\ldots,-E}^{i_1})\big\rangle \subset
  t(\PSL(2n,\R))^* 
\end{equation*}
It follows that
  \set\begin{equation}\label{CVqneq0}
  C_{\pslr{2n}}I\#L(i_0,i_1)=
  \begin{cases}
    \frac{\GL(i_0,\C) \times \GL(i_1,\C)}{\gen{-E,-E}} & \text{$n$ odd}
    \\
    \frac{\GL(i_0,\C) \times \GL(i_1,\C)}{\gen{-E,-E}} \rtimes \gen{c_1
      \cdots c_n} & \text{$n$ even, $i_0 \neq i_1$} \\
    \frac{\GL(i,\C) \times\GL(i,\C)}{\gen{-E,-E}} \rtimes 
    \gen{c_1 \cdots c_n, P}
    & \text{$n$ even, $i_0=i_1$}
  \end{cases}
\end{equation}\add
where $P =  \begin{pmatrix}
        0 & E \\ E & 0
      \end{pmatrix}$ permutes the two identical factors. 

\begin{prop}\label{prop:IL}
   $I\#L(i,n-i) \neq  I\#L(i,n-i)^D$ if and only if $n$ and $i$ are
   even. 
\end{prop}
\begin{proof}
  The auto\m\ group of $\gen{i} \times \gen{\varepsilon} = C_4 \times
  C_2 =  I\#L(i,n-i)^*$ is the dihedral group of order eight
  \begin{equation*}
    \Aut(C_4 \times C_2)=\gen{a,b \mid a^4, b^2, bab=a^3}
  \end{equation*}
  generated by the two auto\m s given by $a(i)=i\varepsilon$,
  $a(\varepsilon)=i^2\varepsilon$ and $b(i)=i$,
  $b(\varepsilon)=i^2\varepsilon$. The auto\m\ $a^2 \in Aut(C_4)
  \subset \Aut(C_4 \times C_2)$ is induced by conjugation with the
  matrix
\begin{equation*}
  \diag(P,\ldots,P), \qquad P=
  \begin{pmatrix}
    0&1\\1&0
  \end{pmatrix},
\end{equation*}
of determinant $(-1)^n$. Thus $A(\SL(2n,\R))(I\#L(i,n-i)^*) \neq
A(\GL(2n,\R))(I\#L(i,n-i)^*)$ and $I\#L(i,n-i)=I\#L(i,n-i)^D$ when $n$
is odd (\ref{lemma:intoSL}).

Assume now that $n$ is even. The group of trace preserving auto\m s
\begin{equation*}
  \A(\GL(2n,\R))(C_4 \times C_2)=
  \begin{cases}
    \gen{a^2,ba} & 2i<n \\
    \Aut(C_4 \times C_2) & 2i=n
  \end{cases}
\end{equation*}
has index $2$ in general but is actually equal to the full auto\m\
group in case $i=n/2$. The conjugating matrix for $ba$ is
\begin{equation*}
  \diag(\overbrace{P,\ldots,P}^i,\overbrace{E,\ldots,E}^{n-i})
\end{equation*}
of determinant $(-1)^i$. Thus  $I\#L(i,n-i)=I\#L(i,n-i)^D$ when $i$ is
odd. If $n=2i$ then the conjugating matrices for the auto\m s $a$ and
$b$ are
\begin{equation*}
  \begin{pmatrix}
    0&E\\E&0
  \end{pmatrix} 
  \diag(\overbrace{P,\ldots,P}^i,\overbrace{E,\ldots,E}^i) \qquad
  \text{and} \qquad 
  \begin{pmatrix}
    0&E\\E&0
  \end{pmatrix}
\end{equation*}
The permutation matrix for $b$ has positive determinant and the matrix
for $a$ has determinant $(-1)^i$. Thus $I\#L(i,n-i)=I\#L(i,n-i)^D$ if
and only if $i$ is odd.
\end{proof}

\section{Centralizers of objects of $\A(\pslr{2n})^{\leq t}_{\leq
    2}$ are LHS}
\label{sec:lhs}

In this section we check that all toral objects of rank $\leq 2$ have
LHS (\ref{subsec:LHS}) centralizers.

\begin{lemma}\label{lemma:lhs}
  The centralizers of the objects of $\A(\pslr{2n})^{\leq t}_{\leq
    2}$, 
  \begin{enumerate}
  \item ${\GL(i,\C)}/\!{\gen{-E}} \rtimes C_2$, $1 \leq i$ 
                              \eqref{CLneq0}\label{lhs1}
  \item ${\SL(2i_0,\R)\circ\SL(2i_1,\R)}\rtimes C_2$, $1 \leq i_0 <
    i_1$ \eqref{CLqeq0}\label{lhs2}
  \item $(\SL(2i,\R) \circ \SL(2i,\R))\rtimes (C_2 \times C_2)$,
    $1 \leq i$ \eqref{CLqeq0}\label{lhs3}
  \item  \label{lhs4} $C_{\pslr{2n}}(V)$, $q(V)=0$ \eqref{CVqeq0}
  \item  \label{lhs5} $C_{\pslr{2n}}(V)$, $q(V) \neq 0$ \eqref{CVqneq0}
 % \item $\big(\prod_{j=0}^3 \slr{2i_j}/\gen{-E}\big) \rtimes C_2^3$, 
%   % $(\slr{2i_0} \times \slr{2i_1} \times \slr{2i_2} \times
%%    \slr{2i_3})/\gen{-E} \rtimes C_2^3$, 
%    $i_0, i_1, i_2, i_3>1$ --
%    what happens when at least one of the factors is $\slr{2}$?
%    \label{lhs4} 
  \end{enumerate}
  are LHS.
\end{lemma}

The cases of interest here are summarized in the following charts,
obtained by use of a computer, for rank one centralizers with quadratic
form $q=0$ (\ref{CLqeq0})
%% /home/moller/manus/dfam/magma/lhs/newlhs/sl2sl2i.prg
%% /home/moller/manus/dfam/magma/lhs/newlhs/sl4sl2i.prg
%% /home/moller/manus/dfam/magma/lhs/newlhs/sl6sl2i.prg
  \begin{center}
    \begin{tabular}[c]{|c||c|c|c|c|} \hline
      $\slr{2i_0} \circ \slr{2i_1}$ & $\ker \theta$ &
      $\Hom(W,\ch{T}^W)$ & $H^1(W;\ch{T})$ & $\theta$ \\ \hline
      $1=i_0,2=i_1$ & $(\Z/2)^2$ & $(\Z/2)²$ & $\Z/2$ & $0$ \\ \hline
      $1=i_0,3=i_1$ & $0$ & $\Z/2$ & $(\Z/2)^2$ & mono \\ \hline
      $1=i_0,4\leq i_1$ & $0$ & $\Z/2$ & $\Z/2$ & iso\\ \hline
      $2=i_0<i_1$ & $(\Z/2)^2$ & $(\Z/2)^3$ & $\Z/2$ &epi\\ \hline
      $3 \leq i_0 < i_1$ & $0$ & $(\Z/2)^2$ &  $(\Z/2)^2$ &iso \\ \hline
    \end{tabular}
  \end{center}

\begin{center}
  \begin{tabular}[c]{|c||c|c|c|c|} \hline
    $\slr{2i} \circ \slr{2i}$ & $\ker \theta$ & $\Hom(W,\ch{T}^W)$ 
    & $H^1(W;\ch{T})$ & $\theta$ \\ \hline
    $i=2$ & $(\Z/2)^4$ & $(\Z/2)^4$ & $(\Z/2)^3$ & $0$ \\ \hline
    $i \geq 3$ & $0$ &  $(\Z/2)^2$ &  $(\Z/2)^2$ & iso\\ \hline
  \end{tabular}
\end{center}

and $q \neq 0$ (\ref{CLneq0})

\begin{center}
    \begin{tabular}[c]{|c||c|c|c|c|} \hline
      $\GL(i,\C)/\gen{-E}$ 
      & $\ker \theta$ & $\Hom(W;\ch{T}^W)$ & $H^1(W;\ch{T})$ &
                   $\theta$  \\ \hline \hline
      $i=2$ & $\Z/2$ & $(\Z/2)²$ & $\Z/2$ &epi \\ \hline 
      $i=3$ & $0$ & $\Z/2$ & $\Z/2$ &iso\\ \hline
      $i=4$ & $0$ & $\Z/2$ & $(\Z/2)^2$ &mono\\ \hline
      $i>4$ & $0$ & $\Z/2$ & $\Z/2$ &iso\\ \hline
    \end{tabular}
\end{center}

and for rank two centralizers with quadratic form $q=0$ (\ref{CVqeq0})

\begin{center}
  \begin{tabular}[c]{|c||c|c|c|c|} \hline
 $\slr{2i_0}^2 \circ \slr{2i_1}^2$ &
    $\ker \theta$ & $\Hom(W;\ch{T}^W)$ & $H^1(W;\ch{T})$ &
                $\theta$ \\ \hline \hline
   $1=i_0, 2=i_1$ & $(\Z/2)^4$  & $(\Z/2)^{12}$  & $(\Z/2)^8$ &epi \\
   \hline  
   $1=i_0, 2 < i_1$ & $0$  & $(\Z/2)^{6}$  & $(\Z/2)^6$ &iso \\
   \hline 
   $2=i_0 < i_1$ & $(\Z/2)^4$  & $(\Z/2)^{18}$  & $(\Z/2)^{14}$ &epi \\
   \hline 
    $2 < i_0 < i_1$ & $0$  & $(\Z/2)^{12}$  & $(\Z/2)^{12}$ &iso \\
   \hline 
  \end{tabular}
\end{center}

\begin{center}
  \begin{tabular}[c]{|c||c|c|c|c|} \hline  
  $\prod_{j=0}^2 \SL(2i_j,\R)/\gen{-E}$ &
  $\ker \theta$ & $\Hom(W;\ch{T}^W)$ & $H^1(W;\ch{T})$ &
                $\theta$ \\ \hline \hline
  $1=i_0, 2=i_1<i_2$ & 
  $(\Z/2)^2$ & $(\Z/2)^6$  & $(\Z/2)^4$ &epi \\ \hline              
  $1=i_0, 2<i_1<i_2$ & 
  $0$ & $(\Z/2)^4$  & $(\Z/2)^4$ &iso \\ \hline  
  $2=i_0<i_1<i_2$ & 
  $(\Z/2)^2$ & $(\Z/2)^8$  & $(\Z/2)^6$ &epi \\ \hline
  $2<i_0<i_1<i_2$ & 
  $0$ & $(\Z/2)^6$  & $(\Z/2)^6$ &iso \\ \hline                         
  \end{tabular}
\end{center}

\begin{center}
  \begin{tabular}[c]{|c||c|c|c|c|} \hline  
  $\prod_{j=0}^3 \SL(2i_j,\R)/\gen{-E}$ &
  $\ker \theta$ & $\Hom(W;\ch{T}^W)$ & $H^1(W;\ch{T})$ &
                $\theta$ \\ \hline \hline
  $1=i_0, 2=i_1<i_2<i_3$ & 
  $(\Z/2)^2$ & $(\Z/2)^{12}$  & $(\Z/2)^{10}$ &epi \\ \hline              
  $1=i_0, 2<i_1<i_2<i_3$ & 
  $0$ & $(\Z/2)^9$  & $(\Z/2)^9$ &iso \\ \hline  
  $2=i_0<i_1<i_2<i_3$ & 
  $(\Z/2)^2$ & $(\Z/2)^{15}$  & $(\Z/2)^{13}$ &epi \\ \hline
  $2<i_0<i_1<i_2<i_3$ & 
  $0$ & $(\Z/2)^{12}$  & $(\Z/2)^{12}$ &iso \\ \hline                         
  \end{tabular}
\end{center}

\begin{center}
  \begin{tabular}[c]{|c||c|c|c|c|} \hline
 $\slr{2i}^4/\gen{-E}$ &
    $\ker \theta$ & $\Hom(W;\ch{T}^W)$ & $H^1(W;\ch{T})$ &
                $\theta$ \\ \hline \hline
 $2=i$ & $(\Z/2)^8$ & $(\Z/2)^{24}$ & $(\Z/2)^{16}$ &epi \\ \hline
 $3 \leq i$ & $0$ & $(\Z/2)^{12}$ & $(\Z/2)^{12}$ &iso \\ \hline
  \end{tabular}
\end{center}

and with quadratic form $q \neq 0$ (\ref{CVqneq0})

\begin{center}
  \begin{tabular}[c]{|c||c|c|c|c|} \hline
    $\GL(i_0,\C) \circ \GL(i_1,\C)$ &
    $\ker \theta$ & $\Hom(W;\ch{T}^W)$ & $H^1(W;\ch{T})$ &
                $\theta$ \\ \hline\hline
    $1=i_0,2=i_1$ & $\Z/2$ & $(\Z/2)^2$ & $\Z/2$ &epi\\ \hline
    $1=i_0, 2<i_1$ & $0$ &  $(\Z/2)^2$ &  $(\Z/2)^2$ &iso\\ \hline
    $2=i_0<i_1$ & $\Z/2$ &  $(\Z/2)^4$ &  $(\Z/2)^3$ &epi\\ \hline
    $2< i_0<i_1$ & $0$ & $(\Z/2)^4$ & $(\Z/2)^4$ &iso\\ \hline
  \end{tabular}
\end{center}

\begin{center}
  \begin{tabular}[c]{|c||c|c|c|c|} \hline
    $\GL(i,\C) \circ \GL(i,\C)$ &
    $\ker \theta$ & $\Hom(W;\ch{T}^W)$ & $H^1(W;\ch{T})$ &
               $\theta$  \\ \hline \hline
    $2=i$ & $(\Z/2)^2$ &  $(\Z/2)^4$ &  $(\Z/2)^3$ & {} \\ \hline
    $3 \leq i$ & $0$ &  $(\Z/2)^4$ &  $(\Z/2)^4$ &iso\\ \hline
  \end{tabular}
\end{center}

%Let $Y$ be a connected \twocg\ containing $Z$ as a central subgroup.
%Put $X=Y/Z$. We consider a \m\ of \twocg s of the form $Y_{h\pi} \to
%X_{h\pi}$ where $\pi$ is a finite $2$-group. The commutative diagram
%\begin{equation*}
%  \xymatrix{
%    {\ch{Z}} \ar[r] \ar@{=}[d] &
%    {\ch{Z}(Y)} \ar[r] \ar[d] &
%    {\ch{Z}(X)} \ar[d] \\
%    {\ch{Z}} \ar[r] &
%    {\ch{T}(Y)} \ar[r] &
%    {\ch{T}(X)} } 
%\end{equation*}
%relating centers and maximal tori, induces a commutative diagram in
%cohomology
%\begin{equation*}
%  \xymatrix{
%    H^i(\pi;\ch{Z}) \ar[r] \ar@{=}[d] &
%    H^i(\pi;\ch{Z}(Y)) \ar[r] \ar[d] &
%    H^i(\pi;\ch{Z}(X)) \ar[r] \ar[d] &
%    H^{i+1}(\pi;\ch{Z}) \ar@{=}[d] \\
%    H^i(\pi;\ch{Z}) \ar[r] &
%    H^i(\pi;\ch{T}(Y)) \ar[r] &
%    H^i(\pi;\ch{T}(X)) \ar[r] &
%    H^{i+1}(\pi;\ch{Z}) }
%\end{equation*}
%with exact rows.  A diagram chase reveals that
%\set\begin{equation}\label{upstairs}
%  \begin{split}
%  & \text{$H^i(\pi;\ch{Z}(X)) \to H^i(\pi;\ch{T}(X)) $ is injective}  \\
%  & \quad \quad \iff 
%  \ker \left( H^i(\pi;\ch{Z}(Y)) \to H^i(\pi;\ch{T}(Y)) \right) \subset 
%  \im \left( H^i(\pi; \ch{Z}) \to  H^i(\pi; \ch{Z}(Y)) \right)  \\ 
%  & \quad \quad \Leftarrow 
%   \text{$H^i(\pi;\ch{Z}(Y)) \to H^i(\pi;\ch{T}(Y)) $ is injective}  
%  \end{split}
%\end{equation}\add
%We shall apply this small observation in the proof of the next lemma.

\begin{proof}[Proof of Lemma~\ref{lemma:lhs}]
  \eqref{lhs1} Let $X=\GL(i,\C)/\!{\gen{-E}} \rtimes C_2$ for $i \geq
1$. Since the Weyl group for $X$ is a direct product $W=W_0 \times
C_2$, $X$ is LHS. 
%We record the situation for $i=4$  %%why should we do that?
% \begin{center}
%    \begin{tabular}[c]{|c||c|c|c|c|} \hline
%      $\GL(i,\C)/\!{\gen{-E}} \rtimes C_2$  
%      & $H^1(\pi;\ch{T}^{W_0})$ & $H^1(W;\ch{T})$ & $H^1(W_0;\ch{T})$
%      & $H^1(W_0;\ch{T})^{\pi}$ \\ \hline \hline
%      $i=4$ & $0$ & $(\Z/2)^2$ & $(\Z/2)^2$ & $(\Z/2)^2$ \\
%      \hline
%    \end{tabular}
%  \end{center} 
% as obtained by explicit computer computations.

  \noindent 
  \eqref{lhs2} Let  $X=(\SL(2i_0,\R)\circ\SL(2i_1,\R))\rtimes
  C_2$ for $1 \leq i_0 < i_1$. The first problematic case is when $i_0=1$
  and $i_1=2$ or $3$. In this case, explicit computer
  computation results in the chart
  \begin{center}
    \begin{tabular}[c]{|c||c|c|c|c|} \hline
      $\SL(2,\R)\circ\SL(2i_1,\R))\rtimes C_2$
      & $H^1(\pi;\ch{T}^{W_0})$ & $H^1(W;\ch{T})$ & $H^1(W_0;\ch{T})$
      & $H^1(W_0;\ch{T})^{\pi}$ \\ \hline \hline
      $i_1=2$ & $0$ & $\Z/2$ & $\Z/2$ & $\Z/2$ \\
      \hline
      $i_1=3$ & $0$ & $(\Z/2)^2$ & $(\Z/2)^2$ & $(\Z/2)^2$ \\ \hline
    \end{tabular}
  \end{center}
  showing that $X$ is LHS. %%% for these two values of $i_1$. 
  The second problematic case is $2=i_0<i_1$ where $\theta(X_0)$ is
  epimorphic. Since $H^1(W_0;\ch{T})=\Z/2$, also $\theta(X_0)^{\pi}$ is
  epimorphic so that $X$ is LHS (\ref{lhscrit1}).

  \noindent
  \eqref{lhs3} Let $X=(\SL(2i,\R) \circ \SL(2i,\R))\rtimes (C_2 \times
  C_2)$ for $i \geq 1$. $X$ is a \twoctg\ when $i=1$ and hence
  obviously LHS\@. 
  For $i \geq 2$ explicit computer computation gives
   \begin{center}
    \begin{tabular}[c]{|c||c|c|c|c|} \hline
      $(\slr{2i} \circ \slr{2i}) \rtimes (C_2 \times C_2)$
      & $H^1(\pi;\ch{T}^{W_0})$ & $H^1(W;\ch{T})$ & $H^1(W_0;\ch{T})$
      & $H^1(W_0;\ch{T})^{\pi}$ \\ \hline \hline
      $i=2$ & $(\Z/2)^2$ & $(\Z/2)^4$ & $(\Z/2)^{3}$ & $(\Z/2)^2$ 
          \\ \hline
      $i \geq 3$ &  $(\Z/2)^2$ & $(\Z/2)^3$ & $(\Z/2)^{2}$ & $\Z/2$ 
          \\ \hline
    \end{tabular}
  \end{center}
  so $X$ is manifestly LHS for $i=2$. For $i>2$, $\theta(X_0)$ is
  bijective so $X$ is LHS (\ref{lhscrit1}).

\noindent
   \eqref{lhs4} Let $X=(\slr{2i}^4/\gen{-E}) \rtimes (D_8 \circ
   D_8)$ for $i \geq 1$. When $i=1$, $X$ is a \twoctg\ which are all
   LHS\@. When $i=2$ explicit computer computation gives
   \begin{center}
    \begin{tabular}[c]{|c||c|c|c|c|} \hline
      $(\slr{2i}^4/\gen{-E}) \rtimes (D_8 \circ D_8)$
      & $H^1(\pi;\ch{T}^{W_0})$ & $H^1(W;\ch{T})$ & $H^1(W_0;\ch{T})$
      & $H^1(W_0;\ch{T})^{\pi}$ \\ \hline \hline
      $i=2$ & $(\Z/2)^7$ & $(\Z/2)^9$ & $(\Z/2)^{16}$ & $(\Z/2)^2$ 
          \\ \hline
    \end{tabular}
  \end{center}
  so $X$ is LHS by definition. For $i>2$, $\theta(X_0)$ is bijective.

  \noindent
  (5) Let $X=(\SL(2i_0,\R)² \circ \SL(2i_1,\R)²) \rtimes (C_2 \times
  D_8)$ for $1 \leq i_0 < i_1$. The problematic cases are
  $i_0=1,i_1=2$ and $2=i_0<i_1$ where $\theta(X_0)$ is surjective but
  not bijective. With the help of computer computations we obtain the table
  %% sl2^2sl4^2.prg,  sl4^2sl6^2.prg, sl4^2sl8^2.prg 
   \begin{center}
    \begin{tabular}[c]{|c||c|c|c|c|} \hline
      $(\slr{2i_0}^2 \circ \slr{2i_1}^2) \rtimes (D_8 \times C_2)$
      & $H^1(\pi;\ch{T}^{W_0})$ & $H^1(W;\ch{T})$ & $H^1(W_0;\ch{T})$
      & $H^1(W_0;\ch{T})^{\pi}$ \\ \hline \hline
      $i_0=1,i_1=2$ & $(\Z/2)^5$ & $(\Z/2)^7$ & $(\Z/2)^{8}$ & $(\Z/2)^2$ 
          \\ \hline
      $i_0=2,3 \leq i_1$ & 
      $(\Z/2)^6$ & $(\Z/2)^{11}$ & $(\Z/2)^{14}$ & $(\Z/2)^{5}$ 
          \\ \hline    
    \end{tabular}
  \end{center}
  showing that $X$ is LHS in these cases also.

  \noindent (6) 
  Let $X= \frac{\prod_{j=0}^2 \SL(2i_j,\R)}
  {\gen{-E}} \rtimes C_2^2$. The problematic cases are $1=i_0,
  2=i_1<i_2$ and $2=i_0<i_1<i_2$. With the help of computer
  computations we obtain the table
   \begin{center}
    \begin{tabular}[c]{|c||c|c|c|c|} \hline
      $\frac{\prod_{j=0}^3 \SL(2i_j,\R)} {\gen{-E}} \rtimes C_2^2$ &
      $H^1(\pi;\ch{T}^{W_0})$ & $H^1(W;\ch{T})$ & $H^1(W_0;\ch{T})$
      & $H^1(W_0;\ch{T})^{\pi}$ \\ \hline \hline
      $i_0=1,2=i_1<i_2$ & $(\Z/2)^3$ & $(\Z/2)^6$ & $(\Z/2)^4$ &
      $(\Z/2)^3$  
          \\ \hline
      $2=i_0<i_1<i_2$ & 
      $(\Z/2)^4$ & $(\Z/2)^{9}$ & $(\Z/2)^{6}$ & $(\Z/2)^{5}$ 
          \\ \hline    
    \end{tabular}
   \end{center}
   showing that $X$ is LHS in these cases also.

  \noindent (7) 
  Let $X= \frac{\prod_{j=0}^3 \SL(2i_j,\R)} {\gen{-E}} \rtimes C_2^3$.
  The problematic cases are $1=i_0, 2=i_1<i_2<i_3$ and
  $2=i_0<i_1<i_2<i_3$.  
  With the help of computer computations we obtain the
  table
  \begin{center}
    \begin{tabular}[c]{|c||c|c|c|c|} \hline
      $\frac{\prod_{j=0}^3 \SL(2i_j,\R)} {\gen{-E}} \rtimes C_2^3$ &
      $H^1(\pi;\ch{T}^{W_0})$ & $H^1(W;\ch{T})$ & $H^1(W_0;\ch{T})$
      & $H^1(W_0;\ch{T})^{\pi}$ \\ \hline \hline
      $i_0=1,2=i_1<i_2<i_3$ & $(\Z/2)^8$ & $(\Z/2)^{16}$ & $(\Z/2)^{10}$ &
      $(\Z/2)^8$  
          \\ \hline
      $2=i_0<i_1<i_2<i_3$ & 
      $(\Z/2)^9$ & $(\Z/2)^{20}$ & $(\Z/2)^{13}$ & $(\Z/2)^{11}$ 
          \\ \hline    
    \end{tabular}
   \end{center}
   showing that $X$ is LHS in these cases also.

  \noindent (8)
  The \twocg\ $(\GL(i,\C) \circ \GL(i,\C)) \rtimes (\Z/2 \times \Z/2)$
  is LHS for $i>2$ where $\theta$ is bijective. When $i=2$ we find 
\begin{center}
    \begin{tabular}[c]{|c||c|c|c|c|} \hline
      $(\GL(i,\C) \circ \GL(i,\C)) \rtimes (C_2 \times C_2)$
      & $H^1(\pi;\ch{T}^{W_0})$ & $H^1(W;\ch{T})$ & $H^1(W_0;\ch{T})$
      & $H^1(W_0;\ch{T})^{\pi}$ \\ \hline \hline
      $i=2$ & $(\Z/2)^2$ & $(\Z/2)^4$ & $(\Z/2)^3$ & $(\Z/2)^2$ 
          \\ \hline
    \end{tabular}
\end{center}
so $X$ is also LHS in this case. 
%(It could be that $\theta$ for $X_0$
%maps onto the invariants  $H^1(W_0;\ch{T})^{\pi}$?)

  \noindent (9)
  Let $X=\GL(i_0,\C) \circ \GL(i_1,\C) \rtimes C_2$, $1 \leq i_0 <
  i_1$. Since the identity component has surjective $\theta$-homo\m\ and
  the component group $\pi=C_2$ acts trivially on $H^1(W_0;\ch{T})$,
  $X$ is LHS by \ref{lhscrit1}. The values of the relevant
  cohomology groups are %%checklhsglC2.prg
  \begin{center}
    \begin{tabular}[c]{|c||c|c|c|c|} \hline
      $(\GL(i_0,\C) \circ \GL(i_1,\C)) \rtimes C_2$
      & $H^1(\pi;\ch{T}^{W_0})$ & $H^1(W;\ch{T})$ & $H^1(W_0;\ch{T})$
      & $H^1(W_0;\ch{T})^{\pi}$ \\ \hline \hline
      $1=i_0, 2=i_1$ & 
      $0$ & $\Z/2$ & $\Z/2$ & $\Z/2$ 
          \\ \hline
      $1=i_0, 2<i_1$ & 
      $0$ & $(\Z/2)^2$ & $(\Z/2)^2$ & $(\Z/2)^2$ 
          \\ \hline
      $2=i_0<i_1$ & 
      $0$ & $(\Z/2)^3$ & $(\Z/2)^3$ & $(\Z/2)^3$ 
          \\ \hline
      $2<i_0<i_1$ & 
      $0$ & $(\Z/2)^4$ & $(\Z/2)^4$ & $(\Z/2)^4$ 
          \\ \hline
    \end{tabular}
\end{center}
   according to computer computations.
\end{proof}

\section{The limit of the functor $H^1(W_0;\ch{T})^{W/W_0}$ on
  $\A(\pslr{2n})^{\leq t}_{\leq 2}$}
\label{sec:lim0}

Let \func{H^1(W_0;\ch{T})}{\A(\pslr{2n})^{\leq t}_{\leq t}}{\Ab} be
the functor that takes the toral \lmntwo\ $V \subset t(\pslr{2n})$ to
the abelian group $H^1(W_0(C_{\pslr{2n}}(V);\ch{T}))$, and 
$H^1(W_0;\ch{T})^{W/W_0}$ the functor that takes $V$ to the the
  invariants for the action of the component group
  $\pi_0C_{\pslr{2n}}(V)$ on this first cohomology group
  (\ref{sec:cond2}). 
%% \marginpar[ind6.prg etc]{ind6.prg etc}
%% /home/moller/manus/dfam/magma/lhs/ind6.prg etc
%% /home/moller/manus/dfam/magma/lhs/indgenn.prg for the general case
  \begin{prop}\label{prop:lim0H1}
    The restriction map
    \begin{equation*}
      H^1(W(\pslr{2n};\ch{T}) \to \lim^0(\A(\pslr{2n})^{\leq t}_{\leq
        2},H^1(W_0;\ch{T})^{W/W_0}) 
    \end{equation*}
    is an iso\m\ for all $n \geq 4$.
  \end{prop}
  \begin{proof}
    Consider first that case where $n=4$.  The \twocg\ $X=\pslr{8}$
    contains (\ref{prop:Apslr2nleqt}) the four rank one \lmntwo s
    $L(2,6), L(4,4), I , I^D$ with centralizers $\slr{2} \circ \slr{6}
    \rtimes C_2$, $\slr{4} \circ \slr{4} \rtimes (C_2 \times C_2)$,
    $\GL(4,\C)/\gen{-E}$ (twice).  The claim of the proposition
    follows from the fact, verifiable by computer computations, that
    in all four cases, the restriction $H¹(W;\ch{T}(X)) \to
    H^1(W_0(C_X(L));\ch{T})^{W/W_0}$ happens to be iso\m .

    For $n>4$, the claim is that the limit of the functor
    $H^1(W_0;\ch{T})^{W/W_0}$ is trivial. In fact, even the limit of
    the functor  $H^1(W_0;\ch{T})$ is trivial. To see this, recall
    (\ref{prop:Apslr2nleqt}) 
    that $\pslr{2n}$ contains the toral lines $L(2i,2n-2i)$, $1 \leq i
    \leq [n/2]$, $I$ and also $I^D$ when $n$ is even. Computer
    computations show that the \m s
    \begin{equation*}
      H^1(W_0;\ch{T})(L(2,2n-2)) \hookrightarrow
      H^1(W_0;\ch{T})(I\#L(1,n-1)) \hookleftarrow
      H^1(W_0;\ch{T})(I)
    \end{equation*}
    are injective with and that their images 
    intersect trivially. When $n \geq 6$ is even, also the images of
    the injective \m s
    \begin{equation*}
      H^1(W_0;\ch{T})(L(4,2n-4)) \hookrightarrow
      H^1(W_0;\ch{T})(I\#L(2,n-2)^D) \hookleftarrow
      H^1(W_0;\ch{T})(I^D)
    \end{equation*}
   intersect trivially. More computer computations show that,
   similarly, the \m s
   \begin{equation*}
      H^1(W_0;\ch{T})(L(2i,2n-2i)) \hookrightarrow 
      H^1(W_0;\ch{T})(I\#L(i,n-i)) \hookleftarrow
      H^1(W_0;\ch{T})(I), \qquad 1 \leq i \leq [n/2],
   \end{equation*}
   are injective with and that their images intersect trivially.
  \end{proof}

\section{Rank two nontoral objects of $\A(\pslr{2n})$}
\label{sec:NdetD}

%When $n$ is odd, there are no nontoral planes in $\pslr{2n}$ for there
%are no mono\m s $2^{1+2}_{\pm} \to \slr{2n}$ with central $\mho_1$
%(\ref{exmp:IDQ}). Well, there could be a nontoral $\Z/2 \times V \to
%\slr{2n}$? See (\ref{cor:toralnontoral}). Yes, but then $q=0$ so that
%is plays no role for the computations of the higher limits, however,
%we still need to check uniqueness of the lift $\nu_L'$ and
%$f_{\nu,L}$.  

In this section we take a closer look at the nontoral rank two objects
of $\A(\pslr{2n})$ in order to verify the conditions of
\ref{lemma:uniquenu}.

Nontoral rank two objects $P$ of $\pslr{2n}$ satisfy either $q(P)=0$
or $[P,P] \neq 0$ (\ref{cor:toralnontoral}) and the latter case only
occurs if $n$ is even.

\noindent
\underline{$q(P)=0$}: 
For any partition $i_0 \geq i_i \geq i_2 \geq i_3 \geq 1$, let
\begin{align*}
  & P[i_0,i_1,i_2,i_3]^*= \gen{(+1)^{i_0}(-1)^{i_1}(+1)^{i_2}(-1)^{i_3},
                           (+1)^{i_0}(+1)^{i_1}(-1)^{i_2}(-1)^{i_3},
                           -E} \subset \Delta_{2n}, \\
  & P[i_0,i_1,i_2,i_3] = P[i_0,i_1,i_2,i_3]^*/\gen{-E} \subset P\Delta_{2n}
\end{align*}
where we apply the notation from notation from
\ref{cha:afam}.\S\ref{sec:pglnc}.  Then $ P[i_0,i_1,i_2,i_3]^* \subset
S\Delta_{2n}$ if and only if $i_0$, $i_1$, $i_2$, and $i_3$ all have
the same parity and $ P[i_0,i_1,i_2,i_3]^*$ is nontoral iff this
parity is odd. It follows that the set of iso\m\ classes of nontoral
rank two objects of $\A(\pslr{2n})^{q=0}$ corresponds bijectively to
the $P(n+2,4)$ partitions $i=[i_0,i_1,i_2,i_3]$ of $n+2$ into sums of
$4$ natural numbers, $n+2=i_0+i_1+i_2+i_3$, $i_0 \geq i_1 \geq i_2
\geq i_3 \geq 1$. The correspondence is via the map
\begin{equation*}
  i=[i_0,i_1,i_2,i_3] \to P[i]=P[2i_0-1,2i_1-1,2i_2-1,2i_3-1] 
\end{equation*}
that to the partition $i=[i_0,i_1,i_2,i_3]$ associates the quotient
$P[i] \subset PS\Delta_{2n}$ of $P[i]^* \subset S\Delta_{2n}$ generated
by the three elements
\begin{align*}
  v_1&=
  \diag(\overbrace{+1,\ldots,+1}^{2i_0-1},\overbrace{-1,\ldots,-1}^{2i_1-1}, 
        \overbrace{+1,\ldots,+1}^{2i_2-1},\overbrace{-1,\ldots,-1}^{2i_3-1}) \\
  v_2&=
  \diag(\overbrace{+1,\ldots,+1}^{2i_0-1},\overbrace{+1,\ldots,+1}^{2i_1-1},
        \overbrace{-1,\ldots,-1}^{2i_2-1},\overbrace{-1,\ldots,-1}^{2i_3-1}),\\
  v_3&=\diag(-1,\ldots,-1)
\end{align*}
The centralizer of $P[i]^*$ in $\SL(2n,\R)$ is
\begin{multline*}
  C_{\SL(2n,\R)}(P[i]^*)=\SL(2n,\R) \cap  C_{\SL(2n,\R)}(P[i]^*) \\
  = \SL(2n,\R) \cap \Big(\prod_{j=0}^3\GL(2i_j-1,\R)\Big) 
  = P[i]^* \times \prod_{j=0}^3\SL(2i_j-1,\R)  
\end{multline*}
and centralizer of $P[i]$ in $\PSL(2n,\R)$ is therefore
\cite[5.11]{jmm:ndet}
\set\begin{equation}\label{eq:CPi1}
   C_{\PSL(2n,\R)}(P[i])=P[i] \times
   \Big(\prod_{j=0}^3\SL(2i_j-1,\R)\Big) \rtimes P[i]^{\vee}_i
\end{equation}\add
where $P[i]^{\vee}_i$ is a group of permutation matrices isomorphic
to $C_2$ if $i=[i_0,i_0,i_2,i_2]$, to $C_2 \times C_2$ if
$i=[i_0,i_0,i_0,i_0]$, and  trivial in all other cases. Note that
$P[i]^*$ is contained in $N(\SL(2n,\R))=\SL(2n,\R) \cap
\GL(2,\R)\wr\Sigma_n$ because we may write
\set\begin{align}
  v_1&=\diag(\overbrace{E,\ldots,E}^{i_0-1},R,\overbrace{-E,\ldots,-E}^{i_1-1},
       \overbrace{E,\ldots,E}^{i_2-1},R,\overbrace{-E,\ldots,-E}^{i_3-1}),
       \qquad R=
       \begin{pmatrix}
         1&0\\0&-1 
       \end{pmatrix} \label{eq:e1} \\
  v_2&=\diag(\overbrace{E,\ldots,E}^{i_0-1},E,\overbrace{E,\ldots,E}^{i_1-1},
       \overbrace{-E,\ldots,-E}^{i_2-1},-E,\overbrace{-E,\ldots,-E}^{i_3-1}) 
       \label{eq:e2}
\end{align}\add\add
and that the centralizer of $P[i]^*$ there is
\begin{multline*}
  C_{N(\SL(2n,\R))}(P[i]^*)=\SL(2n,\R) \cap
  C_{\GL(2,\R)\wr\Sigma_n}(v_1) \cap C_{\GL(2,\R)\wr\Sigma_n}(v_2)\\
  \stackrel{\text{(\ref{lemma:centrgw})}}{=} 
  \SL(2n,\R) \cap 
  C_{\GL(2,\R)\wr(\Sigma_{i_0+i_1-1} \times \Sigma_{i_2+i_3-1})}(v_1) 
  \stackrel{\text{(\ref{lemma:centrgw})}}{=} P[i]^* \times
  \Big(\prod_{j=0}^3 \GL(2,\R)\wr\Sigma_{i_j-1}\Big)
\end{multline*}
It follows that the centralizer of $P[i]$ in $N(\pslr{2n})$ is
\begin{equation*}
  C_{N(\PSL(2n,\R))}(P[i])= P[i] \times  \Big(\prod_{j=0}^3
  \GL(2,\R)\wr\Sigma_{i_j-1}\Big) \rtimes P[i]^{\vee}_i =
  N(C_{\pslr{2n}}(P[i]))  
\end{equation*}
For instance, if $i=[i_0,i_0,i_2,i_2]$, then $ P[i]^{\vee}_i$ is the
group of order two generated by $\diag(C_0,C_2) \in N(\pslr{2n})$
where $C_0$ is the $(4i_0-2) \times (4i_0-2)$ matrix
\begin{equation*}
  \begin{pmatrix}
    0&0&E\\0&T&0\\E&0&0
  \end{pmatrix}, \qquad T=
  \begin{pmatrix}
    0&1\\1&0
  \end{pmatrix}
\end{equation*}
and $C_2$ is a similar $(4i_2-2) \times (4i_2-2)$ matrix. Thus $P[i]
\subset N(\pslr{2n})$ is a \pl\ \cite{jmm:normax} of $P[i] \subset
\pslr{2n}$. The two other \pl s \cite[6.2]{jmm:deter} of $P[i] \subset
\pslr{2n}$ are obtained by composing the inclusion $P[i] \subset
N(\pslr{2n})$ with the inner auto\m\ given by the permutation matrices
$(1,2)(2i_0,2n-2i_3+1)$ 
\begin{center}
\begin{picture}(100,80)(0,-55)
\put(-50,0){$(\overbrace{+1,+1,\ldots,+1}^{2i_0-1},
               \overbrace{-1,\ldots,-1}^{2i_1-1}, 
               \overbrace{+1,\ldots,+1}^{2i_2-1},
               \overbrace{-1,\ldots,-1}^{2i_3-1})$}
\multiput(32,-15)(88,0){2}{\vector(0,1){10}}
\put(32,-15){\line(1,0){88}}
\multiput(-40,-15)(20,0){2}{\vector(0,1){10}}
\put(-40,-15){\line(1,0){20}}
\put(-50,-40){$(\overbrace{+1,+1,\ldots,+1}^{2i_0-1},
               \overbrace{+1,\ldots,+1}^{2i_1-1},
               \overbrace{-1,\ldots,-1}^{2i_2-1},
               \overbrace{-1,\ldots,-1}^{2i_3-1})$}
\multiput(-40,-55)(20,0){2}{\vector(0,1){10}}
\put(-40,-55){\line(1,0){20}}
\multiput(32,-55)(88,0){2}{\vector(0,1){10}}
\put(32,-55){\line(1,0){88}}
\end{picture} 
\end{center}
and $(1,2)(2i_0,2n-2i_3+2)$ 
\begin{center}
\begin{picture}(100,80)(0,-55)
\put(-50,0){$(\overbrace{+1,+1,\ldots,+1}^{2i_0-1},
               \overbrace{-1,\ldots,-1}^{2i_1-1}, 
               \overbrace{+1,\ldots,+1}^{2i_2-1},
               \overbrace{-1,\ldots,-1}^{2i_3-1})$}
\multiput(32,-15)(105,0){2}{\vector(0,1){10}}
\put(32,-15){\line(1,0){105}}
\multiput(-40,-15)(20,0){2}{\vector(0,1){10}}
\put(-40,-15){\line(1,0){20}}
\put(-50,-40){$(\overbrace{+1,+1,\ldots,+1}^{2i_0-1},
               \overbrace{+1,\ldots,+1}^{2i_1-1},
               \overbrace{-1,\ldots,-1}^{2i_2-1},
               \overbrace{-1,\ldots,-1}^{2i_3-1})$}
\multiput(-40,-55)(20,0){2}{\vector(0,1){10}}
\put(-40,-55){\line(1,0){20}}
\multiput(32,-55)(105,0){2}{\vector(0,1){10}}
\put(32,-55){\line(1,0){105}}
\end{picture}
\end{center}
taking $v_1$ and
$v_2$ as in (\ref{eq:e1}, \ref{eq:e2}) to
\begin{align*}
  \ \ \ \ v_1&=
   \diag(\overbrace{E,\ldots,E}^{i_0-1},E,\overbrace{-E,\ldots,-E}^{i_1-1},
       \overbrace{E,\ldots,E}^{i_2-1},-E,\overbrace{-E,\ldots,-E}^{i_3-1}),
        \\
  \ \ \ \  
 v_2&=\diag(\overbrace{E,\ldots,E}^{i_0-1},R,\overbrace{E,\ldots,E}^{i_1-1},
       \overbrace{-E,\ldots,-E}^{i_2-1},R,\overbrace{-E,\ldots,-E}^{i_3-1}) 
\intertext{respectively to}
  \ \ \ \ v_1&=
  \diag(\overbrace{E,\ldots,E}^{i_0-1},R,\overbrace{-E,\ldots,-E}^{i_1-1},
       \overbrace{E,\ldots,E}^{i_2-1},R,\overbrace{-E,\ldots,-E}^{i_3-1}),
        \\
  \ \ \ \ v_2&=
      \diag(\overbrace{E,\ldots,E}^{i_0-1},R,\overbrace{E,\ldots,E}^{i_1-1},
       \overbrace{-E,\ldots,-E}^{i_2-1},R,\overbrace{-E,\ldots,-E}^{i_3-1}) 
\end{align*}
In the same way as above, we see that these are really \pl s of
$P[i]$. The three lifts are not conjugate in $N(\PSL(2n,\R))$ because
the intersection with the \mt\ is $v_2$ in case (\ref{eq:e1},
\ref{eq:e2}) but $v_1$ and $v_1+v_2$ in the two other cases.  
Observe that all three \pl s of $P[i]$ have the same image in
$W(\pslr{2n})=\pi_0N(\pslr{2n}) \subset \pi_0\GL(2,\R) \wr
\Sigma_n$. Observe also that the inclusion $P[i] \times P[i]^{\vee}_i
\to C_{\pslr{2n}}(P[i])$ induces an iso\m\ on component groups and
that the centralizer
\begin{multline*}
  C_{\pslr{2n}}(P[i] \times P[i]^{\vee}_i) =
  C_{C_{\pslr{2n}}(P[i])}(P[i]^{\vee}_i) \\
  =
  \begin{cases}
    P[i] \times \SL(2i_0-1,\R) & i=[i_0,i_0,i_0,i_0] \\
    P[i] \times \SL(2i_0-1,\R)\times\SL(2i_2-1,\R) &
    i=[i_0,i_0,i_2,i_2]\\
    C_{\pslr{2n}}(P[i]) & \text{otherwise}
  \end{cases}
\end{multline*}
has nontrivial identity component when $n>2$.

\noindent
\underline{$[P,P] \neq 0$}: $\A(\pslr{4n})$ contains (up to iso\m )
four rank two objects with nontrivial inner product, namely $H_+$,
$H_+^D$, $H_-$, and $H_-^D$ where $H_{\pm}$ is the image of
$2^{1+2}_{\pm} \subset \SL(2n,\R)$ (\ref{sec:non0innerprodrank2}).

The extraspecial $2$-group $2^{1+2}_+ \subset \SL(4n,\R)$ is described
in \ref{exmp:IDQ}.(\ref{exmp:IDQ6}) or, alternatively, in
\ref{sec:glnCtosl2nR} as 
\begin{equation*}
  2^{1+2}_+=\Big\langle\diag\Big(\overbrace{
      \begin{pmatrix}
      E&0\\0&-E
      \end{pmatrix}, \ldots,
      \begin{pmatrix}
        E&0\\0&-E
      \end{pmatrix}}^n\Big),
      \diag\Big(\overbrace{
        \begin{pmatrix}
        0&E\\E&0
        \end{pmatrix}, \ldots,
        \begin{pmatrix}
          0&E\\E&0
        \end{pmatrix}}^n\Big)\Big\rangle = \gen{g_1,g_2}
\end{equation*}
Note that $2^{1+2}_+$ is contained in $N(\SL(4n,\R)) = \SL(4n,\R) \cap
(\GL(2,\R) \wr \Sigma_{2n})$ and that its centralizer there is
\begin{multline*}
  C_{N(\SL(4n,\R))}(2^{1+2}_+) = 
  \SL(4n,\R) \cap C_{\GL(4n,\R)}(v_1) \cap
  C_{\GL(2,\R)\wr\Sigma_{2n}}(v_2) \\
  \stackrel{\text{(\ref{lemma:centrgw})}}{=}
  \SL(4n,\R) \cap  C_{\GL(2,\R)^n \rtimes(C_2\wr\Sigma_{n})}(v_1) 
  = \GL(2,\R) \wr \Sigma_n 
  =N(\GL(2n,\R))
\end{multline*}
%\begin{align*}
%  C_{N(\SL(4n,\R))}(2^{1+2}_+) &=
%  \SL(4n,\R) \cap C_{\GL(2,\R) \wr \Sigma_{2n}}(2^{1+2}_+) \\
%  &= \SL(4n,\R) \cap C_{\GL(2,\R) \wr \Sigma_{2n}}(g_1) \cap
%      C_{\GL(2,\R) \wr \Sigma_{2n}}(g_2) \\
%  &\stackrel{\text{(\ref{lemma:centrgw})}}{=}  
%    (\GL(2,\R) \times \GL(2,\R))\rtimes\Sigma_n \cap 
%      \GL(2,\R)^n\rtimes (C_2 \wr \Sigma_n) \\ 
%  &= \GL(2,\R) \wr \Sigma_n
%\end{align*}
It follows, as in \ref{sec:non0innerprodrank2}, that the centralizer
of $H_+$ in $N(\pslr{4n})$ is
\begin{equation*}
   C_{N(\SL(4n,\R))}(H_+) = H_+ \times  C_{N(\SL(4n,\R))}(2^{1+2}_+)/\gen{-E}
   = N(H_+ \times \PGL(2n,\R))  
\end{equation*}
which means that $H_+ \subset N(\pslr{4n})$ is a \pl\ of $H_+ \subset
\pslr{4n}$. Another \pl\ can be obtained by pre-composing the inclusion
$H_+ \subset N(\pslr{4n})$ with the nontrivial auto\m\ in
$\A(\pslr{4n})(H_+)=O^+(2,\F_2)$. The final \pl\ is
\begin{multline*}
  (2^{1+2}_+)^{\diag(\overbrace{B, \ldots,B}^n)} = 
  \gen{-(g_1g_2)^{\diag(B,\ldots,B)},g_2^{\diag(B,\ldots,B)}}, 
   \qquad B=
  \frac{1}{\sqrt{2}}\begin{pmatrix}
    E&I\\I&E
  \end{pmatrix}, \\
  -(g_1g_2)^{\diag(B,\ldots,B)}=  \diag\left(
    \begin{pmatrix}
      I&0\\0&-I
    \end{pmatrix},\ldots,
    \begin{pmatrix}
     I&0\\0&-I
    \end{pmatrix}\right), \qquad g_2^B=g_2
\end{multline*}
Also this subgroup is actually contained in the \mtn\ with
centralizer 
\begin{multline*}
  C_{N(\SL(4n,\R))}(2^{1+2}_+)^{\diag(B,\ldots,B)}
  =\SL(4n,\R) \cap
  C_{\GL(2,\R)\wr\Sigma_{2n}}\big(-(g_1g_2)^{\diag(B,\ldots,B)}\big) \cap
  C_{\GL(4n,\R)}(g_2) \\
  \stackrel{\text{(\ref{lemma:centrgw})}}{=} 
  (\GL(1,\C)^2\rtimes C_2 ) \wr \Sigma_n \cap C_{\GL(4n,\R)}(g_2) 
  = (\GL(1,\C) \rtimes C_2) \wr \Sigma_n =N(\GL(2n,\R))
\end{multline*}
where
\begin{equation*}
  C_2=\gen{
    \begin{pmatrix}
      0&T\\T&0
    \end{pmatrix}}, \qquad T=
  \begin{pmatrix}
    0&1\\1&0
  \end{pmatrix}
\end{equation*}
%%\marginpar{{$N \to G$ mono $\Rightarrow$ \\ $N$ is mtn?}}  
Observe that, for all three \pl s of $H_+$, the image in the Weyl
group $W(\pslr{4n})=\pi_0N(\pslr{4n}) \subset
\pi_0\GL(2,\R)\wr\Sigma_{2n}$ is the order $2$ subgroup of
$\Sigma_{2n}$ generated by the permutation $(1,2)(3,4) \cdots
(2n-1,2n)$. Observe also the inclusion $H_+ \#L(1,2n-1) \to
C_{\pslr{4n}}(H_+)$ (\ref{prop:H+L}) induces an iso\m\ on component
groups and that the centralizer $C_{\pslr{4n}}(H_+\#L(1,2n-1))$ has
nontrivial identity component (according to the proof of
\ref{prop:ZC}) when $n \geq 2$.

The extraspecial $2$-group $2^{1+2}_- \subset \SL(4n,\R)$ is described
in \ref{exmp:IDQ}.(\ref{exmp:IDQ7}) or, alternatively, in
\ref{sec:glnCtosl2nR} as 
\begin{equation*}
  2^{1+2}_-=\Big\langle\diag\Big(\overbrace{
      \begin{pmatrix}
      I&0\\0&-I
      \end{pmatrix}, \ldots,
      \begin{pmatrix}
        I&0\\0&-I
      \end{pmatrix}}^n\Big),
      \diag\Big(\overbrace{
        \begin{pmatrix}
        0&I\\I&0
        \end{pmatrix}, \ldots,
        \begin{pmatrix}
          0&I\\I&0
        \end{pmatrix}}^n\Big)\Big\rangle = \gen{g_1,g_2}
\end{equation*}
Note that $2^{1+2}_-$ is contained in $N(\SL(4n,\R)) = \SL(4n,\R) \cap
(\GL(2,\R) \wr \Sigma_{2n})$ and that its centralizer there is
\begin{multline*}
   C_{N(\SL(4n,\R))}(2^{1+2}_-)
  =\SL(4n,\R) \cap
  C_{\GL(2,\R)\wr\Sigma_{2n}}(g_1) \cap
  C_{\GL(4n,\R)}(g_2) \\
  \stackrel{\text{(\ref{lemma:centrgw})}}{=} 
  (\GL(1,\C)^2\rtimes C_2 ) \wr \Sigma_n \cap C_{\GL(4n,\R)}(g_2) 
  = \Big\langle\GL(1,\C),
    \begin{pmatrix}
      0&T\\-T&0
    \end{pmatrix}\Big\rangle\wr\Sigma_n 
  \stackrel{\text{(\ref{sec:NpglnH})}}{=} N(\GL(n,\Ha))
\end{multline*}
It follows, as in \ref{sec:non0innerprodrank2}, that the centralizer
of $H_-$ in $N(\pslr{4n})$ is
\begin{equation*}
   C_{N(\SL(4n,\R))}(H_-) = H_- \times C_{N(\SL(4n,\R))}(2^{1+2}_+)/\gen{-E}
   = N(H_- \times \PGL(n,\Ha))  
\end{equation*}
which means that $H_- \subset N(\pslr{4n})$ is a \pl\ of $H_- \subset
\pslr{4n}$. The two other \pl s can be obtained by pre-composing the
inclusion $H_- \subset N(\pslr{4n})$ with the nontrivial auto\m s in
$\A(\pslr{4n})(H_-)=O^-(2,\F_2)=\Symp(2,\F_2)=\GL(2,\F_2)$.  Observe
that, for all three \pl s of $H_-$, the image in the Weyl group
$W(\pslr{4n})=\pi_0N(\pslr{4n}) \subset \pi_0\GL(2,\R)\wr\Sigma_{2n}$
is the order $2$ subgroup of $\Sigma_{2n}$ generated by $(1,2)(3,4)
\cdots (2n-1,2n)$. Observe also that $H_-$ is contained in the rank
three subgroup $H_+\#L(1,n-1)$ (\ref{prop:H-L}) whose centralizer has
a nontrivial identity component when $n \geq 2$ (according to the
proof of \ref{prop:ZC}).

We conclude that for every nontoral rank two object $P$ of
$\A(\pslr{2n})$ the identity component $C_{\pslr{2n}}(P)_0$ of the
centralizer is center-less. By (part of) \cite[5.2]{jmm:ext}, the
homo\m\ 
\begin{equation*}
  \Aut(C_{\pslr{2n}}(P)) \to \Aut(\pi_0C_{\pslr{2n}}(P)) \times 
                             \Aut(C_{\pslr{2n}}(P)_0),
\end{equation*}
obtained by applying the functors $\pi_0$ and $(\ )_0$, is injective.
Under the inductive assumption that $C_{\pslr{2n}}(P)_0$ (see
(\ref{eq:CPi1}) and (\ref{sec:non0innerprodrank2})) has
$\pi_*(N)$-determined 
%\footnote{For $P[i]$ we here use that
%  $\PSL(\text{odd})$ ($B$-family) has $\pi_*(N)$-determined auto\m s.
%  This is true by JMO so there is at present no independent proof
%  here. For $H_{\pm}$ we use that $\PGL(n,\Ha)$ has
%  $\pi_*(N)$-determined auto\m s} 
auto\m s it then follows from
Lemma~\ref{lemma:uniquenu} and diagram
(\ref{dia:Utriangle}) that condition (3) of
Theorem~\ref{indstepalt} is satisfied.

\section{Limits over the Quillen category of $\pslr{2n}$}
\label{sec:quillen}

In this section we show that the problem of computing the higher
limits of the functors $\pi_i(BZC_{\pslr{2n}})$, $i=1,2$,
(\ref{eq:defnpiiBZC}) is concentrated on objects of the Quillen
category with $q \neq 0$.

%The full subcategory $\A(\pslr{2n})^{\leq t}$ of toral objects is
%equivalent to the category $\A(W,t)$ where $W=W(\pslr{2n})$ and
%$t=t(\pslr{2n})$ is the maximal \lmntwo\ of the \mt\ of $\pslr{2n}$
%\cite[2.8]{jmm:ndet}.  The full subcategory $\A(\pslr{2n})^{q=0}$ of
%objects with trivial quadratic function $q$ is equivalent
%(\ref{cor:equivcat}) to the category $\A(\Sigma_{2n}, PS\Delta_{2n})$.
%The full subcategory $\A(\pslr{2n})^{\leq t} \cap \A(\pslr{2n})^{q=0}$
%of toral objects with trivial quadratic function is thus equivalent to
%the category $\A(W,t \cap PS\Delta_{2n})$ for the Weyl group action on
%$t \cap PS\Delta_{2n} = \gen{e_1,\ldots ,e_n}/\gen{e_1\cdots e_n}$.

\begin{lemma}\label{lemma:H0V}
  Let $V \subset PS\Delta_{2n}$ (\ref{eq:dfamt}) be a nontrivial
  subgroup representing an object of the category
  $\A(\pslr{2n})^{q=0}=\A(\Sigma_{2n},PS\Delta_{2n})$
  (\ref{cor:equivcat}). Then
  \begin{equation*}
    ZC_{\pslr{2n}}(V) = PS\Delta_{2n}^{\Sigma_{2n}(V)}
  \end{equation*}
  where $\Sigma_{2n}(V) \subset \Sigma_{2n}$ is the point-wise
  stabilizer subgroup (\ref{defn:AWt}).
%Then
%\begin{equation*}
%  \pi_1(BZC_{\pslr{2n}})(V)=H^0(\Sigma_{2n}(V); PS\Delta_{2n}), \quad
%  \pi_2(BZC_{\pslr{2n}})(V)=0,
%\end{equation*}
%where $\Sigma_{2n}(V) \subset \Sigma_{2n}$ is the point-wise
%stabilizer subgroup (\ref{defn:AWt}). For $n \geq 2$,
%the cohomology group $H^0(\Sigma_{2n}; PS\Delta_{2n})$ is trivial.
%% $PS\Delta_{2n}^{\Sigma_{2n}}=0$ 
\end{lemma}
\begin{proof}
  Let \func{\nu^*}{V}{S\Delta_{2n}} be a lift to $\SL(2n,\R)$ of the
  inclusion homo\m\ of $V$ into $\PSL(2n,\R)$.  Then 
\begin{equation*}
  C_{\SL(2n,\R)}(\nu^*V)= 
  \SL(2n,\R) \cap \prod_{\rho \in V^{\vee}} \GL(i_{\rho},\R),
  \qquad \Sigma_{2n}(\nu^*V) = \prod_{\rho \in
    V^{\vee}}\Sigma_{i_{\rho}} 
\end{equation*}
where \func{i}{V^{\vee}}{\Z} records the multiplicity of each $\rho
\in V^{\vee}$ in the representation $\nu^*$. 
Write
  \begin{equation*}
    \nu^*(v)=\diag(\overbrace{\rho_1(v),\ldots ,\rho_1(v)}^{i_1},
    \ldots ,\overbrace{\rho_m(v),\ldots ,\rho_m(v)}^{i_m})
  \end{equation*}
  where $\rho_1,\ldots ,\rho_m \in V^{\vee}=\Hom(V,C_2)$ are pairwise
  distinct homo\m s $V \to C_2=\gen{\pm 1}$ and $i_1+\cdots+i_m=2n$.
  There is a corresponding decomposition $\{1,\ldots ,2n\}=I_1 \cup
  \cdots \cup I_m$ of the set $\{1,\ldots ,2n\}$ into $k$ disjoint
  subsets $I_j$ containing $i_j$ elements.

Using
\cite[5.11]{jmm:ndet} and \ref{lemma:WVWVast} we get that
\begin{equation*}
  C_{\PSL(2n,\R)}(V) =  \frac{ C_{\SL(2n,\R)}(\nu^*V) }{\gen{-E}}
  \rtimes V^{\vee}_{\nu^*}, \qquad
  \Sigma_{2n}(V) = \Sigma_{2n}(\nu^*V) \rtimes V^{\vee}_{\nu^*}
\end{equation*}
where $V^{\vee}_{\nu^*} = \{ \zeta \in \Hom(V,\GL(1,\R)) \mid
\forall \rho \in V^{\vee} \colon i_{\zeta\rho} = i_{\rho}\}$.
Suppose that $\zeta$ is a
  nontrivial element of $V^{\vee}_{\nu^*}$. Choose a vector $v \in V$
  such that $\zeta(v)=-1$. Then the determinant of $\nu^*(v)$ is
  $(-1)^n$ for $\nu^*(v)$ consists of an equal number of $+1$ and
  $-1$. Thus $n$ is even. Let $\sigma$ be the permutation associated
  to $\zeta$ that moves the subset $I_j$ monotonically to $I_k$ where
  $\zeta\rho_j=\rho_k$. Then $\sigma$ is even for it is a product of
  $n$ transpositions. In this way we imbed $V^{\vee}_{\nu^*}$ as a
  subgroup of the alternating group $A_{2n} \subset \PSL(2n,\R)$ to
  obtain the semi-direct products.

The center of the centralizer is 
\begin{align*}
  ZC_{\PSL(2n,\R)}(V) 
  &= Z \left( \frac{\prod \SL(2n,\R) \cap \GL(i_{\rho},\R)}{\gen{-E}} 
  \rtimes V^{\vee}_{\nu^*} \right) 
  \stackrel{\textmd{(\ref{cor:cent})}}{=} 
  Z \left( \frac{\prod \SL(2n,\R) \cap \GL(i_{\rho},\R)}{\gen{-E}}
  \right)^{V^{\vee}_{\nu^*}}  \\
  &\stackrel{\textmd{(\ref{lemma:ZGLprod})}}{=}
  \left( \frac{\SL(2n,\R) \cap \prod Z\GL(i_{\rho},\R)}{\gen{-E}} 
  \right)^{V^{\vee}_{\nu^*}}
  =\left( \frac{S\Delta_{2n}^{\Sigma_{2n}(\nu^*V)}}{\gen{-E}} 
  \right)^{V^{\vee}_{\nu^*}} \\
  &=\left( PS\Delta_{2n}^{\Sigma_{2n}(\nu^*V)}
  \right)^{V^{\vee}_{\nu^*}} 
  = PS\Delta_{2n}^{\Sigma_{2n}(V)}
\end{align*}
where the penultimate equlity sign is justified by observing that the
coefficient group homo\m\ 
$H^1(\Sigma_{2n}(\nu^*V);\gen{-E}) \to
 H^1(\Sigma_{2n}(\nu^*V);S\Delta_{2n}) \to
 H^1(\Sigma_{2n}(\nu^*V);\Delta_{2n})$ is injective.
\end{proof}

Let $\pi_i(BZC)=\pi_i(BZC_{\PSL(2n,\R)})$ (\ref{eq:defnpiiBZC}).

\begin{cor}\label{cor:dfamlimqeq0}
  $\lim^*(\A(\PSL(2n,\R))^{q=0},\pi_i(BZC))=0$ for $n \geq 2$ and $i=1,2$.
\end{cor}
\begin{proof}
  This is obvious for $i=2$ as $\pi_2(BZC)=0$. For $i=1$, use
  \ref{dw:limits} to compute the limits of the functor
  $\pi_1(BZC)(V)=PS\Delta_{2n}^{\Sigma_{2n}(V)}$ (\ref{lemma:H0V}). The
  fixed-point group $PSD_{2n}^{\Sigma_{2n}}=0$ since $PSD_{2n}$ is an
  irreducible $\F_2\Sigma_{2n}$-module of dimension $2n-2$ for $n \geq
  2$.
%% /home/moller/manus/dfam/magma/lhs/S2n.prg
\end{proof}

\begin{lemma}\label{lemma:dfamZCtoralqeq0}
  Let $V \subset Pt(\SL)=Pt(\SL(2n,\R))$ (\ref{eq:dfamt}) be a nontrivial
  subgroup representing an object of the category
  $\A(A_{2n} \cap (\Sigma_2\wr\Sigma_n),Pt(\SL)) =
  \A(\pslr{2n})^{\leq t, q=0}$
  (\ref{cor:equivcat}). Then
  \begin{equation*}
    ZC_{\pslr{2n}}(V) = Pt(\SL)^{(A_{2n} \cap (\Sigma_2\wr\Sigma_n))(V)}
  \end{equation*}
  where $(A_{2n} \cap (\Sigma_2\wr\Sigma_n))(V)$ is the point-wise
  stabilizer subgroup (\ref{defn:AWt}).
\end{lemma}
\begin{proof}
  The  point-wise
  stabilizer subgroups are\set
  \begin{equation}\label{eq:dfamstabgroups}
    (A_{2n} \cap (\Sigma_2\wr\Sigma_n))(V)=A_{2n}
    \cap \Sigma_{2n}(V),
    \qquad
    \Sigma_{2n}(V) = \Sigma_2^n \rtimes \Sigma_n(V)
  \end{equation}\add
  Because these stabilizer subgroup have these particular forms and
  $PS\Delta_{2n}^{\Sigma_2^n}=Pt(\SL)$, we get
  \begin{equation*}
    ZC_{\PSL(2n,\R)}(V) 
    =  PS\Delta_{2n}^{\Sigma_{2n}(V)} 
    =  PS\Delta_{2n}^{\Sigma_2^n \rtimes \Sigma_{n}(V)} 
    = Pt(\SL)^{\Sigma_{n}(V)}
    =  Pt(\SL)^{A_{2n} \cap (\Sigma_2^n \cap \Sigma_{n}(V))}
     \end{equation*}
  from \ref{lemma:H0V}
\end{proof}

\begin{cor}\label{cor:dfamlimtoralqeq0}
  $\lim^*(\A(\PSL(2n,\R))^{\leq t, q=0},\pi_i(BZC))=0$ for $n \geq 2$
  and $i=1,2$.
\end{cor}
\begin{proof}
  Similar to \ref{cor:dfamlimqeq0} but using $H^0(A_{2n} \cap
  (\Sigma_2 \wr \Sigma_n);Pt(\SL))= H^0(\Sigma_n;Pt(\SL)) = 0$.
\end{proof}

\begin{lemma}\label{lemma:dfamZCtoral}
  Let $V \subset t(\PSL)=t(\PSL(2n,\R))$ (\ref{eq:dfamt}) be a nontrivial
  subgroup representing an object of the category
  $\A(A_{2n} \cap (\Sigma_2\wr\Sigma_n),t(\PSL)) =
  \A(\pslr{2n})^{\leq t}$
  (\ref{cor:equivcat}) where $n \geq 32$. Then
  \begin{equation*}
    \ch{Z}C_{\pslr{2n}}(V) 
         = \ch{T}^{(A_{2n} \cap (\Sigma_2\wr\Sigma_n))(V)}
  \end{equation*}
  where $\ch{T} = \ch{T}(\PSL(2n,\R))$ is the discrete approximation
  \cite[\S3]{dw:center} to the \mt\ of $\PSL(2n,\R)$ and $(A_{2n} \cap
  (\Sigma_2\wr\Sigma_n))(V)$ is the point-wise stabilizer subgroup
  of $V$ (\ref{defn:AWt}).
\end{lemma}
\begin{proof}
  Consider first the case where $V \subset Pt(\SL) \subset t(\PSL)$. 
  One checks that $\ch{T}^{A_{2n} \cap \Sigma_2^n}= Pt(\SL)$ for $n
  > 2$. 
  Since
  $(A_{2n} \cap (\Sigma_2\wr\Sigma_n))(V) \supset A_{2n} \cap
  \Sigma_2^n$ we get 
  \begin{equation*}
    ZC_{\PSL(2n,\R)}(V) 
   \stackrel{\textmd{\ref{lemma:dfamZCtoralqeq0}}}{=} Pt(\SL)^{(A_{2n} \cap
      (\Sigma_2\wr\Sigma_n))(V)} = \ch{T}^{(A_{2n} \cap
      (\Sigma_2\wr\Sigma_n))(V)}
  \end{equation*}
  in this case.
  
  Consider next the case where $V^*$, the preimage of $V$ in
  $\SL(2n,\R)$, contains $I$ (\ref{eq:dfamIc}) so that $V^*=\gen{I,U}$
  (\ref{lemma:Vast}) for some (possibly trivial) \lmntwo\ $U \subset
  t(\SL) \subset C_{\SL(2n,\R)}(C_4) = \GL(n,\C)$. Then
  \begin{equation*}
    C_{\SL(2n,\R)}(V^*) = \prod_{\rho \in U^{\vee}} \GL(i_{\rho},\C),
    \qquad (\Sigma_2 \wr \Sigma_n)(V^*)=\Sigma_n(U) \subset A_{2n}
  \end{equation*}
  where \func{i}{U^{\vee}}{\Z} records the multiplicity of the linear
  character $\rho\in U^{\vee}$ in the representation $\nu^* \colon U
  \to \GL(n,\C)$ and $\Sigma_n(U)$ is point-wise stabilizer subgroup
  for the action of $\Sigma_n=W(\GL(n,\C))$ on
  $t(\SL)=t(\GL(n,\C))=\gen{e_1, \ldots e_n}$. It now follows
  \cite[5.11]{jmm:ndet} and \ref{lemma:WVWVast}, as in (\ref{CLneq0})
  and (\ref{CVqneq0}), that
  \begin{equation*}
   \begin{split}
    C_{\PSL(2n,\R)}(V) &= 
    \begin{cases}
       \frac{C_{\SL(2n,\R)}(V^*)}{\gen{-E}}  & \textmd{ $n$ odd}\\
     \frac{C_{\SL(2n,\R)}(V^*)}{\gen{-E}} \rtimes
     (U^{\vee}_{\nu^*} \times \gen{c_1 \cdots c_n}) & \textmd{ $n$ even}      
    \end{cases}
    \\
     A_{2n} \supset (\Sigma_2\wr\Sigma_n)(V) &= 
     \begin{cases}
        \Sigma_n(U)  & \textmd{ $n$ odd}\\
       \Sigma_n(U) \rtimes (U^{\vee}_{\nu^*} \times  \gen{c_1 \cdots
         c_n})& \textmd{ $n$ even}  
     \end{cases}  
   \end{split}
  \end{equation*}
  where $U^{\vee}_{\nu^*} = \{ \zeta \in U^{\vee} = \Hom(U,\gen{-E}) 
  \mid \forall \rho \in U^{\vee} \colon i_{\zeta\rho} = i_{\rho}\}$
  can be realized as a subgroup of $\Sigma_n$ and
  $\gen{c_1 \cdots c_n}$ is the diagonal order two subgroup of
  $\Sigma_2^n$. Consequently, if $n$ is odd,
  \begin{multline*}
    \ch{Z}C_{\PSL(2n,\R)}(V) =
    \ch{Z} \left( \frac{\prod \GL(i_{\rho},\C)}{\gen{-E}}\right) =
    \frac{\prod \ch{Z}\GL(i_{\rho},\C)}{\gen{-E}} =
    \frac{\ch{T}(\SL(2n,\R))^{\Sigma_n(U)}}{\gen{-E}} \\ =
    \ch{T}^{\Sigma_n(U)} =
    \ch{T}^{(A_{2n} \cap (\Sigma_2\wr\Sigma_n))(V)} 
  \end{multline*}
  and if $n$ is even,
  \begin{align*}
     \ch{Z}C_{\PSL(2n,\R)}(V) &=
     \ch{Z} \left( \frac{\prod \GL(i_{\rho},\C)}{\gen{-E}} 
        \rtimes (U_{\nu^*} \times \gen{c_1 \cdots c_n}) \right) 
      \stackrel{\textmd{\ref{cor:cent}}}{=} 
       \left( \frac{\prod \ch{Z}\GL(i_{\rho},\C)}{\gen{-E}}\right)^{
      U_{\nu^*} \times \gen{c_1 \cdots c_n}} \\
     &= \left(  \frac{\ch{T}(\SL(2n,\R))^{\Sigma_n(U)}}{\gen{-E}}
     \right)^{U_{\nu^*} \times \gen{c_1 \cdots c_n}} 
     = \left( \ch{T}^{\Sigma_n(U)} \right) ^{U_{\nu^*} \times
       \gen{c_1 \cdots c_n}} \\
     &= \ch{T}^{U_{\nu^*} \times \gen{c_1 \cdots c_n}}
     = \ch{T}^{(A_{2n} \cap (\Sigma_2\wr\Sigma_n))(V)}
  \end{align*}
  where we use that $H^1(\Sigma_n(U);\gen{-E}) \to
  H^1(\Sigma_n(U);\ch{T}(\SL(2n,\R)))$ is injective. (In fact, the
  center of the centralizer,
  $\ch{Z}C_{\PSL(2n,\R)}(V)$, is a product, $\ch{T}^{\Sigma_n(U)}$, of
    $2$-compact tori when $n$ is odd, and a finite abelian group,
    $\ch{T}^{\Sigma_n(U) \rtimes (U_{\nu^*}^{\vee} \rtimes \gen{c_1 \cdots
        c_n})} = (\ch{T}^{\gen{c_1 \cdots c_n}})^{\Sigma_n(U) \rtimes
      U_{\nu^*}^{\vee}} = t(\PSL)^{\Sigma_n(U) \rtimes
      U_{\nu^*}^{\vee}}$ when $n$ is even.)
\end{proof}

Lemma \ref{lemma:dfamZCtoralqeq0} can also be proved along the lines
of \cite[2.8]{jmm:ndet} using \ref{prop:centerofN}.

\begin{cor}\label{cor:dfamtorallim}
  $\lim^*(\A(\PSL(2n,\R))^{\leq t},\pi_i(BZC))=0$ for $n \geq 3$
  and $i=1,2$.
\end{cor}
\begin{proof}
  Similar to \ref{cor:dfamlimqeq0} but using that
  that $H^0(W;\ch{T})(\PSL(2n,\R))=0$ for $n \geq 3$
  (\ref{eq:dfamH01WT}).  
\end{proof}

As we shall next see, Corollaries \ref{cor:dfamlimqeq0},
\ref{cor:dfamlimtoralqeq0} and \ref{cor:dfamtorallim} reduce the
problem of computing the graded abelian group
$\lim^*(\A(\pslr{2n}),\pi_i(BZC))$ considerably.

Let $\A$ be a category containing two full subcategories, $\A_j$,
$j=1,2$, such that any object of $\A$ with a \m\ to an object of
$\A_j$ is an object of $\A_j$. Write $\A_1 \cap \A_2$ for the full
subcategory with objects $\Ob(\A_1 \cap \A_2)=\Ob(\A_1) \cap \Ob(\A_2)$
and $\A_1 \cup \A_2$ for the full subcategory with objects $\Ob(\A_1
\cup \A_2)= \Ob(\A_1) \cup \Ob(\A_2)$.  Let \func{M}{\A}{\Ab} be a
functor taking values in abelian groups.  Consider the subfunctor
$M_{12}$ of $M$ given by
\begin{equation*}
  M_{12}(a)=
  \begin{cases}
    0 & a \in \Ob(\A_1 \cup \A_2) \\
    M(a) & a \not\in  \Ob(\A_1 \cup \A_2)
  \end{cases}
\end{equation*}
We now state a kind of Mayer--Vietoris sequence argument for
cohomology of categories.
\begin{lemma}
  If the graded abelian groups $\lim^*(\A_1,M)$, $\lim^*(\A_2,M)$, and
  $\lim^*(\A_1 \cap \A_2,M)$ are trivial, then $\lim^*(\A,M_{12})
  \cong \lim^*(\A;M)$.
\end{lemma}
\begin{proof}
  Consider also the subfunctor $M_{1}$ of $M$ given by
\begin{equation*}
  M_{1}(a)=
  \begin{cases}
    0 & a \in \Ob(\A_1) \\
    M(a) & a \not\in  \Ob(\A_1)
  \end{cases}
\end{equation*}
Then there are natural transformations $M_{12} \to M_1 \to M$ of
functors. The induced long exact sequences imply that it suffices to
show $\lim^*(\A;M/M_1)=0=\lim^*(\A;M_1/M_{12})$.

The quotient functor $M/M_1$ vanishes outside $\A_1$ where it agrees
with $M$ and therefore \cite[13.12]{jmm:ndet} $\lim^*(\A;M/M_1) \cong
\lim(\A_1; M)$ which is trivial by assumption. 

The same argument applied to $\A_2$ instead of $\A$ gives that
$\lim^*(\A_2;M/M_1) \cong \lim(\A_1 \cap \A_2; M)$.  Since this
abelian group is trivial by assumption, we have that $\lim^*(\A_2;M_1)
\cong \lim^*(\A_2;M)$. Also this abelian group is trivial by
assumption.

The quotient functor $M_1/M_{12}$ vanishes outside $\A_1 \cup \A_2$
where it agrees with $M_1$ and therefore $\lim^*(\A;M_1/M_{12}) \cong
\lim(\A_1 \cup \A_2; M_1)$. Here, the functor $M_1$ vanishes outside
$\A_2$ and hence $\lim(\A_1 \cup \A_2; M_1) \cong \lim^*(\A_2;M_1)$.
Since we just showed that this abelian group is trivial, we have that
$\lim^*(\A;M_1/M_{12})=0$.
\end{proof}

We conclude that 
\begin{equation*}
  \lim^*\left(\A(\pslr{2n}),\pi_j(BZC_{\pslr{2n}})_{12}\right) =
  \lim^*\left(\A(\pslr{2n}),\pi_j(BZC_{\pslr{2n}})\right)  
\end{equation*}
where $\pi_j(BZC_{\pslr{2n}})_{12}$ is the subfunctor of
$\pi_j(BZC_{\pslr{2n}})$ given by
\begin{equation*}
 \pi_j(BZC_{\pslr{2n}})_{12}(V)=
 \begin{cases}
   0 & \text{$V$ is toral or $q(V)=0$} \\
   \pi_j(BZC_{\pslr{2n}}(V)) & \text{$V$ is nontoral and $q(V) \neq 0$}
 \end{cases}
\end{equation*}

According to \ref{cor:toralnontoral} we have
\begin{equation*}
  \text{$V$ is nontoral and $q(V) \neq 0$} \iff [V,V] \neq 0
\end{equation*}
for all \lmntwo s $V$ in $\pslr{2n}$. Thus the problem of computing
the higher limits of the functors $\pi_i(BZC_{\PSL(2n,\R)})$ is
concentrated on the full subcategory $\A(\pslr{2n})^{[\; ,\; ] \neq
  0}$ of $\A(\pslr{2n})$ generated by all \lmntwo s $V \subset
\pslr{2n}$ with nontrivial inner product. Note that if $\pslr{2n}$
contains an \lmntwo\ $V$ with $[V,V] \neq 0$ then $\pslr{2n}$ in
particular contains such a subgroup of rank two. The preimage in
$\SL(2n,\R)$ of rank two $V \subset \pslr{2n}$ with nontrivial inner
product is an extraspecial $2$-group $2^{1+2}_{\pm}$ with central
$\mho_1$ (\ref{lemma:Vast}) so that, by real representation theory
(\ref{sec:xtraspecreps}), $n$ must be even.

\section{Higher limits of the functors $\pi_iBZC_{\pslr{4n}}$ on
  $\A(\pslr{4n})^{[\; ,\; ] \neq 0}$}
\label{sec:lim}

%\begin{tabular}[l]{|l||c|c|c|c|c|}\hline
%rank & \multicolumn{2}{|c|}{2} &
%\multicolumn{3}{c|}{3} \\ \hline
%$V^*$ & $2^{1+2}_+$ & $2^{1+2}_-$ &
%$2^{1+2}_+ \times 2$ & $2^{1+2}_- \times 2$ & $2^{1+2}_{\pm} \circ 4$
% \\ \hline
%$C(V)$ & $V \times \PGL(2m,\R)$ & $V \times \PGL(m,\Ha)$ &
%{} & {} & $V \times \PGL(n,\C)/\{\pm E \}$ \\ \hline
%$q$ & $v_1v_2$ & $v_1^2+v_1v_1+v_2^2$ & 
%$v_1v_2$ & $v_1^2+v_1v_1+v_2^2$ & $v_1v_2+v_3^2$ \\ \hline
%$\#q^{-1}(1)$ & $1$ & $3$ & $2$ & $6$ & $4$ \\ \hline
%$\A(V)$ & $O(q)=O^+(2,\F_2)$ &  $O(q)=O^-(2,\F_2)$ & 
%$
%\begin{pmatrix}
%  O^+(2,\F_2) & 0 \\
%  * & 1
%\end{pmatrix}
%$ & 
%$\begin{pmatrix}
%  O^-(2,\F_2) & 0 \\
%  * & 1
%\end{pmatrix}$ & $O(q)
%\cong \GL(2,\F_2)$
%\\ \hline 
%\end{tabular}

In this section we compute the first higher limits of the functors
$\pi_iBZC_{\pslr{4n}}$, $i=1,2$, by means of Oliver's cochain complex
\cite{bob:steinberg}.  

\begin{lemma}\label{lemma:lim=0}
  $\lim^1\pi_1BZC_{\pslr{4n}} = 0 = \lim^2\pi_1BZC_{\pslr{4n}}$ and 
   $\lim^2\pi_2BZC_{\pslr{4n}} = 0 = \lim^3\pi_2BZC_{\pslr{4n}}$.
\end{lemma}

The case $i=2$ is easy. Since $\pi_2BZC_{\pslr{4n}}$ has value $0$ on
all objects of $\A(\pslr{4n})^{[\; ,\; ] \neq 0}$ of rank $\leq 4$
(\ref{prop:ZC}) it is immediate from Oliver's cochain complex that
$\lim^2$ and $\lim^3$ of this functor are trivial.

We shall therefore now concentrate on the case $i=1$. 
The claim of the above lemma is that
Oliver's cochain complex \cite{bob:steinberg}
\set\begin{equation}\label{eq:Olccc} 0 \to \prod_{|P|=2^2}[P]
  \xrightarrow{d^1} \prod_{|V|=2^3}[V] \xrightarrow{d^2}
  \prod_{|E|=2^4}[E] \xrightarrow{d^3} \cdots
\end{equation}\add
is exact at objects of rank $\leq 3$. Here, as a matter of notational
convention, 
\begin{equation*}
[E]=\Hom_{\A(\pslr{4n})(E)}(\St(E),E)  
\end{equation*}
stands for the $\F_2$-vector space of $\F_2\A(\pslr{4n})(E)$-module
homo\m s from the Steinberg module $\St(E)$ to $E$. The Steinberg
module is the $\F_2\GL(E)$-module obtained in the following way.

Let $P=\F_2e_1+\F_2e_2$ be a $2$-dimensional vector space over $\F_2$
with basis vectors $e_1$, $e_2$. Let $\F_2[0]$ be the $3$-dimensional
$\F_2$-vector space on length zero flags, $[L]$, of nontrivial and
proper subspaces $L$ of $P$. The Steinberg module $\St(P)$ is the
$2$-dimensional kernel of the augmentation map \func{d}{\F_2[0]}{\F_2}
given by $d[L]=1$.

Let $V=\F_2e_1+\F_2e_2+\F_2e_3$ be a $3$-dimensional vector space over
$\F_2$ with basis vectors $e_1$, $e_2$, $e_3$. Let $\F_2[1]$ be the
$21$-dimensional $\F_2$-vector space on length one flags $[P>L]$ of
nontrivial and proper subspaces of $V$ and $\F_2[0]$ the
$14$-dimensional $\F_2$-vector space on all length $0$ flags, $[P]$ or
$[L]$, of nontrivial and proper subspaces of $V$.  The Steinberg
module $\St(V)$ over $\F_2$ for $V$ is the $2^3$-dimensional kernel of
the linear map \func{d}{\F_2[1]}{\F_2[0]} given by $d[P>L]=[P]+[L]$.

\begin{prop}\label{prop:H+H-}
  $H_+ \neq H_+^D$ and $H_- \neq H_-^D$ in $\A(\pslr{4n}$.
  The auto\m\ groups 
  of the objects $H_+$ and $H_-$ (\ref{sec:non0innerprodrank2}) are
  \begin{equation*}
    \A(\pslr{4n})(H_{+})=O^{+}(2,\F_2) \cong C_2, \qquad 
     \A(\pslr{4n})(H_{-})=O^{-}(2,\F_2) = \GL(2,\F_2),
  \end{equation*}
  and the dimensions of the spaces of equivariant maps are
  \begin{equation*}
    \dim [H_+] = 2, \qquad \dim [H_-]=1
  \end{equation*}
\end{prop}
\begin{proof}
  The first part was proved in \ref{sec:non0innerprodrank2}. 
  The Quillen auto\m group $\A(\SL(4n,\R))(2^{1+2}_{\pm}) =
  \A(\GL(4n,\R))(2^{1+2}_{\pm}) = \Out(2^{1+2}_{\pm}) \cong
  O^{\pm}(2,\F_2)$ where the iso\m\ is induced by abelianization
  $2^{1+2}_{\pm} \to H_{\pm}$ (\ref{exmp:IDQ}.(\ref{exmp:IDQ2}),
  \ref{exmp:IDQ}.(\ref{exmp:IDQ3}), \ref{sec:xtraspecreps}).
  According to {\em magma}, $\dim[H_+]=2$ and $\dim[H_-]=1$.
%  Abelianization $2^{1+2}_{\pm} \to H_{\pm}$ induces an iso\m\
%  $\Out(2^{1+2}_{\pm}) \to O^{\pm}(2,\F_2)$ between the outer auto\m\
%  group of the extra special group  $2^{1+2}_{\pm}$ and the auto\m\
%  group of the quadratic function on the $\F_2$-vector space $H_{\pm}$.
%  Since all auto\m s of  $2^{1+2}_{\pm}$ are trace preserving,
%  $\A(\GL(2,\R))(2^{1+2}_{\pm})=\Out(2^{1+2}_{\pm})$ for the faithful
%  representation \func{\phi}{2^{1+2}_{\pm}}{\GL(2,\R)}. But then also 
%   $\A(\SL(4n,\R))(2^{1+2}_{\pm})=\Out(2^{1+2}_{\pm})$ for the
%   oriented representation $2n\phi$. Since
%   $\A(\SL(4n,\R))(2^{1+2}_{\pm}) \to \A(\pslr{4n})(H_{\pm})$ is
%   surjective, we conclude that
%   $\A(\pslr{4n})(H_{\pm})=O^{\pm}(2,\F_2)$. It follows
%   (\ref{lemma:psialpha}) that $H_{\pm}\alpha=H_{\pm}$ for all
%   $\alpha\in\Aut(H_{\pm})$. 
\end{proof}

The $\F_2\A(\pslr{4n})(H_+)$-equivariant maps given by 
\set\begin{equation}\label{eq:deff+}
  f_+[L]=L, \qquad 
   f_0[L]=
    \begin{cases}
      H_+^{\A(H_+)} & q(L)=0 \\
      0 & \text{otherwise}
    \end{cases}
\end{equation}\add
form a basis for the $2$-dimensional space $[H_+]$. The
$\F_2\A(\pslr{4n})(H_-)$-equivariant map given by
\set\begin{equation}\label{eq:deff-}
  f_-[L] = L
\end{equation}\add
is a basis for the  $1$-dimensional space $[H_-]$.

The quadratic function (\ref{sec:xtraspecreps})
$q(v_1,v_2,v_3)=v_1^2+v_2v_3$ on $V_0$ (\ref{sec:non0innerprodrank3})
has auto\m\ group
\begin{equation*}
  O(q) \cong \Symp(2,\F_2)= \gen{
    \begin{pmatrix}
      1&0&0 \\ 1&1&1 \\ 0&1&0
    \end{pmatrix},
    \begin{pmatrix}
      1&0&0 \\ 1&1&1 \\ 0&0&1
    \end{pmatrix}} \subset \GL(3,\F_2)
\end{equation*}
of order $6$.

\begin{prop}
  \label{prop:V0} $V_0 \neq V_0^D$ in $\A(\pslr{4n})$. The auto\m\
  group $\A(\pslr{4n})(V_0)=O(q)$ and $\dim [V_0]=4$. 
\end{prop}
\begin{proof}
  See \ref{exmp:IDQ}.(\ref{exmp:IDQ5}) for the first part. According
  to {\em magma}, $\dim [V_0]=4$.
\end{proof}

The four $\F_2\A(\pslr{4n})(V_0)$-module homo\m s
\set\begin{equation}
  \label{eq:basisV0}
  \{df_+, df_0, df_-,f_0\}
\end{equation}\add
given by
\begin{alignat*}{2}
  df_+[P>L] &=
  \begin{cases}
    L & P=H_+ \\
    0 & \text{otherwise}
  \end{cases} \qquad &
  df_0[P>L] &= 
  \begin{cases}
    P^{\A(P)} & P=H_+,q(L)=0\\
    0 & \text{otherwise}
  \end{cases} \\
  df_-[P>L] &=
  \begin{cases}
    L & P=H_- \\
    0 & \text{otherwise}
  \end{cases} \qquad &
  f_0[P>L] &= 
  \begin{cases}
    V_0^{\A(V_0)} & [P,P]=0,q(L)=0\\
    0 & \text{otherwise}
  \end{cases}
\end{alignat*}
is a basis for $[V_0]$. 

The quadratic function on $H_+ \# L(i,2n-i) \in \Ob(\A(\PSL(4n,\R)))$,
$0 \leq i \leq n$, $q(v_1,v_2,v_3)=v_1v_2$, has auto\m\ group
\begin{equation*}
  O(q) = 
  \begin{pmatrix}
    O^+(2,\F_2) & 0 \\ \ast & 1
  \end{pmatrix} = \gen{
    \begin{pmatrix}
      0&1&0 \\ 1&0&0 \\ 0&0&1 
    \end{pmatrix},
    \begin{pmatrix}
      1&0&0 \\ 0&1&0 \\ 1&0&1
    \end{pmatrix},
    \begin{pmatrix}
      1&0&0 \\ 0&1&0 \\ 0&1&1
    \end{pmatrix} }
\end{equation*}
of order $\vert O^+(2,\F_2)\vert \cdot 2^2 =8$. 

\begin{prop}\label{prop:H+L}
   $H_+\#L(i,2n-i) \neq (H_+\#L(i,2n-i))^D \iff \text{$i$ is
     even}$. The Quillen auto\m\ group is
   \begin{equation*}
     \A(\pslr{4n})(H_+ \# L(i,2n-i))=
     \begin{cases}
         \gen{
    \begin{pmatrix}
      0&1&0 \\ 1&0&0 \\ 0&0&1 
    \end{pmatrix},
    \begin{pmatrix}
      1&0&0 \\ 0&1&0 \\ 1&1&1
    \end{pmatrix} } & \text{$i$ odd} \\
    O(q) & \text{$i$ even} 
     \end{cases}
   \end{equation*}
   and the dimension of the space of equivariant maps is 
   \begin{equation*}
     \dim [H_+\#L(i,2n-i)] = 
     \begin{cases}
       6 & \text{$i$ odd} \\
       3 & \text{$i$ even} \\
     \end{cases}
   \end{equation*}
\end{prop}
\begin{proof}
  $H_+\#L(i,2n-i) \subset \pslr{4n}$ is (\ref{sec:non0innerprodrank3})
  the quotient of
  \begin{equation*}
    G=\langle\diag(R,\ldots ,R), \diag(T,\ldots ,T), 
   \diag(\overbrace{-E,\ldots,-E}^i,\overbrace{E,\ldots
     ,E}^{2n-i})\rangle = \gen{g_1,g_2,g_3} \subset \SL(4n,\R)
  \end{equation*}
  The centralizer of $G$ is $\GL(4n,\R)$ is contained in the
  centralizer of its subgroup $2^{1+2}_+$ which is contained in
  $\SL(4n,\R)$ (\ref{exmp:IDQ}). Observe that
  \begin{itemize}
  \item $R$ and $T$ are conjugate in $\GL(2,\R)$
  \item Conjugation with
    $\diag(\overbrace{T,\ldots,T}^i,\overbrace{E,\ldots,E}^{2n-i})$
    induces $(g_1,g_2,g_3) \xrightarrow{\phi_1} (g_1g_3,g_2,g_3)$
  \item  Conjugation with
    $\diag(\overbrace{R,\ldots,R}^i,\overbrace{E,\ldots,E}^{2n-i})$
    induces $(g_1,g_2,g_3) \xrightarrow{\phi_2} (g_1,g_2g_3,g_3)$
  \item When $i=n$, conjugation with  $
  \begin{pmatrix}
    0 & E \\ E & 0
  \end{pmatrix}$ induces $(g_1,g_2,g_3) \xrightarrow{\phi}
  (g_1,g_2,-g_3)$
  \end{itemize}  
  Consider the auto\m\ groups
  \begin{equation*}
    \A(\SL(4n,\R))(G) \subset \A(\GL(4n,\R))(G) \subset \Out(G) \to
    O(q) \subset \Aut(H_+\# L(i,2n-i))
  \end{equation*}
  where the outer auto\m\ group has order $16$.  Note that the auto\m\ 
  $\phi$ is in the kernel of the homo\m\ $\Out(G) \to O(q)$ induced by
  abelianization $G \to H_+\#L(i,2n-i)$. Using the above observations
  we see that $\A(\GL(4n,\R))(G)$, even $\A(\SL(4n,\R))(G)$ for even
  $i$, maps onto $O(q)$. Thus the Quillen auto\m\ group
  $\A(\GL(4n,\R))(G)$ has order $8$ or $16$. When $i=n$ the auto\m\ 
  $\phi$ is in $\A(\GL(4n,\R))(G)$, even in $\A(\SL(4n,\R))(G)$, and
  when $i \neq n$, $\phi \not\in \A(\GL(4n,\R))(G)$ as it does not
  preserve trace. Thus
  \begin{equation*}
    \vert \A(\GL(4n,\R))(G) \vert = 
    \begin{cases}
      16 & i=n \\
      8  & i \neq n
    \end{cases}
  \end{equation*}
  In any case the group $\A(\SL(4n,\R))(G)$ equals the group
  $\A(\GL(4n,\R))(G)$ if and only if $i$ is even. When $i$ is odd, the
  auto\m\ $\phi_1$ is induced from a matrix of negative determinant so
  that $N_{\GL(4n,\R)}(G) \not\subset \SL(4n,\R)$.
 % and from the iso\m\ 
%  (\ref{eq:NOut}) we see that all matrices inducing $\phi_1$ (they
%  form a $GC_{\GL(4n,\R)}(G)$-coset where $GC_{\GL(4n,\R)}(G) \subset
%  \SL(4n,\R)$) have negative determinant.
  According to {\em magma},
  $\dim[H_+ \# L(i,2n-i)]$ is $3$ when $i$ is even and $6$ when $i$ is
  odd. 
\end{proof}

The six $\F_2\A(\pslr{4n})(H_+\#L(i,2n-i))$-linear maps
\set\begin{equation}
  \label{eq:basisH+L}
  \{df_+,df_0,f_0,df_+^D,df_0^D,f_0^D\}
\end{equation}\add
given by
  \begin{alignat*}{2}
      {df}_+[P>L]&=
    \begin{cases}
      L & P=H_+ \\
      0 & \text{otherwise}
    \end{cases} \qquad &
    {df}_0[P>L]&=
    \begin{cases}
      P^{\A(P)} & P=H_+, q(L)=0 \\
      0 & \text{otherwise}
    \end{cases} \\
    {f}_0[P>L] &=
    \begin{cases}
      v_1 & [P,P]=0,q(L)=0 \\
      0 & \text{otherwise}
    \end{cases}
  \end{alignat*}
  is a basis for the $6$-dimensional $\F_2$-vector space
  $[H_+\#L(i,2n-i)]$ for $i$ odd and $[H_+\#L(i,2n-i)] \times
  [(H_+\#L(i,2n-i))^D]$ for $i$ even. Here, $v_1$ is one of the two
  non-zero vectors of $V^{\A(V)}$ that are not $D$-invariant when $i$
  is odd and the nonzero vector of $V^{\A(V)}$ when $i$ is even where
  $V=H_+\#L(i,2n-i)$.

The quadratic function on $H_- \# L(i,n-i) \in \Ob(\A(\PSL(4n,\R)))$,
$1 \leq i \leq [n/2]$, $q(v_1,v_2,v_3)=v_1^2+v_1v_2+v_2^2$, has 
auto\m\ group
\begin{equation*}
  O(q) = 
  \begin{pmatrix}
     O^-(2,\F_2) & \ast \\ 0 & 1
  \end{pmatrix}
\end{equation*}
of order $\vert O^-(2,\F_2) \vert \cdot 2^2=24$.

\begin{prop}\label{prop:H-L}
   $H_- \# L(i,n-i) \neq  (H_- \# L(i,n-i))^D$ for all $n \geq 2$. The
   Quillen auto\m\ group $\A(\pslr{4n})(H_- \# L(i,n-i))=O(q)$ has
   order $24$ and the dimension of the space of equivariant maps is
   $\dim [H_- \# L(i,n-i)]=1$. 
\end{prop}
\begin{proof}
  $H_- \# L(i,n-i) \subset \pslr{4n}$ is the quotient of
  \begin{multline*}
    G=2^{1+2}_- \times 2 = \langle 
    \diag\left(
      \begin{pmatrix}
        0 & -R \\ R & 0
      \end{pmatrix}, \ldots, 
      \begin{pmatrix}
        0 & -R \\ R & 0
      \end{pmatrix} \right), 
     \diag\left(\begin{pmatrix}
        0 & -T \\ T & 0
      \end{pmatrix}, \ldots, 
      \begin{pmatrix}
        0 & -T \\ T & 0
      \end{pmatrix}\right), \\
    \diag\Big(
      \overbrace{
        \begin{pmatrix}
        -E & 0 \\ 0 & -E
        \end{pmatrix}, \ldots,
        \begin{pmatrix}
          -E & 0 \\ 0 &-E
        \end{pmatrix}}^i,\overbrace{
        \begin{pmatrix}
        E & 0 \\ 0 & E
        \end{pmatrix},\ldots
        \begin{pmatrix}
          E & 0 \\ 0 &E
        \end{pmatrix}}^{n-i}\Big)\rangle = \langle g_1,g_2,g_3 \rangle
      \subset \SL(4n,\R)
  \end{multline*}
The centralizer of $G$ in $\GL(4n,\R)$ is contained in the centralizer
of its subgroup $2^{1+2}_-$ which is contained in
  $\SL(4n,\R)$ (\ref{exmp:IDQ}). Observe that
\begin{itemize}
  \item $\A(\SL(4,\R))(2^{1+2}_-) \cong O(q)$
  \item Conjugation with
    $\diag\Big(\overbrace{ 
      \begin{pmatrix}
        T & 0 \\ 0 &T
      \end{pmatrix},\ldots,
      \begin{pmatrix}
        T & 0 \\ 0 &T
      \end{pmatrix}}^i,\overbrace{
      \begin{pmatrix}
        E & 0 \\ 0 & E
      \end{pmatrix},\ldots,
      \begin{pmatrix}
        E & 0 \\ 0 & E
      \end{pmatrix}}^{n-i}\Big)$
    induces the auto\m\ $(g_1,g_2,g_3) \xrightarrow{\phi_1}
    (g_1g_3,g_2,g_3)$
  \item Conjugation with
    $\diag\Big(\overbrace{ 
      \begin{pmatrix}
        R & 0 \\ 0 & R
      \end{pmatrix},\ldots,
      \begin{pmatrix}
        R & 0 \\ 0 & R
      \end{pmatrix}}^i,\overbrace{
      \begin{pmatrix}
        E & 0 \\ 0 & E
      \end{pmatrix},\ldots,
      \begin{pmatrix}
        E & 0 \\ 0 & E
      \end{pmatrix}}^{n-i}\Big)$ 
    induces the auto\m\ $(g_1,g_2,g_3) \xrightarrow{\phi_2}
    (g_1,g_2g_3,g_3)$
  \item When $i=n/2$, conjugation with  $
  \begin{pmatrix}
    0 & E \\ E & 0
  \end{pmatrix}$ induces $(g_1,g_2,g_3) \xrightarrow{\phi}
  (g_1,g_2,-g_3)$
  \end{itemize}
  Consider the auto\m\ groups
  \begin{equation*}
    \A(\SL(4n,\R))(G) \subset \A(\GL(4n,\R))(G) \subset \Out(G) \to
    O(q) \subset \Aut(H_-\# L(i,n-i))
  \end{equation*}
  where the outer auto\m\ group has order
  $48$. Note that the auto\m\ $\phi$ is in the kernel of the homo\m\
  $\Out(G) \to O(q)$ induced by abelianization $G \to H_-\#L(i,n-i)$.
  Using the above observations
  we see that $\A(\SL(4n,\R))(G)$ maps onto $O(q)$. Thus the Quillen
  auto\m\ group $\A(\GL(4n,\R))(G)$ has order $48$ or $24$. When
  $n$ is even and $i=n/2$ the auto\m\ $\phi$ is in
  $\A(\SL(4n,\R))(G)$ and when $i < n/2$, $\phi$ is not in
  $\A(\GL(4n,\R))(G)$ as it does not preserve trace. Thus  
  \begin{equation*}
    \vert \A(\GL(4n,\R))(G) \vert = 
    \begin{cases}
      48 & i=n/2 \\
      24  & i < n/2
    \end{cases}
  \end{equation*}
  In any case, the group $\A(\SL(4n,\R))(G)$ equals
  $\A(\GL(4n,\R))(G)$ so that $H_- \# L(i,n-i) \neq (H_- \#
  L(i,n-i))^D$ (\ref{lemma:intoSL}). According to {\em magma},
  $\dim[H_- \# L(i,n-i)]=1$.
\end{proof}

The $\F_2\A(\pslr{4n})(H_-\#L(i,n-i))$-linear map
\set\begin{equation}
  \label{eq:basisH-L}
  \{df_-\}
\end{equation}\add
given by
\begin{equation*}
  df_-[P>L]=
  \begin{cases}
    L & P=H_- \\
    0 & \text{otherwise}
  \end{cases}
\end{equation*}
is a basis for the $1$-dimensional $\F_2$-vector space
$[H_-\#L(i,n-i)]$.

The quadratic function $q(v_1,v_2,v_3,v_4)=v_1^2 + v_2v_3$ has auto\m\
group
\begin{equation*}
  O(q) = 
  \begin{pmatrix}
    \Symp(2,\F_2) & \ast \\
    0 & 1
  \end{pmatrix}, \qquad \Symp(2,\F_2) \cong \gen{
    \begin{pmatrix}
      1&0&0 \\ 1&1&1 \\ 0&1&0
    \end{pmatrix},
    \begin{pmatrix}
      1&0&0 \\ 1&1&1 \\ 0&0&1
    \end{pmatrix}} \subset \GL(3,\F_2)
\end{equation*}
of order $48$.

\begin{prop}\label{lemma:intoV0L} %\label{prop:V0L}
  The $4$-dimensional object $V_0\#L(i,n-i)$, $1 \leq i \leq [n/2]$,
  of the category $\A(\pslr{4n})$ satisfies $V_0\#L(i,n-i) \neq
  (V_0\#L(i,n-i))^D$. It contains the objects objects $V_0$,
  $H_+\#L(2i,2n-2i)$, and $H_-\#L(i,n-i)$.  The auto\m\ group
  $ \A(\pslr{4n})(V_0\#L(i,n-i)) = O(q)$ and the dimension of the
  space of equivariant maps is $\dim [V_0\#L(i,n-i)] = 5$.
\end{prop}
\begin{proof}
  $V_0 \# L(i,n-i) \subset \pslr{4n}$ is (\ref{sec:glnCtosl2nR}) the
  quotient of
  \begin{multline*}
    G=2^{1+2}_{\pm} \circ 4 \times 2=  \\ 
    \langle 
    \diag\left(
      \begin{pmatrix}
        0 & -E \\ E & 0
      \end{pmatrix}, \ldots,
      \begin{pmatrix}
        0 & -E \\ E & 0
      \end{pmatrix}\right),
      \diag\left(\begin{pmatrix}
        R & 0 \\ 0 & R
      \end{pmatrix}, \ldots, 
      \begin{pmatrix}
        R & 0 \\ 0 & R
      \end{pmatrix} \right), 
     \diag\left(\begin{pmatrix}
        T & 0 \\ 0 & T
      \end{pmatrix}, \ldots, 
      \begin{pmatrix}
        T & 0 \\ 0 & T
      \end{pmatrix}\right), \\
    \diag\Big(
      \overbrace{
        \begin{pmatrix}
        -E & 0 \\ 0 & -E
        \end{pmatrix}, \ldots,
        \begin{pmatrix}
          -E & 0 \\ 0 &-E
        \end{pmatrix}}^i,\overbrace{
        \begin{pmatrix}
        E & 0 \\ 0 & E
        \end{pmatrix},\ldots
        \begin{pmatrix}
          E & 0 \\ 0 &E
        \end{pmatrix}}^{n-i}\Big)\rangle = \langle g_1,g_2,g_3,g_4 \rangle
      \subset \SL(4n,\R)
  \end{multline*}
  The centralizer of $G$ in $\GL(4n,\R)$ is contained in the
  centralizer of its subgroup $2^{1+2}_-$ which is contained in
  $\SL(4n,\R)$ (\ref{exmp:IDQ}). Observe that
\begin{itemize}
  \item $\A(\SL(4,\R))(2^{1+2}_{\pm}\circ 4) = \Out(G) \cong \Out(C_4)
    \times \Symp(2,\F_2)$ (\ref{exmp:IDQ})
  \item Conjugation with
    $\diag\Big(\overbrace{ 
      \begin{pmatrix}
        0 & E \\ E & 0
      \end{pmatrix},\ldots,
      \begin{pmatrix}
        0 & E \\ E & 0
      \end{pmatrix}}^i,\overbrace{
      \begin{pmatrix}
        E & 0 \\ 0 & E
      \end{pmatrix},\ldots,
      \begin{pmatrix}
        E & 0 \\ 0 & E
      \end{pmatrix}}^{n-i}\Big)$
    induces the auto\m\ $(g_1,g_2,g_3,g_4) \xrightarrow{\phi_1}
    (g_1g_4,g_2,g_3,g_4)$
  \item Conjugation with
    $\diag\Big(\overbrace{ 
      \begin{pmatrix}
        T & 0 \\ 0 &T
      \end{pmatrix},\ldots,
      \begin{pmatrix}
        T & 0 \\ 0 &T
      \end{pmatrix}}^i,\overbrace{
      \begin{pmatrix}
        E & 0 \\ 0 & E
      \end{pmatrix},\ldots,
      \begin{pmatrix}
        E & 0 \\ 0 & E
      \end{pmatrix}}^{n-i}\Big)$
    induces the auto\m\ $(g_1,g_2,g_3,g_4) \xrightarrow{\phi_2}
    (g_1,g_2g_4,g_3,g_4)$
  \item Conjugation with
    $\diag\Big(\overbrace{ 
      \begin{pmatrix}
        R & 0 \\ 0 & R
      \end{pmatrix},\ldots,
      \begin{pmatrix}
        R & 0 \\ 0 & R
      \end{pmatrix}}^i,\overbrace{
      \begin{pmatrix}
        E & 0 \\ 0 & E
      \end{pmatrix},\ldots,
      \begin{pmatrix}
        E & 0 \\ 0 & E
      \end{pmatrix}}^{n-i}\Big)$ 
    induces the auto\m\ $(g_1,g_2,g_3,g_4) \xrightarrow{\phi_3}
    (g_1,g_2,g_3g_4,g_4)$
  \item Conjugation with $\diag\left(
      \begin{pmatrix}
        0 & E \\ E & 0
      \end{pmatrix},\ldots,
      \begin{pmatrix}
         0 & E \\ E & 0
      \end{pmatrix}\right)$  induces the auto\m\ 
    $(g_1,g_2,g_3,g_4) \xrightarrow{\phi_4} (-g_1,g_2,g_3g_4)$
  \item When $i=n/2$, conjugation with  $
  \begin{pmatrix}
    0 & E \\ E & 0
  \end{pmatrix} \in \SL(4n,\R)$ induces the auto\m\
  $(g_1,g_2,g_3,g_4) \xrightarrow{\phi_5} (g_1,g_2,g_3,-g_4)$
  \end{itemize}
  Consider the auto\m\ groups
  \begin{equation*}
    \A(\SL(4n,\R))(G) \subset \A(\GL(4n,\R))(G) \subset \Out(G) \to
    O(q) \subset \Aut(V_0\# L(i,n-i))
  \end{equation*}
  where the outer auto\m\ group has order $196$ and $O(q)$ order $48$.
  Note that the auto\m\ $\phi_4$ of order $2$ is in the kernel of the
  homo\m\ $\Out(G) \to O(q)$ induced by abelianization $G \to
  V_0\#L(i,n-i)$.  Using the above observations we see that
  $\A(\SL(4n,\R))(G)$ maps onto $O(q)$ with a kernel of order at least
  $2$. Thus the Quillen auto\m\ group
  $\A(\GL(4n,\R))(G)$ has order $192$ or $96$. When $n$ is even and
  $i=n/2$ the auto\m\ $\phi_5$ is in $\A(\SL(4n,\R))(G)$ and when $i <
  n/2$, $\phi_5$ is not in $\A(\GL(4n,\R))(G)$ as it does not preserve
  trace. Thus
  \begin{equation*}
    \vert \A(\GL(4n,\R))(G) \vert = 
    \begin{cases}
      192 & i=n/2 \\
      96 & i < n/2
    \end{cases}
  \end{equation*}
  In any case, the group $\A(\SL(4n,\R))(G)$ equals
  $\A(\GL(4n,\R))(G)$ so that $V_0 \# L(i,n-i) \neq (V_0 \#
  L(i,n-i))^D$ (\ref{lemma:intoSL}). According to {\em magma},
  $\dim[V_0 \# L(i,n-i)]=5$.
\end{proof}

The five $\F_2\A(\pslr{4n})(V_0\#L(i,n-i))$-linear maps 
\set\begin{equation}\label{eq:basisV0L} 
  \{ ddf_{+L(2i,2n-2i)}, ddf_{0L(2i,2n-2i)}, df_{0L(2i,2n-2i)}, 
     ddf_{-L(i,n-i)},df_{0V_0} \}  
\end{equation}\add
given by
\begin{align*}
  ddf_{+L(2i,2n-2i)}[V>P>L] &= 
    \begin{cases}
      L & V=H_+\#L(2i,2n-2i),P=H_+ \\
      0 & \text{otherwise}
    \end{cases} \\
  ddf_{0L(2i,2n-2i)}[V>P>L] &= 
    \begin{cases}
      P^{\A(P)} & V=H_+\#L(2i,2n-2i),P=H_+,q(L)=0 \\
      0 & \text{otherwise}
    \end{cases} \\
    df_{0L(2i,2n-2i)}[V>P>L] &= 
    \begin{cases}
      V^{\A(V)} & V=H_+\#L(2i,2n-2i),[P,P]=0,q(L)=0 \\
      0 & \text{otherwise}
    \end{cases} \\
  ddf_{-L(i,n-i)}[V>P>L] &= 
    \begin{cases}
      L & V=H_-\#L(i,n-i),P=H_- \\
      0 & \text{otherwise}
    \end{cases} \\
  df_{0V_0}[V>P>L] &= 
    \begin{cases}
      V^{\A(V)} & V=V_0,[P,P]=0,q(L)=0 \\
      0 & \text{otherwise}
    \end{cases}
\end{align*}
constitute a basis for $[V_0\#L(i,n-i)]$.

%\begin{lemma}\label{lemma:intoV0L}
%  The $4$-dimensional objects $V_0\#L(i,n-i)$, $1 \leq i \leq [n/2]$,
%  contain the $3$-dimensional objects $H_+\#L(2i,2n-2i)$,
%  $H_-\#L(i,n-i)$, and $V_0$. Its Quillen auto\m\ group
%  $\A(V_0\#L(i,n-i))$ equals the auto\m\ group of order $48=6 \cdot
%  2^3$ 
%  \begin{equation*}
%    O(q)=
%    \begin{pmatrix}
%      \GL(2,\F_2) & 0 \\
%      \ast        & 1
%    \end{pmatrix} \qquad (\text{where $\GL(2,\F_2)$ acts on $F_2^3$})
%  \end{equation*}
%  of the quadratic function $q(v_1,v_2,v_3,v_4)=v_1^2+v_2v_3$, and
%  $[V_0\#L(i,n-i)]$ has dimension $5$ with
%  $(ddf_{+(2i,2n-2i)},ddf_{0(2i,2n-2i)},df_{0(2i,2n-2i)},
%  ddf_{-(i,n-i)},df_{0V_0})$ as basis. Thus the matrix for
%  $[H_+\#L(2i,2n-2i),H_-\#L(i,n-i),\gen{f_{0V_0}}]
%  \xrightarrow{d^2}[V_0\#L(i,n-i)]$ is the $5 \times 5$ identity matrix.
%\end{lemma}

%It is a little problematic how much of Lemma \ref{lemma:intoH+P} and
%Lemma \ref{lemma:intoV0L} really is justified by the computer
%computations. 

%%SL8.prg,SL12.prg,SL16.prg
%%/home/moller/manus/dfam/magma/lhs/steinberg/SL16.prg  (probably)

\begin{lemma}\label{lemma:intoH+P}
  The $4$-dimensional object $H_+\#P(1,i-1,2n-i,0)$, $2 < i \leq n$,
  of the category $\A(\pslr{4n})$, $n>2$, satisfies
  $H_+\#P(1,i-1,2n-i,0) = (H_+\#P(1,i-1,2n-i,0))^D$. It contains the
  $3$-dimensional objects
  \begin{equation*}
    H_+\#
    \begin{cases}
      L(1,2n-1), L(i-1,2n-i+1),L(i-1,2n-i+1)^D,L(i,2n-i) & \text{$i$
        odd}\\
      L(1,2n-1), L(i-1,2n-i+1),L(i,2n-i),L(i,2n-i)^D & \text{$i$
        even}
    \end{cases}
  \end{equation*}
  Its Quillen auto\m\ group %% ({\bf wrt a certain basis}!) 
  is
  \begin{multline*}
    \A(\pslr{4n})(H_+\#P(1,i-1,2n-i,0)) = \\ 
    \begin{cases}
    \gen{
    \begin{pmatrix}
      0&1&0&0 \\ 1&0&0&0 \\ 0&0&1&0 \\ 0&0&0&1 
    \end{pmatrix},
    \begin{pmatrix}
      1&0&1&0 \\ 0&1&0&0 \\ 0&0&1&0 \\ 0&0&0&1 
    \end{pmatrix},
    \begin{pmatrix}
      1&0&0&0 \\ 0&1&1&0 \\ 0&0&1&0 \\ 0&0&0&1 
    \end{pmatrix},
    \begin{pmatrix}
      1&0&1&1 \\ 0&1&1&1 \\ 0&0&1&0 \\ 0&0&0&1 
    \end{pmatrix} } 
      & \text{$i>2$ odd} \\
     \gen{
    \begin{pmatrix}
      0&1&0&0 \\ 1&0&0&0 \\ 0&0&1&0 \\ 0&0&0&1 
    \end{pmatrix},
    \begin{pmatrix}
      1&0&0&1 \\ 0&1&0&0 \\ 0&0&1&0 \\ 0&0&0&1 
    \end{pmatrix},
    \begin{pmatrix}
      1&0&0&0 \\ 0&1&0&1 \\ 0&0&1&0 \\ 0&0&0&1 
    \end{pmatrix},
    \begin{pmatrix}
      1&0&1&1 \\ 0&1&1&1 \\ 0&0&1&0 \\ 0&0&0&1 
    \end{pmatrix} } 
      & \text{$i>2$ even}      
    \end{cases}
  \end{multline*}
  of order $16$.  The space of equivariant maps has dimension
  $\dim [H_+\#P(1,i-1,2n-i,0)] = 16$.
  %The ordered set
%  $(df_{0(1,2n-1)},ddf_{+(i,2n-i)},ddf_{0(i,2n-i)},df_{0(i,2n-i)})$ is
%  a basis when $i=2$ and the ordered set
%  $(df_{0(1,2n-1)},ddf_{+(i-1,2n-i+1)},ddf_{0(i-1,2n-i+1)},df_{0(i-1,2n-i+1)},
%  ddf_{+(i,2n-i)},ddf_{0(i,2n-i)},df_{0(i,2n-i)})$ is linearly
%  independent when $i>2$. The matrix for the map
%  $[H_+\#L(1,2n-1),H_+\#L(2,2n-1)]
%  \xrightarrow{d^2}[H_+\#P(1,1,2n-2,0)]$ (where $i=2$) is
%  \begin{equation*}
%    \begin{pmatrix}
%      0 & 1 & 0 & 0 \\
%      0 & 0 & 1 & 0 \\
%      1 & 0 & 0 & 0 \\
%      0 & 1 & 0 & 0 \\
%      0 & 0 & 1 & 0 \\
%      0 & 0 & 0 & 1 \\
%    \end{pmatrix}
%  \end{equation*}
%  and the matrix for  $[H_+\#L(1,2n-1),H_+\#L(i-1,2n-i+1),H_+\#L(i,2n-i)]
%  \xrightarrow{d^2}[H_+\#P(1,1,2n-2,0)]$ (where $i>2$) is 
%  \begin{equation*}
%    \begin{pmatrix}
%      0 & 1 & 0 & 0 & 1 & 0 & 0 & 0 \\
%      0 & 0 & 1 & 0 & 0 & 1 & 0 & 0 \\
%      1 & 0 & 0 & 0 & 0 & 0 & 0 & 0 \\
%      0 & 1 & 0 & 0 & 0 & 0 & 0 & 0 \\
%      0 & 0 & 1 & 0 & 0 & 0 & 0 & 0 \\
%      0 & 0 & 0 & 1 & 0 & 0 & 0 & 0 \\
%      0 & 0 & 0 & 0 & 1 & 0 & 0 & 0 \\
%      0 & 0 & 0 & 0 & 0 & 1 & 0 & 0 \\
%      0 & 0 & 0 & 0 & 0 & 0 & 1 & 0 \\
%    \end{pmatrix}
%  \end{equation*}
%  with respect to a basis that contains the above linearly independent
%  set.
\end{lemma} 
\begin{proof}
  $H_+\#P(1,i-1,2n-i) \subset \pslr{4n}$ is
  (\ref{sec:non0innerprodrank4}) the quotient of
  \begin{multline*}
    G=\langle\diag(R,\ldots ,R), 
             \diag(T,\ldots ,T), \\ 
\diag(E,\overbrace{-E,\ldots,-E}^{i-1},\overbrace{E,\ldots ,E}^{2n-i}),
\diag(E,\overbrace{E,\ldots,E}^{i-1},\overbrace{-E,\ldots ,-E}^{2n-i})
    \rangle = \gen{g_1,g_2,g_3,g_4} \subset \SL(4n,\R)
  \end{multline*}
  The centralizer of $G$ in $\GL(4n,\R)$ is contained in the
  centralizer of its subgroup $2^{1+2}_+$ which is contained in
  $\SL(4n,\R)$ (\ref{exmp:IDQ}). This means (\ref{eq:NOut}) that the
  elements of $\A(\GL(4n,\R))(G)$ and
  $\A(\PGL(4n,\R))(H_+\#P(1,i-1,2n-i,0))$ have a well-defined sign. 
  The Quillen auto\m\ group is
  contained in the group $
  \begin{pmatrix}
    O^+(2,\F_2) & \ast \\ 0 & E
  \end{pmatrix}$ of order $2^5=32$.
  Observe that
  \begin{itemize}
  \item $R$ and $T$ are conjugate in $\GL(2,\R)$ so that
    $(g_1,g_2,g_3,g_4) \xrightarrow{\phi_1} (g_2,g_1,g_3,g_4)$ is in
    the Quillen auto\m\ group and has sign $+1$
  \item Conjugation with
    $\diag(E,\overbrace{E,\ldots,E}^{i-1},\overbrace{T,\ldots,T}^{2n-i})$
    induces $(g_1,g_2,g_3,g_4) \xrightarrow{\phi_2}
    (g_1g_4,g_2,g_3,g_4)$ of sign $(-1)^i$
  \item Conjugation with
    $\diag(E,\overbrace{T,\ldots,T}^{i-1},\overbrace{E,\ldots,E}^{2n-i})$
    induces $(g_1,g_2,g_3,g_4) \xrightarrow{\phi_3}
    (g_1g_3,g_2,g_3,g_4)$ of sign $-(-1)^i$
  \item  Conjugation with
    $\diag(E,\overbrace{E,\ldots,E}^{i-1},\overbrace{R,\ldots,R}^{2n-i})$
    induces $(g_1,g_2,g_3,g_4) \xrightarrow{\phi_4}
    (g_1,g_2g_4,g_3,g_4)$ of sign $(-1)^i$
  \item Conjugation with
    $\diag(E,\overbrace{R,\ldots,R}^{i-1},\overbrace{E,\ldots,E}^{2n-i})$
    induces $(g_1g_3,g_2,g_3,g_4) \xrightarrow{\phi_5}
    (g_1,g_2g_3,g_3,g_4)$ of sign $-(-1)^i$
  \item Conjugation with
    $\diag(E,\overbrace{RT,\ldots,RT}^{i-1},\overbrace{RT,\ldots,RT}^{2n-i})$
    induces the auto\m\ given by $(g_1,g_2,g_3,g_4) \xrightarrow{\phi_6}
    (g_1g_3g_4,g_2g_3g_4,g_3,g_4)$ of sign $+1$.
  \end{itemize}
  It follows that $N_{\GL(4n,\R)}(G) \not\subset\SL(4n,\R)$ as this
  normalizer contains elements of negative determinant regardless of
  the parity of $i$. Also, $\A(\pslr{4n})(H_+\#P(1,i-1,2n-i,0))$ is
  generated by (the auto\m s induced by) $\phi_1$, $\phi_2$, $\phi_4$,
  and $\phi_6$ when $i$ is even, and $\phi_1$, $\phi_3$, $\phi_5$, and
  $\phi_6$ when $i$ is odd.
\end{proof}

%% COMMENT
%It may be possible to get a more invariant form of this proof, less
%dependent on computer computations, like this: Let the Quillen auto\m\
%group $\A(\pslr{4n})(E)$ act on the $3$-flags or on the subset of
%$3$-flags of the form $[V>P>L$ where $V=H_+\#L(i-1,2n-i)$, select one
%of the orbits and define a map as we defined $ddf_+$, $ddf_0$, or
%$df_0$. These maps might be independent elements of
%$[H_+\#P(1,i-1,2n-i,0)$; see the experiments in SL12.add.
%% END COMMENT

The fourteen $\F_2\A(\pslr{4n})(H_+\#P(1,i-1,2n-1))$-linear maps
\set\begin{multline}\label{eq:basispH+P}
  \{ddf_{+L(i-1,2n-i+1)},ddf_{+L(i-1,2n-i+1)}^D,ddf_{0L(i-1,2n-i+1)}, 
    ddf_{0L(i-1,2n-i+10)}^D, \\
    df_{0L(i-1,2n-i+1)},df_{0L(i-1,2n-i+1)}^D, \\
   ddf_{+L(i,2n-i)},ddf_{+L(i,2n-i)}^D,ddf_{0L(i,2n-i)}, 
   ddf_{0L(i,2n-i)}^D,df_{0L(i,2n-i)},df_{0L(i,2n-i)}^D, \\
   df_{0L(1,2n-1)},df_{0L(1,2n-1)}^D \} 
\end{multline}\add
form a partial basis for the $16$-dimensional vector space
$[H_+\#P(1,i-1,2n-1)]$, $2<i\leq n$.  
 For $1 < i \leq n$ and $i$ odd,
   \begin{align*}
     ddf_{+L(i,2n-i)}[V>P>L] &=
     \begin{cases}
       L & V=H_+\#L(i,2n-i),P=H_+ \\
       0 & \text{otherwise}
     \end{cases} \\
     ddf_{+L(i,2n-i)}^D[V>P>L] &=
     \begin{cases}
       L & V=H_+\#L(i,2n-i),P=H_+^D \\
       0 & \text{otherwise}
     \end{cases} \\
     ddf_{0L(i,2n-i)}[V>P>L] &= 
    \begin{cases}
      P^{\A(P)} & V=H_+\#L(i,2n-i),P=H_+,q(L)=0 \\
      0 & \text{otherwise}
    \end{cases} \\
    ddf_{0L(i,2n-i)}^D[V>P>L] &= 
    \begin{cases}
      P^{\A(P)} & V=H_+\#L(i,2n-i),P=H_+^D,q(L)=0 \\
      0 & \text{otherwise}
    \end{cases} \\
    df_{01L(i,2n-i)}[V>P>L] &= 
    \begin{cases}
      V \cap O_1 & V=H_+\#L(i,2n-i),[P,P]=0,q(L)=0 \\
      0 & \text{otherwise}
    \end{cases}  \\
    df_{01L(i,2n-i)}^D[V>P>L] &= 
    \begin{cases}
      V \cap O_2 & V=H_+\#L(i,2n-i),[P,P]=0,q(L)=0 \\
      0 & \text{otherwise}
    \end{cases}
   \end{align*}
   where (in the last two formulas), $O_1$ and $O_2$ are the two
   orbits of length $2$ for the action of
   $\A(\pslr{4n})(H_+\#P(1,i-1,2n-i,0))$ on $H_+\#P(1,i-1,2n-i,0)$.
   Each of the hyperplanes isomorphic to $V=H_+\#L(i,2n-i)$ contains
   precisely one vector $v_1$ from $O_1$ and one vector $v_2$ from
   $O_2$ and $\{v_1,v_2\}$ is a basis for the fixed point group
   $V^{\A(\pslr{4n})(V)}$. 
   For $1<i \leq n$ and $i$ even,
  \begin{align*}
    ddf_{+L(i,2n-i)}[V>P>L] &= 
    \begin{cases}
      L & V=H_+\#L(i,2n-i),P=H_+ \\
      0 & \text{otherwise}
    \end{cases} \\
        ddf_{+L(i,2n-i)}^D[V>P>L] &= 
    \begin{cases}
      L & V=(H_+\#L(i,2n-i))^D,P=(H_+)^D \\
      0 & \text{otherwise}
    \end{cases} \\
    ddf_{0L(i,2n-i)}[V>P>L] &= 
    \begin{cases}
      P^{\A(P)} & V=H_+\#L(i,2n-i),P=H_+,q(L)=0 \\
      0 & \text{otherwise}
    \end{cases} \\
        ddf_{0L(i,2n-i)}^D[V>P>L] &= 
    \begin{cases}
      P^{\A(P)} & V=(H_+\#L(i,2n-i))^D,P=(H_+)^D,q(L)=0 \\
      0 & \text{otherwise}
    \end{cases} \\
    df_{0L(i,2n-i)}[V>P>L] &= 
    \begin{cases}
      V^{\A(V)} & V=H_+\#L(i,2n-i),[P,P]=0,q(L)=0 \\
      0 & \text{otherwise}
    \end{cases} \\
        df_{0L(i,2n-i)}^D[V>P>L] &= 
    \begin{cases}
      V^{\A(V)} & V=(H_+\#L(i,2n-i))^D,[P,P]=0,q(L)=0 \\
      0 & \text{otherwise}
    \end{cases}
  \end{align*}

\begin{lemma}\label{lemma:intoH+P1120}
  The $4$-dimensional object $H_+\#P(1,1,2,0)$ of the category
  $\A(\pslr{8})$ satisfies $H_+\#P(1,1,2,0) = (H_+\#P(1,1,2,0))^D$. It
  contains the $3$-dimensional objects
  \begin{equation*}
    H_+\#L(1,3),  H_+\#L(2,2), (H_+\#L(2,2))^D
  \end{equation*}
  Its Quillen auto\m\ group %% ({\bf wrt a certain basis}!) 
  is
  \begin{multline*}
    \A(\pslr{8})(H_+\#P(1,1,2,0)) = \\ 
     \gen{
    \begin{pmatrix}
      0&1&0&0 \\ 1&0&0&0 \\ 0&0&1&0 \\ 0&0&0&1 
    \end{pmatrix},
    \begin{pmatrix}
      1&0&0&1 \\ 0&1&0&0 \\ 0&0&1&0 \\ 0&0&0&1 
    \end{pmatrix},
    \begin{pmatrix}
      1&0&0&0 \\ 0&1&0&1 \\ 0&0&1&0 \\ 0&0&0&1 
    \end{pmatrix},
    \begin{pmatrix}
      1&0&1&1 \\ 0&1&1&1 \\ 0&0&1&0 \\ 0&0&0&1 
    \end{pmatrix},
    \begin{pmatrix}
     1&0&0&0 \\ 0&1&0&0 \\ 0&0&1&1 \\ 0&0&0&1 
    \end{pmatrix} }
  \end{multline*}
  of order $32$, and
    $\dim [H_+\#P(1,i-1,2n-i,0)] = 8$.
  %The ordered set
%  $(df_{0(1,2n-1)},ddf_{+(i,2n-i)},ddf_{0(i,2n-i)},df_{0(i,2n-i)})$ is
%  a basis when $i=2$ and the ordered set
%  $(df_{0(1,2n-1)},ddf_{+(i-1,2n-i+1)},ddf_{0(i-1,2n-i+1)},df_{0(i-1,2n-i+1)},
%  ddf_{+(i,2n-i)},ddf_{0(i,2n-i)},df_{0(i,2n-i)})$ is linearly
%  independent when $i>2$. The matrix for the map
%  $[H_+\#L(1,2n-1),H_+\#L(2,2n-1)]
%  \xrightarrow{d^2}[H_+\#P(1,1,2n-2,0)]$ (where $i=2$) is
%  \begin{equation*}
%    \begin{pmatrix}
%      0 & 1 & 0 & 0 \\
%      0 & 0 & 1 & 0 \\
%      1 & 0 & 0 & 0 \\
%      0 & 1 & 0 & 0 \\
%      0 & 0 & 1 & 0 \\
%      0 & 0 & 0 & 1 \\
%    \end{pmatrix}
%  \end{equation*}
%  and the matrix for  $[H_+\#L(1,2n-1),H_+\#L(i-1,2n-i+1),H_+\#L(i,2n-i)]
%  \xrightarrow{d^2}[H_+\#P(1,1,2n-2,0)]$ (where $i>2$) is 
%  \begin{equation*}
%    \begin{pmatrix}
%      0 & 1 & 0 & 0 & 1 & 0 & 0 & 0 \\
%      0 & 0 & 1 & 0 & 0 & 1 & 0 & 0 \\
%      1 & 0 & 0 & 0 & 0 & 0 & 0 & 0 \\
%      0 & 1 & 0 & 0 & 0 & 0 & 0 & 0 \\
%      0 & 0 & 1 & 0 & 0 & 0 & 0 & 0 \\
%      0 & 0 & 0 & 1 & 0 & 0 & 0 & 0 \\
%      0 & 0 & 0 & 0 & 1 & 0 & 0 & 0 \\
%      0 & 0 & 0 & 0 & 0 & 1 & 0 & 0 \\
%      0 & 0 & 0 & 0 & 0 & 0 & 1 & 0 \\
%    \end{pmatrix}
%  \end{equation*}
%  with respect to a basis that contains the above linearly independent
%  set.
\end{lemma}
\begin{proof} The proof is similar to that of
  \ref{lemma:intoH+P}. $H_+\#P(1,1,2,0) \subset \pslr{8}$ is the
  quotient of
  \begin{equation*}
    G=\gen{\diag(R,R,R,R),\diag(T,T,T,T),\diag(E,-E,E,E),\diag(E,E,-E,-E)}  
    \subset \SL(8,\R)
  \end{equation*}
   The extra element of $\A(\pslr{8})(H_+\#P(1,1,2,0))$ is induced by
   conjugation with the matrix 
   $\diag\left(
     \begin{pmatrix}
       0&E\\E&0
     \end{pmatrix},
     \begin{pmatrix}
       E&0\\0&E
     \end{pmatrix}\right) \in \SL(8,\R)$. According to {\em magma},
   $\dim [H_+\#P(1,1,2,0)] = 8$.
\end{proof}
The eight $\F_2\A(\pslr{8})(H_+\#P(1,1,2,0))$-linear maps
\set\begin{equation}
  \label{eq:basisH+P112}
  \{ ddf_{+L(2,2)},ddf_{+L(2,2)}^D, ddf_{0L(2,2)},ddf_{0L(2,2)}^D, 
     df_{0L(2,2)},df_{0L(2,2)}^D, df_{01L(1,3)},df_{01L(1,3)}^D \}        
\add\end{equation}
is a basis for the vector space $[H_+\#P(1,1,2,0)]$.

We are now ready to describe the differentials $d^1$ and $d^2$ in
Oliver's cochain complex (\ref{eq:Olccc}) for computing the limits of
the functor $\pi_1(BZC_{\pslr{4n}}(V))=V$ on the category
$\A(\pslr{4n})$.   
The $6 \times
(6n+2[n/2]+8)$ matrix for $d^1$ is of the following form (shown here
for $n=3$) 
\begin{center}
  \begin{tabular}[c]{l|cccc}
    {} & $[H_+\#L(1,5)]$ & $[H_+\#L(2,4)] \times [H_+\#L(2,4)]^D$ &
     $H_+\#L(3,3)$  \\ \hline
    $[H_+]$ & $
    \begin{pmatrix}
      A & 0
    \end{pmatrix}$  & $
    \begin{pmatrix}
      A & 0
    \end{pmatrix}$ & $
    \begin{pmatrix}
      A & 0
    \end{pmatrix}$  \\
     $[H_+]^D$ & $
    \begin{pmatrix}
      0 & A
    \end{pmatrix}$  & $
    \begin{pmatrix}
      0 & A
    \end{pmatrix}$ & $
    \begin{pmatrix}
      0 & A
    \end{pmatrix}$ \\
    $[H_-]$  \\
    $[H_-]^D$  \\
  \end{tabular}
 \begin{tabular}[c]{cc|l}
    $[H_-\#L(1,1)] \times [H_-\#L(1,1)]^D$ & $[V_0] \times [V_0]^D$ & {} \\
    \hline 
     {} & $
    \begin{pmatrix} 
      H & 0
    \end{pmatrix}$ & $[H_+]$  \\
      {} & $
    \begin{pmatrix}
      0 & H
    \end{pmatrix}$ & $[H_+]^D$  \\
     $
    \begin{pmatrix}
      1&0
    \end{pmatrix}$ & $
    \begin{pmatrix}
      B&0
    \end{pmatrix}$ & $[H_-]$  \\
       $
    \begin{pmatrix}
      0&1
    \end{pmatrix}$ & $
    \begin{pmatrix}
      0&B
    \end{pmatrix}$ & $[H_-]^D$
  \end{tabular}
  %\begin{tabular}[c]{l|ccc}
%    {} & $[H_-\#L(1,1)] \times [H_-\#L(1,1)]^D$ & $[V_0] \times [V_0]^D$ \\
%    \hline 
%    $[H_+]$ & {} & $
%    \begin{pmatrix} 
%      H & 0
%    \end{pmatrix}$ \\
%     $[H_+]^D$ & {} & $
%    \begin{pmatrix}
%      0 & H
%    \end{pmatrix}$ \\
%    $[H_-]$ & $
%    \begin{pmatrix}
%      1&0
%    \end{pmatrix}$ & $
%    \begin{pmatrix}
%      B&0
%    \end{pmatrix}$ \\
%    $[H_-]^D$  & $
%    \begin{pmatrix}
%      0&1
%    \end{pmatrix}$ & $
%    \begin{pmatrix}
%      0&B
%    \end{pmatrix}$
%  \end{tabular}
\end{center}
 where
\begin{equation*}
  A=
  \begin{pmatrix}
    1&0&0 \\ 0&1&0
  \end{pmatrix}, \qquad
  H=
  \begin{pmatrix}
    1&0&0&0 \\ 0&1&0&0
  \end{pmatrix}, \qquad
  B=
  \begin{pmatrix}
    0&0&1&0
  \end{pmatrix}
\end{equation*}
is injective so $\lim^1=0$. Exactness is thus equivalent to
\begin{equation*}
  \dim(\im d^2) \geq 6n + 2[n/2] +2
\end{equation*}
We shall show this by  mapping the
$n+[n/2]+2[n/2]+2$ objects of dimension $3$, 
\begin{multline*}
  H_+\#L(i,2n-i), (H_+\#L(i,2n-i))^D \text{($i$ even)}, \quad 
  1 \leq i \leq n,  \\
  H_-\#L(i,n-i),\; (H_-\#L(i,n-i))^D, \quad 1 \leq i \leq [n/2],
  \qquad V_0,\; V_0^D,
\end{multline*}
of $\A(\pslr{4n})$ to the $n-2+2[n/2]$  objects of dimension $4$,
\begin{equation*}
   H_+\#P(1,i-1,2n-i,0),\quad 2 < i \leq n, \qquad 
   V_0\#L(n-i,i),\, (V_0\#L(n-i,i))^D, \quad 1 \leq i \leq [n/2],
\end{equation*}
for $n>2$ and to 
\begin{equation*}
   H_+\#P(1,1,2,0), \quad V_0\#L(1,1),  \quad (V_0\#L(1,1))^D
\end{equation*}
when $n=2$.  The $(6n+2[n/2]+8) \times
(16(n-2)+10[n/2])$-matrix for $d^2$ (shown here for $n=5$) is
\begin{center}
\begin{tabular}[c]{l|ccc}
& $[H_+\#P(1,2,7)]$ & $[H_+\#P(1,3,6)]$ & $[H_+\#P(1,4,5)]$  \\ \hline
$[H_+\#L(1,9)]$ & $
\begin{pmatrix}
  A & A & B
\end{pmatrix}$ &
 $
\begin{pmatrix}
  A & A & B
\end{pmatrix}$ &
 $
\begin{pmatrix}
  A & A & B
\end{pmatrix}$ \\
$[H_+\#L(2,8)] \times [H_+\#L(2,8)]^D$ &
$
\begin{pmatrix}
  E&0&0
\end{pmatrix}$ & {} & {} \\
$[H_+\#L(3,7)]$ & $
\begin{pmatrix}
  0&E&0
\end{pmatrix}$ & $
\begin{pmatrix}
  E&0&0
\end{pmatrix}$ & {} \\
$[H_+\#L(4,6)] \times [H_+\#L(4,6)]^D$ & {} & $
\begin{pmatrix}
  0&E&0
\end{pmatrix}$ & $
\begin{pmatrix}
  E&0&0
\end{pmatrix}$ \\
$[H_+\#L(5,5)]$ & {} & {} & {} $
\begin{pmatrix}
  0&E&0
\end{pmatrix}$ \\
$[H_-\#L(1,4)] \times [H_-\#L(1,4)]^D$ \\
    $[H_-\#L(2,3)] \times [H_-\#L(2,3)]^D$\\
    $[V_0] \times [V_0]^D$
\end{tabular}

  \begin{tabular}[c]{cccc|l}
     $[V_0\#L(1,4)]$ & $[V_0\#L(1,4)]^D$ 
    & $[V_0\#L(2,3)]$ & $[V_0\#L(2,3)]^D$ & {} \\ \hline
    {} & {} & {} & {} & $[H_+\#L(1,9)]$ \\
     $
    \begin{pmatrix}
      H \\ 0
    \end{pmatrix}$ & $
    \begin{pmatrix}
      0 \\ H
    \end{pmatrix}$  &{} &{} & $[H_+\#L(2,8)] \times [H_+\#L(2,8)]^D$  \\
    {} & {} & {} & {} & $[H_+\#L(3,7)]$ \\
     {} & {} & $
    \begin{pmatrix}
      H \\ 0
    \end{pmatrix}$ & $
    \begin{pmatrix}
      0 \\ H
    \end{pmatrix}$ & $[H_+\#L(4,6)] \times [H_+\#L(4,6)]^D$  \\
     {} & {} & {} & {} & $[H_+\#L(5,5)]$ \\
     $
    \begin{pmatrix}
      L \\ 0
    \end{pmatrix}$ & $
    \begin{pmatrix}
      0 \\ L 
    \end{pmatrix}$ & {} & {} & $[H_-\#L(1,4)] \times [H_-\#L(1,4)]^D$  \\
     {} & {} $
    \begin{pmatrix}
      L \\ 0
    \end{pmatrix}$ & $
    \begin{pmatrix}
      0 \\ L
    \end{pmatrix}$ & {} & $[H_-\#L(2,3)] \times [H_-\#L(2,3)]^D$  \\
     $
    \begin{pmatrix}
      K \\ 0
    \end{pmatrix}$ & $
    \begin{pmatrix}
      K \\ 0
    \end{pmatrix}$ & $
    \begin{pmatrix}
      K \\ 0
    \end{pmatrix}$ & $
    \begin{pmatrix}
      K \\ 0
    \end{pmatrix}$ & $[V_0] \times [V_0]^D$  
   \end{tabular}
%%%%%%%%%%%%%%%%%%%%%%%%%%%%%%%%%%%%%
%  \begin{tabular}[c]{l|cccc}
%    {} & $[V_0\#L(1,4)]$ & $[V_0\#L(1,4)]^D$ 
%    & $[V_0\#L(2,3)]$ & $[V_0\#L(2,3)]^D$ \\ \hline
%    $[H_+\#L(1,9)]$ \\
%    $[H_+\#L(2,8)] \times [H_+\#L(2,8)]^D$ & $
%    \begin{pmatrix}
%      H \\ 0
%    \end{pmatrix}$ & $
%    \begin{pmatrix}
%      0 \\ H
%    \end{pmatrix}$ \\
%    $[H_+\#L(3,7)]$ \\
%    $[H_+\#L(4,6)] \times [H_+\#L(4,6)]^D$ & {} & {} & $
%    \begin{pmatrix}
%      H \\ 0
%    \end{pmatrix}$ & $
%    \begin{pmatrix}
%      0 \\ H
%    \end{pmatrix}$ \\
%    $[H_+\#L(5,5)]$ \\
%    $[H_-\#L(1,4)] \times [H_-\#L(1,4)]^D$ & $
%    \begin{pmatrix}
%      L \\ 0
%    \end{pmatrix}$ & $
%    \begin{pmatrix}
%      0 \\ L 
%    \end{pmatrix}$ \\
%    $[H_-\#L(2,3)] \times [H_-\#L(2,3)]^D$ & {} & {} $
%    \begin{pmatrix}
%      L \\ 0
%    \end{pmatrix}$ & $
%    \begin{pmatrix}
%      0 \\ L
%    \end{pmatrix}$ \\
%    $[V_0] \times [V_0]^D$ & $
%    \begin{pmatrix}
%      K \\ 0
%    \end{pmatrix}$ & $
%    \begin{pmatrix}
%      K \\ 0
%    \end{pmatrix}$ & $
%    \begin{pmatrix}
%      K \\ 0
%    \end{pmatrix}$ & $
%    \begin{pmatrix}
%      K \\ 0
%    \end{pmatrix}$ 
%   \end{tabular}
\end{center}
where %each entrance contains a $6 \times (6+6+4)$-matrix and 
\begin{multline*}
  A=
  \begin{pmatrix}
    1&0&0&0&0&0 \\ 0&1&0&0&0&0 \\  0&0&1&0&0&0 \\  0&0&0&1&0&0 \\
    0&0&0&0&0&0 \\  0&0&0&0&0&0
  \end{pmatrix}, \qquad
  B=
  \begin{pmatrix}
    0&0&0&0 \\  0&0&0&0 \\  0&0&0&0 \\  0&0&0&0 \\ 1&0&0&0 \\ 0&1&0&0
  \end{pmatrix}, \qquad
  K=
  \begin{pmatrix}
       1&0&0&0&0 \\ 0&0&0&0&0 \\ 0&0&0&1&0 \\ 0&0&0&0&1 
  \end{pmatrix} \\
  H=
  \begin{pmatrix}
      1&0&0&0&0 \\ 0&1&0&0&0 \\ 0&0&1&0&0  
  \end{pmatrix}, \qquad
  L=
  \begin{pmatrix}
    0&0&0&1&0
  \end{pmatrix}
\end{multline*}
while $E$ is $(6 \times 6)$ unit matrix and $0$ a zero
matrix. 
%%%%%%%%%%%%%%%%%%%%%%%%%%%%%%%%%%%
%\begin{center}
%  \begin{tabular}[c]{l|cccc}
%    {} & $[V_0\#L(1,4)]$ & $[V_0\#L(1,4)]^D$ 
%    & $[V_0\#L(2,3)]$ & $[V_0\#L(2,3)]^D$ \\ \hline
%    $[V_0]$ & {$
%      \begin{pmatrix}
%       1&0&0&0&0 \\  0&0&0&0&0 \\ 0&0&0&1&0 \\ 0&0&0&0&1 
%      \end{pmatrix}$ } & {$
%      \begin{pmatrix}
%       1&0&0&0&0 \\  0&0&0&0&0 \\ 0&0&0&1&0 \\ 0&0&0&0&1
%      \end{pmatrix}$}\\
%    $[H_+\#L(2,8)] \times [H_-\#L(1,4)]$ & {$
%      \begin{pmatrix}
%       1&0&0&0&0 \\  0&1&0&0&0 \\ 0&0&1&0&0 \\ 0&0&0&1&0 
%      \end{pmatrix}$} & {} \\
%    $[H_+\#L(4,6)] \times [H_-\#L(2,3)]$ & {} & {$
%      \begin{pmatrix}
%       1&0&0&0&0 \\  0&1&0&0&0 \\ 0&0&1&0&0 \\ 0&0&0&1&0 
%      \end{pmatrix}$} 
%  \end{tabular}
%\end{center}
%%%%%%%%%%%%%%%%%%%%%%%%%%%%
%\begin{center}
%  \begin{tabular}[c]{l|cc}
%    {} & $[V_0\#L(1,4)]$ & $[V_0\#L(2,3)]$ \\ \hline
%    $[V_0]$ & {$
%      \begin{pmatrix}
%      1&0&0&0&0 \\ 0&1&0&0&0 \\ 0&0&1&0&0 \\ 0&0&0&0&1 
%      \end{pmatrix}$ } & {$
%      \begin{pmatrix}
%      1&0&0&0&0 \\ 0&1&0&0&0 \\ 0&0&1&0&0 \\ 0&0&0&0&1 
%      \end{pmatrix}$}\\
%    $[H_+\#L(2,8)] \times [H_-\#L(1,4)]$ & {$
%      \begin{pmatrix}
%      1&0&0&0&0 \\  0&1&0&0&0 \\  0&0&1&0&0 \\ 0&0&0&1&0
%      \end{pmatrix}$} & {} \\
%    $[H_+\#L(4,6)] \times [H_-\#L(2,3)]$ & {} & {$
%      \begin{pmatrix}
%      1&0&0&0&0 \\  0&1&0&0&0 \\  0&0&1&0&0 \\ 0&0&0&1&0
%      \end{pmatrix}$} 
%  \end{tabular}
%\end{center}
These matrices are given with respect to the bases (\ref{eq:basisV0},
\ref{eq:basisH+L}, \ref{eq:basisH-L}, \ref{eq:basispH+P},
\ref{eq:basisV0L}).

The case $n=2$ of $\pslr{8}$ is special. Part of the matrix
for $d^2$ is the $(22 \times 18)$-matrix
%%%%%%%%%%%%%%%%%%%%%%%%%%%%%%%%%%%%%%%%%%%%%%
%\begin{center}
%  \begin{tabular}[c]{l|ccc}
%    {} & $[H_+\#P(1,1,2,0)]$ & $[V_0\#L(1,1)]$ & $[V_0\#L(1,1)]^D$  \\ \hline
%    $[H_+\#L(1,3)]$ & $
%    \begin{pmatrix}
%      A&B
%    \end{pmatrix}$ & {} & {} \\
%    $[H_+\#L(2,2)] \times [H_+\#L(2,2)]^D$ & $
%    \begin{pmatrix}
%      E&0
%    \end{pmatrix}$ & $
%    \begin{pmatrix}
%      1&0&0&0&0 \\ 0&1&0&0&0 \\ 0&0&1&0&0 \\ 
%      0&0&0&0&0 \\ 0&0&0&0&0 \\ 0&0&0&0&0 \\ 
%    \end{pmatrix}$ & $
%    \begin{pmatrix}
%      0&0&0&0&0 \\ 0&0&0&0&0 \\ 0&0&0&0&0 \\ 
%      1&0&0&0&0 \\ 0&1&0&0&0 \\ 0&0&1&0&0 \\ 
%    \end{pmatrix}$ \\
%    $[H_-\#L(1,1)] \times [H_-\#L(1,1)]^D$ & {} & $
%    \begin{pmatrix}
%       0&0&0&1&0 \\ 0&0&0&0&0
%    \end{pmatrix}$ & $
%    \begin{pmatrix}
%       0&0&0&0&0 \\ 0&0&0&1&0 
%    \end{pmatrix}$ \\
%    $[V_0] \times [V_0]^D$ & {} & $
%    \begin{pmatrix}
%       1&0&0&0&0 \\ 0&0&0&0&0 \\ 0&0&0&1&0 \\ 0&0&0&0&1 \\
%       0&0&0&0&0 \\ 0&0&0&0&0 \\ 0&0&0&0&0 \\ 0&0&0&0&0 
%    \end{pmatrix}$ & $
%    \begin{pmatrix}
%       0&0&0&0&0 \\ 0&0&0&0&0 \\ 0&0&0&0&0 \\ 0&0&0&0&0 \\
%       1&0&0&0&0 \\ 0&0&0&0&0 \\ 0&0&0&1&0 \\ 0&0&0&0&1 \\
%    \end{pmatrix}$
%  \end{tabular}
%\end{center}
%%%%%%%%%%%%%%%%%%%%%%%%%%%%%%%%%%%%%%%%%%%%%
\begin{center}
  \begin{tabular}[c]{l|ccc}
    {} & $[H_+\#P(1,1,2,0)]$ & $[V_0\#L(1,1)]$ & $[V_0\#L(1,1)]^D$  \\ \hline
    $[H_+\#L(1,3)]$ & $
    \begin{pmatrix}
      A&B
    \end{pmatrix}$ & {} & {} \\
    $[H_+\#L(2,2)] \times [H_+\#L(2,2)]^D$ & $
    \begin{pmatrix}
      E&0
    \end{pmatrix}$ & $
    \begin{pmatrix}
      H \\ 0
    \end{pmatrix}$ & $
    \begin{pmatrix}
     0 \\ H
    \end{pmatrix}$ \\
    $[H_-\#L(1,1)] \times [H_-\#L(1,1)]^D$ & {} & $
    \begin{pmatrix}
       L \\ 0
    \end{pmatrix}$ & $
    \begin{pmatrix}
       0 \\ L
    \end{pmatrix}$ \\
    $[V_0] \times [V_0]^D$ & {} & $
    \begin{pmatrix}
      K \\ 0
    \end{pmatrix}$ & $
    \begin{pmatrix}
      0 \\ K
    \end{pmatrix}$
  \end{tabular}
\end{center}
where now 
%% SL8.add, SL16.prg
\begin{equation*}
%   A=
%  \begin{pmatrix}
%    1&0&0&0&0&0 \\ 0&1&0&0&0&0 \\  0&0&1&0&0&0 \\  0&0&0&1&0&0 \\
%    0&0&0&0&0&0 \\  0&0&0&0&0&0
%  \end{pmatrix}, \quad
  B=
  \begin{pmatrix}
    0&0 \\  0&0 \\  0&0 \\  0&0 \\ 1&0 \\ 0&1
  \end{pmatrix}
\end{equation*}
while $E$ is $6 \times 6$ unit matrix and $0$ a zero matrix. 
As (partial)
bases we use
the ordered sets (\ref{eq:basisH+L}, \ref{eq:basisH-L},
\ref{eq:basisV0}, \ref{eq:basisH+P112}, \ref{eq:basisV0L},
\ref{eq:basisH+L}).  This
matrix has rank $16$.

\begin{cor}
  The partial differential
  \begin{multline*}
    \prod_{
      \begin{array}[c]{c}
       1\leq i \leq n \\ \text{$i$ odd}
      \end{array}}[H_+\#L(i,2n-i)] \times
     \prod_{
      \begin{array}[c]{c}
       1\leq i \leq n \\ \text{$i$ even}
      \end{array}}[H_+\#L(i,2n-i)] \times [H_+\#L(i,2n-i)]^D \\
     \times 
    \prod_{1\leq i \leq [n/2]}[H_-\#L(i,n-i)] \times [H_-\#L(i,n-i)]^D
   \times [V_0] \times [V_0]^D\\
    \xrightarrow{d^2} 
    \prod_{2 < i \leq n}[H_+\#P(1,i-1,2n-i,0)] \times 
     \prod_{1\leq i \leq [n/2]}[V_0\#L(i,n-i)]
  \end{multline*}
  has rank $6n+2[n/2]+2$.
\end{cor}
\begin{proof}
  By now we know a matrix for this linear map so we simply check 
  its rank.
\end{proof}

%%%%%%%%%%%%%  Thu Aug 12 13:27:10 CEST 2004 %%%%%%%%%%%%%%%%%
\begin{proof}[Proof of Lemma~\ref{lemma:lim=0}]
  For $\pi_2$ use that it is trivial on the objects with $[,] \neq 0$.
\end{proof}

%\begin{tabular}[c]{|l||c|c|c|c|c|c|c|c|c|c|}\hline
%rank & \multicolumn{2}{|c|}{2} &
%\multicolumn{3}{c|}{3} &
%\multicolumn{5}{c|}{4} \\ \hline
%$V^*$ & $2^{1+2}_+$ & $2^{1+2}_-$ &
%$2^{1+2}_+ \times 2$ & $2^{1+2}_- \times 2$ & $2^{1+2}_{\pm} \circ 4$
%&
%$2^{1+2}_+ \times 2^2$ & $2^{1+2}_- \times 2^2$ & $2^{1+2}_{\pm} \circ
%4 \times 2$ &
%$2^{1+4}_+$ & $2^{1+4}_-$  \\ \hline
%$C(V)$ & $V \times \PGL(2m,\R)$ & $V \times \PGL(m,\Ha)$ & 
%\end{tabular}

\section{The category $\A(\pslr{4n})^{[\; ,\; ] \neq 0}_{\leq 4}$}  
\label{sec:nontoral}

%\noindent
%$\underline{q(V)=0}$: Iso\m\ classes of nontoral objects of
%$\A(\pslr{2n})$ with $q=0$ correspond bijectively (\ref{lemma:Vast},
%\ref{lemma:V*conj}) to the set
%\begin{equation*}
%  V^{\vee} \backslash \{ (V,i) \in \Ob(\A(\SL(2n,\R))_{\not\leq t}) \mid
%  \Eq(S(i))=0 \} / \GL(V)
%\end{equation*}
%In particular, iso\m\ classes of nontoral rank two objects $V$ with
%$q(V)=0$ correspond (\ref{exmp:nontoralqeq0}) to partitions
%$2n=i_0+i_1+i_2+i_3$, $i_0 \geq i_1 \geq i_2 \geq i_3 \geq 1$,
%$i_0,i_1,i_2,i_3$ odd, of $2n$ into sums of four odd positive
%integers.

%\noindent
%$\underline{q(V) \neq 0}$: What can we say about nontoral objects with
%$q \neq 0$?

We shall need information about all objects of $\A(\pslr{4n})^{[\; ,\;
  ] \neq 0}$ of rank $\leq 3$ and some objects of rank $4$. If $V
\subset \pslr{4n}$ is a nontoral \lmntwo\ with nontrivial inner
product then its preimage $V^* \subset \SL(4n,\R)$ is $P \times R(V)$
or $(C_4 \circ P) \times R(V)$ where $P$ is an extraspecial $2$-group,
$C_4 \circ P$ a generalized extraspecial $2$-group, and
$\mho_1(V^*)=\gen{-E}$ (\ref{lemma:Vast}).  We manufacture all oriented
real representations of these product groups as direct sums of tensor
products of irreducible representations of the factors
(\ref{sec:tensor}).

\setcounter{subsection}{\value{thm}}
\subsection{Rank two objects with nontrivial inner product}
\label{sec:non0innerprodrank2}\add

The category $\A(\pslr{4n})$ contains up to iso\m\ four rank two
objects with nontrivial inner product, $H_{\pm}$ and $H_{\pm}^D$.  The
\lmntwo\ $H_{\pm} \subset \pslr{4n}$ is the quotient of the
extraspecial $2$-group $2^{1+2}_{\pm} \subset \SL(4n,\R)$ with
$\mho_1(2^{1+2}_{\pm})=\gen{-E}$ described in
\ref{exmp:IDQ}.(\ref{exmp:IDQ6}) and
\ref{exmp:IDQ}.(\ref{exmp:IDQ7}).
%\begin{exmp}\label{D8Q8}
Their centralizers \cite[Proposition 4]{bob:stubborn} in $\SL(4n,\R)$
and $\PSL(4n,\R)$ are
 \begin{alignat*}{2}
    & C_{\SL(4n,\R)}(2^{1+2}_+) = \GL(2n,\R),\quad & &
    C_{\PSL(4n,\R)}(H_+) = H_+
    \times \PGL(2n,\R) = V_+ \times (\pslr{2n} \rtimes C_2) \\
    & C_{\SL(4n,\R)}(2^{1+2}_-) = \GL(n,\Ha),\quad & &
    C_{\PSL(4n,\R)}(H_-) = H_-
    \times \PGL(n,\Ha) 
  \end{alignat*}
  where $H_+$ and $H_-$ are hyperbolic planes with quadratic functions
  $q_+(v_1,v_2)=v_1v_2$ and $q_-(v_1,v_2)=v_1^2+v_1v_2+v_2^2$
  (\ref{sec:xtraspecreps}), respectively. In the first case, for
  instance, the commutative diagram
  \begin{equation*}
    \xymatrix{
      1 \ar[r] & {\PGL(2n,\R)} \ar[r] & 
      C_{\pslr{4n}}(H_+) \ar[r] & H_+^{\vee} \ar[r] & 0 \\
               && H_+ \ar[u] \ar[ur]_{[\cdot,\cdot]}^{\cong} }
  \end{equation*}
  gives a central section of the \ses\ from \cite[5.11]{jmm:ndet}.
  %In the second case, it seems that
%  $C_2=\gen{e_1}$.  The Quillen auto\m\ groups are the orthogonal
%  groups $\A(\pslr{4n})(V_+)=\mathrm{O}^+(2,\F_2)$ and
%  $\A(\pslr{4n})(V_-)=\mathrm{O}^-(2,\F_2) = \GL(2,\F_2)$ of order $2$
%  and $6$, respectively, as we see from the commutative diagram
%  \begin{equation*}
%    \xymatrix{
%    {\Out(2^{1+2}_{\pm})} \ar@{=}[r] &
%    {\A(\GL(4,\C))(2^{1+2}_{\pm})} \ar@{^(->}[r] &
%    {\A(\SL(8,\R))(2^{1+2}_{\pm})}
%    \ar@{^(->}[r] 
%                            \ar@{->>}[d] &
%    {\Out(2^{1+2}_{\pm})} \ar@{^(->}[d] \\
%    &&   {\A(\PSL(8,\R))(V_{\pm})}   \ar@{^(->}[r] &
%    {\Aut(V_{\pm})} }
%  \end{equation*}
%  where the image of the right vertical injection is
%  \cite[13.9]{huppert:I} $\mathrm{O}^{\pm}(2,\F_2)$. We have from
%  \ref{exmp:IDQ} that the Quillen category auto\m\ group
%  $\A(\GL(4,\C))(2^{1+2}_{\pm})$ \cite[5.8]{jmm:ndet} is the full
%  outer auto\m\ group of the extraspecial group $2^{1+2}_{\pm}$ for
%  all auto\m s are trace preserving.

\setcounter{subsection}{\value{thm}}
\subsection{Rank three objects with nontrivial inner product}
\label{sec:non0innerprodrank3}\add

Let $V$ be a rank three object of $\A(\pslr{4n})$ with nontrivial
inner product. Then $V$ or $V^D$ is isomorphic to $H_+\#L(i,2n-i)$ ($1
\leq i \leq n$), $H_-\#L(i,n-i)$ ($1 \leq i \leq [n/2]$) or $V_0$.
$H_+\#L(i,2n-i) \subset \pslr{4n}$ is defined to be the quotient of
 \begin{multline*}
    \big\langle\diag(R,\ldots ,R), \diag(T,\ldots ,T), 
   \diag(\overbrace{-E,\ldots,-E}^i,\overbrace{E,\ldots
     ,E}^{2n-i})\big\rangle \subset \SL(4n,\R), \\
   R=
   \begin{pmatrix}
     1&0\\0&-1
   \end{pmatrix}, \qquad
   T=
   \begin{pmatrix}
     0&1\\1&0
   \end{pmatrix}
  \end{multline*}
  isomorphic to $2^{1+2}_+ \times C_2$ and $H_-\#L(i,n-i) \subset
  \pslr{4n}$ to be the quotient of
\begin{multline*}
    \Big\langle 
    \diag\left(
      \begin{pmatrix}
        0 & -R \\ R & 0
      \end{pmatrix}, \ldots, 
      \begin{pmatrix}
        0 & -R \\ R & 0
      \end{pmatrix} \right), 
     \diag\left(\begin{pmatrix}
        0 & -T \\ T & 0
      \end{pmatrix}, \ldots, 
      \begin{pmatrix}
        0 & -T \\ T & 0
      \end{pmatrix}\right), \\
    \diag\Big(
      \overbrace{
        \begin{pmatrix}
        -E & 0 \\ 0 & -E
        \end{pmatrix}, \ldots,
        \begin{pmatrix}
          -E & 0 \\ 0 &-E
        \end{pmatrix}}^i,\overbrace{
        \begin{pmatrix}
        E & 0 \\ 0 & E
        \end{pmatrix},\ldots
        \begin{pmatrix}
          E & 0 \\ 0 &E
        \end{pmatrix}}^{n-i}\Big)\Big\rangle 
      \subset \SL(4n,\R)
  \end{multline*}
  isomorphic to $2^{1+2}_- \times C_2$.  
  The \lmntwo\ $V_0 \subset \pslr{4n}$ is the quotient of 
  \begin{multline*}
    \Big\langle \diag\left(
      \begin{pmatrix}
        R&0\\0&R
      \end{pmatrix},\ldots,
      \begin{pmatrix}
        R&0\\0&R
      \end{pmatrix}\right), \diag\left(
      \begin{pmatrix}
        T&0\\0&T
      \end{pmatrix},\ldots,
      \begin{pmatrix}
        T&0\\0&T
      \end{pmatrix}\right), \\
      \diag\left(
        \begin{pmatrix}
          0&-E\\E&0
        \end{pmatrix},\ldots,
        \begin{pmatrix}
          0&-E\\E&0
        \end{pmatrix}\right)\Big\rangle
  \end{multline*}
  isomorphic to the generalized extraspecial $2$-group $C_4 \circ
  2^{1+2}_{\pm} \subset \SL(4n,\R)$ as described in
  \ref{exmp:IDQ}.(\ref{exmp:IDQ5}).

\setcounter{subsection}{\value{thm}}
\subsection{Rank four objects with nontrivial inner product}
\label{sec:non0innerprodrank4}\add
The following partial census of rank four objects with nontrivial
inner product suffices for our purposes.  Define the \lmntwo\ 
$H_+\#P(1,i-1,2n-i) \subset \pslr{4n}$, $2 \leq i \leq n$, to be the
quotient of
 \begin{multline*}
    \big\langle\diag(R,\ldots ,R), 
             \diag(T,\ldots ,T), \\ 
\diag(E,\overbrace{-E,\ldots,-E}^{i-1},\overbrace{E,\ldots ,E}^{2n-i}),
\diag(E,\overbrace{E,\ldots,E}^{i-1},\overbrace{-E,\ldots ,-E}^{2n-i})
    \big\rangle \subset \SL(4n,\R)
  \end{multline*}
Define $V_0\#L(i,n-i) \subset \pslr{4n}$, $1 \leq i \leq [n/2]$, to be
the quotient of 
 \begin{multline*} 
    \Big\langle 
    \diag\left(
      \begin{pmatrix}
        0 & -E \\ E & 0
      \end{pmatrix}, \ldots,
      \begin{pmatrix}
        0 & -E \\ E & 0
      \end{pmatrix}\right),
      \diag\left(\begin{pmatrix}
        R & 0 \\ 0 & R
      \end{pmatrix}, \ldots, 
      \begin{pmatrix}
        R & 0 \\ 0 & R
      \end{pmatrix} \right),\\ 
     \diag\left(\begin{pmatrix}
        T & 0 \\ 0 & T
      \end{pmatrix}, \ldots, 
      \begin{pmatrix}
        T & 0 \\ 0 & T
      \end{pmatrix}\right), \\
    \diag\Big(
      \overbrace{
        \begin{pmatrix}
        -E & 0 \\ 0 & -E
        \end{pmatrix}, \ldots,
        \begin{pmatrix}
          -E & 0 \\ 0 &-E
        \end{pmatrix}}^i,\overbrace{
        \begin{pmatrix}
        E & 0 \\ 0 & E
        \end{pmatrix},\ldots
        \begin{pmatrix}
          E & 0 \\ 0 &E
        \end{pmatrix}}^{n-i}\Big)\Big\rangle
      \subset \SL(4n,\R)
  \end{multline*}
isomorphic to $C_4 \circ 2^{1+2}_{\pm} \times C_2$.

\setcounter{subsection}{\value{thm}}
\subsection{Centers of centralizers}
\label{sec:ZC}\add

For the computations in \S\ref{sec:lim} we need to know the centers of
the centralizers for some of the low dimensional objects of
$\A(\pslr{4n})^{[\; ,\; ] \neq 0}$.

\begin{prop}\label{prop:ZC}
  Let $V \in \Ob(\A(\pslr{4n})^{[\; ,\; ] \neq 0})$ be one of the
  objects
\begin{itemize}
\item $H_+$, $H_-$,
\item $H_+\# L(i,2n-i)$ ($1 \leq i \leq n$), $H_-\# L(i,n-i)$ ($1 \leq
  i \leq [n/2]$), $V_0$, or
\item $H_+\# P(1,i-1,2n-i,0)$ ($1 < i \leq n$), $V_0 \# L(i,n-i)$ ($1
  \leq i \leq [n/2]$)
\end{itemize}
introduced in \ref{sec:non0innerprodrank2}--\ref{sec:non0innerprodrank4}.
Then $ZC_{\pslr{4n}}(V)=V$.
\end{prop}

\begin{proof}
The proof is a case-by-case checking.  

\noindent\underline{$H_+$ and $H_-$}
Since the centralizers of the rank two objects $H_+$ and $H_-$ are
$C_{\pslr{4n}}(H_+)=H_+ \times \PGL(2n,\R)$ and
$C_{\pslr{4n}}(H_-)=H_- \times \PGL(n,\Ha)$, Proposition~\ref{prop:ZC}
is immediate in these cases.

\noindent\underline{$H_+\# L(i,2n-i)$ ($1 \leq i \leq n$) and  $H_+\#
  P(1,i-1,2n-i,0)$ ($1 < i \leq n$)} 
%Let $L(i,2n-i) \subset
%\PGL(2n,\R)$ be the $1$-dimensonal object that is the image of the
%representation $\rho = i\psi_0+(2n-i)\psi_1 \colon L \to \GL(2n,\R)$,
%$1 \leq i \leq n$, and let $P(1,i-1,2n-i,0) \subset \PGL(2n,\R)$ be
%the $2$-dimensional object that is the image of the representation
%$\rho = \psi_0+(i-1)\psi_1+(2n-i)\psi_2 \colon P \to \GL(2n,\R)$, $1 <
%i \leq n$.
We shall only prove the $2$-dimensional case since the $3$-dimensional
case is similar. The centralizer of $H_+\#L(i,2n-i)$ is isomorphic to
the product of $H_+$ with the centralizer of $L=L(i,2n-i)$ in
$\PGL(2n,\R)$. There is \cite[5.11]{jmm:ndet} a \ses\ 
  \begin{equation*}
    1 \to \frac{\GL(i,\R) \times \GL(2n-i,\R)}{\gen{-E}} \to 
    C_{\PGL(2n,\R)}(L) \to \Hom(L,\gen{-E})_{\rho} \to 1
  \end{equation*}
  where the group to the right consists of all homo\m s
  \func{\phi}{L}{\gen{-E}} such that $\rho$ and $\phi\cdot\rho$ are
  conjugate representations in $\GL(2n,\R)$. By trace considerations,
  this group is trivial if $i<n$ and of order two if
  $i=n$. Hence 
  \begin{equation*}
     C_{\PGL(2n,\R)}(L)=
     \begin{cases}
       \frac{\GL(i,\R) \times \GL(2n-i,\R)}{\gen{-E}} & i < n \\
       \frac{\GL(n,\R)^ 2}{\gen{-E}} \rtimes \gen{C_1} & i=n
     \end{cases}
  \end{equation*}
  where $C_1=
  \begin{pmatrix}
    0 & E \\ E & 0
  \end{pmatrix}$
  is the $(2n \times 2n)$-matrix that interchanges the two
  $\GL(n,\R)$-factors. In case $i<n$, use \ref{lemma:ZGLprod}. In case
  $i=n$, the center is (\ref{semicenter}) the pull-back
  of the group homo\m s
  \begin{equation*}
    \frac{\GL(n,\R) \times \gen{(E,-E)}}{\gen{-E}} =
    \left( \frac{\GL(n,\R)^ 2}{\gen{-E}} \right)^{\gen{C_1}} \to
    \Aut\left(\frac{\GL(n,\R)^ 2}{\gen{-E}} \right) \leftarrow
    \gen{C_1}
  \end{equation*}
  which is $\frac{\GL(1,\R) \times \gen{(-E,E)}}{\gen{-E}} = L$
  again.

\noindent\underline{$V_0$ and $V_0\# L(i,n-i)$}
The object $V_0 \subset \pslr{4n}$ is the quotient of $G=4 \circ
2^{1+2}_{\pm} \subset \SL(4n,\R)$ as described in
\ref{exmp:IDQ}.(\ref{exmp:IDQ5}). 
%%
%generated (\ref{sec:xtraspecreps},
%\ref{sec:glnCtosl2nR}) by the matrices 
%\begin{equation*}
%  \diag\left( 
%    \begin{pmatrix}
%      R & 0 \\ 0 & R
%    \end{pmatrix}, \ldots,
%    \begin{pmatrix}
%      R & 0 \\ 0 & R
%    \end{pmatrix}\right), 
%   \diag\left(\begin{pmatrix}
%      T & 0 \\ 0 & T
%    \end{pmatrix}, \ldots,
%    \begin{pmatrix}
%      T & 0 \\ 0 & T
%    \end{pmatrix}\right), 
%    \diag\left(\begin{pmatrix}
%      0 & -E \\ E & 0
%    \end{pmatrix}, \ldots,
%    \begin{pmatrix}
%      0 & -E \\ E & 0
%    \end{pmatrix}\right)
%\end{equation*}
%where $R=
%\begin{pmatrix}
%  1 & 0 \\ 0 & 1
%\end{pmatrix}$, $T=
%\begin{pmatrix}
%  0 & 1 \\ 1 & 0
%\end{pmatrix}$, and $E=
%\begin{pmatrix}
%  1 & 0 \\ 0 & 1
%\end{pmatrix}$ is the $(2 \times 2)$-identity matrix. 
As this representation $\rho=n(\chi+\overline{\chi})$ is the $n$-fold
sum of an irreducible representation of complex type there are exact
sequences
\begin{equation*}
  \xymatrix{
    1 \ar[r] &
    {\frac{\GL(n,\C)}{\gen{-E}}} \ar[r] &
    C_{\pslr{4n}}(V_0) \ar[r] &
    {\Hom(G,\gen{-E})_{\rho}} \ar[r] &
    1 \\
    1 \ar[r] &
    Z(G)/G' \ar[r] \ar@{^(->}[u] &
    G/G' \ar[r] \ar@{^(->}[u] &
    G/Z(G)  \ar[r] \ar@{^(->}[u] &
    1 }
\end{equation*}
where the top row is \cite[5.11]{jmm:ndet}. The elementary abelian
group  $\Hom(G,\gen{-E})_{\rho}$, consisting of all homo\m s
\func{\phi}{G}{\gen{-E}} such that $\rho$ and $\phi\cdot\rho$ are
conjugate in $\SL(4n,\R)$, equals all of  $\Hom(G,\gen{-E})=2^3$ since
conjugation with the 
first two of the generators from \ref{sec:non0innerprodrank3} and with
\begin{equation*}
  C_2=\diag\left(
    \begin{pmatrix}
      E & 0 \\ 0 & -E
    \end{pmatrix}, \ldots,
    \begin{pmatrix}
      E & 0 \\ 0 & -E
    \end{pmatrix}\right) 
\end{equation*}
induce three independent generators. Hence
\begin{equation*}
  C_{\pslr{4n}}(V_0) = 
  \left( \frac{\GL(n,\C)}{\gen{-E}} \times V_0/V_0^{\perp} \right)
  \rtimes \gen{C_2}
\end{equation*}
Note that conjugation with the matrix $C_2$ induces complex
conjugation on $\GL(n,\C)$.  The center of this semi-direct product is
(\ref{semicenter}) the pull-back of the group homo\m s
\begin{equation*}
  \frac{\GL(n,\R) \circ \gen{i}}{\gen{-E}} \times V_0/V_0^{\perp} =
  \left( \frac{\GL(n,\C)}{\gen{-E}} \times V_0/V_0^{\perp}
  \right)^{\gen{C_2}} \to
  \Aut \left( \frac{\GL(n,\C)}{\gen{-E}} \times V_0/V_0^{\perp}
  \right) \leftarrow
  \gen{C_2}
\end{equation*}
which is $\frac{\GL(1,\R) \circ \gen{i}}{\gen{-E}} \times
V_0/V_0^{\perp} = \frac{\gen{i}}{\gen{-E}} \times V_0/V_0^{\perp} =
V_0$. 

The case of $V_0 \# L(i,n-i)$, $1 \leq i < [n/2]$, is quite
similar. The centralizer is 
\begin{equation*}
  C_{\pslr{4n}}(V_0 \# L(i,n-i)) = 
  \left( \frac{\GL(i,\C) \times \GL(n-i,\C)}{\gen{-E}} \times
    V_0/V_0^{\perp} \right) \rtimes \gen{C_2}
\end{equation*} 
and its center is the pull-back of the homo\m s
\begin{multline*}
  \frac{(\GL(i,\R) \times \GL(n-i,\R)) \circ \gen{i}}{\gen{-E}} \times
  V_0/V_0^{\perp} =
  \left( \frac{\GL(i,\C) \times \GL(n-i,\C)}{\gen{-E}} \times  
    V_0/V_0^{\perp} \right)^{\gen{C_2}} \\  \to
  \Aut \left( \frac{\GL(i,\C) \times \GL(n-i,\C)}{\gen{-E}} \times  
    V_0/V_0^{\perp} \right) \leftarrow \gen{C_2}
\end{multline*} 
which is $ZC_{\pslr{4n}}(V_0 \# L(i,n-i)) = 
\frac{(\GL(1,\R) \times \GL(1,\R)) \circ \gen{i}}{\gen{-E}} \times
V_0/V_0^{\perp} = 2^2 \times V_0/V_0^{\perp} = V_0 \times L$. 

If $n$ is even and $i=n/2$, there is a \ses\
\begin{equation*}
  1 \to \frac{\GL(n,\C)^2}{\gen{-E}} \to C_{\pslr{4n}}(V_0 \times L) \to
  \Hom(G \times L, \gen{-E})_{\rho} \to 1
\end{equation*} 
where the elementary abelian group to the right is all of $\Hom(G
\times L, \gen{-E}) = 2^4$. Hence the centralizer
\begin{equation*}
   C_{\pslr{4n}}(V_0 \times L) = 
   \left( \frac{\GL(n,\C)^2}{\gen{-E}} \times V_0/V_0^{\perp} \right)
   \rtimes \gen{C_1,C_2}
\end{equation*}
where $C_2$ is as above and $C_1$ is the $(4n \times 4n)$-matrix $
\begin{pmatrix}
  0 & E \\ E & 0
\end{pmatrix}$. The matrix $C_2$ commutes with $ V_0/V_0^{\perp}$ and
acts as complex conjugation on $\frac{\GL(n,\C)^2}{\gen{-E}}$. The
matrix $C_1$ commutes with $ V_0/V_0^{\perp}$ and switches the two
factors of $\GL(n,\C)^2$. The center of the centralizer is the
pull-back of the group homo\m s
\begin{multline*}
  \frac{\GL(n,\R) \circ \gen{i} \times \gen{(E,-E)}}{\gen{-E}}  
  \times V_0/V_0^{\perp} =
  \left( \frac{\GL(n,\C)^2}{\gen{-E}} \times V_0/V_0^{\perp}
  \right)^{\gen{C_1,C_2}} \\ \to
  \Aut \left( \frac{\GL(n,\C)^2}{\gen{-E}} \times V_0/V_0^{\perp}
 \right) \leftarrow \gen{C_1,C_2} 
\end{multline*}
which is $ZC_{\pslr{4n}}(V_0 \times L) = 
\frac{\GL(1,\R) \circ \gen{i} \times \gen{(E,-E)}}{\gen{-E}} \times
V_0/V_0^{\perp} = 2^2 \times V_0/V_0^{\perp} = V_0 \times L$.

\noindent\underline{$H_-\#L(i,n-i)$}
%Consider the representation $\rho=\psi \# (i\psi_0 + (n-i)\psi_1)$ of
%$2^{1+2}_- \times L$, $1 \leq i \leq [n/2]$. 
As above we have that
\begin{equation*}
  C_{\pslr{4n}}(H_- \times L) = 
  \begin{cases}
  \frac{\GL(i,\Ha) \times \GL(n-i,\Ha)}{\gen{-E}} \times H_- & i <
  [n/2] \\
  \frac{\GL(i,\Ha)^2}{\gen{-E}}\rtimes \gen{C_1} \times H_- 
  & \text{$n$ even and $i=n/2$}  
  \end{cases}
\end{equation*}
with center $ZC_{\pslr{4n}}(H_- \times L) = \frac{\GL(1,\R) \times
  \GL(1,\R)}{\gen{-E}} = 2 \times H_- = H_- \times L$ in case $i \neq
n-i$.  If $n$ is even and $i=n/2$, then the center is the pull-back of
the group homo\m s
\begin{equation*}
  \frac{\GL(i,\Ha) \times \gen{(-E,E)}}{\gen{-E}} \times H_- =
     \frac{\GL(i,\Ha)^2}{\gen{-E}} \times H_- \to
  \Aut \left(  \frac{\GL(i,\Ha)^2}{\gen{-E}} \times H_- \right)
  \leftarrow \gen{C_1} 
\end{equation*}
which is $Z C_{\pslr{4n}}(H_- \times L) = 
\frac{\GL(1,\R) \times \gen{(-E,E)}}{\gen{-E}} \times H_- = 2 \times
  H_- = H_- \times L$.
\end{proof}

\chapter{The $\mathrm{B}$-family}
\label{sec:slodd}

The $B$-family consists of the matrix groups 
\begin{equation*}
  \SL(2n+1,\R), \quad n \geq 2,   
\end{equation*}
of real $(2n+1) \times (2n+1)$ matrices of determinant $+1$. When
$n=1$ we obtain the \twocg\ $\SL(3,\R)=\PGL(2,\C)$ considered in
Chapter~\ref{cha:afam}.  The embedding
\begin{equation*}
  \GL(2n,\R) \to \SL(2n+1,\R) \colon A \to 
  \begin{pmatrix}
    A & 0 \\
    0 & \det A
  \end{pmatrix}
\end{equation*}
permits us to consider $\GL(2n,\R)$ as a maximal rank subgroup of
$\SL(2n+1,\R)$.  The \mtn\ for the subgroup $\GL(2n,\R)$ is also the
\mtn\ for $\SL(2n+1,\R)$, $N(\SL(2n+1,\R))=N(\GL(2n,\R)$
(\ref{eq:dfamTNW}), so that, in particular, the Weyl group
$W(\SL(2n+1,\R)) = W(\GL(2n,\R)) = \Sigma_2 \wr \Sigma_n$
(\ref{eq:dfamW}). 

%The maximal torus normalizer, and
%the Weyl group of $\SL(2n+1,\R)$ and $\GL(2n,\R)$ agree,
%\begin{equation*}
%  T=\SL(2,\R)^n, \quad
%  N=\GL(2,\R) \wr \Sigma_n, \quad
%  W=\Sigma_2 \wr \Sigma_n
%\end{equation*}
%The \mtn\ for $\SL(2n+1,\R)$ sits in the extension
%\begin{equation*}
%  1 \to \SL(2,\R)^n \to \GL(2,\R) \wr \Sigma_n \to
%  \Sigma_2 \wr \Sigma_n \to 1
%\end{equation*}
%which splits for all $n$ because the \twoctg\ $\GL(2,\R)$ is the
%semidirect product of $\SL(2,\R)$ and its Weyl group, $\Sigma_2$.  

It is known that \cite[1.6]{matthey:normalizers} \cite[Main
Theorem]{hms:first}
\begin{equation*}
  H^0(W;\ch{T})=\Z/2, \quad H^1(W;\ch{T})=
  \begin{cases}
    \Z/2 & n=2 \\
    \Z/2 \times \Z/2 & n > 2
  \end{cases}
\end{equation*}
for these groups.

The full general linear group $\GL(2n+1,\R) = \SL(2n+1,\R) \times
\gen{-E}$ is the direct product of $\SL(2n+1,\R)$ with the opposite of
the identity matrix so that $\PGL(2n+1,\R)=\SL(2n+1,\R)$.

\section{The structure of $\SL(2n+1,\R)$}
\label{sec:structslodd}

Consider the \lmntwo s
\begin{align*}
  &\Delta_{2n+1} = 
  \gen{\diag(\pm 1, \ldots , \pm 1)} \subset \GL(2n+1,\R) \\
  &S\Delta_{2n+1} = \SL(2n+1,\R) \cap \Delta_{2n+1} \subset \SL(2n+1,\R)
  \\
  &t=t(\SL(2n+1,\R))=\Delta_{2n+1} \cap T(\SL(2n+1,\R)) =
  \gen{e_1,\ldots ,e_n} \subset T(\SL(2n+1,\R))
\end{align*}
in $\GL(2n+1,\R)$ and $\SL(2n+1,\R)$. 
%%%%%%%%%
%Note that $S\Delta_{2n+1} =
%\Delta_{2n} \subset \GL(2n,\R) \subset \SL(2n+1,\R)$ and
%$\Delta_{2n+1} = S\Delta_{2n} \times \gen {-E}$.  
%%%%%%%%%%%%
%In particular, the matrices $e_j$, $1 \leq j \leq n$,
%$I$, and $c_j$,  $1 \leq j \leq n$, defined in \S\ref{sec:structure}
%as elements of $\GL(2n,\R)$ are also elements of $\SL(2n+1,\R)$. 
%%%%%%%%%%%%%%%%%%%%%%%%%
%Since $\GL(2n+1,\R)$ is the direct product of $\SL(2n+1,\R)$ with the
%central subgroup $\gen{-E}$, 
%\begin{equation*}
%  \A(\SL(2n+1,\R))(V) = \A(\GL(2n+1,\R))(V) = \A(\GL(2n+1,\R))(V
%  \times \gen{-E})
%\end{equation*}
%for all \lmntwo s $V$ in  $\SL(2n+1,\R)$. In particular
%\begin{equation*}
%   \A(\SL(2n+1,\R))(\Delta_{2n}) =  \A(\GL(2n+1,\R))(\Delta_{2n+1}) =
%   \Sigma_{2n+1} 
%\end{equation*}
%by representation theory \cite[Lemma 3, Proposition
%4]{bob:stubborn}. 
%%%%%%%%%%%%%%%%%%%%%%%%%%%%%%%%%%%%%%%%%%%

\begin{lemma}\label{lemma:sloddASL}
  The inclusion functors
  \begin{align*}
    &\A(\Sigma_{2n+1},\Delta_{2n+1}) \to \A(\GL(2n+1,\R)), &
    &\A(\Sigma_{2n+1},S\Delta_{2n+1}) \to \A(\SL(2n+1,\R)), \\
    & {} & &\A(\Sigma_2\wr\Sigma_{n},t) \to \A(\SL(2n+1,\R))^{\leq t}
  \end{align*}
  are equivalences of categories.
\end{lemma}
\begin{proof}
  Similar to \ref{cor:equivcat}.
  $\A(\SL(2n+1,\R)$ is a full subcategory of $\A(\GL(2n+1,\R)$ since
  conjugation with the central element $-E$ of negative determinant is
  the identity.
\end{proof}

The Quillen categories $\A(\GL(2n,\R))=\A(\Sigma_{2n},\Delta_{2n})$
$\A(\SL(2n+1,\R))=\A(\Sigma_{2n+1},\Delta_{2n+1})$
(\ref{cor:equivcat}, \ref{lemma:sloddASL}) are not equivalent.

%Any rank one object of $\A(\GL(2n+1,\R))$ is isomorphic to the image
%of \func{i_0\rho_0+i_1\rho_1}{\Z/2}{\Delta_{2n+1}} where
%$i_0+i_1=2n+1$, $i_1 >0$, and
%$\Hom(\Z/2,\R^{\times})=\{\rho_0,\rho_1\}$ with $\rho_0$ trivial. The
%image is an object of the full subcategory $\A(\SL(2n+1,\R))$ if and
%only if $i_1>0$ is even. (See \S\ref{sec:structure} in
%Chapter~\ref{sec:dfam} for the meaning of $\rho_i$.)

%Any rank two object of $\A(\GL(2n+1,\R))$ is isomorphic to the image
%of the homo\m\ 
%\func{i_0\rho_0+i_1\rho_1+i_2\rho_2+i_3\rho_3}{\Z/2}{\Delta_{2n+1}}
%where $i_0+i_1+i_2+i_3=2n+1$ and at least two of $i_1,i_2,i_3$ are
%positive. The image is an object of the full subcategory
%$\A(\SL(2n+1,\R))$ if and only if $i_1+i_3$, $i_2+i_3$, and $i_1+i_2$
%are even, ie if and only if $i_1,i_2,i_3$ have the same parity (the
%opposite parity if $i_0$). The image is toral if $i_0$ is odd and
%nontoral if $i_0$ is even \cite[5.2]{griess:elem}.  The action of
%$\Aut(\Z/2 \oplus \Z/2) \cong \Sigma_3$ on the set
%$\{\rho_0,\rho_1,\rho_2,\rho_3\}$ fixes $\rho_0$ and permutes the
%other three elements.

For any partition $i=(i_0,i_1)$, $i_0 \geq 0$, $i_1 >0$, of $2n+1$, let
$L[i_0,i_1] \subset \Delta_{2n+1}$ be the subgroup generated by
\begin{equation*}
  \diag(\overbrace{+1, \ldots, +1}^{i_0},\overbrace{-1, \ldots
    ,-1}^{i_1}) = (i_0\rho_0+i_1\rho_1)(e_1)
\end{equation*}
For any partition $(i_0,i_1,i_2,i_3)$ of $2n+1$ where at least two
of $i_1,i_2,i_3$ are positive, let $P[i_0,i_1,i_2,i_3] \subset
\Delta_{2n+1}$ be the subgroup generated by 
\begin{align*}
 &\diag(\overbrace{+1, \ldots, +1}^{i_0},
       \overbrace{-1, \ldots ,-1}^{i_1},
       \overbrace{+1, \ldots ,+1}^{i_2},
       \overbrace{-1, \ldots ,-1}^{i_3}) =
       (i_0\rho_0+i_1\rho_1+i_2\rho_2+i_3\rho_3)(e_1)  \\
  &\diag(\overbrace{+1, \ldots, +1}^{i_0},
        \overbrace{+1, \ldots ,+1}^{i_1},
        \overbrace{-1, \ldots ,-1}^{i_2},
        \overbrace{-1, \ldots ,-1}^{i_3}) =
       (i_0\rho_0+i_1\rho_1+i_2\rho_2+i_3\rho_3)(e_2)
% &\diag(\overbrace{+1, \ldots, +1}^{i_0},
%       \overbrace{+1, \ldots ,+1}^{i_1},
%       \overbrace{-1, \ldots ,-1}^{i_2},
%       \overbrace{-1, \ldots ,-1}^{i_3}) =
%       (i_0\rho_0+i_1\rho_1+i_2\rho_2+i_3\rho_3)(e_2)
\end{align*}
Note that $L[i_0,i_1]$ is a subgroup of $S\Delta_{2n+1}$ if and only
if $i_1$ is even, and that $P[i_0,i_1,i_2,i_3]$ is a subbgroup of
$S\Delta_{2n+1}$ if and only of $i_1,i_2,i_3$ have the same parity,
the opposite parity of $i_0$.

Let $P(k,r)$ denote the number of partitions of $k=i_0+\cdots+i_{r-1}$
into sums of $r$ positive integers $1\leq i_0 \leq \cdots \leq
i_{r-1}$. From the above discussion we conclude 

\begin{prop}\label{prop:sloddrankleq2}
The category $\A(\SL(2n+1,\R))$  contains precisely
\begin{itemize}
\item $n$ iso\m\ classes of rank one objects represented by the lines
  $L[2i_0+1,2i_1]$ where $0 \leq i_0 \leq n-1$ and $i_1=n-i_0$. 
\item $\sum_{j=2}^n P(j,2) + \sum_{j=3}^n P(j,3)$ iso\m\ classes of
  toral rank two objects. They are represented by the subgroups
  $P[2i_0+1,2i_1,2i_2,0]$, where $0 \leq i_0 \leq n-2$ and $(i_1,i_2)$
  is a partition of $n-i_0$, together with the subgroups
  $P[2i_0+1,2i_1,2i_2,2i_3]$, where $0 \leq i_0 \leq n-3$ and
  $(i_1,i_2,i_3)$ is a partition of $n-i_0$.
\item $\sum_{j=3}^{n+2} P(j,3)$ iso\m\ classes of nontoral
  rank two objects represented by the subgroups
  $P[2i_0,2i_1-1,2i_2-1,2i_3-1]$ where $0 \leq i_0 \leq n-1$ and
  $(i_1,i_2,i_3)$ is a partition of $n-i_0+2$.
\end{itemize}
\end{prop}

The centralizers of these objects are
\set\begin{gather}
  \begin{split}\label{eq:Lcent}
     C_{\SL(2n+1,\R)}L[2i_0+1,2i_1]& = \SL(2n+1,\R) \cap 
        \big(\GL(2i_0+1,\R) \times \GL(2i_1,\R)\big) \\
     &= \SL(2i_0+1,\R) \times  \GL(2i_1,\R) 
  \end{split}\\
  \begin{split}\label{eq:Pcent}  
  &C_{\SL(2n+1,\R)}P[i] = \SL(2n+1,\R) \cap
    \prod_{j \in i} \GL(j,\R)\\
   &= 
   \begin{cases}
    \SL(2i_0+1,\R) \times \GL(2i_1,\R) \times \GL(2i_2,\R) \times 
       \GL(2i_3,\R) & \text{$P[i]$ toral}\\ %i=[2i_0,2i_1-1,2i_2-1,2i_3-1] \\
    \GL(2i_0,\R) \times \GL(2i_1-1,\R) \times \GL(2i_2-1,\R) \times 
       \SL(2i_3-1,\R) & \text{$P[i]$ nontoral}%i=[2i_0-1,2i_1,2i_2,2i_3]   
   \end{cases}
\end{split}
\end{gather}\add\add
as, for instance,
\begin{multline*}
  \SL(2n+1,\R) \cap \big(\GL(2i_0+1,\R) \times \GL(2i_1,\R)\big) \\ =
  \SL(2n+1,\R) \cap \big(\SL(2i_0+1,\R) \times \gen{-E} 
           \times \SL(2i_1,\R) \rtimes \gen{D}\big) \\ =
   \SL(2i_0+1,\R) 
           \times \SL(2i_1,\R) \rtimes \gen{-D} =
    \SL(2i_0+1,\R) 
           \times \GL(2i_1,\R), 
\end{multline*}
and the centers of the centralizers are
\set\begin{align}
   \label{eq:ZLcent}
   &ZC_{\SL(2n+1,\R)}L[2i_0+1,2i_1]= L[2i_0+1,2i_1], \\
   \label{eq:CPcent}
   &ZC_{\SL(2n+1,\R)}P[i]= \SL(2n+1,\R) \cap \prod_{i_j>0} Z\GL(i_j,\R)
   = 
   \begin{cases}
     P[i] & \#\{j \mid i_j>0\}=3 \\
     P[i] \times \Z/2 & \#\{j \mid i_j>0\}=4 
   \end{cases}
\end{align}\add\add
%%%%%%%%%%%%%%%%%%%%%%%%%%
%Thus there are $P(n+2,3)+\sum_{k=2}^nP(k,2)$ planes $P$ with
%$ZC_{\SL(2n+1,\R)}(P)=P$ and $P(n+1,3)+2\sum_{k=3}^nP(k,3)$ planes $P$
%with $ZC_{\SL(2n+1,\R)}(P)=P \times \Z/2$.
%%%%%%%%%%%%%%%%%%%%%
 
%%%%%%%%%%%%%%%%%%%%%%%%%
%\begin{align*}
%  C_{\SL(2n+1,\R)}L[2n+1-2i,2i] 
%  &= \SL(2n+1,\R) \cap \big(\GL(2n+1-2i,\R) \times \GL(2i,\R)\big) \\
%  &= \SL(2n+1,\R) \cap \big(\SL(2n+1-2i,\R) \times \gen{-E} \times
%  \SL(2i,\R) \rtimes \gen{D}\big) \\
%  &=\big(\SL(2n+1-2i,\R) \times \SL(2i,\R)\big) \rtimes \gen{-E\cdot D_1} \\
%  &= \SL(2n+1-2i,\R) \times \GL(2i,\R)  
%\end{align*}
%so that $Z_{\SL(2n+1,\R)}L[2n+1-2i,2i] = L[2n+1-2i,2i]$. Similarly, 
%\begin{multline*}
%      C_{\SL(2n+1,\R)}P[i] = \\
%  \begin{cases}
%    \SL(2i_0-1,\R) \times \GL(2i_1,\R) \times \GL(2i_2,\R) &
%    i=[2i_0-1,2i_1,2i_2,0] \\
%    \SL(2i_0-1,\R) \times \GL(2i_1,\R) \times \GL(2i_2,\R) \times
%    \GL(2i_3,\R) &
%    i=[2i_0-1,2i_1,2i_2,2i_3] \\
%    \GL(2i_1-1,\R) \times \GL(2i_2-1,\R) \times
%    \SL(2i_3-1,\R) &
%    i=[0,2i_1-1,2i_2-1,2i_3-1] \\
%    \GL(2i_0,\R) \times \GL(2i_1-1,\R) \times \GL(2i_2-1,\R) \times
%    \SL(2i_3-1,\R) &
%    i=[2i_0,2i_1-1,2i_2-1,2i_3-1] \\
%  \end{cases}
%\end{multline*}
%so that
%\begin{equation*}
%  ZC_{\SL(2n+1,\R)}P[i] = 
%  \begin{cases}
%    P[i] & i=[2i_0-1,2i_1,2i_2,0] \\
%    P[i] \times \Z/2 & i=[2i_0-1,2i_1,2i_2,2i_3] \\
%    P[i]  & i=[0,2i_1-1,2i_2-1,2i_3-1] \\
%    P[i] \times \Z/2 & i=[2i_0,2i_1-1,2i_2-1,2i_3-1] \\
%  \end{cases}
%\end{equation*}
%%%%%%%%%%%%%%%%%%%%%%%%%%%%%%%%%%%%%%%%%%%%%%%%%%%%%%%

\begin{lemma}\label{lemma:H0Vslodd}
  For any nontrivial subgroup $V \subset S\Delta_{2n+1}$ there is a
  natural iso\m\
  \begin{equation*}
    ZC_{\SL(2n+1,\R)}(V)=H^0(\Sigma_{2n+1}(V);S\Delta_{2n+1})
  \end{equation*}
  where $\Sigma_{2n+1}(V)$ is the point-wise stabilizer subgroup
  (\ref{defn:AWt}). 
\end{lemma}
\begin{proof}
  Let $V \subset S\Delta_{2n+1}$ be any nontrivial subgroup of rank
  $r$.  Then $V=V[i]$ is the image of $\sum_{\rho \in
    V^{\vee}}i_{\rho}\rho$ for some function
  \func{i}{\Hom((\Z/2)^r,\R^{\times})}{\Z} where $\sum_{\rho \in
    V^{\vee}}i_{\rho}=2n+1$ and
\begin{multline*}
  ZC_{\SL(2n+1,\R)}V[i]=Z \big(\SL(2n+1,\R) \cap
  \prod_{i_{\rho}>0}\GL(i_{\rho},\R) \big) \\ =
  \SL(2n+1,\R) \cap
  \prod_{i_{\rho}>0}Z\GL(i_{\rho},\R) 
  = S\Delta_{2n+1} \cap
    \Delta_{2n+1}^{\prod\Sigma_{i_{\rho}}} 
  = S\Delta_{2n+1}^{\Sigma_{2n+1}(V[i])} 
\end{multline*}
where the second equality can be proved by using that
$C_{\GL(i,\R)}\SL(i,\R)=Z\GL(i,\R)$ and the final equality follows
from the observation that the stabilizer subgroup
$\Sigma_{2n+1}(V[i])=\prod_{i_{\rho>0}}\Sigma_{i_{\rho}}$.
\end{proof}

\begin{cor}\label{cor:limislodd}
  $\lim^i(\A(\SL(2n+1,\R),\pi_1(BZC_{\SL(2n+1,\R)}))=0$ for all $i>0$.
\end{cor}
\begin{proof}
  Immediate from the general exactness theorem
  (\ref{dw:limits}) for functors of the form as in
  \ref{lemma:H0Vslodd}.
  \end{proof}

  \begin{prop}
    Centralizers of objects of $\A(\SL(2n+1,\R))^{\leq t}_{\leq 2}$
    are LHS.
  \end{prop}
  \begin{proof}
    Let $X_1$ and $X_2$ be connected Lie groups and $\pi_1$ and
    $\pi_2$ finite $2$-groups acting on them. Suppose that the homo\m
    s $\theta(X_1)^{\pi_1}$ and $\theta(X_1)^{\pi_1}$
    (\ref{eq:theta}) are surjective. Then also $\theta(X_1 \times
    X_2)^{\pi_1 \times \pi_2}$ is surjective and so the product $X_1
    \rtimes \pi_1 \times X_2 \rtimes \pi_2$ is LHS
    (\ref{lhscrit1}). This observation applies to the
    products (\ref{eq:Lcent}, \ref{eq:Pcent}) since the
    $\theta$-homo\m s are surjective \cite[5.4]{hms:first}
    (\ref{exmp:sl2C}) for $\SL(2i+1,\R)$, $i \geq 0$, and $\SL(2i,\R)$,
    $i \geq 1$.
  \end{proof}

%Note that centralizers of objects of $\A(\SL(2n+1,\R))^{\leq t}_{\leq
%  2}$ are LHS simply because their Weyl groups split as the direct
%product $W=W_0 \times \pi$ of the Weyl group of the identity component
%$W_0$ with the component group $\pi$.
%\marginpar{\bf Split??}

\setcounter{subsection}{\value{thm}}
\section{The limit of the functor
  $H^1(W;\ch{T})/H^1(\pi_0;\check{Z}( \; )_0)$ on
  $\A(\pslr{2n+1})^{\leq t}_{\leq 2}$}
\label{sec:sloddlim0}\add

In this subsection we check, using a modification of
\ref{sec:cond2}, that conditions (1) and (2) of
\ref{indstepalt} with $X=\SL(2n+1,\R)$ are satisfied
under the inductive assumptions that the connected \twocg s
$\SL(2i+1,\R)$, $0 \leq i < n$, and $\SL(2i,\R)$, $1 \leq i \leq n$,
are uniquely $N$-determined.

The objects $V \subset \SL(2n+1,\R)$ of the category
$\A(\pslr{2n+1})^{\leq t}_{\leq 2}$ are the rank one objects
$L[i_0,i_1]$ and the rank two objects $P[2i_0+1,2i_1,2i_2,0]$ and
$P[2i_0+1,2i_1,2i_2,2i_3]$ as described in \ref{prop:sloddrankleq2}.
The rank two object $P[2i_0+1,2i_1,2i_2,2i_3]$, $i_3 \geq 0$, contains
the three lines $L[2i_0+2i_1+1,2i_2+2i_3]$,
$L[2i_0+2i_2+1,2i_1+2i_3]$, and $L[2i_0+2i_3+1,2i_1+2i_2]$.  Their
centralizers are described in (\ref{eq:Lcent}) and (\ref{eq:Pcent}).
Note that there are functorial iso\m s \set\begin{equation}
  \label{eq:sloddsplit}
  \check{T}^{W_0(C_{\SL(2n+1,\R)}(V))} = 
   (\Z/2)^{\min\{i_0,1\}} \times \check{Z}(C_{\SL(2n+1,\R)}(V)_0) 
\end{equation}\add
as modules over $\pi_0C_{\SL(2n+1,\R)}(V)$.

Condition (1) of \ref{indstepalt} is satisfied as
$C_{X}(V)$ has $N$-determined auto\m s and is $N$-determined for
general reasons (\ref{prodauto}, \ref{lemma:autononcon},
\ref{redtoconnected}). This means that there are iso\m
s, $\alpha_V$ and $f_V$, such that the diagrams
\begin{equation*}
  \xymatrix{
    C_N(V) \ar[r]^{\alpha_V}_{\cong} \ar[d] & C_N(V) \ar[d] \\
    C_X(V) \ar[r]_{f_V}^{\cong} & C_{X'}(V) }
\end{equation*}
commute and $\alpha_V \in H^1(W;\check{T})(C_X(V))$. There may be more
than choice for $\alpha_V$ but for each $\alpha_V$ there is just one
possibility for $f_V$ (\ref{defcons}). The set of
possible $\alpha_V$ for a given $V$ is a $H^1(\pi_0;\check{Z}((\;
)_0))(C_X(V))$-coset in $H^1(W;\ch{T})(C_X(V))$
(\ref{idcomp}). The collection of the $\alpha_V$ for
various $V$ represents an element of the inverse limit
\set\begin{equation}
  \label{eq:sloddlim0}
  \lim^0 \big( \A(\SL(2n+1,\R))^{\leq t}_{\leq 2}, 
               \frac{H^1(W;\ch{T})}{H^1(\pi_0;\check{Z}((\;)_0))} \big)
\end{equation}\add
of the quotient functor over the category $\A(\SL(2n+1,\R))^{\leq
  t}_{\leq 2}$. Condition (2) of \ref{indstepalt} is
satisfied if the restriction map from $H^1(W;\ch{T})(\SL(2n+1,\R))$ to
(\ref{eq:sloddlim0}) is surjective. Because of the natural splitting
(\ref{eq:sloddsplit}) and because the centralizers
$C_{\SL(2n+1,\R)}(V)$ are LHS there is a \ses\ 
\begin{equation*}
  0 \to \Hom(\pi_0,(\Z/2)^{\min\{i_0,1\}}) \to
  \frac{H^1(W;\ch{T})}{H^1(\pi_0;\check{Z}((\;)_0))} \to
  H^1(W_0;\ch{T})^{\pi_0} \to 0
\end{equation*}
of functors on $\A(\SL(2n+1,\R))^{\leq t}_{\leq 2}$. If we apply the
functor $\Hom(\pi_0,(\Z/2)^{\min\{i_0,1\}})$ to the \m s
\set\begin{equation}
  \label{eq:sloddmorphisms}
  L[2i_0+1,2i_1+2i_2] \to P[2i_0+1,2i_1,2i_2,0] \leftarrow  
  L[2i_0+2i_1+1,2i_2]
\end{equation}\add
we see that the induced \m s are injective and that their images
intersect trivially. Thus the inverse limit of this functor is trivial
and from the above \ses\ we obtain an injective map
\begin{equation*}
  \lim^0\big(\A(\SL(2n+1,\R))^{\leq t}_{\leq 2}, 
               \frac{H^1(W;\ch{T})}{H^1(\pi_0;\check{Z}((\;)_0))}
               \big) \to
  \lim^0\big(\A(\SL(2n+1,\R))^{\leq t}_{\leq 2},
  H^1(W_0;\ch{T})^{\pi_0} \big)
\end{equation*}
between the inverse limits. As the inverse
limit to the right is a subgroup of the inverse limit of the functor
$H^1(W_0;\ch{T})$ we conclude that if the restriction map
\set\begin{equation}\label{eq:sloddres} H^1(W_0;\ch{T})(\SL(2n+1,\R))
  \to \lim^0\big(\A(\SL(2n+1,\R))^{\leq t}_{\leq 2}, H^1(W_0;\ch{T})
  \big)
\end{equation}\add
is surjective, then condition (2) of \ref{indstepalt} is
satisfied. 

\begin{lemma}\label{lemma:sloddlimH1W0}
  The restriction homo\m\  (\ref{eq:sloddres}) is an iso\m\ for all $n
  \geq 2$.
\end{lemma}
\begin{proof}
  %%linux:~/manus/dfam/magma/slodd/toral/toral01.prg
  %%linux:~/manus/dfam/magma/slodd/toral/toral06.prg and slodd.log
  %% and sl9.prg especially for n=4
  For $n=2$, the image under the functor $H^1(W_0;\ch{T})$ of the
  category $L[1,4] \to P[1,2,2,0] \leftarrow L[3,2]$ is $0 \to 0
  \leftarrow \Z/2$ so that the limit of the functor $H^1(W_0;\ch{T})$
  is $\Z/2$. Since $\SL(3,\R) \times \SL(2;\R) \to \SL(5,\R)$ turns
  out to induce an iso\m\ on $H^1(W_0;\ch{T})$ the claim follows in
  this case.
  
  For $n=3$, taking into account only the planes of type
  $P[2i_0-1,2i_1,2i_2,0]$, we should compute the limit of the diagram
  \begin{equation*}
    \xymatrix{
     H^1(W_0C_{\SL(7,\R)}L[1,6]) \ar[dr] \\
     {} & H^1(W_0C_{\SL(7,\R)}P[1,4,2,0]) \\
     H^1(W_0C_{\SL(7,\R)}L[3,4]) \ar[ur] \ar[dr] \\
     {} & H^1(W_0C_{\SL(7,\R)}P[3,2,2,0]) \\
     H^1(W_0C_{\SL(7,\R)}L[5,2]) \ar[uuur] \ar[ur] }  
  \end{equation*}
  of $\F_2$-vector spaces. For each of the planes $P$ take the
  intersections of the images in the cohomology groups
  $H^1(W_0C_{\SL(7,\R)}P;\ch{T})$ of $H^1(W_0C_{\SL(7,\R)}L;\ch{T})$
  for each line $L \subset P$.  Take the intersection of the
  pre-images in each $H^1(W_0C_{\SL(7,\R)}L;\ch{T})$ of these
  subspaces of $H^1(W_0C_{\SL(7,\R)}P;\ch{T})$. Using the computer
  program {\em magma\/} one may see that these subspaces have
  dimensions $1,2,2$ for $L=L[1,6],L[3,4],L[5,2]$, respectively, and
  that they equal the image of the restriction maps from
  $H^1(W_0;\ch{T})(\SL(7,\R))$. This shows that the lemma is true in
  this case.
  
  In general, the above mentioned subspaces of
  $H^1(W_0C_{\SL(7,\R)}L;\ch{T})$ have dimension $1$ for $L=L[1,2n]$
  and dimension $2$ for the lines $L=L[2i+1,2n-2i]$ with $1 \leq i
  \leq n-1$ and these subspaces equal the image of the restriction
  maps from $H^1(W_0;\ch{T})(\SL(2n+1,\R))$.
\end{proof}

\section{Rank two nontoral objects of $\A(\SL(2n+1,\R))$}
\label{sec:sloddnontoral2}\add

The nontoral rank two objects of $\A(\SL(2n+1,\R))$ are represented by
the subgroups $P[i] \subset S\Delta_{2n+1}$ 
generated by the elements
\begin{align*}
  e_1&=
  \diag(\overbrace{+1,\ldots,+1}^{2i_0},\overbrace{-1,\ldots,-1}^{2i_1-1}, 
        \overbrace{+1,\ldots,+1}^{2i_2-1},\overbrace{-1,\ldots,-1}^{2i_3-1}) \\
  e_2&=
  \diag(\overbrace{+1,\ldots,+1}^{2i_0},\overbrace{+1,\ldots,+1}^{2i_1-1},
        \overbrace{-1,\ldots,-1}^{2i_2-1},\overbrace{-1,\ldots,-1}^{2i_3-1})
\end{align*}
where $i=(2i_0,2i_1-1,2i_2-1,2i_3-1)$, $0 \leq i_0 \leq n-1$ and
$(i_1,i_2,i_3)$ is a partition of $n+2-i_0$
(\ref{prop:sloddrankleq2}). The generators of $P[i]$ may also be
written as 
\set\begin{align}
  e_1&=\diag(\overbrace{E,\ldots,E}^{i_0-1},E,
  \overbrace{-E,\ldots,-E}^{i_1-1},-R, \overbrace{E,\ldots,E}^{i_2-1},
  \overbrace{-E,\ldots,-E}^{i_3-1},-1), \qquad R=
       \begin{pmatrix}
         1&0\\0&-1 
       \end{pmatrix} \label{eq:slodde1} \\
  e_2&=\diag(\overbrace{E,\ldots,E}^{i_0-1},E,
             \overbrace{E,\ldots,E}^{i_1-1},R,
             \overbrace{-E,\ldots,-E}^{i_2-1},
             \overbrace{-E,\ldots,-E}^{i_3-1},-1) 
       \label{eq:slodde2}
\end{align}\add\add
The centralizer of $P[i]$ is 
\begin{align*}
  C_{\SL(2n+1,\R)}P[i]&=\SL(2n+1,\R) \cap 
  \big(\GL(2i_0,\R) \times \GL(2i_1-1,\R) \times \GL(2i_2-1,\R) \times
  \GL(2i_3-1,\R) \big) \\
  &= \GL(2i_0,\R) \times \GL(2i_1-1,\R) \times \GL(2i_2-1,\R) \times
  \SL(2i_3-1,\R)
\end{align*}
Note that $P[i]$ is contained in the \mtn\ $N(\SL(2n+1,\R))=\GL(2,\R)
\wr \Sigma_n$. Since the centralizer of $P[i]$ in the \mtn,
\begin{multline*}
  C_{\GL(2,\R) \wr \Sigma_n}P[i] = 
  \GL(2,\R) \wr \Sigma_{i_0} \times  \GL(2,\R) \wr \Sigma_{i_1-1}
  \times \GL(1,\R) \\ \times \GL(1,\R) \times
  \GL(2,\R) \wr \Sigma_{i_2-1} \times \GL(2,\R) \wr \Sigma_{i_3-1},
\end{multline*}
is the \mtn\ for the centralizer of $P[i]$, the lift $P[i] \subset
N(\SL(2n+1,\R))$ is a \pl\ of $P[i] \subset \SL(2n+1,\R)$
\cite{jmm:normax}.  The two other \pl s are given by composing with
the permutation matrices for the permutations $(1,2)(i_0+i_1,2n+1)$
and $(1,2)(i_0+i_1+1,2n+1)$ (assuming $i_0>0$) resulting in the lifts
given by
\begin{align*}
  e_1&=\diag(\overbrace{E,\ldots,E}^{i_0-1},E,
             \overbrace{-E,\ldots,-E}^{i_1-1},-E,
             \overbrace{E,\ldots,E}^{i_2-1},
             \overbrace{-E,\ldots,-E}^{i_3-1},-1),
       \\
  e_2&=\diag(\overbrace{E,\ldots,E}^{i_0-1},E,
             \overbrace{E,\ldots,E}^{i_1-1},R,
             \overbrace{-E,\ldots,-E}^{i_2-1},
             \overbrace{-E,\ldots,-E}^{i_3-1},-1) \\
\intertext{and}
  e_1&=\diag(\overbrace{E,\ldots,E}^{i_0-1},E,
             \overbrace{-E,\ldots,-E}^{i_1-1},R,
             \overbrace{E,\ldots,E}^{i_2-1},
             \overbrace{-E,\ldots,-E}^{i_3-1},-1),
       \\
  e_2&=\diag(\overbrace{E,\ldots,E}^{i_0-1},E,
             \overbrace{E,\ldots,E}^{i_1-1},-E,
             \overbrace{-E,\ldots,-E}^{i_2-1},
             \overbrace{-E,\ldots,-E}^{i_3-1},-1) 
\end{align*}
respectively. These two lifts are also \pl s of $P[i] \subset
\SL(2n+1,\R)$. The three \pl s are not conjugate in $N(\SL(2n+1,\R))$
because the intersection with the \mt\ is generated by $e_1+e_2$ in
the first case and by $e_1$, respectively $e_2$, in the next two
cases. Note that all three \pl s have the same \mt, $\SL(2,\R)^{i_0}
\times \SL(2,\R)^{i_1-1} \times \SL(2,\R)^{i_2-1} \times
\SL(2,\R)^{i_3-1}$. 

Let $U=\gen{e_1,e_2,e_3}$ be \lmntwo\ generated by $e_1$ and $e_2$
as in (\ref{eq:slodde1}, \ref{eq:slodde2}) together with
\begin{equation*}
  e_3=\diag(\overbrace{E,\ldots,E}^{i_0-1},R,
             \overbrace{E,\ldots,E}^{i_1-1},E,
             \overbrace{E,\ldots,E}^{i_2-1},
             \overbrace{E,\ldots,E}^{i_3-1},-1),
%       \\
%  e_4&=\diag(\overbrace{-E,\ldots,-E}^{i_0-1},-E,
%             \overbrace{E,\ldots,E}^{i_1-1},E,
%             \overbrace{E,\ldots,E}^{i_2-1},
%             \overbrace{E,\ldots,E}^{i_3-1},+1) 
\end{equation*}
Note that the centralizer of $U$ has a nontrivial identity component,
that the inclusion $U \subset C_{\SL(2n+1,\R)}P[i]$
induces an iso\m\ on $\pi_0$.

Under the inductive assumption that $\SL(2i,\R)$, $1 \leq i \leq n-1$,
and $\SL(2i-1,\R)$, $1 \leq i \leq n$, have $\pi_*(N)$-determined
auto\m s (or using \cite{jmo:selfho}) we conclude from
\ref{lemma:uniquenu} and
\ref{dia:Utriangle} and (part of) \cite[5.2]{jmm:ext}
that condition (3) of \ref{indstepalt} is satisfied for
$\SL(2n+1,\R)$. (Namely,
\ref{lemma:uniquenu}.(\ref{lemma:uniquenu1}) says
that $\nu_L'$ does not depend on the choice of $L<V$. The difference
$f_{\nu,L_2}^{-1} \circ f_{\nu,L_1}$ between any two of the maps
$f_{\nu,L}$ from \ref{indstepalt}.(\ref{indstepalt4})
is an auto\m of
$C_{\SL(2n+1,\R)}P[i]$ that, by
\ref{lemma:uniquenu}.(\ref{lemma:uniquenu2}), is the
identity on the identity component and by the commutative diagram
(\ref{dia:Utriangle}) 
\set\begin{equation}\label{dia:pi0} \xymatrix{
    & U \ar[dl] \ar[dr] \\
    C_{\SL(2n+1,\R)}P[i] \ar[rr]_{f_{\nu,L_2}^{-1} \circ f_{\nu,L_1}}
    && C_{\SL(2n+1,\R)}P[i]}
\end{equation}\add
also the identity on $\pi_0 C_{\SL(2n+1,\R)}P[i]$. Any such auto\m\ of
$C_{\SL(2n+1,\R)}P[i]$ has \cite[5.2]{jmm:ext} the form $A \to
\varphi(A)A$ where
\begin{multline*}
  \varphi \colon \GL(2i_0,\R) \times \GL(2i_1-1,\R) \times
  \GL(2i_2-1,\R) \times \SL(2i_3-1,\R) \to \\
   \pi_0\big(\GL(2i_0,\R) \times \GL(2i_1-1,\R) \times
  \GL(2i_2-1,\R) \times \SL(2i_3-1,\R) \big) \to Z\GL(2i_0,\R)
\end{multline*}
is some homo\m s. Diagram (\ref{dia:pi0}) thus implies that the
inclusion $U \to \SL(2n+1,\R)$ and the mono\m\ given by $e_i \to
\varphi(e_i)e_i$, $1 \leq i \leq 3$, are conjugate. Since the trace of
$e_i$, $1 \leq i \leq 3$, is odd (nonzero), $\varphi$ must be trivial.
Thus $f_{\nu,L_1}$ and $f_{\nu,L_2}$ are identical iso\m s.)

\chapter{The $\mathrm{C}$-family}
\label{sec:spn}

Let $\Ha=\{a+bj \vert a,b\in\C\}$, where $j^2=-1$ and
$ja=\overline{a}j$ for $a\in\C$, be the quaternion algebra.  The
$C$-family consists of the matrix groups
\begin{equation*}
  \PGL(n,\Ha)=\GL(n,\Ha)/\gen{-E}, \quad n \geq 3,   
\end{equation*}
of quaternion projective $n \times n$ matrices.  (These \twocg s also
exist for $n=1$ or $n=2$. However, $\PGL(1,\Ha)= \SL(3,\R)=\PGL(2,\C)$
and $\PGL(2,\Ha)= \SL(5,\R)$
(\ref{eq:lowdegree}) are already covered.)

The \mtn\ for $\GL(1,\Ha)=\Ha^{\times}$, generated by the \mt\ 
$\GL(1,\C)=\C^{\times}$ and the element $j$, sits in the non-split
extension
\begin{equation*}
  1 \to \GL(1,\C) \to N(\GL(1,\Ha)) \to \gen{j}/\gen{-1} \to 1
\end{equation*}
of $\Sigma_2$ by $\GL(1,\C)=\C^{\times}$.
The \mtn\ for $\GL(n,\Ha)$ is the subgroup
\begin{equation*}
  N(\GL(n,\Ha))= N(\GL(1,\Ha)) \wr \Sigma_n,
\end{equation*}
generated by $N(\GL(1,\Ha))^n \subset \GL(n,\Ha)$ and the permutation
matrices. The \mtn\ for $\PGL(n,\Ha)$, the quotient
$N(\GL(n,\Ha))$ by the order two group $\gen{-E}$, sits in the
extension 
\begin{equation*}
  1 \to \frac{\GL(1,\C)^n}{\gen{-E}} \to 
  \frac{N(\GL(1,\Ha))^n}{\gen{-E}} \to
  \frac{N(\GL(1,\Ha))}{\GL(1,\C)} \wr  \Sigma_n
  \to 1
\end{equation*}
which does not split (for $n \geq 3$).

It is known that \cite[1.6]{matthey:normalizers} \cite[Main
Theorem]{hms:first}
\begin{equation*}
  H^0(W;\ch{T})(\PGL(n,\Ha))=0, \quad H^1(W;\ch{T})(\PGL(n,\Ha))=
  \begin{cases}
    \Z/2 & n=3,4 \\
    0 & n > 4
  \end{cases}
\end{equation*}
for the projective groups.

\setcounter{subsection}{\value{thm}}
\section{The structure of $\PGL(n,\Ha)$}
\label{sec:NpglnH}\add
Let
\begin{equation*}
  \Delta_n=t(\GL(n,\Ha))=
  \gen{\diag(\pm 1, \ldots, \pm 1)} \subset \GL(n,\Ha)
\end{equation*}
be the maximal \lmntwo\ in $\GL(n,\Ha)$ and
$C_4=\gen{I} \subset \GL(n,\Ha)$ the cyclic order four
group generated by $I =\diag(i,\ldots,i)$. The maximal \lmntwo\ in
$\PGL(n,\Ha)$ is the quotient
\begin{equation*}
  t(\PGL(n,\Ha))=\frac{t(\PGL(n,\Ha))^*}{\gen{-E}}, \quad
  t(\PGL(n,\Ha))^*=C_4 \circ t(\GL(n,\Ha))
\end{equation*}
so that the toral part of the Quillen category is equivalent
\begin{equation*}
  \A(\PGL(n,\Ha))^{\leq t}=\A(C_2 \wr \Sigma_n, \frac{C_4 \circ \gen{\diag(\pm
      1,\ldots, \pm 1)}}{\gen{-E}})
\end{equation*}
to the category whose objects are nontrivial subgroups of
$t(\PGL(n,\Ha))$ and whose \m s are induced from the action of the
Weyl group. 

For any partition $i=(i_0,i_1)$ of $n=i_0+i_1$ into a sum of two
positive integers $i_0\geq i_1 \geq 1 >0$ let $L[i]=L[i_0,i_1] \subset
\GL(n,\Ha)$ be the subgroup generated by
\begin{equation*}
  \diag(\overbrace{+1, \ldots, +1}^{i_0},\overbrace{-1, \ldots
    ,-1}^{i_1})
\end{equation*}
Then the centralizer
\set\begin{equation}
  \label{eq:cfamCL}
  C_{\PGL(n,\Ha)}L[i_0,i_1]=
  \begin{cases}
    \frac{\GL(i_0,\Ha) \times \GL(i_1,\Ha)}{\gen{-E}} & i_0 \neq i_1
    \\
    \frac{\GL(i_0,\Ha)^2}{\gen{-E}}\rtimes\Big\langle 
      \begin{pmatrix}
        0&E\\E&0
      \end{pmatrix}\gen{-E}\Big\rangle & i_0=i_1
  \end{cases} 
\end{equation}\add
so that the center $ZC_{\PGL(n,\Ha)}L[i_0,i_1]=L[i_0,i_1]$ as in the
proof of \ref{prop:ZC} and \ref{lemma:ZGLprod}.

Let (also) $I \in \PGL(n,\Ha)$ denote the order two element that
is the image of the order four element $i \in \GL(n,\Ha)$. Then
\set\begin{equation}
  \label{eq:cfamCI}
   C_{\PGL(n,\Ha)}(I)=\frac{\GL(n,\C)}{\gen{-E}} \rtimes \gen{j\gen{-E}}
\end{equation}\add
so that the center $ZC_{\PGL(n,\Ha)}(I)=\gen{I}$ as shown in the proof
of \ref{prop:ZC}.

For any partition $(i_0,i_1,i_2,0)$ of $n=i_0+i_1+i_2$ into a sum of
three positive integers $i_0\geq i_1 \geq i_2 >0$ or any partition
$(i_0,i_1,i_2,i_3)$ of $n=i_0+i_1+i_2+i_3$ into a sum of four positive
integers $i_0\geq i_1 \geq i_2 \geq i_3 >0$ let $P[i_0,i_1,i_2,i_3]
\subset \Delta_{2n+1}$ be the subgroup generated by the two elements
\begin{align*}
 &\diag(\overbrace{+1, \ldots, +1}^{i_0},
       \overbrace{-1, \ldots ,-1}^{i_1},
       \overbrace{+1, \ldots ,+1}^{i_2},
       \overbrace{-1, \ldots ,-1}^{i_3})  \\
  &\diag(\overbrace{+1, \ldots, +1}^{i_0},
        \overbrace{+1, \ldots ,+1}^{i_1},
        \overbrace{-1, \ldots ,-1}^{i_2},
        \overbrace{-1, \ldots ,-1}^{i_3})
\end{align*}
Then the centralizer
\set\begin{equation}
  \label{eq:cfamCP}
  C_{\PGL(n,\Ha)}P[i]=
  \begin{cases}
    \frac{\GL(i_0,\Ha)^4}{\gen{-E}} \rtimes (C_2 \times C_2) 
     & i=(i_0,i_0,i_0,i_0) \\
    \frac{\GL(i_0,\Ha)^2 \times \GL(i_2,\Ha)^2}{\gen{-E}} \rtimes C_2
    &  i=(i_0,i_0,i_2,i_2) \\
    \frac{\GL(i_0,\Ha)\times\GL(i_1,\Ha)\times\GL(i_2,\Ha)\times\GL(i_3,\Ha)}
    {\gen{-E}} & \#i=4
  \end{cases}
\end{equation}\add
where the groups $C_2$ are generated by permutation matrices.

For any partition $i=(i_0,i_1)$ of $n=i_0+i_1$ into a sum of two
positive integers $i_0\geq i_1 >0$ let $I\#L[i_0,i_1] \subset
\PGL(n,\Ha)$ be the \lmntwo\ that is the quotient of
\begin{equation*}
  (I\#L[i_0,i_1])^*=\big\langle I,
  \diag(\overbrace{+1,\ldots,
    +1}^{i_0},\overbrace{-1,\ldots,-1}^{i_1})\big\rangle 
\end{equation*}
Then the centralizer
\begin{equation}
  \label{eq:cfamCIL}
  C_{\PGL(n,\Ha)}I\#L[i_0,i_1]=
  \begin{cases}
    \frac{\GL(i_0,\C) \times \GL(i_1,\C)}{\gen{-E}} \rtimes
    \gen{j\gen{-E}} & i_0 \neq i_1 \\
    \frac{\GL(i_0,\C)^2}{\gen{-E}}\rtimes \Big\langle j\gen{-E}, 
      \begin{pmatrix}
        0&E\\E&0
      \end{pmatrix}\gen{-E}\Big\rangle & i_0=i_1
  \end{cases}
\end{equation}\add
 
\begin{prop}
  The category $\A(\PGL(n,\Ha))$ contains exactly 
  \begin{itemize}
  \item $[n/2]+1$ rank one
  toral objects represented by the lines $L[i,n-i]$, $1\leq i \leq
  [n/2]$ (with $q=0$), and by the line $I$ (with $q\neq 0$).
\item $P(n,3)+P(n,4)+[n/2]$ rank two toral objects represented by the
  $P(n,3)$ planes $P[i_0,i_1,i_2,0]$ (with $q=0$), the $P(n,4)$ planes
  $P[i_0,i_1,i_2,i_3]$ (with $q=0$), and the $[n/2]$ planes
  $I\#L[i,n-i]$, $1 \leq i \leq [n/2]$ (with $q\neq 0$).
  \end{itemize}
\end{prop}
%% ~/manus/dfam/magma/cfam/toral.prg

\begin{prop}
  Let $V \subset \PGL(n,\Ha)$ be a nontrivial \lmntwo . Then
  \begin{equation*}
    \textmd{$V$ is toral} \iff \textmd{$[V,V] \neq 0$}
  \end{equation*}
\end{prop}\add
\begin{proof}
  The proof is similar to \ref{cor:toralnontoral} with the extra input
  that all \lmntwo s in $\GL(n,\Ha)$ are toral by quaternion
  representation theory \cite{adams:lie}.
\end{proof}

%Since, according to quaternion representation theory
%\cite[3.57]{adams:lie}, all \lmntwo s in $\GL(n,\Ha)$ are toral it
%follows that $\A(\GL(n,\Ha)) \cong \A(\Sigma_n,(\Z/2)^n)$ and that
%$\A(\PGL(n,\Ha))^{q=0} \cong \A(\Sigma_n,P((\Z/2)^n))$ where
%$P((\Z/2)^n)$ is the quotient of $(\Z/2)^n$ by its diagonal.  It
%contains $P(n,2)=\left[\frac{n}{2}\right]$ lines $L(i_0,i_1)$,
%$i_0+i_1=n$, and $P(n,3)+P(n,4)$ planes $P(i_0,i_1,i_2,i_3)$,
%$i_0+i_1+i_2+i_3=n$, or $P(i_0,i_1,i_2,i_3,i_4)$,
%$i_0+i_1+i_2+i_3+i_4=n$.

%\setcounter{subsection}{\value{thm}}
%\subsection{Centralizers of objects of $\A(\PGL(n,\Ha)^{\leq t}_{\leq
%    2}$ are LHS}
%\label{sec:cfamlhs}
%\add

\begin{prop}\label{prop:cfamlimH1W}
   Centralizers of objects of $\A(\GL(n,\Ha))^{\leq t}_{\leq 2}$
    are LHS.
\end{prop}
\begin{proof}
  The centralizers $C=C_0 \rtimes \pi$ in question are the
  nonconnected centralizers listed in (\ref{eq:cfamCL}),
  (\ref{eq:cfamCI}), (\ref{eq:cfamCP}), and (\ref{eq:cfamCIL}). In
  fact, we only need to deal with
  \begin{equation*}
    \frac{\GL(i,\Ha)^2}{\gen{-E}}\rtimes C_2, \quad
    \frac{\GL(i_0,\Ha)^2 \times \GL(i_1,\Ha)^2}{\gen{-E}}\rtimes C_2,
    \quad 
    \frac{\GL(i,\Ha)^4}{\gen{-E}}\rtimes (C_2 \times C_2)
  \end{equation*}
  as the other cases are covered by \ref{lemma:lhs}. It suffices
  (\ref{lhscrit1})  to show that $\theta(C_0)^{\pi}$
  (\ref{eq:theta}) is surjective.
  
  %% /home/moller/manus/dfam/magma/cfam/toral/lhs/GLiH2.prg
  Computations with the program {\em magma\/} results in the table
  \begin{center}
    \begin{tabular}[c]{|c||c|c|c|c|c|} \hline
      $\frac{\GL(i,\Ha)^2}{\gen{-E}}\rtimes C_2$ & $\ker\theta$ &
      $\Hom(W,\ch{T}^W)$ & $H^1(W;\ch{T})$ & $\theta$ &
      $H^1(W;\ch{T})^{\pi}$ \\ \hline \hline
      $1=i$ & $\dm{2}$ & $\dm{2}$ & $0$ & epi & 0 \\ \hline 
      $2=i$ & $\dm{2}$ & $\dm{4}$ & $\dm{3}$ & \  & $\dm{2}$ \\ \hline
      $2<i$ & $0$ & $\dm{4}$ & $\dm{4}$ & iso  & $\dm{2}$ \\ \hline 
    \end{tabular}
  \end{center}
  From the table we see that $\theta^{\pi}$ is surjective unless
  $i=2$. In that exceptional case, more compute computations show that
  $H^1(\pi;\ch{T}^W)=\Z/2$ and $H^1(W \rtimes C_2;\ch{T})=\dm{3}$
  which means that also  $\frac{\GL(2,\Ha)^2}{\gen{-E}}\rtimes C_2$ is
  LHS.

  %% /home/moller/manus/dfam/magma/cfam/toral/lhs/GLi0H2GLi1H2.prg
  Computations with the program {\em magma\/} results in the table
  \begin{center}
    \begin{tabular}[c]{|c||c|c|c|c|c|} \hline
      $\frac{\GL(i_0,\Ha)^2 \times \GL(i_1,\Ha)^2}{\gen{-E}}$ &  $\ker\theta$ &
      $\Hom(W,\ch{T}^W)$ & $H^1(W;\ch{T})$ & $\theta$ &
      $H^1(W;\ch{T})^{\pi}$ \\ \hline \hline 
      $1=i_0, 2=i_1$ & $\dm{4}$ & $\dm{18}$ & $\dm{14}$ & epi &
      $\dm{7}$ \\ \hline
      $1=i_0, 2<i_1$ & $\dm{2}$ & $\dm{18}$ & $\dm{16}$ & epi &
      $\dm{8}$ \\ \hline
      $2=i_0<i_1$ & $\dm{2}$ & $\dm{24}$ & $\dm{22}$ & epi &
      $\dm{11}$ \\ \hline
      $3<i_0<i_1$ & $0$ & $\dm{24}$ & $\dm{24}$ & iso &
      $\dm{12}$ \\ \hline
    \end{tabular}
  \end{center}
  Since $\theta$ is surjective and  $H^{>0}(\pi; \ker\theta)=0$
  because the action of $\pi$ on $\ker\theta$ is induced from the
  trivial subgroup, $\theta^{\pi}$ is surjective.

  %% /home/moller/manus/dfam/magma/cfam/toral/lhs/GLi0H2GLi1H2.prg
  Computations with the program {\em magma\/} results in the table
  \begin{center}
    \begin{tabular}[c]{|c||c|c|c|c|c|} \hline
      $\frac{\GL(i,\Ha)^4}{\gen{-E}}\rtimes (C_2 \times C_2)$ & 
      $\ker\theta$ & $\Hom(W,\ch{T}^W)$ & $H^1(W;\ch{T})$ & $\theta$ &
      $H^1(W;\ch{T})^{\pi}$ \\ \hline \hline
      $1=i$ & $\dm{4}$ & $\dm{12}$ & $\dm{8}$ & epi & $\dm{2}$ \\ \hline
      $2=i$ & $\dm{4}$ & $\dm{24}$ & $\dm{20}$ & epi & $\dm{5}$ \\
      \hline
      $2<i$ & $0$ & $\dm{24}$ & $\dm{24}$ & iso & $\dm{6}$ \\
      \hline
    \end{tabular}
  \end{center}
  Since $\theta$ is surjective and $H^{>0}(\pi; \ker\theta)=0$ because
  the action of $\pi$ on $\ker\theta$ is induced from the trivial
  subgroup, $\theta^{\pi}$ is surjective.
\end{proof}

\setcounter{subsection}{\value{thm}}
\section{The limit of the functor $H^1(W_0;\ch{T})^{W/W_0}$ on
  $\A(\PGL(n,\Ha)^{\leq t}_{\leq 2}$}
\label{sec:cfamlim0}\add

Let \func{H^1(W_0;\ch{T})}{\A(\PGL(n,\Ha))^{\leq t}_{\leq t}}{\Ab} be
the functor that takes the toral \lmntwo\ $V \subset t(\PGL(n,\Ha))$ to
the abelian group $H^1(W_0(C_{\PGL(n,\Ha)}(V);\ch{T}))$, and 
$H^1(W_0;\ch{T})^{W/W_0}$ the functor that takes $V$ to the the
  invariants for the action of the component group
  $\pi_0C_{\PGL(n,\Ha)}(V)$ on this first cohomology group.

%% /home/moller/manus/dfam/magma/cfam/toral/toral.prg
  \begin{prop}\label{prop:cfamlim0H1}
    The restriction map
    \begin{equation*}
      H^1(W(\PGL(n,\Ha);\ch{T}) \to \lim^0(\A(\PGL(n,\Ha))^{\leq t}_{\leq
        2},H^1(W_0;\ch{T})^{W/W_0}) 
    \end{equation*}
    is an iso\m\ for all $n > 3$.
  \end{prop}
  \begin{proof}
    %\noindent
%    \underline{$\PGL(3,\Ha)$}: Is there more than one map $L[1,2] \to
%    P[1,1,1,0]=P[1,0,1,1]$? 
%    \begin{equation*}
%      \xymatrix@C=50pt@R=45pt{
%         {(\Z/2)^2 \cong H^1(W_0;\ch{T})^{W/W_0}(L[1,2])} 
%         \ar[r]^-{
%           \begin{pmatrix}
%             1&0&0\\0&1&1
%           \end{pmatrix}}
%          \ar[dr]^(.6){
%            \begin{pmatrix}
%              1\\0
%            \end{pmatrix}} &
%         { H^1(W_0;\ch{T})^{W/W_0}(P[1,1,1,0]) \cong (\Z/2)^3} \\
%         {(\Z/2) \cong H^1(W_0;\ch{T})^{W/W_0}(I)} \ar[r]^(.45){
%           \begin{pmatrix}
%             1
%           \end{pmatrix}} &
%          {H^1(W_0;\ch{T})^{W/W_0}(I\#L[1,2]) \cong (\Z/2)} } 
%    \end{equation*}

    \noindent
    \underline{$\PGL(4,\Ha)$}: Computer computations show that the
    intersection of the
    images of the \m s
    \begin{equation*}
      H^1(W_0;\ch{T})^{W/W_0}(L[1,3]) \rightarrow
      H^1(W_0;\ch{T})^{W/W_0}(I\#L[1,3]) \xleftarrow{\cong}
       H^1(W_0;\ch{T})^{W/W_0}(I) 
    \end{equation*}
    is $1$-dimensional and that its pre-image in $
    H^1(W_0;\ch{T})^{W/W_0}(I)$ equals the image of the restriction
    map from $H^1(W,\ch{T})(\PGL(4,\Ha))$. Similarly, the images of
    the mono\m s
    \begin{equation*}
        H^1(W_0;\ch{T})^{W/W_0}(L[1,3]) \hookrightarrow
      H^1(W_0;\ch{T})^{W/W_0}(P[1,1,2,0]) \hookleftarrow
       H^1(W_0;\ch{T})^{W/W_0}(L[2,2])
    \end{equation*}
    meet in a $1$-dimensional subspace whose inverse images in the
    cohomology groups to the right and to the left agree with the
    images of the restriction maps from $H^1(W,\ch{T})(\PGL(4,\Ha))$.

    \noindent
    \underline{$\PGL(n,\Ha)$, $n>4$}:
    Computer computations show that the images of the \m s
    \begin{equation*}
      H^1(W_0;\ch{T})^{W/W_0}(L[1,n-1]) \rightarrow
      H^1(W_0;\ch{T})^{W/W_0}(I\#L[1,n-1]) \xleftarrow{\cong}
       H^1(W_0;\ch{T})^{W/W_0}(I)
    \end{equation*}
    intersect trivially and that the arrow pointing left is an iso\m
    . Similarly, the images of the injective \m s
    \begin{multline*}
       H^1(W_0;\ch{T})^{W/W_0}(L[i,n-i]) \hookrightarrow
      H^1(W_0;\ch{T})^{W/W_0}(P[i,1,n-i-1,0]) \\ \hookleftarrow
       H^1(W_0;\ch{T})^{W/W_0}(L[i+1,n-i-1]), \qquad 1 \leq i < [n/2],
    \end{multline*}
    intersect trivially. These observations imply that $
    \lim^0(\A(\PGL(n,\Ha))^{\leq t}_{\leq
      2},H^1(W_0;\ch{T})^{W/W_0})=0$.
  \end{proof}

\section{The category $\A(\PGL(n,\Ha))^{[\; , \;]\neq 0}_{\leq
    4}$}
\label{sec:cfamnontoral}

We shall need information about all nontoral objects of
$\A(\PGL(n,\Ha))$ of rank $\leq 3$ and some objects of rank $4$. If $V
\subset \PGL(n,\Ha)$ is an \lmntwo\ with nontrivial inner product then
its preimage $V^* \subset \GL(n,\Ha)$ is $P \times R(V)$ or $(C_4
\circ P) \times R(V)$ where $P$ is an extraspecial $2$-group, $C_4
\circ P$ a generalized extraspecial $2$-group, and
$\mho_1(V^*)=\gen{-E}$ (\ref{lemma:Vast}).  We manufacture all oriented
quaternion representations of these product groups as direct sums of
tensor products of irreducible representations of the factors
(\ref{sec:tensor}) as described in \cite[3.7, 3.65]{adams:lie}.

Note that the degrees of the faithful
irreducible representations over $\Ha$ for the groups $2^{1+2}_+$ and $C_4
\circ 2^{1+2}_{\pm}$ are even and that the quaternion group
$2^{1+2}_-$ has a faithful irreducible representation over $\Ha$,
namely the defining representation. 

\setcounter{subsection}{\value{thm}}
\subsection{The category $\A(\PGL(2n+1,\Ha))^{[\; , \;]\neq 0}_{\leq 4}$}
\label{sec:cfamoddnont}\add

The category $\A(\PGL(2n+1,\Ha))$
contains up to iso\m\ just one nontoral rank two object, $H_-$, whose
inverse image in $\GL(2n+1,\Ha)$ is
\begin{equation*}
  Q_8=2^{1+2}_-=\gen{\diag(i,\ldots,i),\diag(j,\ldots,j)} 
\end{equation*}
As in \ref{sec:non0innerprodrank2}, the centralizers \cite[Proposition
4]{bob:stubborn} of $2^{1+2}_-$ and $H_-$ are
\begin{equation*}
   C_{\GL(2n+1,\Ha)}(2^{1+2}_-)=\GL(2n+1,\R), \qquad 
  C_{\PGL(2n+1,\Ha)}(H_-)=H_- \times \SL(2n+1,\R) 
\end{equation*}
so that $ZC_{\PGL(2n+1,\Ha)}(H_-)=H_-$.

There are $n$ nontoral objects of rank three, $H_-\#L[i,2n+1-i]$,  $1
\leq i \leq n$. The inverse image in $\GL(2n+1,\Ha)$ of
$H_-\#L[i,2n+1-i]$ is
\begin{equation*}
  \big\langle\diag(i,\ldots,i),\diag(j,\ldots,j),
\diag(\overbrace{+1,\ldots,+1}^i,\overbrace{-1,\ldots,-1}^{2n+1-i})\big\rangle 
\end{equation*}
and the center of the centralizer,
 $C_{\PGL(2n+1,\Ha)}(H_-\#L[i,2n+1-i])=H_- \times
C_{\SL(2n+1,\R)}L[i,n-1]$, is
$ZC_{\PGL(2n+1,\Ha)}(H_-\#L[i,2n+1-i])=H_-\#L[i,2n+1-i]$ 
according to (\ref{eq:ZLcent}).

The objects $H_-\#P[i_0,i_1,i_2,i_3]$, where $P[i_0,i_1,i_2,i_3]$ is
as in \ref{sec:NpglnH}, are rank four nontoral objects.

We need to know that the nontoral object $H_-$ satisfies condition
($3$) of \ref{indstepalt}.  Note that the conditions of
\ref{lemma:uniquenu} are satisfied because the identity
component of $C_{\PGL(2n+1,\Ha)}(H_-)$ is nontrivial and because the
Quillen auto\m\ group $\A(\PGL(2n+1,\Ha))(H_-)=\GL(2,\F_2)$ acts
transitively on the set \pl s $H_- \subset N(\PGL(2n+1,\Ha))$ of $H_-
\subset \PGL(2n+1,\Ha)$.  Under the inductive assumption that
$\SL(2n+1,\R)$ has $\pi_*(N)$-determined auto\m s (or using
\cite{jmo:selfho}) we conclude from \ref{indstepalt}
and diagram (\ref{dia:Utriangle}) and (part of)
\cite[5.2]{jmm:ext} that condition (3) of
\ref{indstepalt} is satisfied for the nontoral rank $2$
object $H_-$.  (Namely,
\ref{lemma:uniquenu}.(\ref{lemma:uniquenu1}) says that
$\nu_L'$ does not depend on the choice of $L<V$. The difference
$f_{\nu,L_2}^{-1} \circ f_{\nu,L_1}$ between any two of the maps
$f_{\nu,L}$ from \ref{indstepalt}.(\ref{indstepalt4})
is an auto\m\ of $C_{\PGL(2n+1,\Ha)}(H_-)$ that, by
\ref{lemma:uniquenu}.(\ref{lemma:uniquenu2}), is the
identity on the identity component and by the commutative diagram
(\ref{dia:Utriangle}) 
\set\begin{equation}\label{dia:cfampi0H-}
  \xymatrix{
    & H_- \ar[dl] \ar[dr] \\
    C_{\PGL(2n+1,\Ha)}(H_-) \ar[rr]_{f_{\nu,L_2}^{-1} \circ
      f_{\nu,L_1}} && C_{\PGL(2n+1,\Ha)}(H_-)}
\end{equation}\add
also the identity on $\pi_0C_{\PGL(2n+1,\Ha)}(H_-)$. Since the
identity component $\SL(2n+1,\R)$ of the centralizer
$C_{\PGL(2n+1,\Ha)}(H_-)$ has no center, this shows that
$f_{\nu,L_2}^{-1} \circ f_{\nu,L_1}$ is the identity auto\m .)

\setcounter{subsection}{\value{thm}}
\subsection{Rank two nontoral objects of $\A(\PGL(2n,\Ha))$}
\label{sec:cfamnontrank2}\add

The category $\A(\PGL(2n,\Ha))$ contains up to iso\m\ two nontoral
rank two objects, $H_+$ and $H_-$, whose inverse images in
$\GL(2n,\Ha)$ are
\begin{align*}
  2^{1+2}_+&=\gen{\diag(R,\ldots,R),\diag(T,\ldots,T)}, \qquad 
  R=
  \begin{pmatrix}
    1&0\\0&-1
  \end{pmatrix}, \quad
  T=
  \begin{pmatrix}
    0&1\\1&0
  \end{pmatrix} \\
  2^{1+2}_-&=\gen{\diag(i,\ldots,i),\diag(j,\ldots,j)}
\end{align*}
where the representation of the dihedral group $2^{1+2}_+$ is of real
type and the representation of the quaternion group $2^{1+2}_-$ of
quaternion type.  This follows from \ref{lemma:Vast} because
$2^{1+2}_+$ has one faithful irreducible $\Ha$-representation of
degree $2$ and $2^{1+2}_-$ has one faithful irreducible
$\Ha$-representation of degree $1$.  The centralizers are
\cite[Proposition 4]{bob:stubborn}
\begin{align*}
  C_{\GL(2n,\Ha)}(2^{1+2}_+)&=\GL(n,\Ha), &   
  C_{\PGL(2n,\Ha)}(H_+)&=H_+ \times \PGL(n,\Ha) \\
  C_{\GL(2n,\Ha)}(2^{1+2}_-)&=\GL(2n,\R), &   
  C_{\PGL(2n,\Ha)}(H_-)&=H_- \times \PGL(2n,\R) 
\end{align*}
as we see by an argument similar to that of
\ref{sec:non0innerprodrank2}. This implies (\ref{lemma:ZGLprod}) that
$ZC_{\PGL(2n,\Ha)}(H)=H$ for all nontoral rank two objects $H$ of
$\A(\PGL(2n,\Ha))$.

We need to know that these nontoral objects satisfy condition ($3$) of
\ref{indstepalt}. To see this we use
\ref{lemma:uniquenu}. 

\underline{$H_+$}: Condition (1) of \ref{lemma:uniquenu} is
clearly satisfied since the identity component of
$C_{\PGL(2n,\Ha)}(H_+)$ is nontrivial when $n \geq 3$. The group
$H_+^*=2^{1+2}_+$ is contained in
$N(\GL(2n,\Ha))=N(\GL(1,\Ha))\wr\Sigma_{2n}$ and its centralizer there is
\begin{equation*}
  C_{N(\GL(2n,\Ha))}(2^{1+2}_+)=C_{N(\GL(1,\Ha))\wr\Sigma_{2n}}(2^{1+2}_+)
  =N(\GL(1,\Ha))\wr\Sigma_n=N(\GL(n,\Ha))
\end{equation*}
and therefore $H_-$ is contained in
$N(\GL(2n,\Ha))/\gen{-E}=N(\PGL(2n,\Ha))$ where its centralizer is
\begin{equation*}
   C_{N(\PGL(2n,\Ha))}(H_+)= H_+ \times N(\PGL(n,\Ha))=N(
   C_{\GL(2n,\Ha)}(H_+)) 
\end{equation*}
as in \ref{sec:non0innerprodrank2}. This means that $H_+ \subset
N(\PGL(2n,\Ha))$ is a \pl\ \cite{jmm:normax} of $H_+ \subset
\GL(2n,\Ha)$. Precomposing the inclusion $H_+ \subset N(\PGL(2n,\Ha))$
with the nontrivial element of $\A(\PGL(2n,\Ha))(H_+)=O^+(2,\F_2)
\cong C_2$ (\ref{prop:cfamH-L}) leads to another \pl . The third \pl\ 
is the quotient of
\begin{multline*}
  (2^{1+2}_+)^{\diag(B,\ldots,B)}=
  \big\langle \diag(R^B,\ldots,R^B),
  \diag((RT)^B,\ldots,(RT)^B)\big\rangle, \\
  B=\frac{1}{\sqrt{2}}
  \begin{pmatrix}
    1&i\\i&1
  \end{pmatrix}, \qquad
  R^B=T, \qquad
  (RT)^B=
  \begin{pmatrix}
    i&0\\0&-i
  \end{pmatrix}
\end{multline*}
Note that these three \pl s all have the same image in the Weyl group
$\pi_0N(\GL(2n,\Ha))=\pi_0(N(\GL(1,\Ha)))\wr\Sigma_{2n}$, namely the
subgroup generated by the permutation
$(1,2)(3,4)\cdots(2n-1,2n)\in\Sigma_{2n}$. 

Under the inductive assumption that $\PGL(n,\Ha)$ has
$\pi_*(N)$-determined auto\m s (or using \cite{jmo:selfho}) we
conclude from \ref{lemma:uniquenu} and diagram
(\ref{dia:Utriangle}) 
and (part of)
\cite[5.2]{jmm:ext} that condition (3) of
\ref{indstepalt} is
satisfied for the nontoral rank $2$ object $H_+$ of
$\A(\PGL(2n,\Ha))$. (Namely,
\ref{lemma:uniquenu}.(\ref{lemma:uniquenu1}) 
says that
$\nu_L'$ does not depend on the choice of $L<V$. The difference
$f_{\nu,L_2}^{-1} \circ f_{\nu,L_1}$ between any two of the maps
$f_{\nu,L}$ from \ref{indstepalt}.(\ref{indstepalt4})
is an auto\m\ of
$C_{\PGL(2n,\Ha)}(H_+)$ that, by
\ref{lemma:uniquenu}.(\ref{lemma:uniquenu2}), is the
identity on the identity component and by the commutative diagram
(\ref{dia:Utriangle})  
\set\begin{equation}\label{dia:cfampi0H+}
  \xymatrix{
    & H_+ \ar[dl] \ar[dr] \\
    C_{\PGL(2n,\Ha)}(H_+) \ar[rr]_{f_{\nu,L_2}^{-1} \circ f_{\nu,L_1}}
    && C_{\PGL(2n,\Ha)}(H_+)}
\end{equation}\add
also the identity on $\pi_0C_{\PGL(2n,\Ha)}(H_+)$. Since the identity
component of $C_{\PGL(2n,\Ha)}(H_+)$ has no center, this shows
that  $f_{\nu,L_2}^{-1} \circ f_{\nu,L_1}$ is the identity auto\m .)

\underline{$H_-$}: Condition (1) of
\ref{lemma:uniquenu} is clearly satisfied since the
identity component of $C_{\PGL(2n,\Ha)}(H_-)$ is nontrivial when $n
\geq 3$. The group $H_-^*=2^{1+2}_-$ is contained in
$N(\GL(2n,\Ha))=N(\GL(1,\Ha))\wr\Sigma_{2n}$ and its centralizer there
is
\begin{equation*}
  C_{N(\GL(1,\Ha))\wr\Sigma_{2n}}(2^{1+2}_-)
  \stackrel{\textmd{\ref{lemma:centrgw}}}{=}
  C_{N(\GL(1,\Ha))}(i,j)\wr\Sigma_n=\GL(1,\R)\wr\Sigma_{2n}= N(\GL(2n,\R))
\end{equation*}
and therefore $H_-$ is contained in
$N(\GL(2n,\Ha))/\gen{-E}=N(\PGL(2n,\Ha))$ where its centralizer is
\begin{equation*}
  C_{N(\PGL(2n,\Ha))}(H_-)= H_- \times
  N(\GL(2n,\R))/\gen{-E} = H_- \times N(\PGL(2n,\R))=N(
   C_{\PGL(2n,\Ha)}(H_-))  
\end{equation*}
as in \ref{sec:non0innerprodrank2}. This means that $H_- \subset
N(\PGL(2n,\Ha))$ is a \pl\ \cite{jmm:normax} of $H_- \subset
\GL(2n,\Ha)$. Precomposing the inclusion $H_- \subset N(\PGL(2n,\Ha))$
with elements of $\A(\PGL(2n,\Ha))(H_-)=O^-(2,\F_2)=\GL(2,\F_2)$
(\ref{prop:cfamH-L}) leads to other two \pl s of $H_-$.

Under the inductive assumption that the identity component
$\PSL(2n,\R)$ of $\PGL(2n,\R)$ has $\pi_*(N)$-determined auto\m s (or
using \cite{jmo:selfho}) we conclude from
\ref{lemma:uniquenu} and diagam
(\ref{dia:Utriangle})
and (part of) \cite[5.2]{jmm:ext} that condition (3) of
\ref{indstepalt} 
is satisfied for the nontoral rank $2$ object
$H_-$ of $\A(\PGL(2n,\Ha))$. (The argument for this is the same as in
case of $H_+$ with the little extra complication that
$\pi_0C_{\PGL(2n,\Ha)}(H_-)$ has an extra generator so that we replace
diagram (\ref{dia:cfampi0H+}) by
\set\begin{equation}\label{dia:cfampi0H-extra}
  \xymatrix{
  & \langle H_-, \diag(-1,1,\ldots,1) \ar[dl] \ar[dr] \\
  C_{\PGL(2n,\Ha)}(H_+) \ar[rr]_{f_{\nu,L_2}^{-1} \circ f_{\nu,L_1}} &&
  C_{\PGL(2n,\Ha)}(H_+)}
\end{equation}\add
from (\ref{dia:Utriangle}) where the slanted arrows
induce iso\m s on the component groups.)

\setcounter{subsection}{\value{thm}}
\subsection{Rank three nontoral objects of $\A(\PGL(2n,\Ha))$}
\label{sec:cfamnontrank3}\add

The nontoral rank three objects of the category $\A(\PGL(2n,\Ha))$ are
the quotients of $H_+\#L[i,n-i]$, $1 \leq i \leq [n/2]$, $H_-\#L[i,2n-i]$,
$1 \leq i \leq n$, and $V_0$. These subgroups of $\GL(2n,\Ha)$
are defined to be
\begin{gather*}
  \big\langle \diag(\overbrace{R,\ldots,R}^n),
  \diag(\overbrace{T,\ldots,T}^n),
  \diag(\overbrace{E,\ldots,E}^{i},\overbrace{-E,\ldots,-E}^{n-i})\big\rangle
  \\  
  \big\langle\diag(\overbrace{i,\ldots,i}^{2n}),
  \diag(\overbrace{j,\ldots,j}^{2n}),
  \diag(\overbrace{1,\ldots,1}^{i},\overbrace{-1,\ldots,-1}^{2n-i})\big\rangle
  \\
  \big\langle\diag\big(\overbrace{i,\ldots,i}^{2n}\big),
  \diag(\overbrace{R,\ldots,R}^n),\diag(\overbrace{T,\ldots,T}^n),
  \big\rangle
\end{gather*}
and their centralizers are
\begin{gather*}
  C_{\PGL(2n,\Ha)}(H_+\#L[i,n-i])=H_+ \times C_{\PGL(n,\Ha)}L[i,n-i],
  \\
  C_{\PGL(2n,\Ha)}(H_-\#L[i,n-i])=H_- \times C_{\PGL(2n,\R)}L[i,2n-i],
  \\
   C_{\PGL(2n,\Ha)}(V_0)=H_+ \times C_{\PGL(n,\Ha)}(I)
   \stackrel{\ref{eq:cfamCI}}{=} 
   H_+ \times \frac{\GL(n,\C)}{\gen{-E}} \rtimes \gen{j\gen{-E}}
\end{gather*}
so that (\ref{eq:cfamCL}, \ref{prop:ZC}, \ref{eq:cfamCI})
$ZC_{\PGL(2n,\Ha)}(V)=V$ for all nontoral rank three objects $V$ of
$\A(\PGL(2n,\Ha))$. The elements of $H_+\#L[i,n-i]$, $H_-\#L[i,2n-i]$,
and $V_0$ have traces (computed in $\GL(4n,\C)$) in the sets
$\pm\{0,4n-8i,4n\}$, $\pm\{0,4n-4i,4n\}$, and $\pm\{0,4n\}$.

\setcounter{subsection}{\value{thm}}
\subsection{Rank four nontoral objects of $\A(\PGL(2n,\Ha))$}
\label{sec:cfamnontrank4}\add

%$H_+\#P[1,i-1,n-i,0] \subset \GL(2n,\R) \subset \GL(2n,\Ha)$, $1<
%i \leq [n/2]$, is 
%\begin{equation*}
%  \big\langle \diag(\overbrace{R,\ldots,R}^{n}),
%              \diag(\overbrace{T,\ldots,T}^{n}),
%    \diag(E,\overbrace{-E,\ldots,-E}^{i-1},\overbrace{E,\ldots,E}^{n-i}),
%    \diag(E,\overbrace{E,\ldots,E}^{i-1},\overbrace{-E,\ldots,-E}^{n-i})
%  \big\rangle
%\end{equation*}
%It contains the two nontoral objects $H_+\#L[1,n-1]$ and
%$H_+\#L[2,n-2]$ when $i=2$ and the three nontoral objects
%$H_+\#L[1,n-1]$ and $H_+\#L[i-1,n+1-i]$, and $H_+\#L[i,n-i]$ when
%$i>2$.  Traces of the elements of $P$ are $2\{n+2-2i,-n+2i,-n+2\}$ and
%the only possibility for equality is $n+2-2i=-n+2i$ which happens for
%$2i=n+1$ \marginpar{-- which is excluded!}. 
%In this case the Quillen auto\m\ group has order $2^6$, in
%all other cases it has order $2^5$ (\ref{prop:cfamH+L}).

$H_-\#P[1,i-1,2n-i,0] \subset \GL(2n,\Ha)$, $1< i \leq n$, is
\begin{equation*}
  \langle \diag(\overbrace{i,\ldots,i}^{2n}), 
          \diag(\overbrace{j,\ldots,j}^{2n}),
      \diag(1,\overbrace{-1,\ldots,-1}^{i-1},\overbrace{1,\ldots,1}^{2n-i}),
      \diag(1,\overbrace{1,\ldots,1}^{i-1},\overbrace{-1,\ldots,-1}^{2n-i})
  \rangle
\end{equation*}
The elements of $P$ have traces $\{2n+2-2i,-2n+2i,2n+1\}$ and these
three integers are all distinct so that the Quillen auto\m\ group
(\ref{prop:cfamH-L}) has order $3 \cdot 2^5$. This nontoral rank four
object contains the two nontoral rank three objects
$H_-\#L[1,2n-1]$,$H_-\#L[2,2n-2]$ when $i=2$ and the three nontoral
rank three objects $H_-\#L[1,2n-1]$,$H_-\#L[i-1,2n-i+1]$,
$H_-\#L[i,2n-i]$ when $i>2$.

$V_0\#L[i,n-i] \subset \GL(2n,\C) \subset \GL(2n,\Ha)$, $1 \leq i \leq
[n/2]$, is the subgroup
\begin{equation*}
  \big\langle \diag(\overbrace{i,\ldots,i}^{2n}),
              \diag(\overbrace{R,\ldots,R}^{n}),
              \diag(\overbrace{T,\ldots,T}^{n}),
      \diag\big(\overbrace{E,\ldots,E}^{i},\overbrace{-E,\ldots,-E}^{n-i}\big)
  \big\rangle
\end{equation*}
containing the three rank three objects $H_+\#L[i,n-i]$,
$H_-\#L[2i,2n-2i]$, and $V_0$.

For these nontoral rank four objects $E \subset \GL(2n,\Ha)$, 
the center of the centralizer is finite (\ref{prop:ZC}) and as, of
course, $E \subset ZC_{\PGL(2n,\Ha)}(E)$ we see that
$\Hom_{\A(\PGL(2n,\Ha))}(\St(E),E)$ is a subspace of
$\Hom_{\A(\PGL(2n,\Ha))}(\St(E),\pi_1BZC_{\PGL(2n,\Ha)}(E))$.

\section{Higher limits of the functor $\pi_iBZC_{\PGL(n,\Ha)}$ on
  $\A(\PGL(n,\Ha))^{[\; ,\; ] \neq 0}$}
\label{sec:cfamlim}

In this section we compute the first higher limits of the functors
$\pi_jBZC_{\PGL(n,\Ha)}$, $j=1,2$. 

\begin{lemma}\label{lemma:cfamlim=0}
  $\lim^1\pi_1BZC_{\PGL(n,\Ha)} = 0 = \lim^2\pi_1BZC_{\PGL(n,\Ha)}$ and 
   $\lim^2\pi_2BZC_{\PGL(n,\Ha)} = 0 = \lim^3\pi_2BZC_{\PGL(n,\Ha)}$.
\end{lemma}

The case $j=2$ is easy. Since
$\pi_2BZC_{\PGL(n,\Ha)}$ has value $0$ on all objects of
$\A(\PGL(n,\Ha))^{[\; ,\; ] \neq 0}$  of rank $\leq 4$ it
is immediate from Oliver's cochain complex \cite{bob:steinberg} that
$\lim^2$ and $\lim^3$ of this functor are trivial.
We shall therefore now concentrate on the case $j=1$.

For any \lmntwo\ $E$ in $\PGL(n,\Ha)$ we shall write
\set\begin{equation}\label{defn:cfam[E]}
  [E]=\Hom_{\A(\PGL(n,\Ha)(E))}(\St(E),E)
\end{equation}\add
for the $\F_2$-vector space of $\F_2\A(\PGL(n,\Ha)(E)$-equivariant
maps from the Steinberg representation $\St(E)$ over $\F_2$ of
$\GL(E)$ to $E$. Olivers cochain complex has the form
(\ref{eq:Olccc}).

\begin{prop}\label{prop:cfamH-L}
  Regardless of the parity of $n$, the Quillen auto\m\ groups
  \begin{gather*}
  \A(\PGL(n,\Ha))(H_-)=O^-(2,\F_2) \\
  \A(\PGL(n,\Ha))(H_-\#V)=
  \begin{pmatrix}
    O^-(2,\F_2) & \ast \\ 0 & \A(\GL(n,\R))(V)
  \end{pmatrix}
  \end{gather*}
  and $\dim_{\F_2}[H_-] =1= \dim_{\F_2}[H_-\#L[i,2n+1-i]]$ as described
  in \ref{prop:H+H-} and \ref{prop:H-L}.
\end{prop}
\begin{proof}
  $\A(\GL(n,\Ha))(2^{1+2}_-)=\Out(2^{1+2}_-)$ since all auto\m s of
  $2^{1+2}_-$ preserve the trace. This group maps (isomorphically) to
  the subgroup $O^-(2,\F_2) \subset \GL(H_-)$ of auto\m s that
  preserve the quadratic function $q$ on $H_-$. The Quillen auto\m\
  group of  $H_-\#V$ consists of the auto\m s that lift to trace
  preserving auto\m s of $2^{1+2}_-\#V$. The dimension of the vector
  spaces of equivariant maps was computed by {\em magma}.
\end{proof}

In the odd case of $\GL(2n+1,\Ha)$ the cochain complex
(\ref{eq:Olccc}) takes the form
\set\begin{equation}\label{eq:cfamOlcccodd}
  0 \to [H_-] \xrightarrow{d^1} \prod_{1 \leq i \leq
    n}[H_-\#L[i,2n+1-i]] \xrightarrow{d^2} \prod_{|E|=2^4}[E] 
   \xrightarrow{d^3} \cdots
\end{equation}\add
and we need to show that $d^1$ is injective and that $\dim (\im d^2)
\geq n-1$.

If $E=H_-\#P[i]$, where $P[i]$ is as in (\ref{eq:cfamCP}), then
\begin{equation*}
  \A(\PGL(2n+1,\Ha))(H_-\#P[i])=
  \begin{pmatrix}
    O^-(2,\F_2) & \ast \\ 0 & \A(\SL(2n+1,\R))(P[i])
  \end{pmatrix}
\end{equation*}
where $\A(\SL(2n+1,\R))(P[i])$ is the group of trace preserving auto\m
s of $P[i]$. It turns out that 
\begin{equation*}
  \dim_{\F_2}[H_-\#P[i_0,i_1,i_2,i_3]]=
  \begin{cases}
    2 &  \A(\SL(2n+1,\R))(P[i])=\{E\} \\
    1 &  \A(\SL(2n+1,\R))(P[i])=C_2 \\
    0 &  \A(\SL(2n+1,\R))(P[i])=\GL(2,\F_2) \\
  \end{cases}
\end{equation*}
When $n=1$ or $n=2$, the cochain complex (\ref{eq:cfamOlcccodd}) has
the form
\begin{gather*}
  0 \to [H_-] \xrightarrow{d^1} [H_-\#L[1,2]] \xrightarrow{d^2}
  [H_-\#P[1,1,1,0]] \to \cdots \\
   0 \to [H_-] \xrightarrow{d^1} [H_-\#L[1,4]] \times [H_-\#L[2,3]]  
  \xrightarrow{d^2}  [H_-\#P[1,1,3,0]] \times  [H_-\#P[1,2,2,0]] \to
  \cdots 
\end{gather*}
where all vector spaces are one-dimensional. In the case of $n=1$,
$d^1$ is an iso\m , and in the case $n=2$, $d^1$ has matrix $
\begin{pmatrix}
  1 & 1 
\end{pmatrix}$ and $d^2$ has matrix $
\begin{pmatrix}
  1 & 1 \\ 1 & 1
\end{pmatrix}$. In case $n \geq 3$, it is enough to show that $d^1$
is injective and $d^2$ has rank $n-1$ in the cochain complex
\begin{equation*}
  0 \to [H_-] \xrightarrow{d^1} \prod_{1 \leq i \leq n}
  [H_-\#L[i,2n+1-i]] \xrightarrow{d^2}
  \prod_{2 <i \leq n}[H_-\#P[1,i-1,2n-i+1,0]] 
\end{equation*}
that agrees with (\ref{eq:cfamOlcccodd}) in degrees $1$, a product of
one-dimensional vector spaces, and $2$, a product of two-dimensional
vector spaces. The \lmntwo\ $H_-\#P[1,i-1,2n-i+1,0] \subset
\GL(2n+1,\Ha)$ contains the nontoral subspaces $H_-\#L[1,2n]$,
$H_-\#L[i-1,2n-i+2]$, and $H_-\#L[i,2n-i+1]$. The map $f_-$, defined
exactly as in \ref{eq:deff-}, is the nonzero element
of $[H_-]$ and the maps $df_-$, defined exactly as in
\ref{eq:basisH-L}, are nonzero in $H_-\#L[i,2n+1-i]$.
Thus $d^1$ is injective. A {\em magma\/} computation reveals that $\{
ddf_{-L[i-1,2n-i+2]}, ddf_{-L[i,2n-i+1]} \}$, where these
$\F_2\A(\GL(2n+1,\Ha))(H_-\#P[1,i-1,2n-i+1,0])$-maps are defined as in
\ref{eq:basisV0L}, is a basis for the two-dimensional
space $H_-\#P[1,i-1,2n-i+1,0]$ and that $ ddf_{-L[1,2n]}=
ddf_{-L[i-1,2n-i+2]}+ ddf_{-L[i,2n-i+1]}$. This shows that $d^2$ has
rank $n-1$.

In the even case of $\GL(2n,\Ha)$ the cochain complex (\ref{eq:Olccc})
takes the form
\begin{equation*}
  0 \to [H_-] \times [H_+] \xrightarrow{d^1}
  \prod_{1 \leq i \leq n}[H_-\#L[i,2n-i]] \times
   \prod_{1 \leq i \leq [n/2]}[H_+\#L[i,n-i]] \times [V_0] 
   \xrightarrow{d^2} \prod_{|E|=2^4} [E]
\end{equation*}

\begin{prop}\label{prop:cfamH+L}
  The auto\m\ groups of the low-degree nontoral objects of
  the Quillen category $\A(\PGL(2n,\Ha))$ are:
  \begin{align*}
  &\A(\PGL(2n,\Ha))(H_+)=O^+(2,\F_2) &
  &\A(\PGL(2n,\Ha))(H_+\#V)=
  \begin{pmatrix}
    O^+(2,\F_2) & \ast \\ 0 & \A(\GL(n,\Ha))(V) 
  \end{pmatrix} \\
   &\A(\PGL(2n,\Ha))(V_0) \cong \Symp(2,\F_2)  &
   &\A(\PGL(2n,\Ha))(V_0\#L[i,n-i]) \cong 
     \begin{pmatrix}
    \Symp(2,\F_2) & \ast \\ 0 & 1
  \end{pmatrix}
  \end{align*}
  and $\dim_{\F_2}[H_+] =2$, $\dim_{\F_2}[H_+\#L[i,n-i]]=3$,
  $\dim_{\F_2}[V_0]=4$, and $\dim_{\F_2}[V_0\#L[i,n-i]=5$ as described
  in \ref{prop:H+H-}, \ref{prop:H+L},
  and \ref{prop:V0}, and
  \ref{eq:basisV0L}.
\end{prop}
\begin{proof}
  The Quillen auto\m\ groups of the dihedral group $2^{1+2}_+$ and the
  generalized extraspecial group $4 \circ 2^{1+2}_{\pm}$ are the full
  outer auto\m\ groups because the traces are nonzero only on the
  derived groups which are characteristic. The images in $\GL(H_+)$,
  respectively, $\GL(V_0)$, isomorphic to $O^+(2,\F_2) \cong C_2$ and
  to $\Symp(2,\F_2)=\GL(2,\F_2)$, are the Quillen auto\m\ groups for
  $H_+$ and $V_0$. 
  For the middle formula, recall that the trace of
  $H_{\pm}\#V$ is the product of the traces.
\end{proof}

As in the real case (Chp~\ref{sec:dfam}) we get that $d^1$ embeds
$[H_-] \times [H_+]$ into $[V_0]$. The only problem is to show that
the rank of $d^2$ is $\geq n + 3[n/2]+4-3 = n + 3[n/2]+1$. We have to
show that
\begin{equation*}
  \dim(\im d^2) \geq n + 3[n/2] +1
\end{equation*}
We show this by mapping the $n+[n/2]+1$ nontoral rank three objects
(\ref{sec:cfamnontrank3}),
\begin{itemize}
\item $[H_-\#L[i,2n-i]]$, $1 \leq i \leq n$, with basis $\{df_\}$ as
  in (\ref{eq:basisH-L}),
\item $[H_+\#L[i,n-i]]$, $1 \leq i \leq [n/2]$, with basis
  $\{df_+,df_0,f_0\}$ as in (\ref{eq:basisH+L}), and
\item $[V_0]$ with basis $\{df_+,df_0,df_-,f_0\}$ as in
  (\ref{eq:basisV0})  
\end{itemize}
into the $(n-2)+[n/2]$ nontoral rank four
objects (\ref{sec:cfamnontrank4})
\begin{itemize}
\item $H_-\#P[1,i-1,2n+1-i]$, $2 < i \leq n$, with basis
  $\{ddf_{-L[i-1,2n+1-i]},ddf_{-L[i,2n-i]}\}$ where these maps are
  defined as the similar maps in (\ref{eq:basisV0L}),
\item $V_0\#L[i,n-i]$, $1 \leq i \leq [n/2]$, with basis
  \begin{equation*}
    \{ddf_{+L[i,n-i]}, ddf_{0L[i,n-i]}, df_{0L[i,n-i]},
  ddf_{-L[2i,2n-2i]}, df_{0V_0}\}
  \end{equation*}
  as in (\ref{eq:basisV0L})     
\end{itemize}
Computations with {\em magma\/} shows that the resulting $(n+3[n/2]+4)
\times (2n+5[n/2])$-matrix has rank $n + 3[n/2]+1$. The matrix has the
form (shown here for $n=5$)
\begin{center}
  \begin{tabular}[c]{r|ccc}
& $[H_-\#P[1,2,7]]$ & $[H_-\#P[1,3,6]]$ & $[H_-\#P[1,4,5]]$  \\ \hline
$H_-\#L[1,9]$ & 
$\begin{pmatrix}
  1 & 1
\end{pmatrix}$ &
$\begin{pmatrix}
  1& 1
\end{pmatrix}$ &
$\begin{pmatrix}
  1 & 1
\end{pmatrix}$ \\
$H_-\#L[2,8]$   & 
$\begin{pmatrix}
  1 & 0
\end{pmatrix}$ \\
$H_-\#L[3,7]$   &
$\begin{pmatrix}
  0 & 1
\end{pmatrix}$ &
$\begin{pmatrix}
  1 & 0
\end{pmatrix}$ \\
$H_-\#L[4,6]$   &&
$\begin{pmatrix}
  0 & 1
\end{pmatrix}$ &
$\begin{pmatrix}
  1 & 0
\end{pmatrix}$ \\
$H_-\#L[5,5]$   &&&
$\begin{pmatrix}
  0 & 1
\end{pmatrix}$ \\
$H_+\#L[1,4]$ \\
$H_+\#L[2,3]$ \\
$V_0$
  \end{tabular}
  \begin{tabular}[c]{cc|l}
$V_0\#L[1,4]$ & $V_0\#L[2,3]$ \\ \hline
&& $H_-\#L[1,9]$ \\
$\begin{pmatrix}
  0&0&0&1&0
\end{pmatrix}$ &&
$H_-\#L[2,8]$ \\
&& $H_-\#L[3,7]$ \\
&$\begin{pmatrix}
  0&0&0&1&0
\end{pmatrix}$ &
$H_-\#L[4,5]$ \\
&& $H_-\#L[5,5]$ \\
$A$ &&  $H_+\#L[1,4]$ \\
& $A$ &  $H_+\#L[2,3]$ \\
$B$ & $B$ & $V_0$
  \end{tabular}
\end{center}
where
\begin{equation*}
  A=
  \begin{pmatrix}
    1&0&0&0&0 \\
    0&1&0&0&0 \\
    0&0&1&0&0 \\
  \end{pmatrix}, \qquad
  B=
  \begin{pmatrix}
    1&0&0&0&0 \\
    0&1&0&0&0 \\
    0&0&0&1&0 \\
    0&0&0&0&1 \\
  \end{pmatrix}
\end{equation*}

\chapter{The exceptional $2$-compact groups}
\label{chp:di4andf4}

We use the material of the previous chapters to (re)prove that the
\twocg s $\Gtwo$, $\di$ and $\Ffour$ are uniquely $N$-determined.

\section{The $2$-compact group $\Gtwo$}
\label{sec:G2}

$B\Gtwo$ is a rank two \twocg\ containing a rank three \lmntwo\ $E_3
\subset \Gtwo$ such that $\A(\Gtwo)(E_3)=\GL(3,\F_2)$
\cite[6.1]{griess:elem} \cite[5.3]{dw:diagrams} and 
\begin{equation*}
  H^*(B\Gtwo :\F_2) \cong H^*(BE_3;\F_2)^{\GL(3,\F_2)} 
  \cong \F_2[c_4,c_6,c_7]
\end{equation*}
realizes the mod $2$ rank $3$ Dickson algebra \cite{mimura:char}.  The
Quillen category $\A(\Gtwo))$ contains exactly one isomorphism class
of objects $E_1,E_2,E_3$ of ranks $1,2,3$ as  Lannes theory
\cite{lannes} implies that the inclusion functor $\A(E_3,\GL(3,\F_2))
\to \A(\Gtwo)$ is an equivalence of categories. 
The centralizers of $E_1 \subset E_2 \subset E_3$ are
\begin{equation*}
  \SO(4) \supset T \rtimes \gen{-E} \supset E_3,
\end{equation*}
In all three cases, $ZC_{\Gtwo}(E_i)=E_i$ so that $\pi_2BZC_{\Gtwo}=0$
and $\pi_1BZC_{\Gtwo}=H^0(\GL(3,\F_2)(-);E_3)$.  Thus
$\pi_1BZC_{\Gtwo}$ is an exact functor
(\ref{dw:limits}) with
$\lim^0\pi_1BZC_{\Gtwo}=H^0(\GL(3,\F_2);E_3)=0$.

The Weyl $W(\Gtwo) \subset \GL(2,\Z) \subset \GL(2,\Z_2)$, of order
$12$, is generated by the two matrices  \cite[VI.4.13]{bo4-6}
\begin{equation*}
  \begin{pmatrix}
    -1&0\\\phantom{-}3&1
  \end{pmatrix}, \quad
  \begin{pmatrix}
    1&\phantom{-}1\\0&-1
  \end{pmatrix}
\end{equation*}
and the \mtn\ $N(\Gtwo)$ is the semi-direct product of the \mt\ and
the Weyl group \cite{cww}.

It is known that $H^0(W;\ch{T})(\Gtwo)=0$, $H^1(W;\ch{T})(\Gtwo)=0$,
and  $H^2(W;\ch{T})(\Gtwo)=0$
\cite{hms:first,hammerli:remarks}.

\begin{proof}[Proof of Theorem~\ref{thm:g2}]
  The rank one centralizer, $\SL(4,\R) = \SL(2,\C) \circ \SL(2,\C)$,
  is uniquely ${\N}$-determined (\ref{thm:afam}). Condition
  \ref{indstepalt}.(\ref{indstepalt1}) is satisfied because
  $H^1({W}(X);\ch{{T}}(X))=0$ for $X=\Gtwo$, $\SL(4,\R)$
  \cite{hms:first}, \ref{indstepalt}.(\ref{indstepalt3}) and
  \ref{indstepalt}.(\ref{indstepalt4}) because the only rank two
  object in $\Gtwo$ is toral and its centralizer is a \twoctg .  We
  noted above that the higher limits vanish. Now, \ref{ndetauto} and
  \ref{indstepalt} show that $\Gtwo$ is uniquely ${\N}$-determined.
  
  We have $\Aut(\Gtwo) = {W}(\Gtwo) \backslash
  {\N}_{\GL(2,\Z_2)}({W}(\Gtwo))$ (\ref{lemma:autX}) as the extension
  class $e(\Gtwo)=0$ \cite{cww}. The exact sequence (\ref{outN}) can
  be used to calculate the auto\m\ group.  Using the description of
  the root system from \cite[VI.4.13]{bo4-6} with short root
  $\alpha_1=\varepsilon_1-\varepsilon_2$ and long root
  $\alpha_2=2\varepsilon - \varepsilon_2 - \varepsilon_3$ generating
  the integral lattice in $\Z_2^3$ one finds that
  \begin{equation*}
    {\N}_{GL(2,\Z_2)}({W}(\Gtwo)) = \gen{\Z_2^{\times}, A, {W}(\Gtwo)}, \quad
    A=\sqrt{-3} 
    \begin{pmatrix}
      0 & 3 \\ 1 & 0
    \end{pmatrix}
  \end{equation*}
  and therefore $\Aut(\Gtwo) = \Z_2^{\times} / \Z^{\times} \times C_2$
  where the cyclic group of order two is generated by the exotic
  auto\m\ $A$ interchanging the two roots.
\end{proof}

%\begin{exmp}\cite{antonio:bg2} \cite{dw:diagrams}
%  ($\Gtwo$) $\A(\Gtwo)=\A((\Z/2)^3,\GL(3,\F_2))$ and the centralizers
%  are $\SO(4)=\SL(4,\R)=\SL(2,\C)^2/C_2$, $T \rtimes \Z/2$ and
%  $(\Z/2)^3$.
%\end{exmp}
%\begin{exmp}
%  $\Pin(n)$
%\end{exmp}

\section{The $2$-compact group $\di$}
\label{sec:di4}
%% /home/moller/manus/dfam/magma/di4/di4.prg
$B\di$ is a rank three \twocg\ containing a rank four \lmntwo\ $E_4
\subset \di$ such that $\A(\di)(E_4)=\GL(4,\F_2)$ and \cite{dw:new}
\begin{equation*}
  H^*(B\di ;\F_2) \cong H^*(BE_4;\F_2)^{\GL(4,\F_2)} \cong
  \F_2[c_8, c_{12}, c_{14}, c_{15}]
\end{equation*}
realizes the mod $2$ rank $4$ Dickson algebra.  Lannes theory
\cite{lannes} implies that the Quillen category $\A(\di)$ is
equivalent to $\A(\GL(4,\F_2),E_4)$ with exactly one \lmntwo\ (iso\m\ 
class), $E_1,\ldots, E_4$, of each rank $1,\ldots,4$.  The
centralizers of the toral subgroups $E_1,E_2, E_3$ and the nontoral
subgrouop $E_4$ are, respectively,
\begin{equation*}
  \Spin(7) \supset \SU(2)^3/\gen{(-E,-E,-E)} \supset T \rtimes \gen{-E}
  \supset E_4
\end{equation*}
and $ZC_{\di}(E_i)=E_i$ in all four cases so that the functor
$\pi_jBZC_{\di} \colon \A(\GL(4,\F_2),E_4) \to \Ab$ is the $0$-functor
for $j=2$ and equivalent to the functor $H^0(\GL(4,\F_2)(-);E_4)$ for
$j=1$. This is  an exact functor (\ref{dw:limits}) and
$\lim^0\pi_1BZC_{\di}  = H^0(\GL(4,\F_2);E_4)=0$.

%The first three \lmntwo s are toral and $E_4$ is nontoral so that the
%toral subcategory $\A(\di)^{\leq t}$ is equivalent to
%$\A(\GL(3,\F_2),E_3)$.

As may be seen from \cite{sheptodd}, the Weyl group $W(\di) \subset
\GL(3,\Z_2)$ of order $2|\GL(3,\F_2)|=336$ is generated by the
matrices
%%/home/moller/manus/dfam/magma/di4/di4.prg
\begin{equation*}
  \begin{pmatrix}
    1&0&0\\2&-1&-1\\0&0&1
  \end{pmatrix}, \quad
  \begin{pmatrix}
    1&0&0\\0&0&1\\2&-1&-1
  \end{pmatrix}, \quad
  \begin{pmatrix}
    -1&1&1\\0&1&0\\0&0&1
  \end{pmatrix},\quad
  \begin{pmatrix}
    -v&0&v^2+v\\-1&1&v\\-2v&0&v
  \end{pmatrix}
\end{equation*}
where $v \in \Z_2$ is the unique $2$-adic integer with $2v^2-v+1=0$.
The first three matrices generate $W(\Spin(7))$ \cite[3.9,
3.11]{brockerdieck}.  Since $W(\di)$ is isomorphic to $\GL(3,\F_2)
\times \gen{-E}$,
\begin{equation*}
  H^n(W;\ch{T})(\di) = 
 \bigoplus_{2i \leq n} H^{n-2i}(\GL(3,\F_2);H^{2i}(\gen{-E};\ch{T})=
 \bigoplus_{2i \leq n} H^{n-2i}(\GL(3,\F_2);(\Z/2)^3)
\end{equation*}
and, in particular,
\begin{equation*}
  H^0(W;\ch{T})(\di) =0, \quad  H^1(W;\ch{T})(\di) = \Z/2, \quad 
  H^2(W;\ch{T})(\di) = \Z/2
\end{equation*}
%and since the center of $\di$ is trivial it is not regular in the
%sense of \cite{hms:first}. 
We may characterize the \mtn\ \ses\ for $\di$ as the unique nonsplit
extension of $\ch{T}$ by $W(\di)$; it is nonsplit because the
restriction to $W(\Spin(7)) \subset W(\di)$ is nonsplit \cite{cww}.

We can not use \ref{indstepalt} as it stands because
condition (2) fails: The restriction map
\begin{equation*}
\Z/2=H^1(W;\ch{T})(\di) \to 
H^1(W;\ch{T})(\Spin(7))\stackrel{\textmd{\cite{hms:first}}}{=}(\Z/2)^2  
\end{equation*}
is not surjective.  Note, however, that the proof of
\ref{indstepalt} goes through with only insignificant
changes if we replace hypotheses (1) and (2) by
    \begin{enumerate}
    \item[(1 \& 2)] The centralizer of any toral $(V,\nu) \in
      \Ob(\A(X)^{\leq t}_{\leq 2})$ is uniquely $N$-determined. 
    \end{enumerate}
and leave the other conditions unchanged.

\begin{proof}[Proof of Theorem~\ref{thm:di4}]
  Condition (1 \& 2) is satisfied for $\di$ since the connected \twocg
  s $\Spin(7)$ and $\SU(2)^2/\Delta$ are uniquely $N$-determined
  (\ref{thm:afam}, \ref{mainthm}) and the general results of
  \ref{cha:ndet}.\S\ref{sec:red}. Since also the relevant higher
  limits vanish \cite[2.4]{dw:new}, $\di$ is uniquely $N$-determined
  by \ref{ndetauto} and \ref{indstepalt}. Since $\Out_{\tr}(W(\di))$
  is trivial and $Z(W(\di))=\gen{-E}$ has order two, $\Aut(\di)$ can
  be read off from \ref{lemma:autX} and the \ses\ (\ref{autwNgl}).
\end{proof}

\section{The $2$-compact group $\Ffour$}
\label{sec:f4}
%%/home/moller/manus/dfam/magma/f4/f4.prg
$B\Ffour$ is  a rank four \twocg\ containing a rank
five \lmntwo\ $E_5 \subset \Ffour$ such that \cite[2.1]{antonio:f4p2}
\begin{equation*}
  H^*(B\Ffour;\F_2) \cong H^*(BE_5;\F_2)^{\A(\Ffour)(E_5)} \cong 
  \F_2[y_4, y_6, y_7, y_{16}, y_{24}]
\end{equation*}
where the Quillen auto\m\ group is the parabolic subgroup
\begin{equation*}
  \A(\Ffour)(E_5) =
  \begin{pmatrix}
    \GL(2,\F_2) & \ast \\ 0 & \GL(3,\F_2)
  \end{pmatrix} \subset \GL(5,\F_2)
\end{equation*}
of order $2^6|\GL(2,\F_2)| |\GL(3,\F_2)|$. The inclusion functor
$\A(\A(\Ffour)(E_5),E_5) \to \A(\Ffour)$ is a category equivalence by
Lannes theory \cite{lannes}.  Inspection of the list of centralizers
of \lmntwo s in $\Ffour$ \cite[3.2]{antonio:f4p2} shows that
$ZC_{\Ffour}(V)=V$ for each nontrivial $V \subset \E_5$ so that the
functor $\pi_2BZC_{\Ffour}=0$ and $\pi_1BZC_{\Ffour}=
H^0(\A(\Ffour)(E_5)(-);E_5)$.  Thus $\pi_1BZC_{\Ffour}$ is an exact
functor (\ref{dw:limits}) and
$\lim^0\pi_1BZC_{\Ffour}=H^0(\A(\Ffour)(E_5);E_5)=0$.

It is known that $ H^0(W;\ch{T})(\Ffour)=0$,
$H^1(W;\ch{T})(\Ffour)=0$, and $H^2(W;\ch{T})(\Ffour)=\Z/2$
\cite{hms:first,hammerli:remarks}.

\begin{proof}[Proof of Theorem~\ref{thm:f4}]
  Condition (1 \& 2) is satisfied for $\Ffour$ because
  centralizers of rank one objects \cite[3.2]{antonio:f4p2}
  \begin{equation*}
    \frac{\SU(2) \times \Symp(3)}{\Z/2}, \quad
    \Spin(9)
  \end{equation*}
  and centralizers of rank two objects \cite[3.2]{antonio:f4p2},
  \begin{equation*}
    \frac{\U(1) \times \U(3)}{\Z/2}\rtimes \Z/2, \quad
    \frac{\Spin(4) \times \Spin(5)}{\Z/2}, \quad
    \Spin(8)
  \end{equation*}
  have uniquely $N$-determined centralizers.  This follows from
  \ref{cha:ndet}.\S\ref{sec:red} as the simple factors are uniquely
  $N$-determined (\ref{thm:afam}, \ref{mainthm}). There are no
  nontoral \lmntwo\ of rank two. There is a unique nontoral $E_3$ of
  rank three and a unique nontoral \lmntwo\ $E_4$ of rank four and
  $\A(\Ffour)(E_3)=\GL(3,\F_2)$, $\A(\Ffour)(E_3)=2^3 \cdot
  \GL(3,\F_2)$, $C_{\Ffour}(E_3)=(\Z/2)^2 \times \Or(3)$,
  $C_{\Ffour}(E_4)=(\Z/2)^3 \times \Or(2)$ \cite{antonio:f4p2}. It
  follows that the relevant higher limits vanish, and since there are
  no nontoral \lmntwo s of rank two, $\Ffour$ is uniquely
  $N$-determined by \ref{ndetauto} and \ref{indstepalt}.
  
  The auto\m\ group of the \twocg\ $\Ffour$ is (\ref{lemma:autX}) the
  middle term of the exact sequence (\ref{autwNgl}). ( All auto\m s of
  $\Ffour$ preserve the extension class $e(\Ffour)$ which is the
  nontrivial element of $H^2(W;\ch{T})=\Z/2$ \cite{cww,
    matthey:second}.)  The group $\Out_{\tr}(W(\Ffour))$ of trace
  preserving outer auto\m s is cyclic of order two but the nontrivial
  outer auto\m\ of $W(\Ffour)$ can not be realized as conjugating with
  an element of $N_{\GL(L)}(W)$. The center of $W(\Ffour)$ is
  $C_2=\gen{-E}$. We conclude that $\Aut(\Ffour)=\Z^{\times}
  \backslash \Z_2^{\times}$ consists entirely of unstable Adams
  operations.
\end{proof}

\section{The $E$-family}
\label{sec:efam}

We consider the \twocg s $E_6$, $PE_7$, and $E_8$.

\begin{proof}[Proof of Theorem~\ref{thm:efam}]
  According to \cite{griess:elem}, the Lie group $E_6$ contains two
  conjugacy classes of elements of order two, $2A$ and $2B$. Their
  centralizers are $C_{E_6}(2A)=\SU(2) \circ \SU(6)$ and $C_{E_6}(2B)$
  has shape $T_1D_5$. There is a unique maximal nontoral \lmntwo\ $E_5
  \subset \Ffour \subset E_6$. It follows that $E_6$ contains a unique
  nontoral rank three \lmntwo\ $E_3$ and a a unique nontoral rank four
  \lmntwo\ $E_4$ with auto\m\ groups
  $\A(E_6)(E_3)=\A(\Ffour)(E_3)=\GL(3,\F_2)$ and $\A(E_6)(E_4) \supset
  \A(\Ffour)(E_4)=2^3 \cdot \GL(3,\F_2)$. We may compute the
  centralizer of $E_3$, $C_{E_6}(E_3)=E_3 \times \SU(3)$, as the
  centralizer of rank two \lmntwo\ in $C_{E_6}(2A)$. Consequently, the
  centralizer $C_{E_6}(E_4)=E_3 \times C_{\SU(3)}(\Z/2)=E_3 \times
  \U(2)$. Using Oliver's cochain complex we see that the relevant
  higher limits vanish. We verify the first two conditions of
  \ref{indstepalt} by using \ref{lemma:altcond2}. Since $H^1(W;
  \ch{T})(E_6)=0$ \cite{matthey:normalizers} we need to show that the
  limit from \ref{lemma:altcond2} is trivial. This is a machine
  computation in the toral subcategory of $\E_6$. There are two toral
  subgroups of rank one (as already mentioned), both with
  $H^1(W_0;\ch{T})=\Z/2$, and four toral subgroups of rank two; two of
  these are pure and two are nonpure. By considering the two nonpure
  \lmntwo s of rank two, we see that the limit is indeed trivial.
  Since also there are no rank two nontoral \lmntwo s, it follows
  from \ref{indstepalt}, that $E_6$ is $N$-determined. It has
  $\pi_*(N)$-determined auto\m s by  \ref{ndetauto}. 
  
  According to \cite{griess:elem}, the Lie group $PE_7$ contains three
  conjugacy classes of elements of order two with centralizers
  $C_{PE_7}(2B)=\SL(2,\C) \circ \SL(12,\R)$, $C_{PE_7}(4A)=(\GL(1,\C)
  \circ E_6) \rtimes \Z/2$, and $C_{PE_7}(4H)={\SL(8,\C)}/{\gen{i}}
  \rtimes \Z/2$.  We verify the first two conditions of
  \ref{indstepalt} by using \ref{lemma:altcond2}. Since $H^1(W;
  \ch{T})(PE_7)=0$ \cite{matthey:normalizers} we need to show that the
  limit from \ref{lemma:altcond2} is trivial. This is a machine
  computation in the toral subcategory of $PE_7$. We have
  $H^1(W_0;\ch{T})=\Z/2$ for each of the centralizers just
  mentioned. There is a unique rank two toral object $E_2 \subset
  PE_7$ containing all three conjugacy classes $2B$, $4A$, $4H$. The
  centralizer of $E_2$ has shape $T_1A_1A_5$ and
  $H^1(W_0;\ch{T})=(\Z/2)^3$. 
  %% /home/moller/manus/e7/w7.prg
  By mapping the three rank one objects into this single rank two
  object, we see that the limit is indeed trivial. Still according to
  \cite{griess:elem}, there are three nontoral rnak two \lmntwo s. One
  is $2A$-pure, has centralizer $(\Z/2)^2 \times \PSL(8,\R)$ and
  auto\m\ group $\A(PE_7)(E_2)=\GL(2,\F_2)$, one is nonpure, has
  centralizer $(\Z/2)^2 \times \PGL(4,\Ha)$ and auto\m\ group
  $\A(PE_7)(E_2)=C_2$, and the third one is $4H$-pure, has centralizer
  $(\Z/2)^2 \times \Ffour$ and auto\m\ group
  $\A(PE_7)(E_2)=\GL(2,\F_3)$. We use \ref{lemma:uniquenu} to verify
  the third condition of \ref{indstepalt}. It turns out that $W(PE_7)$
  contains two elements, $v_1$ and $v_2$, of order two with
  $+1$-eigenspace of dimension four. The image of $C_{W(PE_7)}(v_1)$
  in $\GL(\pi_1(T^{v_1}\otimes\Q)=\GL(4,\Q)$ has order $2^5 \cdot 3^1$
  while in case of $v_2$ we get an image of order $2^7 \cdot 3^2$.
  We conclude that if $E_2 \subset N(PE_7)$ is \pl\ of the nonpure
  nontoral rank two object, then the image in $W(PE_7)$ is $v_1$. 
  %% /home/moller/manus/e7/O2w7.prg  
  This observation can be used in connection with \ref{lemma:uniquenu}
  to verify the third condition of \ref{indstepalt}. It remains to be
  seen that the relevant higher limits vanish. This is a computation
  with Oliver's cochain complex very similar to the ones alredy seen.
  We omit the details. Now \ref{ndetauto} and \ref{indstepalt}
  imply that $PE_7$ is uniquely $N$-determined.
  
  The proof for the \twocg\ $E_8$ is similar to the proof for $E_6$
  since there are no nontoral rank two \lmntwo s in the simply
  connected Lie group $E_8$.
\end{proof}

\chapter{Proofs of the main theorems}
\label{cha:proofs}

This chapter contains the proofs the main results stated in the
Introduction.

\section{Proof of Theorem~\ref{thm:afam}}

We show that $\PGL(n+1,\C)$ is uniquely $N$-determined by induction
over $n$. The induction step is provided by Lemma~\ref{vercond} and
the start of the induction by Proposition~\ref{pgl2C}.

\begin{lemma}\label{vercond}
  Suppose that $\pgl{r+1}$ is uniquely ${\N}$-determined for all $0
  \leq r <n$. Then $\pgl{n+1}$, $n \geq 1$, satisfies conditions
  \ref{ndetauto}.(\ref{ndetauto1}) (for
  $\pi_*(\N)$-determined auto\m s),
  \ref{indstepalt}.(\ref{indstepalt3}),
  \ref{indstepalt}.(\ref{indstepalt1}), and
  \ref{indstepalt}.(\ref{indstepalt4}).
\end{lemma}
\begin{proof}
  Condition (1) of \ref{ndetauto} (for
  $\pi_*(\N)$-determined auto\m s) is concerned with centralizers
  $C_{\pgl{n+1}}(L,\lambda)$ of rank one objects (\ref{eq:afamCL}).
  The condition is satisfied for all connected rank one centralizers
  by the induction hypothesis and \ref{ndetcenter},
  \ref{prodauto}.  The condition is satisfied for the
  nonconnected rank one centralizer (when $n+1$ is even) by
  \ref{lemma:autononcon} since
  $H¹(C_2;\Z/2^{\infty})=0$ for the nontrivial action of the cyclic
  group $C_2$ of order two on $\Z/2^{\infty}$.
  
  We use \ref{lemma:altcond2} to verify conditions (1)
  and (2) of \ref{indstepalt}.  Let $(V,\nu)$ be a
  toral elementary abelian $2$-subgroup of $\pgl{n+1}$ of rank $\leq
  2$ and $C(\nu) = C_{\pgl{n+1}}(\nu)$ its centralizer. We have seen
  that $C(\nu)$ is LHS (\S\ref{sec:afamlhs}) and that
  $\ch{{\Ze}}(C(\nu)_0)=\ch{{\Ze}}({\N}_0(C(\nu)))$ as $C(\nu)_0$ does
  not contain a direct factor isomorphic to
  $\GL(2,\C)/\GL(1,\C)=\SO(3)$ (\ref{sec:centermtn},
  (\ref{eq:afamCL})).  The identity component $C(\nu)_0$ has
  $\pi_*({\N})$-determined auto\m s according to
  \ref{autondetcenter} and
  \ref{prodauto}, and $C(\nu)$ has ${\N}$-determined
  auto\m s by \ref{lemma:autononcon}.  The identity
  component $C(\nu)_0$ is ${\N}$-determined according to
  \ref{ndetcenter} and \ref{ndetprod}, and $C(\nu)$ is
  ${\N}$-determined by \ref{redtoconnected}. Thus
  $C(\nu)$ is LHS and totally ${\N}$-determined.
  
  The functor $H^1({\W}/{\W}_0;\ch{{\T}}^{\W}_0)$ is zero on
  $\A(\pgl{n+1})^{\leq t}_{\leq 2}$ except on the object
  $(V,\nu)=(i_0,i_0,i_0,i_0)$, when $n+1=4i_0$, where it has value
  $\Z/2$. However, this object has Quillen auto\m\ group $\GL(V)$ and
  since the only $\GL(V)$-equivariant homo\m\ $\St(V)=V \to \Z/2$ is
  the trivial homo\m , $\lim^1(\A(\pgl{n+1})^{\leq t}_{\leq
    2};H^1({\W}/{\W}_0;\ch{{\T}}^{\W}_0))=0$ follows from Oliver's
  cochain complex \cite{bob:steinberg}.
 
  We now turn to condition (3) of \ref{indstepalt}.
  When $n+1$ is odd there are no nontoral rank two objects
  (\ref{tpglnC}) and so there is nothing to prove.  When $n+1=2m$ is
  even, let $(H,\nu)$ be the unique nontoral rank two object of
  $\A(\pgl{2m})$ (\ref{nonTpglnC}.(1)).  Let $X'$ be a connected
  \twocg\ with \mtn\ \func{j'}{{\N}(\PGL(2m,\C)}{X'}.  We must show
  that $\nu_L'$ and
  \func{f_{\nu,L}}{C_{\pgl{2m}}(H,\nu)}{C_{X'}(H,\nu_L')} as defined
  in \ref{indstepalt}.(\ref{indstepalt4}) are
  independent of the choice of rank one subgroup $L \subset V$.  When
  $m=1$, the claim follows from \ref{lemma:C3cond},
  \ref{help1}, \ref{dia:triangle} (where
  $\overline{\nu}(V)$ and $\overline{\nu}'(V)$ are iso\m s in this
  case) since $\PGL(2,\C)$ does contain a unique rank one \lmntwo\ 
  with nonconnected centralizer (\ref{prop:afamtoral}) and a unique
  nontoral rank two \lmntwo\ (\ref{nonTpglnC}.(1)). When $m>1$, we use
  \ref{lemma:uniquenu} which immediately yields that
  $\nu_L'$ is independent of the choice of $L<V$. There exists a torus
  $T_{\nu} \to C_N(V,\nu_L^N)$ as in
  \ref{lemma:uniquenu}.(\ref{lemma:uniquenu2}) because
  the three \pl s $\nu_L^N$, $L < V$, differ by an auto\m\ of $H$ (the
  Quillen auto\m\ group of $(H,\nu)$ is the full auto\m\ group
  $\Aut(H)$ of $H$ (\ref{nonTpglnC}.(1))). Since the identity
  component of $C_{\pgl{n+1}}(H,\nu)$ is uniquely $\N$-determined by
  induction hypothesis, the restriction $(f_{\nu,L})_0$ of $f_{\nu,L}$
  to the identity components is independent of the choice of $L$
  (\ref{underT}.(\ref{underT2})). Also
  $\pi_0(f_{\nu,L})$ is independent of the choice of $L<V$ by
  \ref{dia:triangle} (where $\pi_0(\overline{\nu}(V))$ and
  $\pi_0(\overline{\nu}'(V))$ are iso\m s).  But since $\pgl{m}$ is
  centerfree $f_{\nu,L}$ is in fact determined (use one half of
  \cite[5.2]{jmm:ext}) by $(f_{\nu,L})_0$ and $\pi_0(f_{\nu,L})$. We
  conclude that $f_{\nu,L}$ is independent of the choice of $L<V$.

\end{proof}

\begin{prop}\label{pgl2C}
  The \twocg\ $\PGL(2,\C)$ is uniquely ${\N}$-determined.
\end{prop}
\begin{proof}
  The centralizer cofunctor $C_{\PGL(2,\C)}$ takes the Quillen category of
  $\PGL(2,\C)$, consisting (\ref{prop:afamtoral}, \ref{nonTpglnC})
  of one toral line, $L$, and one nontoral plane, $H$,
  \set\begin{equation}\label{Qforpgl2C}\SelectTips{cm}{}
    \xymatrix@C=65pt{ *\txt{\makebox[12mm][r]{$L$}} \ar[r] &
      *\txt{\makebox[5mm][l]{$H$}} \ar@(ur,dr)^{\GL(H)} }
  \end{equation}\add
 to the diagram
 \begin{equation}\SelectTips{cm}{}
  {\xymatrix@C=65pt{ *\txt{\makebox[35mm][r]{$\GL(1,\C)^2/\GL(1,\C)
          \rtimes C_2$}} 
        & *\txt{\makebox[5mm][l]{$H$}}
      \ar@(ur,dr)^{\GL(H)^{\op}} \ar[l]  }}
  \end{equation}
  of uniquely ${\N}$-determined \twocg s. The \twoctg\ to the left is
  uniquely ${\N}$-determined  because (\ref{ex:toral})
  $H^1(C_2;\Z/2^{\infty}) = 0$ for the nontrivial action of $C_2$ on
  $\Z/2^{\infty}$.  The center cofunctor takes this diagram back to
  the starting point (\ref{Qforpgl2C}) for which the higher limits
  vanish (\ref{dw:limits}).  $\PGL(2,\C)$ is thus uniquely
  ${\N}$-determined by \ref{ndetauto} and
  \ref{indstepalt}. 
\end{proof}

\begin{proof}[Proof of Theorem~\ref{thm:afam}]
  The proof is by induction over $n \geq 1$. The start of the
  induction is provided by \ref{pgl2C}. The induction step is provided
  by \ref{vercond} and \ref{lim=0} using \ref{ndetauto} and
  \ref{indstepalt}.  
  
  According to \ref{lemma:autX}, the auto\m\ group
  \begin{equation*}
  \Aut(\PGL(n+1,\C)) = W \backslash N_{\GL(L)}(W) = W \backslash
  \gen{\Z_2^{\times},W} = Z(W) \backslash \Z_2^{\times}
  \end{equation*}
  is isomorphic to $\Z^{\times} \backslash \Z_2^{\times}$ for $n=1$
  and to $\Z_2^{\times}$ for $n>1$. Here we use \cite{jmm:norm} or the
  exact sequence (\ref{autwNgl}) where we note that $\Out_{\tr}(W)$,
  is trivial for all $n \geq 1$; $\Out(\Sigma_{n+1})$ is trivial for
  all $n \neq 5$ \cite[II.5.5]{huppert:I} and the nontrivial outer
  auto\m\ of $\Sigma_6$ does not preserve trace.
\end{proof}

\begin{proof}[Proof of Corollary \ref{cor:afam}]
  Let $X=\GL(n,\C)$, $n \geq 1$, and write $\ch{T}$, $W$ and $L$ for
  $\ch{T}(X)$, $W(X)$, and $L(X)$.  Since the adjoint form
  $PX=P\GL(n,\C)$ of $X$ is uniquely $N$-determined (\ref{thm:afam}),
  so is $X$ (\ref{autondetcenter}, \ref{ndetcenter}). The extension
  class $e(X) \in H^2(W;\ch{T})$ (\S\ref{sec:AM}) is the zero class
  since the \mtn\ $N(X)=\GL(1,\C)\wr\Sigma_n$ splits. Therefore,
  $\Aut(X)$ is isomorphic to $W \backslash N_{\GL(L)}(W)$
  (\ref{lemma:autX}).  Using the exact sequence (\ref{autwNgl}), we
  conclude, as in the proof of Theorem~\ref{thm:afam}, that $\Aut(X)
  \cong Z(W) \backslash \Aut_{\Z_2W}(L) = Z(W) \backslash
  \Aut_{\Z_2\Sigma_n}(\Z_2^n)$.
\end{proof}

\section{Proof of Theorem~\ref{mainthm}}

The proof of Theorem~\ref{mainthm} uses induction over $n$
simultaneously applied to the three infinite families $\PSL(2n,\R)$,
$\SL(2n+1,\R)$, and $\PGL(n,\Ha)$.

\begin{proof}[Proof of Theorem~\ref{mainthm}]
  The statement of the theorem means (\ref{defn:det})
  that the \twocg s
  \begin{itemize}
  \item $\pslr{2n}$, $\SL(2n+1,\R)$, and $\PGL(n,\Ha)$ have
    $\pi_*(N)$-determined auto\m s, 
  \item $\pslr{2n}$, $\SL(2n+1,\R)$, and $\PGL(n,\Ha)$ are
    $N$-determined. 
  \end{itemize}
  We may inductively assume that the connected \twocg s $\pslr{2i}$,
  $1 \leq i \leq n-1$, $\SL(2i+1,\R)$, $1 \leq i <n-1$, and
  $\PGL(i,\Ha)$, $ 1 \leq i < n$, are uniquely $N$-determined.  From
  Theorem~\ref{thm:afam} we know that $\PGL(i,\C)$ is uniquely
  $N$-determined for all $i \geq 1$. The plan is now to use
  \ref{ndetauto} and \ref{indstepalt} inductively.

  Consider first the connected, centerless \twocg\ $\pslr{2n}$.
 
  \noindent
  \underline{$\pslr{2n}$ has $N$-determined auto\m s}:
  According to  \ref{ndetauto} it suffices to show that
  \begin{enumerate}
  \item $C_{\pslr{2n}}(L)$ has $N$-determined auto\m\ for any rank one
      \lmntwo\ $L \subset \pslr{2n}$.
    \item $\lim^1(\A(\pslr{2n}), \pi_1BZC_{\pslr{2n}}) = 0 = 
      \lim^2(\A(\pslr{2n}), \pi_2BZC_{\pslr{2n}})$.
  \end{enumerate}
  Item (2) is proved in \ref{lemma:lim=0}. The
  centralizers that occur in item (1) are listed in
  (\ref{CLqeq0}) and (\ref{CLneq0}).
  That the centralizers of (\ref{CLqeq0}) have
  $N$-determined auto\m s follows, under the induction hypothesis that
  the \twocg s $\pslr{2i}$, $1 \leq i \leq n-1$, have $N$-determined
  auto\m s, from general hereditary properties of $N$-determined
  \twocg s (\S\ref{sec:red}). Note here that
  $\ch{Z}(C_0)=\ch{T}(C_0)^{W(C_0)}$ for $C=C_{\pslr{2n}}(L)$ by
  \cite[1.6]{matthey:normalizers}.  Similarly, the centralizers of
  (\ref{CLneq0}) have $N$-determined auto\m s because
  the \twocg s $\PGL(n,\C)$, $1 \leq n < \infty$, have $N$-determined
  auto\m s (\ref{thm:afam}).
  
  \noindent
  \underline{$\pslr{2n}$ is $N$-determined}: We verify the four
  conditions of \ref{indstepalt}. Let $V \subset
  \pslr{2n}$ be a toral \lmntwo\ of rank at most $2$. The centralizer
  $C=C_{\pslr{2n}}(V)$ is one of the \twocg s listed in
  (\ref{CLqeq0}), (\ref{CVqeq0}),
  (\ref{CLneq0}), or (\ref{CVqneq0}), so
  it is LHS (\ref{lemma:lhs}).  The identity component $C_0$ of $C$
  satisfies the equation $\ch{Z}(C_0)=\ch{T}(C_0)^{W(C_0)}$
  \cite[1.6]{matthey:normalizers} and the adjoint form
    \begin{equation*}
       PC_0=
      \begin{cases}
        \PSL(2i_0,\R) \times \PSL(2i_1,\R) & i_0+i_1=n\\
        \PSL(2i_0,\R) \times \PSL(2i_1,\R) \times 
        \PSL(2i_2,\R) \times \PSL(2i_3,\R) & i_0+i_1+i_2+i_3=n\\
        \PGL(n,\C) \\
        \PGL(i_0,\C) \times \PGL(i_1,\C) & i_0+i_1=n
      \end{cases}
    \end{equation*}
    in these four cases. The induction hypothesis and the general
    results of \S\ref{sec:red} imply that $C_0$ is
    uniquely $N$-determined and that $C$ is totally $N$-determined.
    Since also the homo\m\ $H^1(W;\ch{T}) \twoheadrightarrow
    \lim^1(\A(\pslr{2n})^{\leq t}_{\leq 2}; H^1(W_0;\ch{T})^{W/W_0})$
    is surjective (\ref{sec:dfam}.\S\ref{sec:lim0}), we get from
    \ref{lemma:altcond2} that the first two conditions
    of \ref{indstepalt} are satisfied. The third
    condition has been verified in \ref{sec:dfam}.\S\ref{sec:NdetD}
    % \footnote{Using that $\PGL(n,\Ha)$ has
    %  $\pi_*(N)$-determined auto\m s} 
    and the fourth, and final, condition in
    \ref{lemma:lim=0}.

  \noindent
  \underline{$\pslr{2n}$ has $\pi_*(N)$-determined auto\m s}: This
  means that the only auto\m\ of $\pslr{2n}$ that restricts to the
  identity on the \mt\ is the identity, ie that 
  \begin{equation*}
   H^1(W;\ch{T})(\pslr{2n}) \cap \AM(\Aut(\pslr{2n}))=\{0\} 
  \end{equation*}
  where $\AM$ is the Adams--Mahmud homo\m\ 
  (\S\ref{sec:AM}).  For $n>4$,
  $H^1(W;\ch{T})(\pslr{2n})=0$, and there is nothing to prove.
  Consider the case $n=4$. Let $f$ be an auto\m\ of $\pslr{8}$ such
  that $\AM(f) \in H^1(W;\ch{T})$. Let $L \subset \pslr{8}$ be any
  rank one \lmntwo . Since $f$ is the identity on the \mt , $f(L)$ is
  conjugate to $L$ so that $f$ restricts to an auto\m\ of
  $C_{\pslr{8}}(L)$ and to an auto\m\ of the identity component of
  $C_{\pslr{8}}(L)$.  Since $C_{\pslr{8}}(L)_0$ has
  $\pi_*(N)$-determined auto\m s by
  \ref{autondetcenter} and
  \ref{prodauto}, $f \in H^1(W;\ch{T})(\pslr{2n})$
  restricts to $0$ in $H^1(W;\ch{T})(C_{\pslr{8}}(L)_0)$. However, the
  restriction map is injective (see the proof of
  \ref{prop:lim0H1}) so that $f=0$. This shows that
  $\pslr{8}$ has $\pi_*(N)$-determined auto\m s.
  
  Consider next the \twocg\ $\SL(2m+1,\R)$ where $m=n-1$. 
  
  \noindent
  \underline{$\SL(2m+1,\R)$ has $N$-determined auto\m s}: We verify
  the conditions of \ref{ndetauto}. Let $L \subset \SL(2m+1,\R)$
  be an \lmntwo\ of rank $1$. The centralizer $C=C_{\SL(2m+1,\R)}(L)$
  is given in (\ref{eq:Lcent}). According to
  \S\ref{sec:red}, $C$ has $N$-determined auto\m s.  (Use the
  natural splitting of (\ref{eq:sloddsplit}) in
  connection with \ref{lemma:autononcon}.) See
  \ref{cor:limislodd} for the vanishing of the higher
  limits.

  \noindent
  \underline{$\SL(2m+1,\R)$ is $N$-determined}: Conditions (1) and (2)
  of \ref{indstepalt} are verified in
  \ref{sec:slodd}.\S\ref{sec:sloddlim0}, condition (3) in
  \ref{sec:slodd}.\S\ref{sec:sloddnontoral2}, and condition (4) in
  \ref{cor:limislodd}.

  \noindent
  \underline{$\SL(2m+1,\R)$ has $\pi_*(N)$-determined auto\m s}: To
  prove this, it suffices to find a rank one \lmntwo\ $L \subset
  \SL(2m+1,\R)$ such that $C_{\SL(2m+1,\R)}(L)_0$ has
  $\pi_*(N)$-determined auto\m s and such that $C_{\SL(2m+1,\R)}(L)_0
  \to \SL(2m+1,\R)$ induces a mono\m\ on $H^1(W;\ch{T})$. Such a line
  is provided by $L=L[2m-1,2]$ with centralizer identity component
  $C_{\SL(2m+1,\R)}(L)_0=\SL(2m-1,\R) \times \SL(2,\R)$; see the proof
  of \ref{lemma:sloddlimH1W0}.

  Consider finally the \twocg\ $\PGL(n,\Ha)$ for $n \geq 3$.

  \noindent
  \underline{$\PGL(n,\Ha)$ has $N$-determined auto\m s}: We verify the
  conditions of \ref{ndetauto}. Let $L \subset \PGL(n,\Ha)$ be
  an \lmntwo\ of rank $1$. The centralizer $C=C_{\PGL(n,\Ha)}(L)$ is
  given in (\ref{eq:cfamCL}) and
  (\ref{eq:cfamCI}). According to the general results of
  \S\ref{sec:red}, $C$ has $N$-determined auto\m s and
  according to \ref{lemma:cfamlim=0} the higher limits
  vanish.

  \noindent
  \underline{$\PGL(n,\Ha)$ is $N$-determined}: Note that $\PGL(3,\Ha)$
  satisfies condition (1 \& 2) of \ref{chp:di4andf4}.\S\ref{sec:di4}
  so that we may apply the same variant of
  \ref{indstepalt} used for $\di$
  (\ref{chp:di4andf4}.\S\ref{sec:di4}). When $n>3$, conditions (1) and
  (2) of \ref{indstepalt} follow if we can verify that
  the conditions of \ref{lemma:altcond2} are satisfied.
  That the centralizer $C_{\PGL(n,\Ha)}(V)$
  (\ref{eq:cfamCL}, \ref{eq:cfamCI},
  \ref{eq:cfamCP}, \ref{eq:cfamCIL}),
  where $V$ is an \lmntwo\ of rank at most two, satisfies the
  conditions of \ref{lemma:altcond2} is a consequence
  of the general results of \S\ref{sec:red} and
  \ref{prop:cfamlimH1W},
  \ref{prop:cfamlim0H1}. See
  \ref{sec:cfamoddnont} and
  \ref{sec:cfamnontrank2} for condition (3) and
  \ref{lemma:cfamlim=0} for condition (4) of
  \ref{indstepalt}.

  \noindent
  \underline{$\PGL(n,\Ha)$ has $\pi_*(N)$-determined auto\m s}: We
  only need to consider the cases $n=3$ and $n=4$ as
  $H^1(W;\ch{T})(\PGL(n,\Ha))=0$ for $n>4$ \cite{hms:first}.  In those
  two cases, it suffices, as above, to find a rank one \lmntwo\ $L
  \subset \PGL(n,\Ha)$ such that $C_{\PGL(n,\Ha)}(L)_0$ has
  $\pi_*(N)$-determined auto\m s and such that $C_{\PGL(n,\Ha))}(L)_0
  \to \PGL(n,\Ha)$ induces a mono\m\ on $H^1(W;\ch{T})$. Such a line
  is provided by $L=I$ for which
  $C_{\PGL(n,\Ha)}(I)_0=\GL(n,\C)/\gen{-E}$
  (\ref{eq:cfamCI}).
  
  Since $\PSL(2n,\R)$, $n \geq 4$, is uniquely $N$-determined and has
  a split \mtn , its auto\m\ group is isomorphic to $W \backslash
  N_{\GL(L)}(W)$ by \ref{lemma:autX}. When $n=4$, the group,
  $\Out_{\tr}(W)$, to the right in the exact sequence (\ref{autwNgl})
  is the permutation group $\Sigma_3$. There are Lie group outer
  auto\m s inducing $\Sigma_3$. When $n>4$,
  \begin{multline*}
    \Aut(\PSL(2n,\R)) \cong W \backslash N_{\GL(L)}(W)  = W
    \backslash \gen{\Z_2^*, W(\PGL(2n,\R))} = W \backslash
    \gen{\Z_2^{\times},W,c_1} \\ = (W \cap \gen{\Z_2^{\times},c_1})
    \backslash \gen{\Z_2^{\times},c_1} = 
    \begin{cases}
      \gen{-c_1} \backslash  \gen{\Z_2^{\times},c_1} = \Z_2^{\times} &
      \text{$n$ odd} \\
       \gen{-1} \backslash  \gen{\Z_2^{\times},c_1} = \Z^{\times}
       \backslash \Z_2^{\times} \times \gen{c_1} &
      \text{$n$ even}
    \end{cases}
  \end{multline*}  
   %% /home/moller/manus/2cgs/magma/aut/out.prg and outPGLnH
 %%%%%%%%%%%%%%%%%%%%%%%%%%%%%%%%%%%%%%%%%%%%%%
  %(\ref{lemma:autX}) isomorphic to the middle group, $W \backslash
%  N_{\GL(L)}(W)$, of the exact sequence (\ref{autwNgl}). The group
%  $\Out_{\tr}(W)$ is the permutation group $\Sigma_3$ when $n=4$, is
%  of order two when $n>4$ is even, and is trivial when $n>4$ is odd.
%  In the cases where $\Out_{\tr}(W)$ is nontrivial, its elements are
%  induced from Lie group auto\m s of $\PSL(2n,\R)$. The center of $W$
%  is $\{\pm 1\}$ when $n$ is even and is trivial when $n$ is odd. We
%  can now read off the structure of $\Aut(\PSL(2n,\R))$ from
%  (\ref{autwNgl}). (Conjugation with the element $c_1$
%  (\ref{eq:dfamIc}), or any other matrix of negative determinant, is
%  an auto\m\ of $\PSL(2n,\R)$. When $n$ is odd, conjugation with $c_1$
%  in $\Out(N)$ equals conjugation with $c_1c_2 \cdots c_n$ which is
%  the unstable Adams operation $\psi^{-1}$.  When $n$ is even,
%  conjugation with $c_1$ induces the diagonal matrix
%  $\diag(-1,1,\ldots,1)$ on $L=\pi_2(BT)$, and $\psi^{-1}$,
%  conjugation with $c_1c_2 \cdots c_n$, is the trivial element of
%  $\Out(N)$.)

  Similarly, 
  \begin{equation*}
   \Aut(\SL(2n+1,\R)) \cong  W \backslash N_{\GL(L)}(W) = W \backslash
   \gen{\Z_2^{\times},W} = (W \cap \Z_2^{\times}) \backslash \Z_2^{\times}
    = \Z^{\times} \backslash \Z_2^{\times}
  \end{equation*}
  for $n \geq 2$ by \ref{lemma:autX}.
  %isomorphic to the middle group, $W \backslash N_{\GL(L)}(W)$, of the
%  exact sequence (\ref{autwNgl}).  Since there are no trace preserving
%  outer auto\m s of the Weyl group
%  $W(\SL(2n+1,\R))=\Sigma_2\wr\Sigma_n$, we conclude that
%  $\Aut(\SL(2n+1,\R)) \cong \Z^{\times} \backslash \Z_2^{\times}$.
  %% /home/moller/manus/2cgs/magma/aut/out.prg

  The auto\m\ group $\Aut(\PGL(n,\Ha))$, $n \geq 3$, is
  (\ref{lemma:autX}) contained in $W \backslash N_{\GL(L)}(W) \cong
  \Z^{\times} \backslash \Z_2^{\times}$. Since $H^2(W;\ch{T})$ is an
  \lmntwo\ \cite{matthey:second}, it is isomorphic to the second
  cohomology group $H^2(W;t(\PGL(n,\Ha)))$ with coefficient module
  $t(\PGL(n,\Ha))$, the maximal \lmntwo\ in the \mt . The unstable
  Adams operations with index in $\Z_2^{\times}$ act trivially here
  since they act as coefficient group auto\m s. Thus all elements of $
  W \backslash N_{\GL(L)}(W)$ preserve the extension class $e \in
  H^2(W;\ch{T})$ and we conclude that $\Aut(\PGL(n,\Ha)) \cong
  \Z^{\times} \backslash \Z_2^{\times}$. 
\end{proof}

\begin{proof}[Proof of Corollary~\ref{cor:glnR}]
  Note first that $\GL(n,\R)$ is LHS for all $n \geq 1$. If $n$ is
  odd, $\GL(n,\R)=\SL(n,\R) \times \gen{-E}$ is LHS because its Weyl
  group is the direct product of the Weyl group of the identity
  component with the component group. If $n$ is even, see
  \ref{exmp:sl2C}.(5).  According to \ref{lemma:autononcon} and
  \ref{redtoconnected}, $\GL(n,\R)$ is totally $N$-determined.
 
  If $n$ is odd, the identity component has trivial center, so that
  $\Aut(\GL(n,\R))= \Aut(\SL(n,\R)) = \Z^{\times} \backslash
  \Z_2^{\times}$ by the \ses\ \cite[5.2]{jmm:ext}.  
   
  Suppose next that $n=2m$ is even. When $m=1$, $\Aut(\GL(2,\R)) =
  \Aut(\Z/2,\Z/2^{\infty},0)=\Aut(\Z/2^{\infty})=\Z_2^{\times}$
  according to (\ref{eq:H1TW}).  When $m>1$,
  $H^1(\pi;\ch{Z}(\SL(2m,\R))) = H^1(\pi;\gen{-E})$ is the order two
  subgroup $\gen{\delta}$ of $\Aut(\GL(2m,R))$ generated by the group
  iso\m\ $\delta(A) =(\det A)A$, $A \in \GL(2m,\R)$, and
  $H^1(W;\ch{T}) = \Hom(W_{\textmd{ab}},\gen{-E}) = \Z/2 \times \Z/2$
  (for $m>2$) \cite{hms:first,matthey:second} is the middle group of
  an exact sequence
    \begin{equation*}
      0 \to H^1(\pi;\gen{-E}) \to H^1(W;\ch{T}) \to H^1(W_0;\ch{T})
      \to 0
    \end{equation*}
    because $\GL(2m,\R)$ is LHS. (For $m=2$, $H^1(W_0;\ch{T})=0$ and
    $H^1(\pi;\ch{Z}(\SL(2m,\R)))=\Z/2$, though.) In the exact sequence
    (\ref{eq:H1TW}) for the auto\m\ group of
    $N=N(\GL(2m,\R))=N(\SL(2m+1,\R))$, the group on the right hand side
    is $\Aut(W,\ch{T},0)=\gen{W,\Z_2^{\times}}$ as for $\SL(2m+1,\R)$.
    Thus $\Aut(N)$ is generated by $H^1(W;\ch{T})$, $W$, and
    $\Z_2^{\times}$ so that $\Aut(N,N_0)=\Aut(N)$ as $W_0$ is normal
    in $W$.  Note that these three subgroups of $\Aut(N,N_0)$
    commute because of the special form of the elements of
    $H^1(W;\ch{T}) =  \Hom(W_{\textmd{ab}},\gen{-E})$. Hence
    \begin{multline*}
      \frac{\Aut(N,N_0)}{W_0} 
      = \frac{\gen{H^1(W;\ch{T}), W, \Z_2^{\times}}}{W_0} 
      =  \frac{\gen{H^1(W;\ch{T}), W_0, c_1, \Z_2^{\times}}}{W_0}
      =  \frac{\gen{H^1(W;\ch{T}), c_1, \Z_2^{\times}}}{W_0 \cap
        \gen{H^1(W;\ch{T}), c_1, \Z_2^{\times}} } \\
      = 
      \begin{cases}
        \frac{\gen{H^1(W;\ch{T}), c_1, \Z_2^{\times}}}{\gen{-c_1}} 
        = H^1(W;\ch{T}) \times \Z_2^{\times} &
        \text{$m$ odd} \\
         \frac{\gen{H^1(W;\ch{T}), c_1, \Z_2^{\times}}}{\gen{-1}} 
         = H^1(W;\ch{T}) \times \gen{c_1} \times 
           \Z^{\times} \backslash \Z_2^{\times} &
        \text{$m$ even}
      \end{cases}
    \end{multline*}
    According to \ref{lemma:autXnoncon}, the auto\m\ group
    $\Aut(\GL(2m,\R))$ is a subgroup of the above group and
    \begin{equation*}
      \Aut(\GL(2m,\R)) = 
      \begin{cases}
        \gen{\delta} \times \Z_2^{\times} & \text{$m$ odd} \\
        \gen{\delta} \times \gen{c_1} \times 
           \Z^{\times} \backslash\Z_2^{\times} & \text{$m$ even} 
      \end{cases}
    \end{equation*}
    for $m>1$.
%the Lie group auto\m\ 
%  $\delta(A)=\det(A)A$, $A \in \GL(2m,\R)$, represents the nontrivial
%  element of the left term $H^1(\pi;\ch{Z}(\SL(2m,\R)))=\Z/2$ of the
%  exact sequence \cite[5.2]{jmm:ext}. Using \ref{lemma:autXnoncon} we
%  see that
%  \begin{align*}
%    &\Aut(N,N_0) = \Aut(N) = \gen{H^1(W;\ch{T}),\Z_2^{\times},W}
%    = \gen{H^1(W;\ch{T}),\Z_2^{\times},W_0,c_1} \\
%    &\Aut(\GL(2m,\R)) = W_0 \backslash \gen{\delta,
%      \Z_2^{\times}, W_0, c_1} = \gen{\delta} \times \Aut(\SL(2m,\R))  
%     \\
%    &\Out(\GL(n,\R)) = W \backslash \gen{\delta,
%      \Z_2^{\times},W} = \gen{\delta} \times \Z^{\times} \backslash
%    \Z_2^{\times} 
%  \end{align*}
%  commute and $H^1(W;\ch{T})=\Hom(W_{\textmd{ab}};Z(\GL(n,\R))=\Z/2
%  \times \Z/2$.
\end{proof}

\section{Proof of Theorem~\ref{thm:conclusion}}

At this stage we know that $\di$ and $G$, for any compact, connected
simple, centerless Lie group $G$, are uniquely $N$-determined when
considered as \twocg s.

\begin{proof}[Proof of Theorem~\ref{thm:conclusion}]
  Let $X$ be a connected \twocg .  The splitting conjecture is true on
  the level of \mtn s \cite[1.12]{dw:normalizers} in the sense that
  $N(X)=N(G) \times N(\di)^t$.  From \ref{autondetcenter},
  \ref{prodauto}, \ref{ndetcenter}, and \ref{ndetprod} we know that $G
  \times \di^t$ is uniquely $N$-determined. In particular, $X$ and $G
  \times \di^t$ are isomorphic.  Let now $X$ be any \twocg , not
  necessarily connected. The remarks at the beginning of Section
  \ref{sec:red} in Chapter~\ref{cha:ndet} show that the
  $H^i$-injectivity, $i=1,2$, condition holds for $X$. Thus $X$ has
  $N$-determined auto\m s by \ref{lemma:autononcon} and $X$ is
  $N$-determined by \ref{redtoconnected} if we also assume that $X$ is
  LHS (\ref{LHSinit}).
\end{proof}

\chapter{Miscellaneous}
\label{sec:misc}

This chapter contains standard facts used at various places in this
paper. 

\section{Real representation theory}
\label{sec:realreps}

Real representations are semi-simple and determined by their
characters \cite[2.11, 3.12.(c)]{huppert:char}.  Any simple real
representation arises from a simple complex representation in the
following way: Let $\chi$ be the character of a simple complex
representation of a finite group $G$.  Then \cite[13.1, 13.11,
13.12]{huppert:char}
\begin{description}
\item[$\chi \neq \overline{\chi}, \varepsilon_2(\chi)=\phantom{+} 0\,$]
  $\psi=\chi+\overline{\chi}$ is the character of a simple $\R$-module
  of complex type.  
  %$\chi$ is not the character of a real representation.
\item[$\chi=\overline{\chi}, \varepsilon_2(\chi)=+1\,$] $\chi$ is the
  character of a simple $\R$-module of real type.
\item[$\chi=\overline{\chi}, \varepsilon_2(\chi)=-1\,$] $\psi=2\chi$ is the
  character of simple $\R$-module of quaternion type.
 %$\chi$ is not the character of a real representation.
\end{description}
where $\varepsilon_2(\chi)=\frac{1}{\order{G}} \sum_{g \in
  G}\chi(g^2)$. 
\begin{exmp}
  (1) The character table of the cyclic group $C_4$ of order $4$ 
\begin{center}
  \begin{tabular}{|r||r|r|r|r|r|} \hline
   $C_4$ & $\varepsilon_2$ & $1$ & $-1$ & $i$ & $-i$ \\ \hline\hline
   $\chi_1$ & $+$ & $1$ & $1$ & $1$ & $1$ \\ \hline
   $\chi_2$ & $+$ & $1$ & $1$ & $-1$ & $-1$ \\ \hline
   $\chi_3$ & $0$ & $1$ & $-1$ & $i$ & $-i$ \\ \hline
   $\chi_4$ & $0$ & $1$ & $-1$ & $-i$ & $i$ \\ \hline
  \end{tabular}
\end{center}
shows that there are two linear real representations and one
$2$-dimensional simple real faithful representation of complex type with
character $\psi=\chi_3+\chi_4=(2,-2,0,0)$.

\noindent (2)
The character table of the dihedral group $D_8=2^{1+2}_+$
\begin{center}
  \begin{tabular}{|r||r|r|r|r|r|r|} \hline
   $D_8$ & $\varepsilon_2$ & $1$ & $-1$ & $R_1$ & $R_2$ & $i$ \\ \hline\hline
   $\chi_1$ & $+$ & $1$ & $1$ & $1$ & $1$ & $1$ \\ \hline
   $\chi_2$ & $+$ & $1$ & $1$ & $-1$ & $1$ & $-1$\\ \hline
   $\chi_3$ & $+$ & $1$ & $1$ & $1$ & $-1$ & $-1$\\ \hline
   $\chi_4$ & $+$ & $1$ & $1$ & $-1$ & $-1$ & $1$\\ \hline
   $\chi_5$ & $+$ & $2$ & $-2$ & $0$ & $0$ & $0$\\ \hline
  \end{tabular}
\end{center}
shows that there are four linear real representations and one
$2$-dimensional simple real faithful representation of real type with
character $\chi_5=(2,-2,0,0,0)$.
%%\marginpar{Check table}

\noindent (3)
The character table of the quaternion group $Q_8=2^{1+2}_-$ (identical
to the one for $D_8$ except for one value of $\varepsilon_2$)
\begin{center}
  \begin{tabular}{|r||r|r|r|r|r|r|} \hline
   $Q_8$ & $\varepsilon_2$ & $1$ & $-1$ & $k$ & $j$ & $i$ \\ \hline\hline
   $\chi_1$ & $+$ & $1$ & $1$ & $1$ & $1$ & $1$ \\ \hline
   $\chi_2$ & $+$ & $1$ & $1$ & $-1$ & $1$ & $-1$\\ \hline
   $\chi_3$ & $+$ & $1$ & $1$ & $1$ & $-1$ & $-1$\\ \hline
   $\chi_4$ & $+$ & $1$ & $1$ & $-1$ & $-1$ & $1$\\ \hline
   $\chi_5$ & $-$ & $2$ & $-2$ & $0$ & $0$ & $0$\\ \hline
  \end{tabular}
\end{center}
shows that there are four linear real representations and one
$4$-dimensional simple real faithful representation of quaternion type
with character $\psi=2\chi_5=(4,-4,0,0,0)$.
\end{exmp}
We are interested in  real oriented representations, i.e.\ homo\m s of
finite groups into the special linear group $\SL(2n,\R)$ (as opposed
to homo\m s into the general linear group $\GL(2n,\R)$). The outer
auto\m\ of $\SL(2n,\R)$ is conjugation by any orientation reversing
matrix such as $D=\diag(-1,1,\ldots,1)$.
%%\marginpar{$\GL^+(2n,\R)$}

\begin{lemma}\label{lemma:intoSL}
  Let $V \subset \pslr{2n}$ be an object of $\A(\pslr{2n})$ and $G=V^*
  \subset \SL(2n,\R)$ its inverse image in $\SL(2n,\R)$ as in
  \ref{lemma:Vast}.  Then
  \begin{equation*}
    \text{$V$ and $V^D$ are nonisomorphic objects of $\A(\pslr{2n})$} 
    \iff  N_{\GL(2n,\R)}(G) \subset \SL(2n,\R)
  \end{equation*}
%  \begin{align*}
%    \text{$V$ and $V^D$ are nonidentical objects of $\A(\pslr{2n})$}
%    &\iff 
%    C_{\GL(2n,\R)}(\phi) \subset \SL(2n,\R) \\
%    \text{$V$ and $V^D$ are nonisomorphic objects of $\A(\pslr{2n})$}
%    &\iff 
%    N_{\GL(2n,\R)}(\phi) \subset \SL(2n,\R)
%  \end{align*}
 % Let $G$ be a finite group and \func{\phi}{G}{\SL(2n,\R)}  a degree
%  $2n$ real oriented representation of $G$. Then
%  \begin{equation*}
%    \phi^D \not\in \phi^{\SL(2n,\R)} \iff
%    C_{\GL(2n,\R)}(\phi) \subset \SL(2n,\R)
%  \end{equation*}
 % \begin{enumerate}
%  \item If $C_{\GL(2n,\R)}(\phi) \subset \SL(2n,\R)$, then the
%    $\GL(2n,\R)$-conjugacy class of $\phi$ consists of the two
%    $\SL(2n,\R)$-conjugacy classes represented by $\phi$ and $\phi^D$.
%  \item If $C_{\GL(2n,\R)}(\phi) \not\subset \SL(2n,\R)$, then the
%    $\GL(2n,\R)$-conjugacy class of $\phi$ consists of a single
%    $\SL(2n,\R)$-conjugacy class.
%  \end{enumerate}
\end{lemma}
\begin{proof}
We note that
\begin{align*}
  \text{$V$, $V^D$ are isomorphic objects of $\A(\pslr{2n})$} &\iff 
  \text{$G$, $G^D$ are conjugate subgroups of $\SL(2n,\R)$} \\
  &\iff G \in G^{D\SL(2n,\R)} \\
  &\iff N_{\GL(2n,\R)}(G) \cap D\SL(2n,\R) \neq \emptyset \\
  &\iff N_{\GL(2n,\R)}(G) \not\subset \SL(2n,\R)
\end{align*}
for any nontrivial elementary abelian $2$-group $V \subset
\pslr{2n}$. 
%%%%%%%%%%%%%%%%%%%
 % In order to prove the first part note that
%  \begin{align*}
%     \text{$V$ and $V^D$ are identical objects of $\A(\pslr{2n})$}
%     &\iff
%     \text{$\phi$ and $\phi^D$ are conjugate in $\SL(2n,\R)$} \\
%     &\iff \phi^D \in \phi^{\SL(2n,\R)} \\
%     &\iff D \in \SL(2n,\R)C_{\GL(2n,\R)}(\phi) \\
%     &\iff C_{\GL(2n,\R)}(\phi) \not\subset \SL(2n,\R)
%  \end{align*}
%  The proof of the second part  if somewhat similar.
%%%%%%%%%%%%%%%%%%%%%%%%%%%%%%%%%%%%%%%%%%%%%%%%%%%%%%%%%
  %The $\GL(2n,\R)$-conjugacy class of $\phi$ \marginpar{$\SL$ or $\GL^+$?}
%  \begin{equation*}
%    \phi^{\GL(2n,\R)} = \phi^{\SL(2n,\R) \cup D\SL(2n,\R)} =
%     \phi^{\SL(2n,\R)} \cup   \phi^{D\SL(2n,\R)}
%  \end{equation*}
%  consists of at most two $\SL(2n,\R)$-conjugacy classes represented
%  by $\phi$ and $\phi^D$. Now,
%  \begin{multline*}
%    \phi^{\SL(2n,\R)} =  \phi^{D\SL(2n,\R)} \iff
%     \phi^{\SL(2n,\R)} \ni \phi^D \iff
%     C_{\GL(2n,\R)}(\phi) \SL(2n,\R) \ni D \\ \iff
%    C_{\GL(2n,\R)}(\phi)\not\subset \SL(2n,\R) 
%  \end{multline*}
%  meaning that these two $\SL(2n,\R)$-conjugacy classes coincide iff
%  $C_{\GL(2n,\R)}(\phi)\not\subset \SL(2n,\R)$.
\end{proof}
%\begin{cor}\label{cor:realphiphiD}
%  Let $G$ be a finite group and \func{\phi}{G}{\SL(2n,\R)}  a degree
%  $2n$ real oriented representation of $G$. If the character of $\phi$
%  contains a real irreducible representation of real type then $\phi$
%  and $\phi^D$ are conjugate as real oriented representations.
%\end{cor}
%\begin{proof}
%  The centralizer $C_{\GL(2n,\R)}(\phi)$ contains a factor of the form
%  $\GL(i,\R)$, $i \geq 1$, so it is not contained in $\SL(2n,\R)$.
%\end{proof}
For instance, all representations of \lmn s are conjugate in $\SL(2n,\R)$ if
and only if they are conjugate in $\GL(2n,\R)$ - as in
\ref{cor:equivcat}. 

Let $\A(\GL(2n,\R))(G)$ be the subgroup of $\Out(G)$ consisting of all
outer auto\m s of $G$ induced by conjugation with some element of
$\GL(2n,\R)$ \cite[5.8]{jmm:ndet} (ie $\A(\GL(2n,\R))(G)$ is the group
$\Out_{\tr}(G)$ of all trace preserving outer auto\m s of $G$)
and $\A(\SL(2n,\R))(G)$ the subgroup of $\Out(G)$ consisting of all
outer auto\m s of $G$ induced by conjugation with some element of
$\SL(2n,\R)$. Since
\set\begin{equation}\label{eq:NOut}
  N_{\GL(2n,\R)}(G)/G C_{\GL(2n,\R)}(G) \xrightarrow{\cong}
  \A(\GL(2n,\R))(G) 
\end{equation}\add
we conclude from \ref{lemma:intoSL} that 
\begin{multline*}
  \text{$G$, $G^D$ are nonconjugate subgroups of $\SL(2n,\R)$} \iff 
  N_{\GL(2n,\R)}(G) \subset \SL(2n,\R) \\
  \iff \begin{cases} C_{\GL(2n,\R)}(G) \subset \SL(2n,\R) \\
                \A(\SL(2n,\R))(G) = \A(\GL(2n,\R))(G) \end{cases}
\end{multline*}

Let $V$ and $E$ be objects of $\A(\pslr{2n})$ such that $\dim V +1
=\dim E$. If there are \m s $V \to E$ and $V^D \to E$, then (some
representative of) $E=\gen{V,V^D}$ is generated by the images of (some
representatives of) $V$ and $V^D$ so that $E=E^D$. Conversely, if
$E=E^D$ and there is \m\ $V \to E$ then there is also a \m\ $V^D \to
E^D=E$. 

We have $2(\phi^D)=2\phi \in \Rep(G,\SL(4n,\R))$ for any oriented real
degree $2n$ representation $\phi \in \Rep(G,\SL(2n,\R))$ as the
conjugating matrix $2D$ is orientation preserving.

%\begin{lemma}\label{lemma:intoSL}
%  Let \func{\phi_1,\phi_2}{G}{\SL(2n,\R)} be two  real oriented
%  representations of a finite group $G$. Assume that $\phi_1$ and
%  $\phi_2$ are conjugate in $\GL(2n,\R)$. Then $\phi_1$ is conjugate
%  in $\SL(2n,\R)$ to $\phi_2$ or to $\phi_2^D$.
%\end{lemma} 
%\begin{proof} 
%  Since $\phi_1$ and $\phi_2$ are conjugate in $\GL(2n,\R)$, there is
%  a matrix $A \in \GL(2n,\R)$ such that
%  \begin{equation*}
%    \phi_1 = \phi_2^A  = (\phi_2^{D})^{D^{-1}A}
%  \end{equation*}
%  where either $A$ or $D^{-1}A$ preserves orientation.
%\end{proof} 

%\begin{lemma}
%  \label{lemma:phiandphiD}
%  Let \func{\varphi}{G}{\SL(2n,\R)} be a real orientable degree $2n$
%  representation. Then
%  \begin{equation*}
%    \text{$\varphi$ and $\varphi^D$ are conjugate in $\SL(2n,\R)$}
%    \iff 
%    C_{\GL(2n,\R)}(\varphi) \not\subset \SL(2n,\R)
%  \end{equation*}
%\end{lemma}
%\begin{proof}
%  If $C_{\GL(2n,\R)}(\varphi)$ contains an orientation reversing
%  matrix $B$, then $\varphi = \varphi^B = (\varphi^D)^{D^{-1}B}$ where
%  $D^{-1}B$ preserves the orientation. Conversely, if $\varphi$ and
%  $\varphi^D$ are conjugate in $\SL(2n,\R)$, then $\varphi =
%  (\varphi^D)^A = \varphi^{DA}$ for some orientation preserving matrix
%  $A$.
%\end{proof}
%%%%%%%%%%%%%%  Fri Aug 20 14:47:55 CEST 2004 %%%%%%%%%%%%%%%%%%%%%

\begin{exmp}\label{exmp:IDQ}
  \begin{simplelist}
 \item \label{exmp:IDQ1}
Let $G \subset \SL(2d,\R)$ be a finite group making $\R^{2d}$ a simple
$\R G$-module of complex type. Consider the image of $G \subset
\SL(2nd,\R)$ of $G$ under the $n$-fold diagonal $\SL(2d,\R)
\xrightarrow{\Delta_n} \SL(2dn,\R)$. The centralizer
$C_{\GL(2nd,\R)}(G)=\GL(n,\C)$ is connected, hence contained in
$\SL(2dn,\R)$. Since $C_{\GL(2d,\R)}(G)=\GL(1,\C)$, the elements of
$G$ commute with $i \in \GL(1,\C)$ and we may factor the inclusion of
$G$ into $ \SL(2dn,\R)$ as
\begin{equation*}
  G \to C_{\GL(2d,\R)}(i)=\GL(d,\C) 
  \xrightarrow{\Delta_n} \GL(dn,\C) \to \SL(2dn,\R)
\end{equation*}
Let $\chi$ be the character for $G$ in $\GL(d,\C)$ so that the
character for $G$ in $\GL(2d,\R)$ is $\chi+\overline{\chi}$. There are
inclusions
\begin{equation*}
  \xymatrix{
    {\A(\GL(2dn,\R))(G)} \ar@{=}[r] &
    {\A(\GL(2d,\R))(G)} \ar@{=}[r] &
    {\Out_{\chi+\overline{\chi}}(G)} \\
    {\A(\SL(2dn,\R))(G)} \ar@{^(->}[u] &
    {\A(\GL(d,\C))(G)} \ar@{_(->}[l]  \ar@{^(->}[u] \ar@{=}[r] &
    {\Out_{\chi}(G)} \ar@{^(->}[u] }
\end{equation*}
where ${\Out_{\phi}(G)}$ is the group of all outer auto\m s that
respect the function $\phi$.

\item  \label{exmp:IDQ2}
Let $G \subset \GL(d,\R)$ be a finite group making $\R^{d}$ a simple
$\R G$-module of real type. Consider the image of $G \subset
\SL(2nd,\R)$ of $G$ under the $2n$-fold diagonal $\GL(d,\R)
\xrightarrow{\Delta_{2n}} \SL(2dn,\R)$. The centralizer
$C_{\GL(2nd,\R)}(G)=\GL(2n,\R)$ is contained in $\SL(2dn,\R)$ when $d$
is even. We my factor the inclusion of $G$ into $\SL(2dn,\R)$ as
\begin{equation*}
  G \to \GL(d,\R) \to \GL(d,\C) \xrightarrow{\Delta_n} \GL(nd,\C) \to
  \SL(2nd,\R) 
\end{equation*}
and as the trace functions for $G$ in $\GL(d,\C)$ and $\GL(2nd,\R)$
are proportional $\A(\GL(2nd,\R))(G)=\A(\GL(d,\C))(G) \subset
\A(\SL(2dn,\R))(G)$. Hence $G \neq G^D$ in $\SL(2nd,\R)$. 

\item  \label{exmp:IDQ3}
Let $G \subset \SL(4d,\R)$ be a finite group making $\R^{4d}$ a simple
$\R G$-module of quaternion type. Consider the image of $G \subset
\SL(4nd,\R)$ of $G$ under the $n$-fold diagonal $\SL(4d,\R)
\xrightarrow{\Delta_{n}} \SL(4dn,\R)$. The centralizer
$C_{\GL(4dn,\R)}(G)=\GL(n,\Ha)$ is connected so it is contained in
$\SL(4dn,\R)$. Since $C_{\GL(4d,\R)}(G)=\GL(1,\Ha)\subset\GL(2,\C)$
the elements of $G$ commute with $i\in\GL(2,\C)$ and we may factor the
inclusion of $G$ into $\SL(4dn,\R)$ as
\begin{equation*}
  G \to C_{\GL(4d,\R)}(i)=\GL(2d,\C) \xrightarrow{\Delta_n}
  \GL(2nd,\C) \to \SL(4nd,\R)
\end{equation*}
and as the trace functions for $G$ in $\GL(2d,\C)$ and $\GL(4nd,\R)$
are proportional  $\A(\GL(4nd,\R))(G)=\A(\GL(2d,\C))(G) \subset
\A(\SL(4dn,\R))(G)$. Hence $G \neq G^D$ in $\SL(4nd,\R)$.

\item \label{exmp:IDQ4} $\R^2$ is a simple $\R C_4$-module of complex
  type with respect to the group
  \begin{equation*}
    C_4=\gen{I} \subset \SL(2,\R), \qquad 
      I=\begin{pmatrix}
        0 & -1 \\ 1 & 0
      \end{pmatrix} 
  \end{equation*}
  Consider the image of $C_4$ in $\SL(2n,\R)$ under the $n$-fold
  diagonal. The Quillen auto\m\ group
  $\A(\GL(2n,\R))(C_4)=\A(\GL(2,\R))(C_4)=\Out(C_4)$ since the trace
  lives on $\mho_1(C_4)=\gen{-E}$ only. However,
  \begin{equation*}
    \A(\SL(2n,\R))(C_4)=
    \begin{cases}
      \Out(C_4) & \text{$n$ even} \\
      \{1\}     & \text{$n$ odd} 
    \end{cases}
  \end{equation*}
  so that $C_4 \neq C_4^D \iff \text{$n$ even}$.

\item \label{exmp:IDQ5} $\R^4$ is a simple $\R G_{16}$-module of
  complex type with respect to the group
  \begin{equation*}
    G_{16}=4\circ 2^{1+2}_{\pm}=\gen{
      \begin{pmatrix}
        R & 0 \\ 0 & R
      \end{pmatrix},
      \begin{pmatrix}
        T & 0 \\ 0 & T
      \end{pmatrix},
      \begin{pmatrix}
        0 & -E \\ E & 0
      \end{pmatrix} } \subset \SL(4,\R), \quad
  R=
  \begin{pmatrix}
    1 & 0 \\ 0 & -1
  \end{pmatrix}, \quad
  T=
  \begin{pmatrix}
    0 & 1 \\ 1 & 0
  \end{pmatrix}
  \end{equation*}
  Consider the image of $G_{16}$ in $\SL(4n,\R)$ under the $n$-fold
  diagonal. The Quillen auto\m\ group
  $\A(\GL(4n,\R))(G_{16})=\A(\GL(4,\R))(G_{16})=\Out(G_{16}) \cong
  \Out(C_4) \times \Symp(2,\F_2)$ since the trace
  lives on the derived group $[G_{16},G_{16}]=\gen{-E}$ only. In fact,
  $\A(\GL(2,\C))(G_{16})$ is the factor $\Symp(2,\F_2)$ and since the
  generator of the factor $\Out(C_4)$ is induced from conjugation with
  the matrix $
  \begin{pmatrix}
    0 & E \\ E & 0
  \end{pmatrix}$ of $\SL(4,\R)$, we see that also $\A(\SL(4,\R))(G_{16})
  \subset \A(\SL(4n,\R))(G_{16})$ is the full outer auto\m\ group of
  $G_{16}$. Hence $G_{16} \neq G_{16}^D$ in $\SL(4n,\R)$.

  \item \label{exmp:IDQ6} $\R^2$ is a simple $\R G$-module of 
  real type with respect to the group
  \begin{equation*}
    G = 2^{1+2}_+ = \gen{R,T} \subset \GL(2,\R)
  \end{equation*}
  Consider the image of $G$ in $\SL(4n,\R)$ under the $2n$-fold
  diagonal map. Then $G \neq G^D$ in $\SL(4n,\R)$.

 \item \label{exmp:IDQ7} $\R^4$ is a simple $\R G$-module 
 of quaternion type with respect to the group
  \begin{equation*}
    G = 2^{1+2}_- = \gen{
      \begin{pmatrix}
        0 & -R \\ R & 0
      \end{pmatrix},
      \begin{pmatrix}
        0 & -T \\ T & 0
      \end{pmatrix} } \subset \SL(4,\R)
  \end{equation*}
  Consider the image of $G$ in $\SL(4n,\R)$ under the $n$-fold
  diagonal map. Then $G \neq G^D$ in $\SL(4n,\R)$.
 \end{simplelist}
\end{exmp}

%  There are no embeddings of $D_8$ or $Q_8$ into $\slr{4n+2}$ with
%  central $\mho_1$.

  %%$D_8$: Such an embedding into $\GL(4n+2,\R)$ is $(2n+1)\psi$ which
%%does not factor through $\slr{4n+2}$. $Q_8$: Such embeddings are
%%sums of $\psi$ so they have degree divisible by $4$.

\setcounter{subsection}{\value{thm}}
\subsection{Representations of (generalized) extraspecial
  $2$-groups}\label{sec:xtraspecreps}\add
%%%%%%%%%%%%%%%%%%% Sun May 30 10:40:01 CEST 2004 %%%%%%%%%%%%%%%%%%%%%%
The extraspecial $2$-groups $G=2^{1+2d}_{\pm}$ have
\cite[7.5]{huppert:char} $2^d$ linear characters (that vanish on
$\mho_1(G)=G'=Z(G)=C_2$) and one simple complex character
\begin{equation*}
  \chi(g)=
  \begin{cases}
    0 & g \not\in Z(G) \\
    2^d\lambda(g) & g \in Z(G)
  \end{cases}
\end{equation*}
induced from the nontrivial linear character (group iso\m )
\func{\lambda}{Z(G)}{\{\pm 1\}}.

If $G=2^{1+2d}_+$ is of positive type, $\varepsilon_2(\chi)=+1$ and
$\chi\alpha=\chi$ for all $\alpha\in \Out(G)$, isomorphic to
$O^+(2d,2)$ \cite[III.13.9.b]{huppert:I}. This complex character is
also the character of the unique simple real representation $G \to
\GL(2^d,\R)$ which is of real type; when $d$ is even this
representation actually takes values in $\slr{2^d}$ but when $d$ is
odd this representation is not oriented. The unique
faithful real representation $G \to \GL(2 \cdot 2^d,\R)$ with central
$\mho_1$ has character $2\chi$ and it splits into two distinct
oriented real faithful representations $\psi, \psi^D \colon G \to
\SL(2 \cdot 2^d,\R)$ invariant under the action of $\Out(G)$
(\ref{exmp:IDQ}.(\ref{exmp:IDQ2})). 

If $G=2^{1+2d}_-$ is of negative type, $\varepsilon_2(\chi)=-1$ and
$\chi\alpha=\chi$ for all $\alpha\in \Out(G)$, isomorphic to
$O^-(2d,2)$ \cite[III.13.9.b)]{huppert:I}. The unique simple real
representation $G \to \GL(2 \cdot 2^d,\R)$ with character $2\chi$ is
of quaternion type. It splits into two distinct oriented
representations $\psi, \psi^D \colon G \to \SL(2 \cdot 2^d,\R)$
invariant under the action of $\Out(G)$
(\ref{exmp:IDQ}.(\ref{exmp:IDQ3})).

%%%%%%%%%%%%%%%% Sun May 30 10:40:01 CEST 2004 END %%%%%%%%%%%%%%%%%%%

The generalized extraspecial $2$-group $G=4 \circ 2^{1+2d}_{\pm}$ has
\cite[7.5]{huppert:char} $2^{1+d}$ linear characters (that vanish on
$\mho_1(G)=G'=C_2 \subsetneq Z(G)=C_4$) and two simple complex
characters
\begin{equation*}
  \chi(g)=
  \begin{cases}
    0 & g \not\in Z(G) \\
    2^d\lambda(g) & g \in Z(G)
  \end{cases}
\end{equation*}
induced from the two faithful linear characters
\func{\lambda}{Z(G)}{\gen{i}=C_4}. These two degree $2^d$ simple
characters, $\chi$ and $\overline{\chi}$, are interchanged by the
action of $\Out(G)= \Out(C_4) \times \Symp(2d,2)$ \cite{griess:autos}
(interchanged by the first factor $\Out(C_4)$ and preserved by the
second factor $\Symp(2d,2)$).  The unique simple real representation
$G \to \GL(2 \cdot 2^d,\R)$ has character $\chi+\overline{\chi}$ and
is of complex type as $\varepsilon_2(\chi)=0$. It splits up
into two distinct oriented representations
$\psi,\psi^D \colon G \to \SL(2 \cdot 2^d,\R)$ invariant under the
action of $\Out(G)$ (\ref{exmp:IDQ}.(\ref{exmp:IDQ1})). 

These irreducible faithful representations have easy explicit
constructions that we now explain.

Let $E$ be a nontrivial \lmntwo\ of rank $d\geq 1$ and $\R[E]$ its
real group algebra.  For $\zeta \in E^{\vee} = \Hom(E,\R^{\times})$
and $u \in v$, let $R_{\zeta},T_u \in \GL(\R[E])$ be the linear auto\m
s given by $R_{\zeta}(v)=\zeta(v)v$ and $T_u(v)=u+v$ for all $v \in
E$. The computation
\begin{equation*}
  R_{\zeta}T_u(v) = R_{\zeta}(u\cdot v) = \zeta(u)\zeta(v)(u\cdot v) =
  \zeta(u)T_u(\zeta(v)v) = \zeta(u)T_uR_{\zeta}(v)
\end{equation*}
shows that $R_{\zeta}T_u = \zeta(u)T_uR_{\zeta}$ or, equivalently,
$[R_{\zeta},T_u] = \zeta(u)$.

The group $2^{1+2d}_+=\gen{R_{\zeta},T_u} \subset \GL(\R[E]) \subset
\GL(\C[E]) \stackrel{\tau}{\subset} \SL(2^{d+1},\R)$ is extraspecial
and the quadratic form on its abelianization $2^{2d}$ is given by
\begin{equation*}
  q(x_1,\ldots ,x_d,y_1,\ldots ,y_d) = x_1y_1+\cdots+x_dy_d
\end{equation*}
because
\begin{equation*}
  (R_1^{x_1} \cdots R_d^{x_d}T_1^{y_1} \cdots T_d^{y_d})^2 = 
  \prod _{i=1}^d (R_i^{x_i}T_i^{y_i})^2 =\prod _{i=1}^d (-E)^{x_iy_i}  
\end{equation*}
where $T_1,\ldots ,T_d$ correspond to a basis of $E$, $R_1,\ldots,R_d$
correspond to the dual basis, and $x_i,y_i \in \{0,1\}=\F_2$.  This is
the unique faithful complex representation of degree $2^d$. It is also
the character of a simple real representation $G \to \GL(2^d,\R)$,
even $G \to \SL(2^d,\R)$ when $d$ is even, of real type.

The group
$2^{1+2d}_-=\gen{R_1,\ldots,R_{d-1},iR_d,T_1,\ldots,T_{d-1},iT_d}
\subset \GL(\C[E]) \stackrel{\tau}{\subset} \SL(2^{d+1},\R)$ is
extraspecial and the quadratic form on its abelianization $2^{2d}$ is
given by
\begin{equation*}
  q(x_1,\ldots ,x_d,y_1,\ldots ,y_d) = 
        x_1y_1+\cdots+x_{d-1}y_{d-1}+x_d^2+x_dy_d+y_d^2
\end{equation*}
because
\begin{equation*}
  (R_1^{x_1} \cdots R_{d-1}^{x_{d-1}}(iR_d)^{x_d}
           T_1^{y_1} \cdots T_{d-1}^{y_{d-1}}(iT_d)^{y_d})^2
  = (-E)^{x_d^2+y_d^2} (R_1^{x_1} \cdots
  R_d^{x_d}T_1^{y_1} \cdots T_d^{y_d})^2 
\end{equation*}
where $x_i,y_i\in\{0,1\}=\F_2$.  This is
the unique faithful complex representation of degree $2^d$.

The group $4 \circ 2^{1+2d}_{\pm} = 4 \circ 2^{1+2d}_+ =
\gen{i,R_{\zeta},T_u} = \gen{2^{1+2d}_+,2^{1+2d}_-} = 4 \circ
2^{1+2d}_-\subset \GL(\C[E]) \stackrel{\tau}{\subset} \SL(2^{d+1},\R)$
is generalized extraspecial with center $C_4 =\gen{i}$, derived group
$[4 \circ 2^{1+2d}_{\pm}, 4 \circ 2^{1+2d}_{\pm}]=\mho_1(4 \circ
2^{1+2d}_{\pm}) = C_2 \subset C_4=Z(4 \circ 2^{1+2d}_{\pm})$, and
elementary abelian abelianization $2 \times 2^{2d}$.  The quadratic
form on its abelianization is given by
\begin{equation*}
   q(z,x_1,\ldots ,x_d,y_1,\ldots ,y_d) = z^2 + \sum_{i=1}^d x_iy_i
\end{equation*}
because
\begin{equation*}
  (i^zR_1^{x_1} \cdots R_d^{x_d}T_1^{y_1} \cdots T_d^{y_d})^2 = 
  (-E)^{z^2} \prod_{i=1}^d
   (R_i^{x_i}T_i^{y_i})^2 = (-E)^{z^2}(R_1^{x_1} \cdots
  R_d^{x_d}T_1^{y_1} \cdots T_d^{y_d})^2
\end{equation*}
where $z,x_i,y_i \in \{0,1\}$.  This representation and its conjugate
are the two faithful complex representation of degree $2^d$.

In the first two cases the associated symplectic inner product is 
\begin{equation*}
  [(x_1,\ldots,x_d,y_1,\ldots,y_d),(x'_1,\ldots
  ,x'_d,y'_1,\ldots,y'_d)] = \sum_{i=1}^d(x_iy'_i+x'_iy_i)
\end{equation*}
while it is
\begin{equation*}
   [(z,x_1,\ldots,x_d,y_1,\ldots,y_d),(z',x'_1,\ldots
  ,x'_d,y'_1,\ldots,y'_d)] = \sum_{i=1}^d(x_iy'_i+x'_iy_i)
\end{equation*}
in the last case.

\setcounter{subsection}{\value{thm}}
\subsection{Tensor products of real representations}
\label{sec:tensor}\add
Suppose that $\R^m$ is an $\R G$-module with trace $\chi$ and $\R^n$
an $\R H$ module with trace $\rho$. Consider $\R^{mn}=\R^m \otimes
\R^n$ as an $\R (G \times H)$-module in the usual way where
$(g,h)(u\otimes v)=gu\otimes hv$. The trace of this representation is
$\chi\#\rho(g,h)=\chi(g)\rho(h)$ and the determinant is
$\det(g,h)=(\det g)^n (\det h)^m$. This means that if 
\begin{itemize}
\item $m$ and $n$ are both even, or if
\item $m$ is even and $\R^m$ an oriented $G$-representation
\end{itemize}
then $\R^{mn}$ is a real oriented $G\times H$-representation.

\setcounter{subsection}{\value{thm}}
\subsection{Embedding $\GL(n,\C)$ in $\SL(2n,\R)$}
\label{sec:glnCtosl2nR}\add
Here are two embeddings \func{\tau}{\GL(n,\C)}{\GL^+(2n,\R)} with the
property that $\tr(\tau(A))=\tr(A)+\overline{\tr(A)}$ for all
$\A\in\GL(n,\C)$. 

If we write $\C^n=(\R+i\R)^n$, then
\begin{equation*}
  \GL(n,\C) \ni A + i B \xrightarrow{\tau}
  \left(\begin{pmatrix}
    a_{ij} & -b_{ij} \\ b_{ij} & a_{ij}
  \end{pmatrix}\right)_{1 \leq i,j \leq n} \in \GL^+(2n,\R)
\end{equation*}
In particular, $i \in \GL(n,\C)$ is sent to
$\diag(I,\cdots,I)\in\slr{2n}$ and $C_{\slr{2n}}(i)$ consists of
matrices with $2 \times 2$ blocks as above. For $(2 \times
2)$-matrices this embedding has the form
\begin{equation*}
  \GL(2,\C) \ni
  \begin{pmatrix}
    a_1+ia_2 & b_1+ib_2 \\
    c_1+ic_2 & d_1+id_2
  \end{pmatrix} \xrightarrow{\tau}
  \begin{pmatrix}
    \begin{pmatrix}
      a_1 & -a_2 \\ a_2 & a_1
    \end{pmatrix} &
    \begin{pmatrix}
      b_1 & -b_2 \\ b_2 & b_1
    \end{pmatrix} \\
    \begin{pmatrix}
      c_1 & -c_2 \\ c_2 & c_1
    \end{pmatrix} &
    \begin{pmatrix}
      d_1 & -d_2 \\ d_2 & d_1
    \end{pmatrix}
  \end{pmatrix} \in \slr{4}
\end{equation*}
and with this convention the six subgroups of $\slr{4}$ isomorphic
to $D_8$, $Q_8$, $G_{16}=4 \circ 2^{1+2}_{\pm}$ are
\begin{align*}
  D_8 &= \gen{
    \begin{pmatrix}
      E & 0 \\ 0 & -E
    \end{pmatrix},
    \begin{pmatrix}
      0 & E \\ E & 0
    \end{pmatrix} }, \quad &
    D_8^D&=\gen{
    \begin{pmatrix}
      E & 0 \\ 0 & -E
    \end{pmatrix},
    \begin{pmatrix}
      0 & -R \\ -R & 0 
    \end{pmatrix} } \\ 
  Q_8 &= \gen{
    \begin{pmatrix}
      I & 0 \\ 0 & -I
    \end{pmatrix},
    \begin{pmatrix}
      0 & I \\ I & 0
    \end{pmatrix} }, \quad &
  Q_8^D&=\gen{
    \begin{pmatrix}
      -I & 0 \\ 0 & -I
    \end{pmatrix},
    \begin{pmatrix}
      0 & T \\ -T & 0
    \end{pmatrix} } \\
  G_{16} &= \gen{
    \begin{pmatrix}
      E & 0 \\ 0 & -E
    \end{pmatrix},
    \begin{pmatrix}
      0 & E \\ E & 0
    \end{pmatrix},
    \begin{pmatrix}
      I & 0 \\ 0 & I
    \end{pmatrix}}, \quad &
  G_{16}^D &= \gen{
    \begin{pmatrix}
      E & 0 \\ 0 & -E
    \end{pmatrix},
    \begin{pmatrix}
      0 & -R \\ -R & 0
    \end{pmatrix},
    \begin{pmatrix}
      -I & 0 \\ 0 & I
    \end{pmatrix}}
\end{align*}
%where $R=
%\begin{pmatrix}
%  1 & 0 \\ 0 & -1
%\end{pmatrix}$, $T=
%\begin{pmatrix}
%  0 & 1 \\ 1 & 0
%\end{pmatrix}$, and $I=
%\begin{pmatrix}
%  0 & -1 \\ 1 & 0
%\end{pmatrix}=TR=-RT$. 

%%\marginpar{~/manus/2cgs/dfam/magma/nontoral/QD}

If we write $\C^n=\R^n+i\R^n$ then
\begin{equation*}
  \GL(n,\C) \ni A+iB \xrightarrow{\tau}
  \begin{pmatrix}
    A & -B \\ B & \phantom{-}A
  \end{pmatrix}
  \in \GL^+(2n,\R)
\end{equation*}
In particular, $i \in \GL(n,\C)$ is sent to
$
\begin{pmatrix}
  0 & -E \\ E & 0
\end{pmatrix} \in \slr{2n}
$ and $C_{\slr{2n}}(i)$ consists of block matrices of the form as
above.  For $(2 \times 2)$-matrices this embedding has the form
\begin{equation*}
   \GL(2,\C) \ni
  \begin{pmatrix}
    a_1+ia_2 & b_1+ib_2 \\
    c_1+ic_2 & d_1+id_2
  \end{pmatrix} \xrightarrow{\tau}
  \begin{pmatrix}
    \begin{pmatrix}
      a_1 & b_1 \\ c_1 & d_1
    \end{pmatrix} &
    -\begin{pmatrix}
      a_2 & b_2 \\ c_2 & d_2
    \end{pmatrix} \\
    \begin{pmatrix}
       a_2 & b_2 \\ c_2 & d_2
    \end{pmatrix} & 
    \phantom{-}\begin{pmatrix}
       a_1 & b_1 \\ c_1 & d_1
    \end{pmatrix}
  \end{pmatrix} \in \SL(4,\R)
\end{equation*}
and with this convention the six subgroups of $\slr{4}$ isomorphic
to $D_8$, $Q_8$,  $G_{16}=4 \circ 2^{1+2}_{\pm}$ are
\begin{align*}
  D_8&=\gen{
    \begin{pmatrix}
      R & 0 \\ 0 & R
    \end{pmatrix},
    \begin{pmatrix}
      T & 0 \\ 0 & T
    \end{pmatrix} }, \quad &
  D_8^D&=\gen{
    \begin{pmatrix}
      R & 0 \\ 0 & R
    \end{pmatrix},
    \begin{pmatrix}
      -T & 0 \\ 0 & T
    \end{pmatrix} } \\
  Q_8&= \gen{
    \begin{pmatrix}
      0 & -R \\R & 0
    \end{pmatrix},
    \begin{pmatrix}
      0 & -T \\ T & 0
    \end{pmatrix} }, \quad &
  Q_8^D&=\gen{
    \begin{pmatrix}
      0 & E \\ -E & 0
    \end{pmatrix},
    \begin{pmatrix}
      0 & -I \\ -I & 0
    \end{pmatrix} } \\
  G_{16}&=\gen{
    \begin{pmatrix}
      R & 0 \\ 0 & R
    \end{pmatrix},
    \begin{pmatrix}
      T & 0 \\ 0 & T
    \end{pmatrix},
    \begin{pmatrix}
      0 & -E \\ E & 0
    \end{pmatrix}}, \quad &
   G_{16}^D&=\gen{
     \begin{pmatrix}
       R & 0 \\ 0 & R
     \end{pmatrix},
     \begin{pmatrix}
       -T & 0 \\ 0 & T
     \end{pmatrix},
     \begin{pmatrix}
       0 & R \\ -R & 0
     \end{pmatrix}}
\end{align*}

%%% Local Variables: 
%%% mode: latex
%%% TeX-master: "dfam"
%%% End: 

%%%%%%%%%%%%%%%%%%%%%%%%%%%%%%%%%%%%%%%

\setcounter{subsection}{\value{thm}}

\section{Lie group theory}
\label{sec:lietheory}

Some facts from Lie theory are collected here.

\setcounter{subsection}{\value{thm}}
\subsection{Centerings}
\label{sec:centerings}\add

There are centerings \cite{plesken:orders} 
\begin{equation*}
  L(\Pin(2n)) \xrightarrow{P} L(\GL(2n,\R)) \xrightarrow{Q} L(\PGL(2n,\R)) 
\end{equation*}
where
\begin{equation*}
  P(x_1,x_2,\ldots ,x_n) = (2x_1-x_2- \cdots -x_n,x_2,\ldots ,x_n),
  \quad
  Q(x_1,x_2,\ldots ,x_n)  = (2x_1, x_1-x_2, \ldots , x_1-x_n)
\end{equation*}
The expression for $P$ is worked out in \cite[pp.\
174--175]{brockerdieck}. The expression for $Q$ follows from the
commutative diagram 
\begin{equation*}
  \xymatrix{
    {T(\GL(2n,\R)) = \U(1)^n} \ar[r] \ar@{->>}[dr]_{\varphi} & 
    {\U(1)^n/\gen{(-1,\ldots ,-1)} = T(\PGL(2n,\R)} \ar[d]_{\cong} \\
    & {\U(1)^n} }
\end{equation*}
where $\varphi(z_1,z_2,\ldots ,z_n)=(z_1^2,z_1z_2^{-1},\ldots
,z_1z_n^{-1})$ 
is surjective with kernel $C_2=\gen{(-1,\ldots ,-1)}$. Since the
action of the Weyl group $C_2 \wr \Sigma_n$ is known in 
$L(\GL(2n,\R))$, the two other actions can be worked out as well. The
action in $L(\Pin(2n))$ is $P^{-1}(C_2 \wr \Sigma_n) P$ and the
action in $L(\PGL(2n,\R))$ is  $Q (C_2 \wr \Sigma_n) Q^{-1}$. Here,
\begin{multline*}
  P^{-1}(u_1,u_2, \ldots ,u_n)= (\frac{1}{2}(u_1+\cdots +u_n),u_2,
  \ldots ,u_n), \\
   Q^{-1}(u_1,u_2, \ldots ,u_n)= (\frac{1}{2}u_1,\frac{1}{2}u_1-u_2,
   \ldots ,\frac{1}{2}u_1-u_n)
\end{multline*}
are the inverses.

\setcounter{subsection}{\value{thm}}
\subsection{Centralizers in semi-direct products}
\label{sec:semicent}\add

Let $G \rtimes W$ be the semi-direct product for a group action of $W$
on $G$. The following lemma is elementary.
\begin{lemma}\label{lemma:centrgw}
  For any $g \in G$ and $w \in W$,
  \begin{align*}
    &C_{G \rtimes W}(g,w)=\{(h,v) \mid \exists w \in C_W(w) \colon
    g(wh)=h(vg) \}, \\
    &C_{G \rtimes W}(g)=\{(h,v) \in G \rtimes W \mid vg=g^h\}, \qquad
    C_{G \rtimes W}(w)=G^w \rtimes C_W(w)
  \end{align*}
  where $G^w$ is the fixed point group for the action of $w$ on
  $G$. If $G$ is abelian then
  \begin{equation*}
    C_{G \rtimes W}(g) = G \rtimes W(g)
  \end{equation*}
  where $W(g)=\{w \in W \mid wg=g\}$ is is the isotropy subgroup at $g$.
\end{lemma}

Let \func{\mu}{V}{\ch{T} \rtimes W} be a group homo\m\ of an \lmntwo\ 
$V$ into the semi-direct product of a discrete \twoct\ $\ch{T}$ and a
group $W$. Write $\mu=(\ch{T}(\mu),W(\mu))$ for the two coordinates of
$\mu$. Then \func{W(\mu)}{V}{W} is a group homo\m\ and
\func{\ch{T}(\mu)}{V}{\ch{T}} a crossed homo\m\ into the $V$-module $V
\xrightarrow{W(\mu)} W \to \Aut(\ch{T})$. Let $H^1(V;\ch{T})$ be the
first cohomology group for this $V$-module and $[\ch{T}(\mu)] \in
H^1(V;\ch{T})$ the cohomology class represented by the crossed homo\m\ 
$\ch{T}(\mu)$.

\begin{lemma}
  There is a \ses\ 
  \begin{equation*}
    0 \to H^0(V;\ch{T}) \to C_{\ch{T} \rtimes W}(\mu) \to
    C_W(W(\mu))_{[\ch{T}(\mu)]} \to 1
  \end{equation*}
where $ C_W(W(\mu))_{[\ch{T}(\mu)]}$ is the isotopy subgroup at
$[\ch{T}(\mu)]$ for the action of $C_W(W(\mu))$ on $H^1(V;\ch{T})$. 
\end{lemma}
\begin{proof}
  We first determine the kernel of the homo\m\ $C_{\ch{T} \rtimes
    W}(\mu) \to C_W(W(\mu))$. Let $t \in \ch{T}$. Then
  \begin{align*}
    \text{$(t,1)$ commutes with $V$} & \iff
    \forall v \in V \colon (t,1)(\ch{T}(\mu)(v),W(\mu)(v))=
      (\ch{T}(\mu)(v),W(\mu)(v))(t,1) \\
    &\iff \forall v \in V \colon t +\ch{T}(\mu)(v) =
      \ch{T}(\mu)(v)+W(\mu)(v)(t) \\
    &\iff \forall v \in V \colon W(\mu)(v)t=t \\
     &\iff t \in H^0(V;\ch{T})
  \end{align*}
  More generally, for any element $(t,w) \in \ch{T} \rtimes W$ we have
  \begin{align*}
     \text{$(t,w)$ commutes with $V$} & \iff
      \forall v \in V \colon (t,w) (\ch{T}(\mu)(v),W(\mu)(v)) =
      (\ch{T}(\mu)(v),W(\mu)(v)) (t,w) \\
     &\iff  \forall v \in V \colon s+w\ch{T}(\mu)(v) =
     \ch{T}(\mu)(v) + W(\mu)(v)(t), wW(\mu)(v)= W(\mu)(v)w \\
     &\iff w \in C_{W}(W(\mu)) \text{\ and\ } 
       \forall v \in V \colon (1-w)\ch{T}(\mu)(v)=(1-W(\mu)(v))t \\
      &\iff w \in C_{W}(W(\mu)) \text{\ and\ } 
       \forall v \in V \colon w\ch{T}(\mu)(v)= \ch{T}(\mu)(v)-(1-W(\mu)(v))t
  \end{align*}
  It follows that for $w \in C_{W}(W(\mu))$ we have
  \begin{align*}
    w \in\im \left(C_{\ch{T} \rtimes W}(\mu) \to C_W(W(\mu)) \right) 
    & \iff \exists t \in \ch{T} \colon (t,w) \in C_{\ch{T}
      \rtimes W}(\mu) \\
     & \iff w[\ch{T}(\mu)]=[\ch{T}(\mu)]
  \end{align*}
  i.e.\ that $w$ fixes the crossed homo\m\ $\ch{T}(\mu)$ up to a
  principal crossed homo\m .
\end{proof}

\setcounter{subsection}{\value{thm}}
\subsection{Centers of semi-direct products}
\label{sec:centsemi}\add

Let $G\rtimes \Sigma$ be the semi-direct product for the action 
$\Sigma \to \Aut(G)$ of the 
group $\Sigma$ on the group $G$. Let $G^{\Sigma}=\{g\in G \vert \Sigma 
g=g\}$ and $\Sigma_G=\{\sigma\in\Sigma\vert \sigma(g)=g\text{\ for all
    $g\in G$}\}$.
\begin{lemma}\label{semicenter}
  The center ${Z}(G\rtimes\Sigma)=G^{\Sigma}\times_{\Aut(G)}{Z}(\Sigma)$
  of $G\rtimes\Sigma$ is the pull-back
  \begin{equation*}
    \xymatrix{
     {\Ze}(G\rtimes\Sigma) \ar[r] \ar[d] & {\Ze}(\Sigma)\ar[d]\\
     G^{\Sigma} \ar[r] & {\Aut(G)} }
  \end{equation*}
  of the action map restricted to the center of $\Sigma$ along the map 
  $G^{\Sigma} \to \Aut(G)$ given by inner auto\m s.
\end{lemma}
\begin{proof}
  Suppose that $(g,\sigma)\in G \times\Sigma$ is in the center of
  $G\rtimes\Sigma$. Since
  \begin{equation*}
    (g,\sigma) \cdot (1,\tau) = (g,\sigma\tau) = (1,\tau) \cdot
    (g,\sigma) = (\tau(g),\tau\sigma)
  \end{equation*}
  for all $\tau\in\Sigma$, $g$ is fixed by $\Sigma$ and $\sigma$ is
  central in $\Sigma$. Moreover, from
  \begin{equation*}
    (g,\sigma) \cdot (h,1) = (g\cdot\sigma(h),\sigma) = (h,1) \cdot
    (g,\sigma) = (hg,\sigma)
  \end{equation*}
  we see that $\sigma(h) = h^g$ for all $h\in G$. 
\end{proof}
\begin{cor}\label{cor:cent}
  If the center of $\Sigma$
    acts faithfully on $G$ through auto\m s that are not inner, then
    $\Ze(G \rtimes \Sigma)=\Ze(G)^{\Sigma}$.
  If $G$ is abelian, then ${\Ze}(G\rtimes\Sigma)=G^{\Sigma}\times
    {\Ze}(\Sigma)_G$ is a direct product. 
\end{cor}
\begin{proof}
  In the first case, the vertical map $\Ze(\Sigma) \to \Aut(G)$ is
  injective and its image intersects trivially with the image of the
  horizontal map $G^{\Sigma} \to \Aut(G)$. So the pull-back is
  $G^{\Sigma} \cap \Ze(G) = \Ze(G)^{\Sigma}$. In the second case,
  the bottom horizontal homo\m\ $G^{\Sigma} \to \Aut(G)$ is trivial.
\end{proof}
\begin{cor}\label{cycliccenter}
  Let $G$ be a group and ${\Ze} \neq G$ a central subgroup. Let the cyclic
  group $C_p$ of prime order $p$ act on $G^p/{\Ze}$ by cyclic permutation.
  Then
  \begin{equation*}
    {\Ze}(G)/{\Ze} \times \{ z \in {\Ze} \vert z^p=1\} \cong {\Ze}(G^p/{\Ze}\rtimes C_p)
  \end{equation*}
via the iso\m\ that takes the element $z\in {\Ze}$ of order $p$ to
$(1,z,\ldots ,z^{p-1}){\Ze} \in G^p/{\Ze}$ and is the diagonal on
${\Ze}(G)/{\Ze}$. 
\end{cor}
\begin{proof}
  Observe that
  \begin{equation*}
    G/{\Ze} \times \{ z \in {\Ze} \vert z^p=1\} \xrightarrow{\cong}
    \big(G^p/{\Ze}\big)^{C_p} 
  \end{equation*}
  via the iso\m\ that takes $(g{\Ze},z)$ to $g(1,z,\ldots ,z^{p-1}){\Ze}$. To
  see this, consider an element $(g_1,\ldots ,g_p){\Ze}$ which is fixed by 
  $C_p$. Then $(g_1,g_2,\ldots ,g_p){\Ze}=(g_p,g_1,\ldots ,g_{p-1}){\Ze}$ so
  there exists an element $z\in {\Ze}$ so that
  $g_2=g_1z,g_3=g_2z=g_1z^2,\ldots ,
  g_p=g_1z^{p-1},g_1=g_1z^p$. Therefore, $z^p=1$ and $(g_1,g_2,\ldots
  ,g_p)=g_1(1,z,\ldots ,z^{p-1})$. 

  Thus ${\Ze}(G^p/{\Ze}\rtimes C_p)$ is the pull back of the group homo\m s
  \begin{equation*}
    G/{\Ze} \times  \{ z \in {\Ze} \vert z^p=1\} \xrightarrow{\varphi}
    \Aut(G^p/{\Ze}) \leftarrow C_p
  \end{equation*}
  where $\varphi(g{\Ze},z)((g_1,\ldots ,g_p){\Ze})=(g_1^g,\ldots
  ,g_p^g){\Ze}$. Let $((g{\Ze},z),\sigma)$ be an element of the pull
  back. Assume that $\sigma$ is non-trivial. Since $p$ is a prime
  number, $\sigma$ has no fixed points. The equation
  \begin{equation*}
    \forall g_1,\ldots ,g_p\in G \colon (g_1^g,\ldots ,g_p^g){\Ze}=
    (g_{\sigma(1)},\ldots ,g_{\sigma(p)}){\Ze}
  \end{equation*}
  shows that $g_1^g{\Ze}=g_{\sigma(1)}{\Ze}$. This is impossible unless
  $\sigma$ is the identity since otherwise we can
  find a $g_1\in {\Ze}$ and a $g_{\sigma(1)} \not\in {\Ze}$. Thus the
  permutation $\sigma$ must be the identity. The requirement for
  $((g{\Ze},z),1)$ to be in the pull back is that
  \begin{equation*}
    \forall (g_1,\ldots ,g_p)\in G^p \exists u\in {\Ze} \colon
    (g_1^g,g_2^g,\ldots ,g_p^g)=(g_1u,g_2u,\ldots ,g_pu)
  \end{equation*}
  which implies that $[g_1,g]=u=[g_2,g]$ for all $g_1,g_2\in G$. If we
  take $g_1=1$ to be the identity, we see that $g$ must be central.
\end{proof}

\setcounter{subsection}{\value{thm}}
\subsection{Centers of Lie groups and \pcg s}
\label{sec:centLie}\add
Let $Y$ be a compact connected Lie group and $ZY$ its center. Let $BY$
denote the $p$-completed classifying space of $Y$, ie the \pcg\ 
associated to $Y$. Lie group multiplication $ZY \times Y \to Y$
induces a \he\ $BZY \to \map(BY,BY)_{B1}$ \cite[1.4]{dw:center} of the
$p$-completed classifying space $BZY$ to the the mapping space
component containing the identity map. We need a version that holds
for nonconnected Lie groups as well.

Let $G$ be a possible nonconnected Lie group and $ZG$ its center. Let
$BZG$ and $BG$ denote the $\F_p$-localized classifying spaces of $ZG$
and $G$, respectively. The space $\map(BG,BG)_{B1}$ is the center of
the \pcg\ $BG$ \cite[1.3]{dw:center}.

\begin{lemma}\label{lemma:ZG}
  The map
  \begin{equation*}
    BZG \to \map(BG,BG)_{B1}
  \end{equation*}
   induced by Lie group multiplication $ZG \times G \to G$, is a weak
   homotopy equivalence.
\end{lemma}
\begin{proof}
  Let $Y$ be the identity component of $G$ and $\pi=G/Y$ the group of
  components. Note that the group $\pi$ acts on the center $ZY$ of
  $Y$ and that there is an exact sequence of abelian groups
  \begin{equation*}
    1 \to H^0(\pi;ZY) \to ZG \to Z\pi \to H^1(\pi; ZY)
  \end{equation*}
  relating the centers $ZY$, $ZG$, and $Z\pi$, of $Y$, $G$, and
  $\pi$. The abelian Lie group $ZG$, a product of a torus and a
  finite abelian group, is described by the data
  \begin{gather*}
    \pi_1(ZG) \otimes \Q \cong  H^0(\pi;\pi_1(ZY) \otimes \Q) \\
    1 \to H^0(\pi; \pi_0(ZY)) \times H^1(\pi; \pi_1(ZY)) \to \pi_0(ZG)
    \to Z\pi \to H^1(\pi;\pi_0(ZY)) \times  H^2(\pi;\pi_1(ZY))
  \end{gather*}
  where the second line is an exact sequence. 

  Similarly, there is a fibration of mapping spaces 
  \begin{equation*}
    \map(BY,BY)^{h\pi} \to \map(BG,BG) \to \map(BG,B\pi)
  \end{equation*}
  where the fibre over $BG \xrightarrow{B\pi_0} B\pi$ is the space
  $\map(BY,BY)^{h\pi}$ of self-maps of $BG$ over $B\pi$. 
  If we restrict to a single component of the total space, we obtain
  a fibration
  \begin{equation*}
    \map(BG,BG)_{B1} \to \map(BG,B\pi)_{B\pi_0} \simeq
    \map(B\pi,B\pi)_{B1} 
  \end{equation*}
  between path-connected spaces. The base space is $BZ\pi$. The fibre
  consists of some path-components of
  $\map(BY,BY)_{B1}^{h\pi}=(BZY)^{h\pi}$, the space of self-maps of
  $BG$ over $B\pi$ with restriction to $BY$ homotopic to the identity
  map. We have 
  \begin{equation*}
    \pi_i((BZY)^{h\pi})=H^{1-i}(\pi;\pi_0(ZY)) \times
    H^{2-i}(\pi;\pi_1(ZY)) 
  \end{equation*}
  because $BZY=K(\pi_0Y,1) \times K(\pi_1Y,2)$ is a product of
  Eilenberg--MacLane spaces \cite[3.1]{mn:center}
  \cite[1.1]{dw:center}. It follows that $\map(BG,BG)_{B1}$ is an
  abelian \cite[3.5, 8.6]{dw:fixpt} \pctg\ described by the data
  \begin{gather*}
    \pi_2((BZY)^{h\pi}) \otimes \Q \cong \pi_2(\map(BG,BG)_{B1})
    \otimes \Q  \\
    1 \to  \pi_1((BZY)^{h\pi}) \to \pi_1(\map(BG,BG)_{B1}) \to Z\pi
    \to \pi_0((BZY)^{h\pi}) 
  \end{gather*}
  where the second line is an exact sequence is an exact sequence. 

  Finally, the left commutative diagram of Lie groups 
  \begin{equation*}
    \xymatrix{
      (ZY)^{\pi} \times G \ar[d] \ar[r] & G \ar@{=}[d] & 
      B(ZY)^{\pi} \ar[r] \ar[d] & (BZY)^{h\pi} \ar[d]\\
      ZG \times G \ar[d] \ar[r] & G \ar[d]  \ar@{~>}[r]&
      BZG \ar[r] \ar[d] & {\map(BG,BG)_{B1}} \ar[d] \\
      Z\pi \times {\pi} \ar[r] & {\pi} &
      BZ\pi \ar[r] & {\map(B\pi,B\pi)_{B1}}}
  \end{equation*}
  induces the right commutative diagram of mapping space. Compairing
  the homotopy groups, we see that $BZ(G) \to \map(BG,BG)_{B1}$ is
  weak homotopy equivalence.
\end{proof}

\begin{lemma}\label{lemma:ZGLprod}
  We have 
  \begin{gather*}
    Z\left( \frac{\GL(i_0,\R) \times \cdots \GL(i_t,\R)}{\gen{-E}}\right) = 
      \frac{\gen{-E} \times \cdots \times \gen{-E}}{\gen{-E}} \cong
      C_2^{t} \\
    Z\left( \frac{\SL(n,\R) \cap (\GL(i_0,\R) \times \cdots
        \GL(i_t,\R))}{\gen{-E}}\right) =  
      \frac{\SL(n,\R) \cap (\gen{-E} \times \cdots \times
        \gen{-E})}{\gen{-E}}  \\
    Z\left( \frac{\GL(i_0,\Ha) \times \cdots \GL(i_t,\Ha)}{\gen{-E}}\right) = 
      \frac{\gen{-E} \times \cdots \times \gen{-E}}{\gen{-E}} \cong
      C_2^{t}
  \end{gather*}
  for all natural numbers $i_0,\ldots, i_t >0$ with sum $n$ (which, in
  the second line, is even).
\end{lemma}
\begin{proof}(For the case where the field is $\R$.)
  Put $G=\GL(i_0,\R) \times \cdots \GL(i_t,\R)$. There is
  \cite[5.11]{jmm:ndet} a \ses\ 
  \begin{equation*}\label{ses:prodGL}
    1 \to Z(G)/\gen{-E} \to Z(G/\gen{-E}) \to \Hom(G,\gen{-E})_{\mathrm{id}}
    \to 1
  \end{equation*}
  where the group to the right consists of all homo\m s
  \func{\phi}{G}{\gen{-E}} such that the map $g \to \phi(g)g$ is
  conjugate to the identity of $G$.  Let $B \in \GL(i_j,\R)$ be any
  matrix of positive trace. Then $\phi(E,\ldots,E,B,E,\ldots ,E)=E$
  since the map $g \to \phi(g)g$ preserves trace. It follows that
  $\phi(g)=E$ for all $g \in G$ since $\phi$ is constant on the $2^t$
  components of $G$. Thus the group to the right in the above \ses,
  $\Hom(G,\gen{-E})_{\mathrm{id}}$, is trivial.

  Since 
  \begin{equation*}
  Z\big(\SL(n,\R) \cap \prod \GL(i_j,\R)\big) \subset C_{\prod
    \GL(i_j,\R)}(\prod \SL(i_j,\R)) = \prod Z\GL(i_j,\R)  
  \end{equation*}
  we see that $Z\big(\SL(n,\R) \cap \prod \GL(i_j,\R)\big)=\SL(n,\R)
  \cap \prod Z\GL(i_j,\R)$. Suppose that the homo\m\ 
  \func{\phi}{\SL(n,\R) \cap \prod \GL(i_j,\R)}{\gen{-E}} is  
  such that the map $g \to \phi(g)g$ is conjugate to the identity. Let
  $B_1 \in \GL(i_{j_1},\R)$ and $B_2 \in \GL(i_{j_2}),\R)$ be any pair
  of matrices such that $\tr(B_1)+\tr(B_2)>0$. Then
  $\phi(E,\ldots, B_{i_1},\ldots , B_{i_2}, \ldots ,E)=E$ by trace
  considerations. The \ses\ from  \cite[5.11]{jmm:ndet}, similar to
  (\ref{ses:prodGL}), now yields the formula of the second line. 
 
  The formula of the third line has a similar proof. 
\end{proof}

It is not true in general that $Z(G)/Z$ is the center of the quotient
$G/Z$ of the Lie group $G$ by the central subgroup $Z$.

\setcounter{subsection}{\value{thm}}
\subsection{Centralizers in quotients}
\label{sec:centralquot}\add

Let $G$ be a Lie group and $Z \subset G$ a central subgroup. Write
$N(G)$ for the normalizer of the \mt , $T(G)$, and $W=W(G)=N(G)/T(G)$
for the Weyl group. Suppose that $V \subset T(G)/Z$ is a toral
subgroup of the quotient Lie group $G/Z$ and let $V^* \subset T(G)
\subset G$ be the preimage of $V$ in $G$. 

There is an exact sequence
\begin{equation*}
  1 \to W(V^*) \to W(V) \to \Hom(V^*,Z)
\end{equation*}
relating the point-wise stabilizer subgroups for the action of the
Weyl group $W$ on $V^*$ and $V$. The image of homo\m\ to the right
consists of all $\zeta \in \Hom(V^*,Z)$ for which the auto\m\ of $V^*$
given by $v^* \to \zeta(v^*)v^*$, $v^* \in V^*$, is of the form $v^*
\to wv^*$ for some Weyl group element $w \in W$.

Similarly, there is an exact sequence \cite[5.11]{jmm:ndet}
\begin{equation*}
  1 \to C_G(V^*)/Z \to C_{G/Z}(V) \to \Hom(V^*,Z)
\end{equation*}
relating the centraizers of $V^* \subset G$ and $G \subset G/Z$. The
image of homo\m\ to the right consists of all $\zeta \in \Hom(V^*,Z)$
for which the auto\m\ of $V^*$ given by $v^* \to \zeta(v^*)v^*$, $v^*
\in V^*$, is of the form $v^* \to g^{-1}v^*g$ for some $g \in G$.

\begin{lemma}\label{lemma:WVWVast}
  $W(V)/W(V^*) =  C_{G/Z}(V)/C_G(V)$.
\end{lemma}
\begin{proof}
  Any auto\m\ of the toral subgroup $V^*$ that is induced by
  conjugation with an element of $G$ is in fact induced by conjugation
  with an element of $N(G)$ \cite[IV.2.5]{brockerdieck} and hence
  agrees with the action of a Weyl group element.
\end{proof}

\setcounter{subsection}{\value{thm}}
\subsection{Action on centralizers in Lie case}
\label{sec:actLie}\add

Let \func{\nu}{V}{G} be a mono\m\ of a non-trivial \lmn\ to a compact
Lie group $G$. There is a canonical map $\B C_G(\nu(V)) \to \map(\B V,
\B G)_{\B\nu}$ from the classifying space of the Lie theoretic
centralizer of $\nu(V)$ to the mapping space component containing
$\B\nu$. Write $c_g$ for conjugation with $g\in G$.
\begin{lemma}
  Suppose that $\nu\alpha=c_g\nu$ for some element $g \in G$ and some
  auto\m\ $\alpha\in\GL(V)$. Then conjugation by $g$ takes
  $C_G(\nu(V))$ to $C_G(c_g\nu(V)) = C_G(\nu\alpha(V)) = C_G(\nu(V))$
  and the diagram
  \begin{equation*}
    \xymatrix{
    {\B C_G(\nu(V))} \ar[r] &
    {\map(\B V,\B G)_{\B\nu}} \ar[d]^-{(\B\alpha)^*}_-{\cong} \\
    {\B C_G(\nu(V))} \ar[u]^-{\B c_g}_-{\cong} \ar[r] &
    {\map(\B V,\B G)_{\B\nu}} }
  \end{equation*}
  is homotopy commutative.
\end{lemma}
\begin{proof}
  The commutative diagram of Lie group \m s
  \begin{equation*}
    \xymatrix{
      V \times C_G(\nu(V)) \ar[d]_-{\alpha\times c_g}
      \ar[r]^-{\nu\times 1} &
      {\nu(V) \times  C_G(\nu(V))} \ar[r]^-{\mathrm{mult}} & 
      G \ar@{=}[d] \\
      V \times C_G(\nu(V)) \ar[r]_-{\nu\times 1} &
      {\nu(V) \times  C_G(\nu(V))} \ar[r]_-{\mathrm{mult}} & 
      G }
  \end{equation*}
  induces a commutative diagram
  \begin{equation*}
    \xymatrix@C=50pt{
      {\B V \times \B C_G(\nu(V))} 
      \ar[r]^-{\B (\mathrm{mult} \circ (\nu \times 1))}
      \ar[d]_-{\B\alpha\times\B c_g} &
      {\B G} \ar@{=}[d] \\
      {\B V \times \B C_G(\nu(V))} 
      \ar[r]^-{\B (\mathrm{mult} \circ (\nu \times 1))} &
      {\B G} }
  \end{equation*}
 of classifying spaces. Taking adjoints, we obtain the homotopy
 commutative diagram 
 \begin{equation*}
   \xymatrix{
      {\B C_G(\nu(V))} \ar[r] &
      {\map(\B V,\B G)_{\B\nu}} \ar[d]^-{(\B\alpha)^*} \\
      {\B C_G(\nu(V))} \ar[u]^-{\B c_g} \ar[r] &
      {\map(\B V,\B G)_{\B\nu}} }
 \end{equation*}
 as claimed.
\end{proof}
\begin{cor}
  Suppose that \func{\mu}{V}{{\N}(G)} is a mono\m\ and that
  $\mu\alpha=c_n\mu$ for some $\alpha\in\GL(V)$ and $n\in {\N}(G)$. Then 
  \begin{equation*}
    w^{-1}=\pi_2((\B\alpha)^*) \colon \pi_2(\B {\T}(G))^{\pi_0(\mu)(V)}
    \to \pi_2(\B {\T}(G))^{\pi_0(\mu)(V)}
  \end{equation*}
  where $w\in {\W}(G)$ is the image of $n\in {\N}(G)$.
\end{cor}
\begin{proof}
  There is a commutative diagram
  \begin{equation*}
    \xymatrix{
      {\pi_2(\B {\T})} \ar@{=}[r] &
      {\pi_2(\B {\N}(G))} &
      {\pi_2(\B C_{{\N}(G)}(V,\mu))} \ar@{_(->}[l] \ar[r]^-{\cong} &
      {\pi_2(\map(\B V,\B {\N}), \B\mu)} \ar[d]^-{\pi_2((\B\alpha)^*)} \\
      {\pi_2(\B {\T})} \ar@{=}[r] \ar[u]^-{w} &
      {\pi_2(\B {\N}(G))} \ar[u]_-{\pi_2(\B c_n)} &
      {\pi_2(\B C_{{\N}(G)}(V,\mu))} \ar@{_(->}[l] \ar[r]^-{\cong}
      \ar[u]_-{\pi_2(\B c_n)} &
      {\pi_2(\map(\B V,\B {\N}), \B\mu)} }
\end{equation*}
where $\pi_2(\B C_{{\N}(G)}(V,\mu)) = \pi_2(\B {\T}(G))^{\pi_0(\mu)(V)}$ denotes
the fixed point group for the group action \func{\pi_0(\mu)}{V}{{\W}(G)
  \subseteq \Aut(\pi_2(\B {\T}(G)))}. Since \func{\B c_n}{\B {\N}}{\B {\N}} is
freely homotopic to the identity along the loop $w \in \pi_1(\B {\N})$
its effect on the $\Z_p[\pi_1(\B {\N})]$-module $\pi_2(\B {\N})$ is
multiplication by $w$.
\end{proof}

\setcounter{subsection}{\value{thm}}
\subsection{Low degree identifications}
\label{lowdegree}\add

There are the following low degree identifications \cite[pp.\ 61,
292]{brockerdieck} \cite[above def 3.3]{matthey:second}
\set\begin{equation}\label{eq:lowdegree}
  \begin{split}
    &  \Spin(3)= \Symp(1) = \SU(2), 
                \quad \SO(3)=  \mathrm{PSp}(1)=\mathrm{PSU}(2) \\
    &  \Spin(4) = \Spin(3) \times \Spin(3) = \SU(2) \times \SU(2), \quad
      \mathrm{PSO}(4) = \SO(3) \times \SO(3) \\
    & \Spin(5)=\mathrm{Sp}(2), \quad \SO(5)=\mathrm{PSp}(2) \\
      & \Spin(6) = \SU(4), \quad \mathrm{PSO}(6) = \mathrm{PSU}(4)
  \end{split}
\add\end{equation}

\providecommand{\bysame}{\leavevmode\hbox to3em{\hrulefill}\thinspace}
\providecommand{\MR}{\relax\ifhmode\unskip\space\fi MR }
% \MRhref is called by the amsart/book/proc definition of \MR.
\providecommand{\MRhref}[2]{%
  \href{http://www.ams.org/mathscinet-getitem?mr=#1}{#2}
}
\providecommand{\href}[2]{#2}

\end{document}